\documentclass[11pt]{article}
\usepackage{etex}

\usepackage[margin=1in]{geometry}
\usepackage{stackengine}

\usepackage{stmaryrd}

\usepackage{color}
\usepackage[latin1]{inputenc}
\usepackage[T1]{fontenc}
\usepackage[normalem]{ulem}
\usepackage[english]{babel}
\usepackage{verbatim}
\usepackage{graphicx}
\usepackage{enumerate}
\usepackage{amsmath,amssymb,amsfonts,amsthm,amscd,mathrsfs}
\usepackage{array}

\usepackage{dsfont}

\usepackage[T1]{fontenc}
\usepackage{babel}

\usepackage{tocloft}
\usepackage{url}

\usepackage{caption}
\usepackage{subcaption}
\usepackage[bookmarksopen, bookmarksnumbered]{hyperref}

\usepackage{tikz}
\usetikzlibrary{arrows}
\usetikzlibrary{arrows.meta}
\definecolor{wwhhii}{rgb}{1.,1.,1.}
\definecolor{rreedd}{rgb}{1.,0.,0.}
\definecolor{uuuuuu}{rgb}{0.26666666666666666,0.26666666666666666,0.26666666666666666}
\definecolor{darkgreen}{HTML}{0d8513}

\newtheorem{theorem}{Theorem}[section]

\newtheorem{lemma}[theorem]{Lemma}
\newtheorem{cla}[theorem]{Claim}
\newtheorem{prop}[theorem]{Proposition}
\newtheorem{cor}[theorem]{Corollary}

\theoremstyle{definition}
\newtheorem{defn}[theorem]{Definition}
\newtheorem{rem}[theorem]{Remark}

\input xy
\xyoption{all}

\DeclareMathOperator{\argmax}{argmax}
\DeclareMathOperator{\Exp}{Exp}

\DeclareMathOperator{\diam}{diam}
\DeclareMathOperator{\Dist}{Disj}

\DeclareMathOperator{\Graph}{Graph}

\DeclareMathOperator{\ConsDiff}{ConstDiff}
\DeclareMathOperator{\Good}{Good}

\newcommand{\E}{\mathbb E}
\newcommand{\PP}{\mathbb P}
\newcommand{\Z}{\mathbb Z}
\newcommand{\R}{\mathbb R}

\newcommand{\Q}{\mathbb Q}
\newcommand{\N}{\mathbb N}

\newcommand{\LL}{\mathbb L}

\newcommand{\HH}{\mathbb H}

\newcommand{\don}{\mathds{1}}

\newcommand{\cC}{\mathcal C}
\newcommand{\cD}{\mathcal D}
\newcommand{\cK}{\mathcal K}
\newcommand{\cP}{\mathcal P}

\newcommand{\cB}{\mathcal B}

\newcommand{\cF}{\mathcal F}
\newcommand{\cE}{\mathcal E}
\newcommand{\cA}{\mathcal P}
\newcommand{\ocA}{\overline{\mathcal P}}
\newcommand{\cL}{\mathcal L}
\newcommand{\cM}{\mathcal M}

\newcommand{\cS}{\mathcal S}

\newcommand{\cU}{\mathcal U}

\newcommand{\sS}{\mathscr S}
\newcommand{\sZ}{\mathscr Z}

\newcommand{\LE}{\mathscr L}

\newcommand{\LV}{\vartheta}
\newcommand{\KA}{\kappa}

\newcommand{\fG}{\mathsf R}

\newcommand{\fr}{\mathsf r}

\newcommand{\bc}{\mathbf c}

\newcommand{\boo}{\mathbf 0}

\newcommand{\by}{\mathbf y}
\newcommand{\bx}{\mathbf x}

\newcommand{\oK}{\overline{K}}

\newcommand{\hcL}{\hat{\cL}}

\newcommand{\hK}{\hat{K}}

\newcommand{\wdf}{\mathcal{D}}
\newcommand{\NC}{\mathsf{NC}}
\newcommand{\hNC}{\overline{\NC}}
\newcommand{\vNC}{\widehat{\NC}}

\makeatletter
\renewcommand\tableofcontents{%
  \null\hfill\textbf{\contentsname}\hfill\null\par
  \@mkboth{\MakeUppercase\contentsname}{\MakeUppercase\contentsname}%
  \@starttoc{toc}%
}

\g@addto@macro\normalsize{%
  \setlength\abovedisplayskip{5pt}
  \setlength\belowdisplayskip{5pt}
  \setlength\abovedisplayshortskip{3pt}
  \setlength\belowdisplayshortskip{3pt}
}

\makeatother

\makeatletter
\renewenvironment{abstract}{
    \if@twocolumn
      \section*{\abstractname}
    \else
      \begin{center}
        {\bfseries \small\abstractname\vspace{\z@}}
      \end{center}
      \scriptsize\quotation
    \fi}
    {\if@twocolumn\else\endquotation\fi}
\makeatother

\numberwithin{equation}{section}

\begin{document}
\title{Fractal geometry of the space-time difference profile in the directed landscape via construction of geodesic local times}
\author{Shirshendu Ganguly
\thanks{Department of Statistics, UC Berkeley, Berkeley, CA, USA. e-mail: sganguly@berkeley.edu}
\and
Lingfu Zhang
\thanks{Department of Mathematics, Princeton University, NJ, USA. e-mail: lingfuz@math.princeton.edu}
}
\date{}

\maketitle

\begin{abstract}
The Directed Landscape, $\cL$,  a random four parameter energy field, was constructed in the breakthrough work of Dauvergne, Ortmann, and Vir\'ag \cite{DOV}, and has since been shown in \cite{DV21} to be the scaling limit of various integrable models of Last Passage percolation, a central member of the Kardar-Parisi-Zhang universality class, where one considers the maximum passage time between points attained by weight maximizing paths termed as  geodesics. Here $\cL(x,s;y,t)$ denotes the scaled  passage time between the points $(x,s)$ and $(y,t).$
$\cL$ exhibits several scale invariance properties making it a natural source of rich fractal behavior, a topic that has been of constant investigation since its construction.
This was initiated in \cite{BGH}, where a novel object called the difference profile: $\cD(y,t)=\cL(-1,0;y,t)-\cL(-1,0;y,t)$, the difference of passage times from $(-1,0)$ and $(1,0)$ to $(y,t)$ was introduced.  Owing to geodesic geometry, it turns out that this process is almost surely locally constant and the set of non-constancy inherits remarkable fractal properties. In particular, for any fixed $t,$ the ``spatial process'' $\cD(\cdot, t)$ is a monotone function, which can be interpreted as the distribution function of a random fractal measure. Across the articles \cite{BGH,BGHhau,GH21,D21last}, it has been established that this random measure is supported on a fractal set of Hausdorff dimension $1/2$, connected to disjointness of geodesics and bears a rather strong resemblance to Brownian local time. 
However, the arguments in the above works crucially rely on the aforementioned monotonicity property which is absent when the temporal structure of $\cD(\cdot,\cdot)$ is probed, necessitating the development of new methods. 

In this paper, we put forth several new ideas en-route to computing new fractal exponents governing the same. More precisely, we show that the set of non-constancy of  the two dimensional process $\cD(\cdot,\cdot)$ and the one dimensional temporal process $\cD(0,\cdot)$ (a mean zero process owing to symmetry) have Hausdorff dimensions $5/3$ and $2/3$ respectively. 

Beyond the construction of infinite geodesics, Busemann functions, competition interfaces for $\cL$ (also constructed independently and simultaneously in \cite{RV21} and \cite{BSS22}), all objects of much broader interest, a particularly crucial ingredient in our analysis is the novel construction of a local time process for the geodesic akin to Brownian local time, supported on the ``zero set'' of the geodesic. Further, we show that the latter has Hausdorff dimension $1/3$ in contrast to the zero set of Brownian motion which has dimension $1/2.$  

\maketitle
\begin{figure}[h]
\begin{center}
\includegraphics[height=0.27\linewidth]{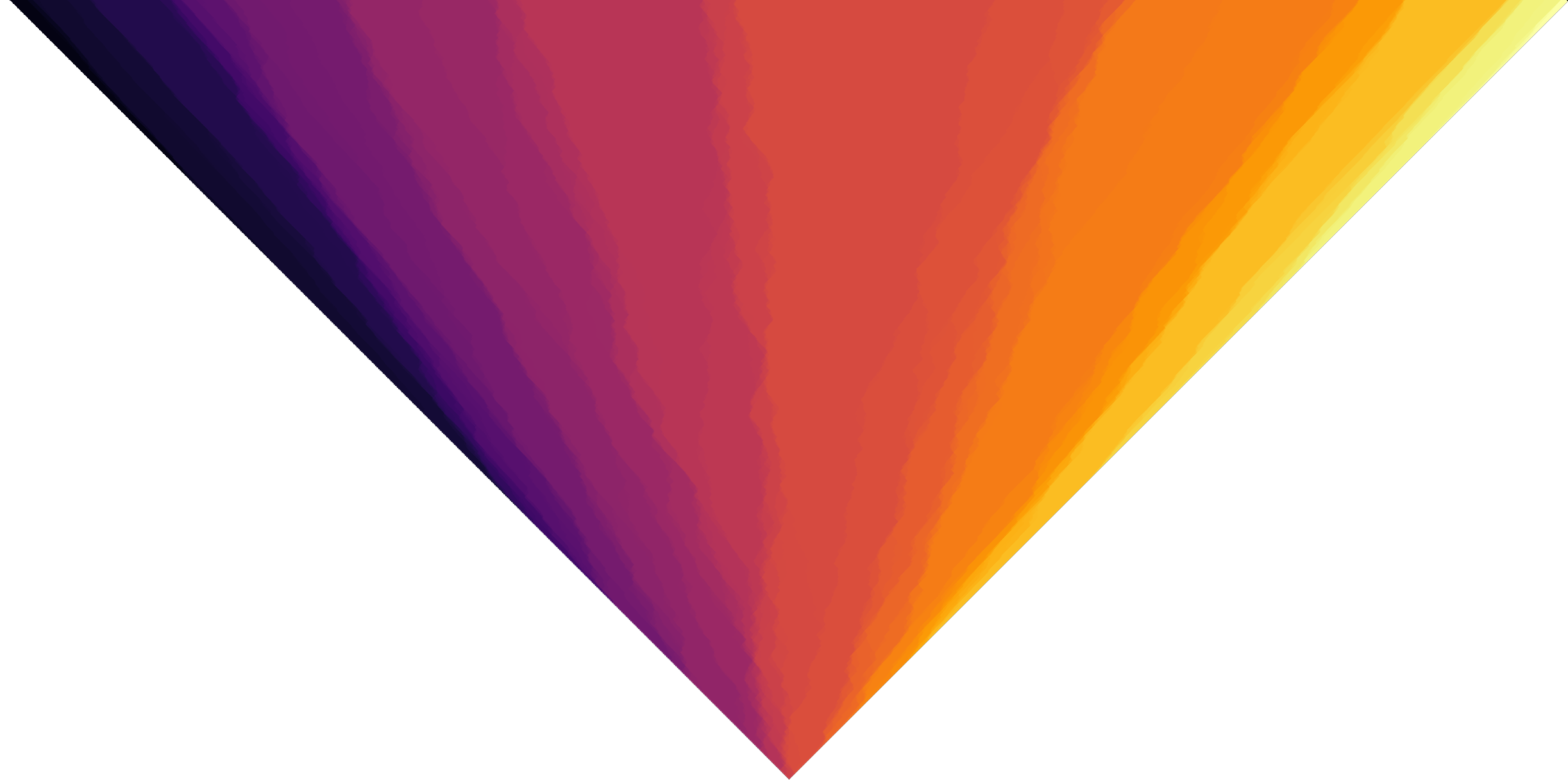}
\end{center}
\caption{\scriptsize{Simulation of the difference profile for the discrete model of Exponential LPP on the positive quadrant $\Z^2_+$, rooted at the points $(0,50)$ and $(50,0)$ and hence defined on the cone emanating from (50,50) (rotated by $45^{\circ}$). Colors denote the level sets. (Courtesy of Milind Hegde.)}}
\label{f.simulations}
\end{figure}

\newpage

\end{abstract}

\setcounter{tocdepth}{2}
\tableofcontents

\section{Introduction}

The Kardar-Parisi-Zhang (KPZ) universality class encompasses  a broad family of models of one-dimensional random growth 
which are believed to exhibit certain common features, such as local growth driven by white noise and a 
smoothing effect from some notion of surface tension and a slope dependent non-linear growth.
 It is expected that the fluctuation theories of such models are dictated by universal scaling exponents and limiting distributions, 
exhibited by solutions to a stochastic PDE known as the KPZ equation \cite{Corwin12,Corwin16,HQ18}.
The models known or believed to belong to this collection, include asymmetric exclusion processes, first and last passage percolation, and directed polymers in random media.
For these $(1+1)$-dimensional models, the resulting  evolution of the growth interface $h(t,x)$ is described by  the characteristic triple $(1,\frac{1}{3},\frac{2}{3})$ of exponents: at time $t$, the value of $h(t,x)$ is of order $t^1$, with order $t^{1/3}$ deviations from its mean, and further, non-trivial correlations are observed when $x$ is varied on the scale of $t^{2/3}$.
Furthermore, $h(t,x)$ on proper centering and scaling  
according to these exponents, is expected to converge to a universal limit $t\to\infty$ \cite{QR14, MQR}. 

Nonetheless, in spite of the breadth of models thought to lie in the class, only a handful have been rigorously shown to do so.
An important subclass of models 
believed to be members of the KPZ universality class is 
known as last passage percolation (LPP). Ignoring microscopic specifications, in general terms they all consist of an environment of random noise
 through which \emph{directed} paths travel, accruing the integral of the noise along it---a quantity known as energy or weight. Given two endpoints, a maximization is done 
over the weights of all paths with these endpoints with the 
optimizing path called a \emph{geodesic}, although, traditionally, the latter is used to denote \emph{shortest} paths in a metric space.

While most of the above picture is conjectural, recent 
developments have 
confirmed the convergence of several exactly solvable models to a well-defined scaling limit.
This object was introduced in \cite{DOV} and named the \textit{directed landscape} which possesses interesting fractal properties which has been the subject of intense recent study since its construction (see e.g., \cite{Ganguly21} for an account of such and related advances).

This paper continues this program and makes several novel contributions. Before elaborating on this further, let us introduce the key object of study formally.

\subsection{Directed landscape}
The directed landscape $\cL$ is a random continuous function from the parameter space 
\[
\R^4_\uparrow = \{u = (p; q) = (x, s; y, t) \in \R^4 : s < t\}
\]
to $\R$ constructed in the breakthrough work \cite{DOV}. It satisfies the following composition law 
\begin{equation}  \label{eq:DL-compo}
\cL(x,r;z,t) = \max_{y\in \R}     \cL(x,r;y,s) + \cL(y,s;z,t)
\end{equation}
for any $x,z \in \R$ and $r<s<t$.
For a continuous path $\pi:[s,t] \to \R$, its weight under the metric $\cL$ is given by
\[
\|\pi\|_\cL = \inf_{k\in \N} \inf_{s=t_0< t_1 < \cdots < t_k = t} \sum_{i=1}^k \cL(\pi(t_{i-1}), t_{i-1}; \pi(t_i), t_i).
\]
Such a path $\pi$ is said to be a geodesic from $(\pi(s), s)$ to $(\pi(t), t)$ if $\|\pi\|_\cL = \cL(\pi(s), s; \pi(t), t)$.
It is shown in \cite{DOV} that almost surely, there exists at least one geodesic from $p$ to $q$ for any $u = (p; q) = (x, s; y, t) \in \R^4_\uparrow$.
We let $\pi_u = \pi_{(p;q)}=\pi_{(x,s;y,t)}$ denote any such geodesic from $p$ to $q$.
Further, for any fixed $u\in \R^4_\uparrow$, almost surely the geodesic $\pi_u$ is unique.\\
\noindent
\textbf{Semi-infinite geodesics:} Going beyond finite geodesics, 
{for each $p=(x,s) \in \R^2$ and $r\in \R$, we denote by $\pi_p^r:[s,\infty)\to \R$ a semi-infinite geodesic started from $p$ in the $r$ direction; i.e. we have $t^{-1}\pi_p^r(t) \to r$ as $t\to\infty$, and for any $s\le s_1 < s_2 <\infty$, the restriction of $\pi_p^r$ on $[s_1,s_2]$ is a geodesic (from $(\pi_p^r(s_1), s_1)$ to $(\pi_p^r(s_2), s_2)$).
When $r=0$, we will drop the $r$ dependence and simply use $\pi_p=\pi_p^0$. The issues regarding the existence and uniqueness of such semi-infinite geodesics are not straightforward and will be addressed in Section \ref{sec:buse-semiinf-dl}.} \\

The landscape admits invariance under appropriate scalings guided by the exponent triple $(1,1/3,2/3)$, much like Brownian motion which is invariant under diffusive scaling. Hence, as in the case of the latter, $\cL$ is expected to exhibit various fractal or self-similar properties.   

The present paper continues a recently initiated program devoted to studying such properties and makes several novel contributions.

While the main motivation of the paper lies in studying an object called the ``\emph{difference profile}'', to be defined shortly, en-route we introduce the \emph{geodesic local time} (GLT), analogous to the classical Brownian local time (BLT), an object we expect to be independently interesting on its own, much akin to its Brownian counterpart which has witnessed intense investigation since its construction by L\'evy (see the monograph \cite{BC21} for a beautiful account of the remarkable properties of BLT, and \cite[Section 6]{morters2010brownian} for a comprehensive treatment).

Thus, one of the main contributions of this paper is to formally construct, for the semi-infinite geodesic $\pi_{(0,0)}$, its local time, a non-decreasing function $L:\R_{\ge 0} \to \R$ recording how much time $\pi_{(0,0)}$ has spent at $0$. Further, while much remains to be understood about $L$,  we do succeed in establishing some of its properties, tailored for our applications, such as H\"older regularity and the fractal dimension of its support.

We now move towards defining the aforementioned difference profile, an object introduced in \cite{BGH} and subsequently the subject of  several recent works \cite{BGHhau,GH21, D21last}.

To motivate its origin, we begin by pointing out that another central object in the KPZ universality class which had been rigorously established much before $\cL$, is known as the parabolic Airy$_2$ process, $\cP_1:\R\to\R$ \cite{prahofer2002PNG} obtained as the projection $\cL(0,0,\cdot,1)$, i.e., when the second spatial coordinate is allowed to vary while the remaining space coordinate is fixed to be $0$ and the two temporal coordinates are frozen to be $0$ and $1$ respectively. In words, the  parabolic Airy$_2$ process traces out the last passage time in $\cL$ from the origin to the points $(\cdot,1).$ 

However, by natural translation invariant properties of $\cL,$ we have the following equality in distribution (among many other symmetries, see e.g., Lemma \ref{lem:DL-symmetry}) for any $y\in \R,$
$$\cL(0,0;\cdot,1)\overset{d}{=}\cL(y,0; y+\cdot,1).$$

Thus the landscape $\cL$ provides a natural coupling of infinitely many parabolic Airy$_2$ processes ``rooted'' at different spatial points, making the coupling structure an intriguing object of study. Such a study was initiated in \cite{BGH}, who introduced the following process, which, as we will soon see, encodes valuable information about the coupling.

\noindent
\textbf{Difference profile:} Fix $y_a,y_b\in\R$ with $y_a < y_b$. Consider the random function $\wdf:\HH\to\R$ (where $\HH=\R\times \R_+$ denotes the upper half plane) given by 
\begin{equation}\label{e.W defn}
\wdf(x,t) = \cL(y_a,0,x,t) - \cL(y_b,0, x,t).
\end{equation}
We will call $\wdf$ the \emph{weight difference profile}; it is the difference in the weights of two geodesics with differing but fixed starting points and a  common ending point as the latter varies.

\begin{center}
\textbf{For concreteness, we will fix $y_a=-1$ and $y_b=1$ respectively. }
\end{center}

Much of this paper is devoted towards computing new fractal exponents associated to $\cD$. Since its introduction in \cite{BGH}, several of its properties have already been recorded in the literature. We list some of them below. \\

Owing to certain key properties of geodesics, it turns out that $\wdf$ is almost everywhere locally constant, with probability one ({while this is a consequence of known arguments, this has not been completely spelled out in the literature. One of our main results, Theorem \ref{thm:2dhaus}, will in fact imply this}). Further, the following is known.

\begin{lemma}\label{l.W increasing} Almost surely, for any fixed $t,$
$\wdf(\cdot,t)$ is a continuous non-increasing function.
\end{lemma}

This has been proved several times in the literature, for example \cite[Lemma 9.1]{DOV} or \cite[Theorem 1.1(2)]{BGH}.
Thus, in particular, Lemma~\ref{l.W increasing} implies that, for any fixed $t,$ at any given point $x$, $\wdf(x,t)$ is either constant in a neighborhood of $x$ or decreasing at $x$ on at least one side. 

Moreover, it is an easy argument that the non-constant points of a continuous non-increasing function must form a perfect set, i.e., be closed and have no isolated points. Perfect sets must necessarily be uncountable. In particular, $-\wdf(\cdot,t)$ can be interpreted as the distribution function of a random measure supported on the uncountable set $\hNC_t$ (the notation $\hNC_t$ is chosen to be in correspondence with the 2D and temporal counterparts $\NC$ and $\vNC_x$, to be defined shortly). Canonical examples of a similar nature include the Cantor function or the non-constant points of BLT. Associated to any such set is its fractal dimension, which quantifies how ``sparse'' the set is. With the dimension of the mentioned two examples being classically known, one is led to ask: what is the fractal dimension of $\hNC_t$?

This question was originally raised and answered in \cite{BGH}, where it was shown that the Hausdorff dimension of $\hNC_t$ is $\frac{1}{2}$ almost surely. (The definition of the Hausdorff dimension of a set is recalled ahead in Definition~\ref{d.hd} for the reader's convenience.)
Further, in \cite{BGHhau}, it was shown that $\hNC_t$ is precisely the set of points $x,$ which admit geodesics to the points $(-1,0)$ and $(1,0)$ respectively, which are disjoint except at their initial point $(x,t).$ 
Finally, more recently, in \cite{GH21} (and subsequently \cite{D21last}), comparisons were made between $\wdf(\cdot,t)$ and BLT, showing, among other things, that the former is absolutely continuous with respect to the latter in a certain precise sense. 

\begin{figure}
     \centering
     \begin{subfigure}[t]{0.45\textwidth}
         \centering
\includegraphics[width=0.8\textwidth]{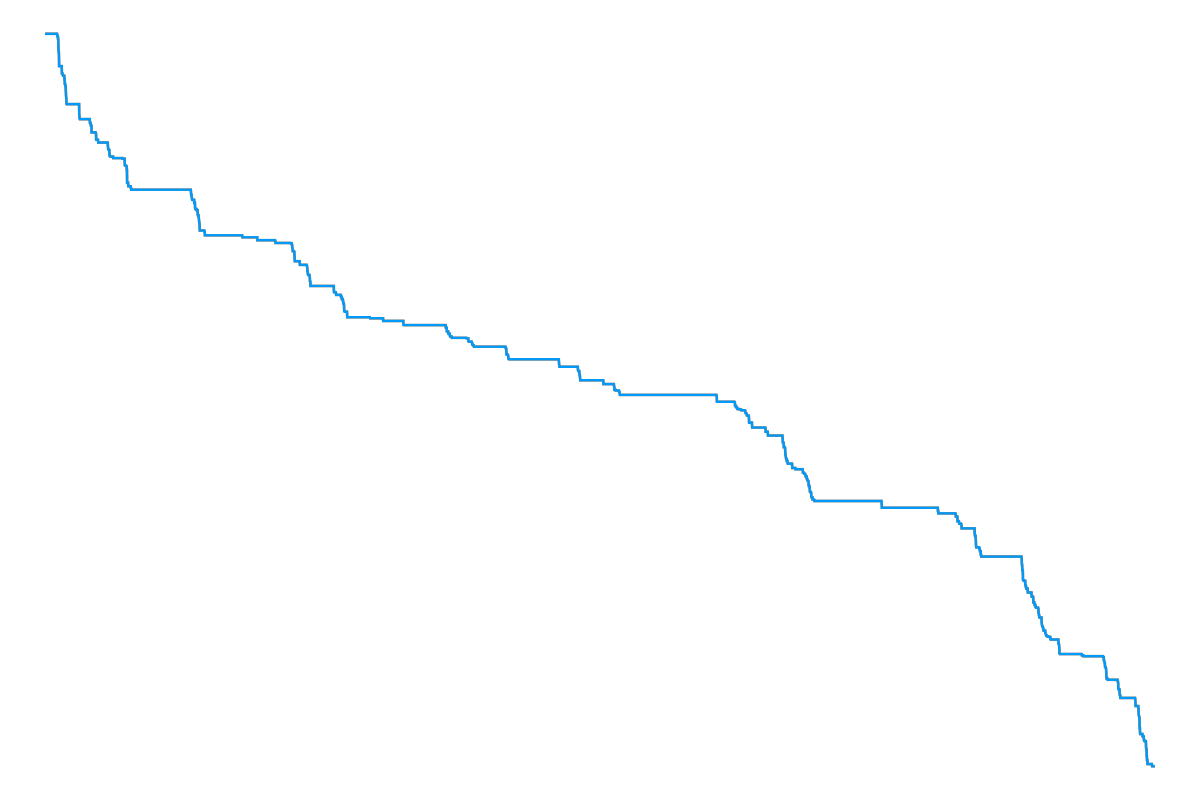}
     \end{subfigure}
     \hfill
     \begin{subfigure}[t]{0.53\textwidth}
     \centering
\includegraphics[width=\textwidth]{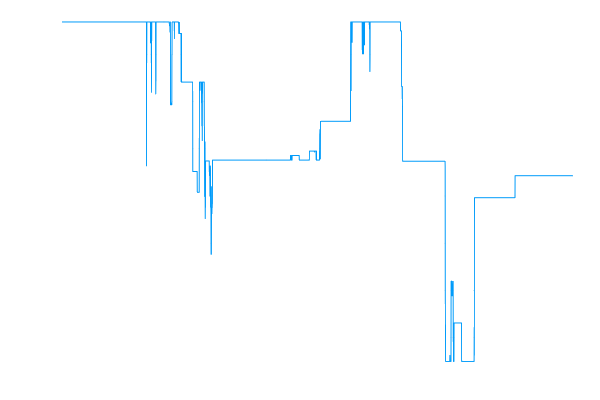}
     \end{subfigure}
        \caption{Simulations of the difference profile $\cD$: the left panel is $\cD(\cdot, 1)$, and the right panel is $\cD(0, \cdot)$. (Courtesy of Milind Hegde.)}
        \label{fig:simu}
\end{figure}

However, in all the above works, the monotonicity of $\wdf(\cdot,t)$ was crucially used with the arguments breaking down if the temporal coordinate was varied instead, in which case the process $\cD(x,\cdot)$ is no longer monotone, necessitating the development of new ideas. (See Figure \ref{fig:simu} for simulations of $\cD$.)\\

\noindent
\textbf{This is the main purpose of this article. Thus, to summarize, we will investigate the spatio-temporal behavior of the process $\wdf$ and compute new fractal dimensions associated to it. As a key ingredient, we introduce and study the notion of GLT, an object we anticipate to be of broader interest and applicability. The connection between the difference profile and GLT is established through a novel level set decomposition which is presented shortly, in Section \ref{s:iop}.} \\

We now proceed to stating our main results. 
\section{Main results}
As indicated already, the  key object of our investigation is the difference profile $\cD(x,t) = \cL(-1,0;x,t)-\cL(1,0;x,t).$

\subsection{Fractal properties of the difference profile}
Our main results pertain to the fractal dimension of the support set of $\cD$ which is almost surely locally constant. 
To this end, we define  $\NC\subset \HH$ to be the set of all $(x,t)$ where the difference profile $\cD$ is not constant in a neighborhood. A one dimensional variant is  $\vNC_x\subset \R_+$ consisting of all $t$ such that the function $t\mapsto \cD(x,t)$ is not constant in a neighborhood of $t$.
We note that $\vNC_x \subset \{t\in \R_+: (x,t)\in \NC\}$. The set $\vNC_x$ is exactly the temporal analogue of the spatial set $\hNC_t$ of all $x$ such that $\cD(x,t)$ is not locally constant, considered in the papers \cite{BGH, BGHhau, GH21, D21last}.
As we will mostly work with the temporal non-constant set at $x=0$, we let $\vNC=\vNC_0$.

With the above preparation, we can now state our two main results. 

\begin{theorem} \label{thm:2dhaus}
Almost surely, the set $\NC$ has Hausdorff dimension $5/3$.
\end{theorem}

\begin{theorem} \label{thm:1dhaus}
Almost surely, the set $\vNC$ has Hausdorff dimension $2/3$.
\end{theorem}

As the reader will notice from the arguments presented in the paper, unlike the process $\cD(\cdot, t)$ which is monotone and decreases from $-\infty$ to $\infty$ as $x$ varies on $\R,$ the temporal process $\cD(0,\cdot)$ eventually becomes a constant owing to geodesic coalescence. Further, the process, almost surely, changes sign infinitely often as $t$ approaches $0$ owing to  symmetry, de-correlation across different values of $t$ as the latter converges to $0$ suitably fast, and ergodicity of the directed landscape.\\

As already indicated in the introduction, the proofs of our main theorems rely on a key ingredient, namely, a novel construction of the geodesic local time (GLT), and analysis of its fractal properties which we include as our final main result. 
Recall that we use $\pi_p^r$ to denote the semi-infinite geodesic started from $p$ in the $r$ direction (the formal treatement will appear in Section \ref{sec:buse-semiinf-dl}), and write $\pi_p=\pi_p^0$.\\

\subsection{The geodesic local time}   \label{ssec:geo-loc}
Throughout this paper we let $\boo$ denote the origin $(0,0)$. Consider the geodesic $\pi_\boo$, and its intersection with the line $\{(0,t):t\ge 0\}$.
Denote by 
\begin{equation}\label{zeroset12}
Z:=\{t\ge 0: \pi_\boo(t)=0\},
\end{equation} the geodesic zero set.
\begin{theorem} \label{thm:geohaus}
Almost surely, the set $Z$ has Hausdorff dimension $1/3$.
\end{theorem}

The proofs of all these three theorems rely on the construction of GLT, a non-decreasing function $L:\R_{\ge 0} \to \R$, which records the time $\pi_{\boo}$ has spent at $0$.
This induces a (random) measure $\mu$ on $\R_{\ge 0}$, supported on the set $Z$ similar to how BLT is supported on the zero set of Brownian motion. An important step in the construction of $L$ is the construction of the geodesic occupation measure.
We define the occupation measure $\nu_h$ for the geodesic $\pi_\boo$ at height $h$, by setting for any measurable $A\subset \R$, 
\begin{equation}\label{occmeasure1}
\nu_h(A) = \int_0^h \don[\pi_\boo(t)\in A] dt.
\end{equation}
Given the above, the construction of $L$ is completed by proving the existence of the density of $\nu_h$ at $0$. 
\begin{prop}  \label{prop:om-limit}
For any height $h>0,$ almost surely the limit $\lim_{w\searrow 0} (2w)^{-1}\nu_h([-w, w])$ exists and is finite.
\end{prop}

Before proceeding further, it is worth pointing out  some other recent advances around the theme of fractal properties of last passage percolation. While the works \cite{JR20, JRS19, SS21a} focussed on pre-limiting models, \cite{CHHM, Dau22N} worked with the evolution induced by the directed landscape started from rather general initial data (an object termed as the KPZ-fixed point, see \cite{MQR}). In another direction, fractal properties often manifests in remarkable behavior of the system pertaining to ``chaos and noise-sensitivity''. Refer  e.g. to the success story of the study of dynamical critical percolation summarized beautifully in \cite{GS14}. In the KPZ universality paradigm, such a study was initiated recently in \cite{GH20Sta}.

\subsection{Key ideas in the proofs}\label{s:iop}
The arguments in this paper involve several new ideas, which we describe in some detail in this section. 

We first setup some notations that will be used throughout this paper.
For any set $A,$ we use $|A|$ to denote its cardinality. We use $\LE$ to denote the Lebesgue measure on $\R$. The notation $\|\cdot\|$ denotes the $l^2$ norm unless otherwise stated. Recall that we let $\boo=(0,0)$. For any $a\le b$, we let $\llbracket a, b\rrbracket = [a, b] \cap \Z$. 
We would need to consider the space of all non-empty compact subsets of $\R^2$ which equipped with the Hausdorff metric, becomes a complete metric space, with the induced topology referred to as the Hausdorff topology.
We will also freely use the notation $\argmax$ even when there isn't a unique maximum (but there must exist at least one) using it to denote an arbitrary one.

We next recall the definitions of the Hausdorff dimension and Hausdorff measure of a set.
\begin{defn}\label{d.hd}
For any $d\ge 0$ and metric space $X$, the $d$-dimensional Hausdorff measure of $X$ is defined as
\[
\lim_{\delta\searrow 0} \inf\left\{\sum_i \diam(U_i)^d : \{U_i\}\text{ is a countable cover of } X \text{ with } 0<\diam(U_i)<\delta\right\}.
\]
The Hausdorff dimension of $X$ is $\inf\{d>0: \text{ the $d$-dimensional Hausdorff measure of $X$ is zero }\}$.
\end{defn}

We start by stating two general strategies to obtain upper and lower bounds for the Hausdorff dimension for $X$ being a subset of $\R^n$ (for some $n\in\N$).
\begin{itemize}
  \item Upper bound: Show that for any $U\subset \R^n$, $\PP[X\cap U \neq \emptyset]<\diam(U)^{n-d}$ (potentially with some extra logarithmic factors).
  Then for any compact $S\subset \R^n$, on dividing $S$ into roughly $\delta^{-n}$ subsets, each with diameter $\le \delta$,  $X$ is covered by the union of roughly $\delta^{-d}$ of them. Thus for any $\epsilon>0$, the $d+\epsilon$-dimensional Hausdorff measure of $X$ is zero. So the Hausdorff dimension of $X$ is at most $d$.
  \item Lower bound: This is usually the harder direction and involves constructing a measure $m$ whose support is contained in $X$, such that $m(X)>0$, and for any measurable $U\subset X$ the H\"older estimate $m(U)<\diam(U)^d$  (potentially with extra sub-polynomial factors) holds. Then for any countable cover $\{U_i\}$ of $X$, one must have $\sum_i \diam(U_i)^d>m(X)$, implying that the Hausdorff dimension of $X$ is at least $d$.
\end{itemize}

Before diving into our arguments, for contrast, let us first briefly sketch the nature of the arguments appearing in the analysis of $\cD(\cdot,t)$ in \cite{BGH} (from here on we will fix $t=1$) and why they fail in the case of $\cD(x,\cdot)$ (we will fix $x=0$).
A key starting observation in \cite{BGH} is the following (also see Figure \ref{fig:00l}): for two points $x_1<x_2,$ $\cD(x,1)$ is constant on $[x_1,x_2]$ if the geodesics $\pi_{(-1,0;x_1,1)}$ $\pi_{(1,0;x_2,1)}$ coalesce.

\begin{figure}
     \centering
     \begin{subfigure}[t]{0.45\textwidth}
         \centering
\begin{tikzpicture}[line cap=round,line join=round,>=triangle 45,x=0.4cm,y=0.3cm]
\clip(-2,-1) rectangle (12,10);

\draw (-10,0) -- (30,0);
\draw [dashed] (-10,8) -- (30,8);
\draw [dashed] [blue] (3,8) -- (6,8);

\draw [line width=1pt] (0,0) -- (2.5,1) -- (1.7,2) -- (4.2,3) -- (2.3,4) -- (4.8,5) -- (3.6,6) -- (4.5,7) -- (3,8);
\draw [line width=1pt] (10,0) -- (8.5,1) -- (8.7,2) -- (5.7,3) -- (5,4) -- (4.8,5) -- (3.6,6) -- (6.5,7) -- (7,8);

\draw [line width=0.5pt,color=blue] (4.5,7) -- (4.8,7.5) -- (4.35,8);
\draw [line width=0.5pt,color=blue] (5.05,6.5) -- (5.1,7) -- (6.2,7.5) -- (6.3,8);

\draw [fill=uuuuuu] (0,0) circle (1.5pt);
\draw [fill=uuuuuu] (10,0) circle (1.5pt);

\draw [fill=uuuuuu] (3,8) circle (1.5pt);
\draw [fill=uuuuuu] (7,8) circle (1.5pt);
\draw [blue] [fill=blue] (4.35,8) circle (1.5pt);
\draw [blue] [fill=blue] (6.3,8) circle (1.5pt);

\begin{scriptsize}
\draw (0,0) node[anchor=north]{$(-1,0)$};
\draw (10,0) node[anchor=north]{$(1,0)$};

\draw (3,8) node[anchor=north east]{$(x_1, 1)$};
\draw (7,8) node[anchor=north west]{$(x_2, 1)$};
\end{scriptsize}

\end{tikzpicture}
 \caption{If the geodesics $\pi_{(-1,0;x_1,1)}$ and $\pi_{(1,0;x_2,1)}$ coalesce, any geodesic from $(-1,0)$ or $(1,0)$ to any $(x, 1)$ with $x\in [x_1, x_2]$ must stay with $\pi_{(-1,0;x_1,1)}$ or $\pi_{(1,0;x_2,1)}$ up to the end of their intersection. This implies that $\cD(x,1)$ is constant on $[x_1,x_2]$.}
\label{fig:00l}
     \end{subfigure}
     \hfill
     \begin{subfigure}[t]{0.5\textwidth}
         \centering
\begin{tikzpicture}[line cap=round,line join=round,>=triangle 45,x=0.4cm,y=0.3cm]
\clip(-2,-1) rectangle (12,11);

\draw (-10,0) -- (30,0);
\draw [dashed] (-10,8) -- (30,8);
\draw [dashed] (5,-10) -- (5,20);

\draw [line width=1pt] (0,0) -- (2.5,1) -- (1.7,2) -- (4.2,3) -- (2.3,4) -- (5.4,5) -- (4.1,6) -- (5.1,7) -- (5,8);
\draw [line width=1pt] (10,0) -- (8.5,1) -- (8.7,2) -- (5.7,3) -- (5,4) -- (5.4,5) -- (4.1,6) -- (6.1,7) -- (6.2,8) -- (4.4,9) -- (5,10);

\draw [fill=uuuuuu] (0,0) circle (1.5pt);
\draw [fill=uuuuuu] (10,0) circle (1.5pt);

\draw [fill=uuuuuu] (5,8) circle (1.5pt);
\draw [fill=uuuuuu] (6.2,8) circle (1.5pt);
\draw [fill=uuuuuu] (5,10) circle (1.5pt);

\begin{scriptsize}
\draw (0,0) node[anchor=north]{$(-1,0)$};
\draw (10,0) node[anchor=north]{$(1,0)$};

\draw (5,8) node[anchor=north east]{$(0, t_1)$};

\draw (6.2,8) node[anchor=north west]{$(z, t_1)$};
\draw (5,10) node[anchor=north west]{$(0, t_2)$};
\end{scriptsize}

\end{tikzpicture}
 \caption{If all geodesics from $(-1,0)$ and $(1,0)$ to $(x,t_1)$ coalescence, for $x$ in an order $(t_2-t_1)^{2/3}$ neighborhood of $0$, the process $\cD(0,\cdot)$ is (likely to be) constant on $[t_1, t_2]$.}
\label{fig:00r}
     \end{subfigure}
        \caption{Illustrations for the arguments in \cite{BGH} (for upper bounding the Hausdorff dimension), and the same strategy in the temporal direction.}
        \label{fig:00}
\end{figure}

Thus to prove an upper bound on the Hausdorff dimension of $\hNC_1,$ \cite{BGH} uses as input the following estimate from \cite{H20Exp}: 
For any $x_1< x_2$, $y_1<y_2$, and $s<t$, if we let $\Dist_{(x_1, x_2)\to (y_1, y_2)}^{[s, t]}$ be the event that $\pi_{(x_1, s;y_1, t)}$ and $\pi_{(x_2, s;y_2, t)}$ are unique and disjoint, then 

\begin{theorem} \label{thm:disj1}
There is a constant $C$, such that for any $\epsilon>0$, we have
\[
\PP[ \Dist_{(-1,1)\to (-\epsilon, \epsilon)}^{[0,1]} ] < C\epsilon^{1/2} e^{C |\log(\epsilon)|^{5/6}}.
\]
\end{theorem}

Thus, for any small $\epsilon,$ ignoring sub-polynomial factors, the probability that $\hNC_1$ intersects any interval of size $\epsilon$ is at most $\epsilon^{1/2}$, which then leads to an upper bound of $1/2$ on the Hausdorff dimension, according to the above general strategy. 

For the lower bound, the argument relies on H\"{o}lder regularity of $\cL.$ Namely that the Parabolic Airy$_2$ process is locally $1/2-$ H\"{o}lder-regular (in fact it is globally parabolic and locally Brownian, see e.g., \cite{CH, Ha16, HHJ}) and so is  $\cD(\cdot,1)$ being the difference of two Airy$_2$ processes rooted at $\pm 1$. 
By interpreting $-\cD(\cdot,1)$ as the distribution function of a measure and the above general strategy, we get a Hausdorff dimension lower bound of $1/2$.

It is instructive to apply now the same strategy to analyze the process $\cD(0,\cdot)$ instead (also see Figure \ref{fig:00r}). For two points $t_1<t_2$, consider the geodesics $\pi_{(-1,0;0,t_1)}$ and $\pi_{(1,0;0,t_2)}.$ Let the latter intersect the line $t=t_1$ at location $z,$ then by a similar argument as before, one can essentially conclude (the precise conclusion is somewhat different which we will not spell out here) that $\cD(0,\cdot)$ is constant on $[t_1,t_2]$ unless the geodesics $\pi_{(-1,0;0,t_1)}$ and $\pi_{(1,0; z,t_1)}$ are disjoint. Now, it is known, and goes back to the work of Johansson \cite{Jo00} (see also \cite{BSS, HS20, DSV}),  that the geodesics in LPP and hence $\cL$ are $2/3-$ H\"{o}lder-regular, which then implies that $z\approx |t_2-t_1|^{2/3}.$ 
This in conjunction with the above estimate on existence of disjoint geodesics imply that the probability that  $\vNC$ intersects $[t_1,t_2]$ is at most, up to sub-polynomial factors, $|t_2-t_1|^{1/3}$ which translates into a Hausdorff dimension upper bound of $2/3$ for $\vNC.$

On the other hand, applying the same reasoning as before, we now seek to use the fact that $\cL$ is $1/3-$ H\"{o}lder-regular in the temporal direction, implying at least as much regularity for $\cD(0,\cdot)$. This then yields, again by the above general strategy, a lower bound of $1/3$ for the Hausdorff dimension which unlike the spatial case, falls short of matching the upper bound. 

The reason for the disparity is that the process $\cD(0,\cdot)$ is a mean zero process and hence not monotone. 
In the spatial case, when applying the above general strategy, the process $\cD(\cdot,1)$ was itself inducing a measure owing to monotonicity, whereas the process $\cD(0, \cdot)$ induces a signed measure instead. Thus the argument which essentially works with the ``absolute value'' of the measure, or the ``un-signed'' measure fails to capture certain cancellations exhibited by $\cD(0, \cdot),$ which we expect leads to a square-root fluctuation effect boosting the lower bound from $1/3$ to $2/3.$

Thus to obtain a sharp lower bound we will seek to construct new measures capturing the appropriate cancellation. 
This involves many new ideas and is quite technical. We will review the key steps in the sequel but first let us digress a bit allowing us to introduce a few new objects which will play central roles in our arguments.\\

\noindent 
\textbf{GLT and the zero set:}  As is the case in Brownian motion, the zero set of the semi-infinite geodesic $\pi_\boo$ (see the beginning of Section \ref{ssec:geo-loc}) denoted by $Z$, will turn out to be a particularly fundamental object.   As indicated in Section \ref{ssec:geo-loc}, we will construct GLT $L$ and show that it is $1/3$ H\"{o}lder-regular, and using that we will show that $Z$ has Hausdorff dimension $1/3$ (as stated in Theorem \ref{thm:geohaus}). This will follow from the following estimate, whose proof will be briefly reviewed at the end of this section:
\begin{lemma}  \label{lem:bound-zb-s}
There are universal constants $C,c>0$ such that the following is true.
For any closed interval $I\subset \R$ and $J\subset \R_{\ge 0}$, we have
\[
\PP[\LE(\{t \in J:\pi_\boo(t)\in I\}) > M\LE(I)\LE(J)^{1/3}] < Ce^{-cM}
\]
for any $M>0$.
\end{lemma}

\noindent
\textbf{Competition interfaces:} The next object that we need to consider pertains to level sets. As is perhaps not surprising from the fact that we are investigating the behavior of the difference profile $\cD$, it would be particularly convenient to consider level sets of the form 

\begin{equation}\label{levelcurve1}
\LV_\ell(t)=\inf\{x\in \R: \cD(x,t) \le \ell\}.
\end{equation}
It turns out however that a remarkable identity holds if we were to take $\ell=0$ and consider the difference profile   not rooted at $\pm1 $, but with Brownian initial data instead. 
Namely, let $\cB$ be a standard two sided Brownian motion. 
For any $x\in\R, t>0$, we let
\begin{align*}
\cL^L(x,t)&=\sup_{y\le 0} \cL(y,0;x,t) + \cB(y),\\
\cL^R(x,t)&=\sup_{y\ge 0} \cL(y,0;x,t) + \cB(y),
\end{align*}
be the passage times induced by the Brownian initial data from the negative and positive parts of the $x-$axis respectively.
One can show that the function $x\mapsto \cL^L(x,t)-\cL^R(x,t)$ is non-increasing in $x$ (see the beginning of Section \ref{dualsec}).
We then let 
\begin{equation} \label{leftrightprofile}
\LV^\cB(t)=\inf\{x\in \R: \cL^L(x,t)-\cL^R(x,t) \le 0\}.
\end{equation}

One then has the remarkable duality result:
\begin{prop} \label{prop:equal-dist-geo-comp}
The processes $\LV^\cB$ and $\pi_{\boo}$ are equal in distribution.
\end{prop}
In the discrete setting (of exponential LPP), this goes back to \cite{FP05, BCS06cube, FMP09}, and we get the above directed landscape version using a limit transition from \cite{DV21}. During the preparation of this paper, two articles \cite{RV21, SS21} were posted on the arXiv where this was independently developed.  While \cite{RV21} worked with the directed landscape as in the present paper, \cite{SS21} considered the prelimiting model of Brownian LPP instead.

Now, how does this help us?  Suppose for the moment that the same were true for the level sets $\LV_\ell(\cdot)$ for all $\ell$, in place of $\LV^\cB(\cdot).$ 
Recalling that we only seek to prove a lower bound of $2/3$ for the Hausdorff dimension of $\vNC,$ we will now construct a measure supported on $\vNC$ which will be $2/3$ H\"{o}lder-regular. At a very high level, this is achieved by taking the average of the local time for each of the level sets. The construction proceeds as in the case of GLT. Now using the above assumed comparison between level sets and the geodesic, it follows that each of the individual level set local times are $1/3$ H\"{o}lder-regular using Lemma \ref{lem:bound-zb-s}. Given this, we obtain the desired $2/3$ H\"{o}lder regularity, using H\"{o}lder regularity of $\cL$. 

To see this, consider an interval $[t_1,t_2].$ Now each $\LV_{\ell}$ induces a $1/3$ H\"{o}lder-regular measure $\kappa_\ell$. To define this we start by defining the corresponding occupation measures analogous to \eqref{occmeasure1}. 
For each $h>0$, we denote by $\lambda_{\ell,h}$ the measure on $\R$, such that
\begin{equation}\label{occmeasure1l}
\lambda_{\ell,h}(A) =\int_0^h \don[\LV_\ell(t) \in A] dt,
\end{equation}
for any measurable $A\subset \R$. Given this, $\kappa_\ell$ on an interval $(g,h)$ is then defined to be the density of $\lambda_{\ell,h}-\lambda_{\ell,g}$ at $0,$ that is, 
\begin{equation}\label{density1}
\KA_{\ell}((g,h))=\lim_{w\searrow 0} (2w)^{-1} (\lambda_{\ell,h}([-w,w]) - \lambda_{\ell,g}([-w,w])).
\end{equation}
(The actual definition proceeds a bit differently owing to certain measure theoretic considerations which we will ignore in this discussion).
Define now the measure $\KA$ by setting 

\begin{equation}\label{onedimmeasure}
\KA(A)=\int \KA_\ell(A) d\ell.
\end{equation}

\begin{figure}
     \centering
     \begin{subfigure}[t]{0.45\textwidth}
         \centering
\begin{tikzpicture}[line cap=round,line join=round,>=triangle 45,x=0.4cm,y=0.32cm]
\clip(-2,-1) rectangle (12,11);

\draw (-10,0) -- (30,0);
\draw [dashed] (5,-10) -- (5,80);

\draw [line width=1pt,color=red] (5,0) -- (4.5,1) -- (5.7,2) -- (4.8,3) -- (6.3,4) -- (6.8,5) -- (5.6,6) -- (4.2,7) -- (4.4,8) -- (3.8,9);
\draw [line width=1pt,color=green] (5,0) -- (4.7,1) -- (6.3,2) -- (5.7,3) -- (6.6,4) -- (7.1,5) -- (6.6,6) -- (4.5,7) -- (6.2,8) -- (5.7,9);

\draw [fill=uuuuuu] (5,0) circle (1.5pt);
\draw [fill=uuuuuu] (5,6) circle (1.5pt);
\draw [fill=uuuuuu] (5,8) circle (1.5pt);

\begin{scriptsize}
\draw (5,0) node[anchor=north]{$(0,0)$};
\draw (5,6) node[anchor=east]{$(0,t_1)$};
\draw (5,8) node[anchor=west]{$(0,t_2)$};

\draw (3.8,9) node[anchor=south]{$\LV_{\ell_1}$};
\draw (5.7,9) node[anchor=south]{$\LV_{\ell_2}$};
\end{scriptsize}

\end{tikzpicture}
 \caption{If both $\LV_{\ell_1}$ and $\LV_{\ell_2}$ intersect the interval $\{0\}\times[t_1,t_2]$, it is likely that $|\ell_2-\ell_1|$ is in the order of at most $|t_2-t_1|^{1/3}$.}
\label{fig:measurebd}
     \end{subfigure}
     \hfill
     \begin{subfigure}[t]{0.5\textwidth}
         \centering
\begin{tikzpicture}[line cap=round,line join=round,>=triangle 45,x=0.6cm,y=.08cm]
\clip(0,-0.8) rectangle (10,58);

\draw [line width=1pt,color=blue] (1.65,25) -- (2,27) -- (2.5,29);
\draw [line width=1pt,color=blue] (2.6,25) -- (2,27) -- (2.5,29);
\draw [line width=1pt,color=blue] (2.5,29) -- (2.4,31) -- (2.7,33) -- (2.6,35);
\draw [line width=1pt,color=blue] (2.6,35) -- (3.3,38) -- (2.8,42) -- (3.1,46) -- (5,50);

\draw [line width=1pt,color=blue] (2.4,31) -- (4.7,33) -- (5,35);
\draw [line width=1pt,color=red] (8.35,25) -- (8.2,27);
\draw [line width=1pt,color=red] (7.4,25) -- (8.2,27) -- (8.6,29) -- (7,31);
\draw [line width=1pt,color=red] (7,31) -- (7.1,33) -- (5,35);

\draw [line width=1pt,color=red] (7,31) -- (8,38) -- (7.3,43) -- (5.3,46) -- (5,50);

\draw (0,25) -- (10,25);
\draw (5,25) -- (5,55);

\draw [fill=uuuuuu] (5,25) circle (1.5pt);
\draw [fill=uuuuuu] (5,50) circle (1.5pt);
\draw [fill=uuuuuu] (5,35) circle (1.5pt);

\draw [line width=0.6pt] (1,5) -- (1.125,3.) -- (1.25,4.5) -- (1.375,4.1) -- (1.5,2) -- (1.75,3) -- (2,19.8) -- (2.25,19.3) -- (2.5,3.2) -- (2.75,3.4) -- (3,4.3) -- (3.25,5.1) -- (3.5,4.6) -- (3.75,4.9) -- (4,3.8) -- (4.25,4.1) -- (4.5,2.9) -- (4.75,3.5) -- (5,5) -- (5.25,4.6) -- (5.5,2.8) -- (5.75,3.1) -- (6,4.2) -- (6.25,3.3) -- (6.5,4.8) -- (6.75,4.6) -- (7,3.1) -- (7.25,4.4) -- (7.5,4) -- (7.75,19.1) -- (8,18.4) -- (8.25,4.1) -- (8.5,5.3) -- (8.75,3.8) -- (9,4.5);

\foreach \i in {1.65, 2.6, 7.4, 8.35}
{
\draw [dashed] (\i,0) -- (\i,27);
}

\begin{tiny}

\draw (5,25) node[anchor=north]{$\boo$};
\draw (0.5,4) node[anchor=east]{$\cB$};
\end{tiny}

\end{tikzpicture}
\caption{An illustration of the event where the Brownian motion $\cB$ has a very spiky behavior with spikes near $\pm 1$, and geodesic coalescence.}
\label{fig:0}
     \end{subfigure}
\caption{}
\end{figure}

This is the measure that we use to lower bound the Hausdorff dimension of $\vNC$ by showing that $\kappa$ is supported in $\vNC$ and further for any interval $[t_1,t_2],$ $\kappa ([t_1,t_2])\approx |t_2-t_1|^{2/3}.$  Very briefly, the latter holds because by definition $\kappa ([t_1,t_2])=\int\kappa_\ell ([t_1,t_2])d\ell.$ Now by the claimed $1/3$ H\"{o}lder regularity of GLT and the assumed comparison between $\kappa_\ell$ and the same, the integrand is at most $|t_2-t_1|^{1/3}.$ Now by the $1/3$ H\"{o}lder regularity of the landscape in the temporal direction (a consequence of $1/3$ being the universal exponent governing weight fluctuations in the KPZ universality class), for $\LV_{\ell_1}$ and $\LV_{\ell_2}$ to both intersect the interval $\{0\}\times[t_1,t_2]$, we must have, again ignoring sub-polynomial factors, $|\ell_2-\ell_1| \le |t_2-t_1|^{1/3}$ (see Figure \ref{fig:measurebd}), which yields the desired estimate $\kappa ([t_1,t_2])=\int\kappa_\ell ([t_1,t_2]) =O(|t_2-t_1|^{2/3})$ since the integrand is non-zero only on a interval of size at most $|t_2-t_1|^{1/3}.$

Note that the above reasoning is contingent on the assumption that the measures $\kappa_\ell$ are $1/3$ H\"{o}lder-regular which was a consequence of the claim that the level curves $\LV_{\ell}$ are ``similar'' to the geodesic $\pi_{\boo}.$ We rely on Proposition \ref{prop:equal-dist-geo-comp} to prove a result of this kind which suffices for our purpose.  At a very high level, the arguments in this part is based on the observation that if the Brownian motion $\cB$ has a very spiky behavior with spikes near $\pm 1$ respectively as in Figure \ref{fig:0} then this resembles the original situation of considering the difference profile rooted at $\pm1$ respectively. In particular, in presence of such spikes, owing to geodesic coalescence, $\LV^\cB$ is the same as $\LV_\ell$ for some random $\ell$, where the law of the latter dominates a multiple of Lebesgue measure owing to the fluctuation of the difference in the height of the spikes at $\pm1$, which can be obtained by resampling suitable parts of $\cB$.

Theorem \ref{thm:2dhaus} is proved in the same way as Theorem \ref{thm:1dhaus} and is in fact somewhat easier technically since instead of working with the ``density'' measures in \eqref{density1}, we work with the occupation measures in \eqref{occmeasure1l} instead and the ``occupation'' counterpart of $\KA$ by taking their integral over different values of $\ell$.  The same set of arguments then concludes the proof. In this case the measure we work with is simply $\zeta$, such that for any $a < b$ and $0\le g<h$, we let
\begin{equation} 
\zeta([a,b]\times (g,h)) = \int_g^h [\cD(a,t)-\cD(b,t)] dt=\int_{\ell}[\lambda_{\ell,h}([a,b])-\lambda_{\ell,g}([a,b])]d\ell.    
\end{equation}

Now when $b-a=h-g\approx \epsilon,$ using the same argument as before one can conclude that $\zeta([a,b]\times (g,h))\approx \epsilon^{5/3}$ leading to the desired lower bound on the dimension. \\

We end the above discussion by pointing out that while the difference profile $\cD$ we consider is rooted at $(\pm 1, 0),$ one may generalize our argument to a rather general class of ``rooting functions'' where for functions $f_1,f_2,$ one considers
\[
(x, t)\mapsto \sup_{y}  [\cL(y,0;x,t) + f_1(y)] - \sup_{y}  [\cL(y,0;x,t) + f_2(y)].
\]

In particular, for the special case $f_1, f_2: \R\to \R\cup\{-\infty\}$, such that (1) they are bounded from above on any compact interval, (2) there exists some $x_*>0$ such that $f_1(x)=-\infty$ for any $x>-x_*$ and $f_2(x)=-\infty$ for any $x<x_*$, (3) $\limsup_{|x|\to\infty} f_1(x)/x^2, f_2(x)/x^2 \le 0$, essentially our arguments can be carried out without any significant changes to prove the counterparts of  Theorems \ref{thm:2dhaus} and \ref{thm:1dhaus} in this setting. However we do not pursue spelling out all the details which will be done in future work where more general $f_1$ and $f_2$ will be treated as well.\\

The remainder of the discussion is devoted to providing a brief overview of the ideas going into the other novel aspects of the paper. This includes the construction of the semi-infinite geodesic $\pi_{\boo}$ as well as the proof of the duality relation with the competition interface \eqref{leftrightprofile}, and finally the construction of GLT and the local times for the level curves \eqref{levelcurve1}. As already mentioned, these objects and similar ideas have recently appeared on the arXiv in the articles \cite{RV21,SS21, BSS22} during the preparation of this paper. \\

\noindent
\textbf{Construction of $\pi_{\boo}$}: The construction relies on geodesic coalescence. To construct a semi-infinite geodesic we show that the finite geodesics $\pi_{(p;z_n,h_n)},$ for any sequences $\{z_n\}_{n=1}^\infty$ and $\{h_n\}_{n=1}^\infty$ satisfying that $z_n\to \infty$ as $n\to\infty$ and $|z_n-h_nr|<h_n^{2/3}$ for each $n\in \N$,  owing to geodesic coalescence, all share the same initial segment. Further, the length of the initial segment itself goes off to infinity as $n\to \infty.$ $\pi_{\boo}$ is then defined to be the path made up of such initial segments. 

While the construction is intuitive, it does involve some amount of work to show that it works. Namely it has the desired direction. We further show that there is a unique such semi-infinite geodesic.  

We will not comment much on how the above is carried out except that we rely on a construction of the Busemann function in the vertical direction (again relying on geodesic coalescence).
 Busemann functions were originally used to study the large-scale geometry of geodesics in Riemannian manifolds, and was introduced to study general first passage percolation (FPP) models by Hoffman \cite{Ho05, Ho08}. It has also been intensively studied and used in the exponential LPP setting; see e.g. \cite{Se17}. 
We then show that in our setting the Busemann function is a Brownian motion and that $\pi_{\boo}(t)$, is the argmax of the sum of the Busemann function and an independent copy of an appropriately scale Parabolic Airy$_2$ process. By the known Brownianity of both the functions and their independence, it follows that $\pi_{\boo}(t)$ is unique, at least for rational $t$ and hence for all $t$ by continuity. Further, it follows by estimates of Brownian motion as well as the Parabolic Airy$_2$ process, that $\pi_{\boo}(t)\approx O(t^{2/3})$ and thereby implying that $\pi_{\boo}$ indeed has the desired direction.\\

\noindent
\textbf{Proof of duality:} The duality is known to hold in the discrete model of Exponential LPP. To prove the directed landscape version, we rely on the recent convergence results proved in \cite{DV21}, which in particular implies that both the semi-infinite geodesic and the competition interfaces in the discrete setting converge to their counterparts $\pi_{\boo}$ and  $\LV^\cB$ which is enough to conclude the proof. To be completely precise, however, the results in \cite{DV21} only can be used to deduce convergence of certain ``compactified'' versions of the above objects. To make use of this, we work with certain truncations of $\pi_{\boo}$ and  $\LV^\cB$, and their discrete versions, show that with high probability the truncated versions agree with the untruncated ones, the discrete ones converge to the continuous counterparts and hence the discrete duality transfers to the desired continuous duality statement.\\

\noindent
\textbf{Construction of GLT:} Akin to the construction of BLT, this proceeds by first constructing the occupation measure for $\pi_{\boo}$ as defined in \eqref{occmeasure1} at height $h$, by setting for any measurable $A\subset \R$, 
\begin{equation}
\nu_h(A) = \int_0^h \don[\pi_\boo(t)\in A] dt.
\end{equation}
All that remains is to show that $\nu_h$ has a density at $0,$ (see the Proposition \ref{prop:om-limit}), i.e., 
almost surely, the limit $\lim_{w\searrow 0} (2w)^{-1}\nu_h([-w, w])$ exists and is finite. While the finiteness is a consequence of the estimate in Lemma \ref{lem:bound-zb-s}, we briefly explain the idea behind the proof of the existence of a density. The starting point is the observation that, by Lebesgue's theorem for the differentiability of monotone functions, for any measure on $\R$ that is locally finite, almost every point is a point of density. However, a priori, this set of density points (which is a random set) could avoid $0$. We now want to show that almost surely that is not the case. Now this would be immediate if the occupation measure  $\nu_h$ was translation invariant in law, i.e., for any real number $x,$  $\nu_h(\cdot)\overset{d}{=}\nu_h(\cdot+x),$ since any translation invariant set must contain a given point, in this case $0$, with probability one. The proof now proceeds by constructing $\nu'_h,$ a translation invariant ``proxy'' for $\nu_h,$ namely considering the occupation measure of a collection of semi-infinite geodesics started from a Poisson process of intensity one. The latter measure does admit a density at $0$ almost surely given a sample of the Poisson process and the landscape $\cL.$ Now owing to geodesic coalescence and that there is a positive probability that the Poisson process consists of a point quite close to the origin and no other points in the vicinity, it follows that the occupation measure of the two measures agree, (at least away from $0$), i.e., 
$$\nu'_h-\nu'_g=\nu_h-\nu_g$$ for any $0<g<h.$ This is enough to finish the proof once we show that $\nu_g$ is negligible if $g$ is small enough.\\

\noindent
We finish with a brief discussion on the proof of Lemma \ref{lem:bound-zb-s}.
A slightly more general statement in the discrete setting of Exponential LPP had already appeared previously in \cite{SSZ}. For completeness we include the arguments adapted to the setting of the Directed Landscape and describe some of the ideas here. 

To begin, a simple argument shows that $\LE(I)\LE(J)^{1/3}$ is the correct order of the random variable of interest, namely $\LE(\{t \in J:\pi_\boo(t)\in I\})$. This is because, by KPZ transversal fluctuation considerations,  for any height $t,$ $\pi_{\boo}(t)$ is distributed roughly uniformly on an interval $t^{2/3}$ and hence, if $J=[a,b],$ 

\[
\E[\LE(\{t \in J:\pi_\boo(t)\in I\})]\approx \int_{a}^{b}\frac{\LE(I)}{t^{2/3}} \lesssim  \LE(I)\LE(J)^{1/3}
\]

The proof of the lemma now relies on the above expectation bound for ``all'' $I$ and $J$  and dominating $\LE(\{t \in J:\pi_\boo(t)\in I\}$ by a sum of independent random variables and appealing to standard concentration of measure results. However the geodesic does not admit a decomposition into independent segments. To get around this, one uses the observation that a similar expectation bound as above holds when one considers the maximum occupation time, say  among all geodesics with endpoints outside  $I\times J.$ Note that such a class of geodesics in particular includes  segments of the geodesic $\pi_{\boo}$ crossing such a box but crucially such a quantity across different boxes is independent. 

Thus the occupation time in a box can be dominated by an independent sum of occupation times in smaller boxes with the same width, which by scale invariance, corresponds to the original box with a larger width.  At this point we use induction, assuming the desired exponential tail for larger widths (which is straightforward to establish for a large enough width) and standard concentration bounds to conclude the same for a smaller width, going down in dyadic scales.

\vspace{15pt}

\noindent \begin{large}\textbf{Organization of the remaining text.}\end{large}
The rest of this paper is organized in the following manner. 
In Section \ref{sec:prelim} we set up formally the directed landscape and the related pre-limit model we will be analyzing, and recall some basic results from the literature.
In Section \ref{sec:buse-semiinf-dl} we construct the Busemann function and semi-infinite geodesics in the directed landscape.
The duality between semi-infinite geodesics and competition interface is given in Section \ref{dualsec}.
The last three sections contain proofs of the main results concerning Hausdorff dimensions: Section \ref{sec:2ddphau} is for the 2D non-constant set (Theorem \ref{thm:2dhaus}), Section \ref{sec:geo-loc} focusses on GLT (Theorem \ref{thm:geohaus}) while Section \ref{sec:1D-diff-prof} is for the non-constant set in the temporal direction (Theorem \ref{thm:1dhaus}).
The appendices contain proofs of several useful intermediate results which could be of independent interest. They are either technical or resemble arguments appearing earlier in the literature.\\

\noindent \begin{large}\textbf{Acknowledgements}\end{large}
SG thanks Alan Hammond and Milind Hegde for useful discussions on the topic of fractal geometry of the directed landscape.  He is partially supported by NSF grant DMS-1855688, NSF Career grant DMS-1945172, and a Sloan Fellowship. The authors also thank Milind Hegde for helping them with the simulations of the difference profile in Figures \ref{f.simulations} and \ref{fig:simu}, and {Evan Sorensen for helpful comments}.

\section{Preliminaries}  \label{sec:prelim}

As indicated in the previous section, our arguments will rely on various existing definitions and results in the last passage percolation literature which we collect in this section. This includes various results about the pre-limiting model of exponential LPP on the lattice (Section \ref{ss:elpp}, a significant part of it will be devoted to reviewing the duality between competition interfaces and infinite geodesics as well as the recent results proving its convergence to the directed landscape), basic properties of the directed landscape (Section \ref{ss:dl}) and the Airy line ensemble (Section \ref{ss:ale}).  The reader can skip this section on first read and come back to it later as needed.

\subsection{Exponential LPP}  \label{ss:elpp}
We consider directed last passage percolation (LPP) on $\Z^2$ with i.i.d.\ exponential weights on the vertices, i.e., we have a random field $\{\omega_p: p\in \Z^2\}$ where $\omega_p$ are i.i.d. $\Exp(1)$ random variables. For any up/right path $\gamma$ from $p$ to $q$ where $p\preceq q$ (i.e., $p$ is co-ordinate wise smaller or equal than $q$)  the weight of $\gamma$, denoted $T(\gamma)$ is defined by 
\[T(\gamma):=\sum_{w\in \gamma} \omega_{w}.\]
For any two points $p$ and $q$ with $p\preceq q$, we shall denote by $T_{p,q}$ the last passage time from $p$ to $q$; i.e., the maximum weight among weights of all directed paths from $p$ to $q$.
The (almost surely unique) directed path with the maximum weight is the geodesic from $p$ to $q$, denoted as $\Gamma_{p,q}$.
When $p\not\preceq q$, we let $T_{p,q}=-\infty$.

In the setting of exponential LPP, the existences and uniqueness of semi-infinite geodesics have been established (see \cite{Cou11, FP05}), and we summarize the known facts here.
For any $p\in\Z^2$, there is an almost surely unique upper-right path from it, denoted as $$\Gamma_p=\{p=(x_0, y_0), (x_1,y_1), \cdots \},$$ such that for any $0\le i \le j$, the path $(x_i,y_i), \cdots, (x_j,y_j)$ is $\Gamma_{(x_i,y_i), (x_j,y_j)}$, and $\lim_{i\to\infty} \frac{x_i}{y_i}=1$.
The path $\Gamma_p$ is also referred to as the $(1,1)$-directional semi-infinite geodesic from $p$.
Almost surely, for any sequence $q_n\in\Z^2$ that goes to infinity in the $(1,1)$ direction, $\Gamma_{p,q_n}$ converges to $\Gamma_p$, in the sense that in any compact set $\Gamma_{p,q_n}$ is the same as $\Gamma_p$ when $n$ is large enough.
The collection of all the semi-infinite geodesics (in the $(1,1)$-direction) exhibits a tree structure. In particular, for any $p, q\in \Z^2$, the paths $\Gamma_p$ and $\Gamma_q$ are the same except for finitely many initial points.

Below we shall always assume (the probability $1$ event) that for any $p\preceq q\in \Z^2$, there is a unique geodesic $\Gamma_{p,q}$, and a unique $(1,1)$-direction semi-infinite geodesic $\Gamma_p$.

For the rest of this paper, we let $d, ad:\Z^2\to \Z$ be the function where $d(x,y)=x+y$ and $ad(x,y)=x-y$. For each $n\in\Z$ the anti-diagonal line $\LL_n=\{p\in\Z^2: d(p)=2n\}$ will be important.

\subsubsection{Convergence to the directed landscape}   \label{ssec:lpptodl}
In \cite{DV21}, it is proved that the exponential LPP converges to the directed landscape, in a sense that we record precisely next.
Define a metric $\cK_n$ on $\R^2$ as follows.
For any $(x,s;y,t) \in \R^4$, we let
\[
\cK_n(x,s;y,t) = 2^{-4/3}n^{-1/3} F(ns+2^{5/3}n^{2/3}x, ns; nt+2^{5/3}n^{2/3}y, nt),
\]
where $F(p;q)= T_{p, q} - 2d(q-p) - \omega_p$ for $p,q \in \Z^2$, 
and $F$ is extended to be a function on $\R^4$ using the following rounding from \cite{DV21}.
Let $\fr:\R^2 \to \Z^2$ be the function where
\[
\fr(x,y) =
\begin{cases}
(x, y)  &  x, y \in \Z;\\
(\lfloor x \rfloor, y) &  x\not\in \Z, y \in \Z;\\
(x, \lfloor y \rfloor) &  x\in \Z, y \not\in \Z;\\
(\lceil x\rceil, \lfloor y \rfloor) &  x\not\in \Z, y \not\in \Z.
\end{cases}
\]
We then let $F(x,y;z,w)= T_{(\lceil x\rceil, \lceil y\rceil), \fr(z,w)} - 2d(\fr(z,w)-\fr(x,y)) - \omega_{(x,y)}\don[x, y\in \Z]$.

Given this, we define the geodesic sets in $\cK_n$ as follows.
For any $p, q\in \R^2$, a set $A\subset \R^2$ is called a \emph{geodesic set} in $\cK_n$ from $p$ to $q$, if there is a total ordering $\preceq_A$ on $A$, such that $p\preceq_A p' \preceq_A q$ for any $p'\in A$; 
and for any $p_1\preceq_A p_2 \preceq_A p_3 \in A$, there is $\cK_n(p_1;p_2) + \cK_n(p_2;p_3) = \cK_n(p_1;p_3)$.
Such a set $A$ is called a \emph{maximal geodesic set} in $\cK_n$ if it is not contained in any other geodesic set from $p$ to $q$.

\begin{defn}\label{rangegraph}
For any $p=(x,y)\in\R^2$, let $\fG_n (p) = ( 2^{-5/3}n^{-2/3}(x-y), n^{-1}y)$ and
for any set $A\subset \R^2$ and $n\in \N$, let $\fG_n(A) = \{( \fG_n(p) : p \in A\}$. Finally, for any continuous path $\pi:I\to\R$, where $I\subset \R$ is any subset,  define its \emph{graph} $\Graph(\pi)$ as the set $\{(\pi(t),t): t\in I\}$.
\end{defn}
Given the above, it is straightforward to check that for any $p, q\in\Z^2$, the set $\fG_n(\Gamma_{p,q})$ is a maximal geodesic set from $\fG_n(p)$ to $\fG_n(q)$.
We end with the following recent convergence result.

\begin{theorem}[\protect{\cite[Theorem 13.8(2)]{DV21}}]  \label{thm:exp-to-dl}
There is a coupling of $\cL$ with $\cK_n$ for all $n\in\N$, such that the following holds almost surely.
First, for any compact $K\subset \R^4_\uparrow$, we have $\cK_n \to \cL$ uniformly in $K$.
Second, for any $(p_n; q_n) \to (p; q) \in \R^4_\uparrow$,
and any $\cK_n$ maximal geodesic set from $p_n$ to $q_n$, denoted as $\pi_n$, they are pre-compact in Hausdorff topology, and any subsequential limit is $\Graph(\pi)$, for some $\cL$ geodesic $\pi$ from $p$ to $q$.
\end{theorem}

\subsubsection{Transversal fluctuation}

Estimates on transversal fluctuation of geodesics have appeared and been used in several works recently (see e.g. \cite{BGZ,MSZ,BSS,BSS19,GH20,BGHH20}. In this paper we would use the following estimate for semi-infinite geodesics.
\begin{lemma}[\protect{\cite[Corollary 2.11]{MSZ}}] \label{lem:trans-dis-inf}
There exist constants $C,c>0$ such that the following is true.
Take any $n\in\N$ large enough and any $x>0$, and denote $m_{n,x}=xn^{2/3}$.
Consider the rectangle whose four vertices are
$(-m_{n,x}, m_{n,x})$, $(m_{n,x}, -m_{n,x})$, and $(n-m_{n,x}, n+m_{n,x})$, $(n+m_{n,x}, n-m_{n,x})$.
Then with probability $1-Ce^{-cx^3}$, the part of the geodesic $\Gamma_\boo$ below $\LL_n$ is contained in that rectangle.
\end{lemma}
We also use the following verison in this paper, where the rectangle is replaced by a parallelogram.
\begin{lemma} \label{lem:trans-semi-inf}
There exist constants $C,c>0$ such that the following is true.
Take any $n\in\N$ large enough and any $x>0$, and denote $m_{n,x}=xn^{2/3}$.
Consider the parallelogram whose four vertices are
$(-m_{n,x}, 0)$, $(m_{n,x}, 0)$, and $(n-m_{n,x}, n)$, $(n+m_{n,x}, n)$.
Then with probability $1-Ce^{-c(x^3\wedge n^{2/3}x)}$, we have $\Gamma_\boo \cap (\Z\times \llbracket 0,n\rrbracket)$ is contained in this parallelogram.
\end{lemma}
\begin{proof}
We consider the rectangle whose four vertices are
$(-m_{n,x}/2, m_{n,x}/2)$ and $(m_{n,x}/2, -m_{n,x}/2)$, and $(n+m_{n,x}, n)$ and $(n+m_{n,x}, n)$.
If the part of the geodesic $\Gamma_\boo$ below $\LL_n$ is contained in this rectangle, then $\Gamma_\boo \cap (\Z\times \llbracket 0,n\rrbracket)$ is contained in the parallelogram. 
Thus by Lemma \ref{lem:trans-dis-inf}, the probability is at least $1-C'e^{-c'm_{n,x}^3(n+m_{n,x})^{-2}}$ for some constants $C,c>0$; and this implies the conclusion.
\end{proof}

\subsubsection{Duality}  \label{ssec:dualityexp}
We next state the already alluded to duality between semi-infinite geodesics and competition interfaces in exponential LPP.
Towards this we begin by setting up the appropriate Busemann function (recall that we had indicated how this will serve as a valuable device in our constructions in Section \ref{s:iop}) in exponential LPP.
\begin{defn}\label{expbusemann}
By the already stated tree structure of semi-infinite geodesics, for each $p\in\Z^2$ we let $G(p)= T_{p,\bc} - T_{\boo,\bc}$, where $\bc\in\Z^2$ is the coalescing point of $\Gamma_p$ and $\Gamma_\boo$;
i.e. $\bc$ is the vertex in $\Gamma_p \cap \Gamma_\boo$ with the smallest $d(\bc)$.
\end{defn}
Then $G$ satisfies the following properties.
\begin{enumerate}
\item For any $p=(x,y)\in\Z^2$ and $n>y$, if we take $m_*=\argmax_{m\ge x} T_{p,(m,n)}+G(m,n+1)$, then $\Gamma_p = \Gamma_{p,(m_*,n)}\cup \Gamma_{(m_*,n+1)}$, and $G(p)=T_{p,(m_*,n)}+G((m_*,n+1))$.
\item For any down-right path $\cU=\{u_k\}_{k\in\Z}$, the random variables $G(u_k)-G(u_{k-1})$ are independent. The law of $G(u_k)-G(u_{k-1})$ is $\Exp(1/2)$ if $u_k=u_{k-1}-(0,1)$, and is $-\Exp(1/2)$ if $u_k=u_{k-1}+(1,0)$. They are also independent of $\omega(p)$ for all $p\in \cU_-$, where $\cU_-$ is the lower part of $\Z^2\setminus \cU$,
\end{enumerate}
The first property is a consequence of the definition, and a proof of the second property can be found in \cite{Sep20}.

Now, for each $p\in\Z^2$, let $\Gamma_p^\vee$ be the semi-infinite geodesic in the $(-1,-1)$ direction. {We shall always assume the almost sure event where $\Gamma_p^\vee$ exists and is unique for every $p\in\Z^2$.}

We let $G^\vee$ be the $(-1,-1)$ direction Busemann function; i.e. for any $p\in\Z^2$ we let $G^\vee(p)= T_{\bc, p} - T_{\bc, \boo}$, where $\bc\in\Z^2$ is the coalescing point of $\Gamma_p^\vee$ and $\Gamma_\boo^\vee$.
By symmetry, $G^\vee$ has the same distribution and hence similar properties as $G$.

We next split $\Z_{\ge 0}^2\setminus\{\boo\}$ into two parts: let $P_1$ consist of $p\in \Z_{\ge 0}^2\setminus\{\boo\}$ such that $\Gamma_p^\vee \cap \{(x,0):x\in\N\}\neq\emptyset$, and $P_2$ consist of $p\in \Z_{\ge 0}^2\setminus\{\boo\}$ such that $\Gamma_p^\vee \cap \{(0,x):x\in\N\}\neq\emptyset$.
Equivalently, let $\N\times \{0\} \subset P_1$, and for any $p\in \Z_{\ge 0}\times \N$, if we let $m_*=\argmax_{m\in\Z} G^\vee(m,0)+T_{(m,1),p}$, then $p\in P_1$ if and only if $m_*>0$.

One can think of the boundary of $P_1$ and $P_2$ as a path in the dual lattice $\Z^2+(1/2, 1/2)$,
Let $\Delta \subset \Z^2$ be the boundary shifted by $(-1/2, -1/2)$; or equivalently, let $\Delta$ consist of all $p\in\Z_{\ge 0}^2$, such that $p+(1,0)\in P_1$, and $p+(0,1)\in P_2$.
See the left panel in Figure \ref{fig:49} below for an illustration.

As developed in \cite{FP05, FMP09} and also \cite{BCS06cube}, the following  remarkable duality between competition interfaces and semi-infinite geodesics holds.
\begin{lemma}[\protect{\cite[Lemma 4.3]{FMP09}}]  \label{lem:dual-prelim}
$\Delta$ and $\Gamma_\boo$ have the same distribution.
\end{lemma}

While we have already introduced and defined the Directed Landscape, in this section we record many of its properties that we will be using throughout the article. 

\subsection{Directed landscape}  \label{ss:dl}

\subsubsection{Basic properties}
There are several symmetries of the directed landscape.
\begin{lemma}[\protect{\cite[Lemma 10.2]{DOV}, \cite[Proposition 1.23]{DV21}}]  \label{lem:DL-symmetry}
As continuous functions on $\R^4_\uparrow$, $\cL$ is equal in distribution as
\[
(x,s;y,t) \mapsto \cL(-y,-t;-x,-s),\quad (x,s;y,t) \mapsto \cL(-x,s,-y,t);
\]
and
\[
(x,s;y,t) \mapsto \cL(x+cs+z,s+r;y+ct+z,t+r) + (t-s)^{-1}((x-y-c(t-s))^2 - (x-y)^2)
\]
for any $c,z,r\in \R$; and
\[
(x,s;y,t) \mapsto w\cL(w^{-2}x,w^{-3}s;w^{-2}y,w^{-3}t);
\]
for any $w>0$.
We shall refer to these invariances as flip invariance, skew-shift invariance (or translation invariance when $c=0$), and scaling invariance, respectively.
\end{lemma}

A basic property of the directed landscape is the following quadrangle inequality, which is equivalent to Lemma \ref{l.W increasing}. This and versions of it in different settings has been proved and widely used in many previous studies \cite{DOV, DZ, DV21, BGH}.
\begin{lemma}[\protect{\cite[Lemma 9.1]{DOV}}]  \label{lem:DL-quad}
For any $s<t$, $x_1<x_2$, $y_1<y_2$, we have
\[
\cL(x_1,s;y_1,t) + \cL(x_2,s;y_2,t) \ge \cL(x_1,s;y_2,t) + \cL(x_2,s;y_1,t).
\]
\end{lemma}

We next state the following basic regularity estimates on $\cL$. Recall (from Section \ref{s:iop}) that $\|\cdot\|$ denotes the usual $l^2$ norm.

\begin{lemma}[\protect{\cite[Corollary 10.7]{DOV}}]  \label{lem:DLbound}
There is a random number $R$ such that the following is true.
First, for any $M>0$ we have $\PP[R>M] < Ce^{-cM^{3/2}}$ for some universal constants $c,C>0$.
Second, for any $u=(x,s;y,t) \in \R_\uparrow^4$, we have
\[
\left| \cL(x,s;y,t) + \frac{(x-y)^2}{t-s}\right| \le R(t-s)^{1/3} \log^{4/3}(2(\|u\|+2)^{3/2}/(t-s)) \log^{2/3}(\|u\|+2). 
\]
\end{lemma}
We remark that the logarithmic factors come from the exponential tails of $\cL(u)$ for any fixed $u \in \R_\uparrow^4$, continuity estimates of $\cL$, and a union bound and one cannot hope to obtain a uniform estimate without such factors due to the translation invariance and space-time decorrelation of $\cL$.
We will also use the following continuity estimate of $\cL$.

\begin{lemma}[\protect{\cite[Proposition 1.6]{DOV}}]  \label{lem:modcont}
Let $K$ be a compact subset of $\R_\uparrow^4$.
There is a random number $R$ depending only on $K$, such that the following is true.
First, for any $M>0$ we have $\PP[R>M] < Ce^{-cM^{3/2}}$ for some universal constants $c,C>0$.
Second, for any $(x,s;y,t), (x',s';y',t') \in K$, we have
\begin{multline*}
\left| \left(\cL(x,s;y,t) + \frac{(x-y)^2}{t-s} \right)- \left(\cL(x',s';y',t') + \frac{(x'-y')^2}{t'-s'} \right)\right|
\\ 
\le R\left(\tau^{1/3}\log^{2/3}(1+\tau^{-1}) + \xi^{1/2}\log^{1/2}(1+\xi^{-1}) \right),
\end{multline*}
where $\tau = |s-s'| \vee |t-t'|$ and $\tau = |x-x'| \vee |y-y'|$.
\end{lemma}

Like in the exponential LPP setting, (using these regularity estimates above) in the directed landscape we get the following estimates on transversal fluctuation of geodesics, which is uniform in the endpoints.
\begin{lemma} \label{lem:dl-trans-f}
There is a random number $R$ such that the following is true.
First, for any $M>0$ we have $\PP[R>M]< Ce^{-cM^{9/4}\log(M)^{-4}}$ for some constants $c,C>0$.
Second, for any $u=(x,s;y,t) \in \R_\uparrow^4$, any geodesic $\pi_u$, and $(s+t)/2\le r<t$, we have
\begin{equation}  \label{eq:dl-trans-bd-0}
\left|\pi_u(r) - \frac{x(t-r)+y(r-s)}{t-s}\right| < R(t-r)^{2/3} \log^3(1+\|u\|/(t-r)) .
\end{equation}
Similar bound holds when $s<r<(s+t)/2$ by symmetry.
\end{lemma}
Estimates similar to this have appeared in the literature, see e.g. \cite[Proposition 12.3]{DOV}, and \cite[Theorem 1.4]{GH20} for a uniform in endpoints estimate (in the Brownian LPP setting).
The proof of Lemma \ref{lem:dl-trans-f} does not contain particularly new ideas, and relies on computations using Lemma \ref{lem:DLbound} and Lemma \ref{lem:modcont}.
For completeness, we give the details in Appendix \ref{sec:apptrans}.

The next subsection records several estimates concerning the geometry of geodesics in $\cL.$

\subsubsection{Intersections of geodesics}

We begin with an estimate on the number of intersections that geodesics can have with a line.

\begin{lemma}  \label{lem:number-of-inter}
For any $\alpha \in (0,1),$ there are constants $c,C>0$ such that the following is true.
For any $M>\ell^3$, we have $\PP[|\{\pi_{(x,0;y,1)}(\alpha): |x|, |y| \le \ell, \pi_{(x,0;y,1)} \text{is unique} \}| > M] < Ce^{-cM^{1/384}}$.
\end{lemma}

This can be deduced from a similar result in the exponential LPP setting, proved in \cite{BHS}. 
The condition ``$\pi_{(x,0;y,1)}$ is unique'' is assumed due to technical reasons, and is used in passing the geodesics in exponential LPP to the limit.
In fact, while adapting the arguments of \cite{BHS} to the directed landscape setting, we can prove a stronger version of Lemma \ref{lem:number-of-inter} by removing this extra condition, we choose to not pursue this for the  simplicity of the proof. The current formulation would suffice for applications in this paper.
\begin{prop}[\protect{\cite[Proposition 3.10]{BHS}}] \label{prop:numnoncoalgeo}
Let $A_n$ and $B_n$ be line segment along the lines $\LL_0$ and $\LL_n$, with midpoints $(xn^{2/3}, -xn^{2/3})$ and $(n,n)$, respectively, and each has length $2n^{2/3}$.
Let $F_{n,m}=\left(\bigcup_{u\in A_n, v\in B_n} \Gamma_{u,v}\right) \cap \LL_m$, for $0<m<n$.
For any $\psi \in (0,1)$, $\alpha\in (0, 1/2)$,  there is a constant $c>0$, such that when $|x|<\psi n^{1/3}$ and $\alpha n < m < (1-\alpha)n$, we have
\[
\PP[|F_{n,m}| > M] < e^{-cM^{1/128}}
\]
for any sufficiently large $M<n^{0.01}$ and sufficiently large $n$.
\end{prop}
We remark that while the statement of \cite[Proposition 3.10]{BHS} is for $\alpha = 1/3$, the proof therein goes through essentially verbatim for a general $\alpha \in (0,1)$.

Using this and the limit transition from Theorem \ref{thm:exp-to-dl}, we can prove Lemma \ref{lem:number-of-inter}.
\begin{proof}[Proof of Lemma \ref{lem:number-of-inter}]
Below we let $C, c>0$ be constants depending on $\alpha$, whose values may change from line to line.

Let $A_{n,\ell}$ and $B_{n,\ell}$ be the line segments along the lines $\LL_0$ and $\LL_n$,  with midpoints $\boo$ and $(n,n)$, respectively, and length $2^{5/3}\ell n^{2/3}$.
Let $F_{n,m,\ell}=\left(\bigcup_{u\in A_{n,\ell}, v\in B_{n,\ell}} \Gamma_{u,v}\right) \cap \LL_m$, for $0<m<n$.
Then by splitting $A_{n,\ell}$ and $B_{n,\ell}$ into segments of length $2n^{2/3}$, and applying Proposition \ref{prop:numnoncoalgeo} to each pair, we get the following result. For any $M>0$, and any $n$ large enough (depending on $\ell$ and $M$), we have
\[
\PP[|F_{n,\lfloor \alpha n\rfloor,\ell}| > \ell^2 M] < C\ell^2 e^{-cM^{1/128}}.
\]
Thus, when $M>\ell^3$, we have
\[
\PP[|F_{n,\lfloor \alpha n\rfloor,\ell}| > M] < C e^{-cM^{1/384}}.
\]
Finally, consider the coupling from Theorem \ref{thm:exp-to-dl}.
Suppose that in the directed landscape, we have $|\{\pi_{(x,0;y,1)}(\alpha): |x|, |y| \le \ell, \pi_{(x,0;y,1)} \text{is unique} \}| > M$ for some $M>\ell^3$.
Then for any $n$ large enough, one would have $|F_{n,\lfloor \alpha n\rfloor,\ell}| > M$, since any unique geodesic $\pi_{(x,0;y,1)},$ is the Hausdorff limit of a sequence of pre-limiting geodesics. Thus the conclusion follows.
\end{proof}

\subsubsection{Ordering of geodesics}

We recall the notion of leftmost geodesics from \cite{DOV}: for any $(x,s;y,t)\in\R^4_\uparrow$, a geodesic $\pi$ from $(x,s)$ to $(y, t)$ is called the leftmost geodesic, if for any other geodesic  $\tau$ from $(x,s)$ to $(y, t)$, there is $\tau(t')\ge \pi(t')$ for any $t'\in [s, t]$.
By \cite[Lemma 13.2]{DOV}, almost surely, for any $(x,s;y,t)\in\R^4_\uparrow$ there is a leftmost geodesic from $(x,s)$ to $(y, t)$.
We also have the following ordering property.

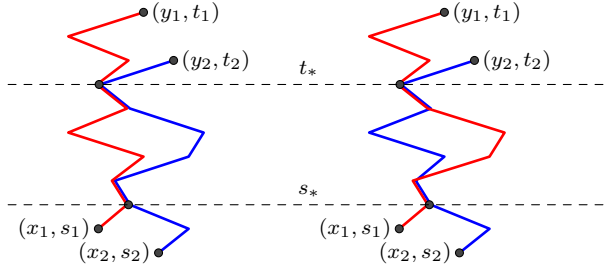
\begin{figure}[hbt!]
    \centering
\begin{tikzpicture}[line cap=round,line join=round,>=triangle 45,x=0.2cm,y=0.32cm]
\clip(-11,-2) rectangle (31,11);

\draw [line width=1pt,color=blue] (0,0) -- (2,1) -- (-2+0.1,2) -- (-3+0.1,3) -- (2,4) -- (3,5) -- (-2+0.1,6) -- (-4+0.1,7) -- (1,8);
\draw [line width=1pt,color=red] (-4,1) -- (-2.1,2) -- (-3.1,3) -- (-1,4) -- (-6,5) -- (-2.1,6) -- (-4.1,7) -- (-2,8) -- (-6,9) -- (-1,10);

\draw [line width=1pt,color=blue] (20,0) -- (22,1) -- (18+0.1,2) -- (17+0.1,3) -- (19,4) -- (14,5) -- (18+0.1,6) -- (16+0.1,7) -- (21,8);
\draw [line width=1pt,color=red] (16,1) -- (18-0.1,2) -- (17-0.1,3) -- (22,4) -- (23,5) -- (18-0.1,6) -- (16-0.1,7) -- (18,8) -- (14,9) -- (19,10);

\draw [fill=uuuuuu] (0,0) circle (1.5pt);
\draw [fill=uuuuuu] (1,8) circle (1.5pt);
\draw [fill=uuuuuu] (-4,1) circle (1.5pt);
\draw [fill=uuuuuu] (-1,10) circle (1.5pt);

\draw [fill=uuuuuu] (20,0) circle (1.5pt);
\draw [fill=uuuuuu] (21,8) circle (1.5pt);
\draw [fill=uuuuuu] (16,1) circle (1.5pt);
\draw [fill=uuuuuu] (19,10) circle (1.5pt);

\draw [fill=uuuuuu] (-2,2) circle (1.5pt);
\draw [fill=uuuuuu] (18,2) circle (1.5pt);
\draw [fill=uuuuuu] (-4+0.05,7) circle (1.5pt);
\draw [fill=uuuuuu] (16+0.05,7) circle (1.5pt);

\begin{scriptsize}
\draw (0,0) node[anchor=east]{$(x_2,s_2)$};
\draw (-4,1) node[anchor=east]{$(x_1,s_1)$};
\draw (20,0) node[anchor=east]{$(x_2,s_2)$};
\draw (16,1) node[anchor=east]{$(x_1,s_1)$};

\draw (-1,10) node[anchor=west]{$(y_1,t_1)$};
\draw (1,8) node[anchor=west]{$(y_2,t_2)$};
\draw (19,10) node[anchor=west]{$(y_1,t_1)$};
\draw (21,8) node[anchor=west]{$(y_2,t_2)$};

\draw (10,2) node[anchor=south]{$s_*$};
\draw (10,7) node[anchor=south]{$t_*$};

\draw [dashed] (-10,2) -- (30,2);
\draw [dashed] (-10,7) -- (30,7);
\end{scriptsize}

\end{tikzpicture}
\caption{An illustration of the proof of Lemma \ref{lem:ordering-geo}: the left figure depicts $\pi_{(x_1,s_1;y_1,t_1)}$ and $\pi_{(x_2,s_2;y_2,t_2)}$, and the right figure depicts $\pi_1'$ (the red path) and $\pi_2'$ (the blue path), which are alternative candidates for the geodesics from $(x_1,s_1)$ to $(y_1,t_1)$ and from $(x_2,s_2)$ to $(y_2,t_2)$, respectively.}
\label{fig:211}
\end{figure}
\begin{lemma}  \label{lem:ordering-geo}
Almost surely the following is true.
For any $(x_1,s_1;y_1,t_1), (x_2,s_2;y_2,t_2)\in\R^4_\uparrow$ with $s_1\vee s_2 \le t_1 \wedge t_2$, the set $O=\{t'\in\R: \pi_{(x_1,s_1;y_1,t_1)}(t') = \pi_{(x_2,s_2;y_2,t_2)}(t')\}$ is either empty or a closed interval, if one of the following two conditions hold:
\begin{enumerate}
    \item both $\pi_{(x_1,s_1;y_1,t_1)}$ and $\pi_{(x_2,s_2;y_2,t_2)}$ are the leftmost geodesics;
    \item at least one of $\pi_{(x_1,s_1;y_1,t_1)}$ and $\pi_{(x_2,s_2;y_2,t_2)}$ is the unique geodesic (between its endpoints).
\end{enumerate}
\end{lemma}
\begin{proof}
We assume that $O\neq \emptyset$.
The denote $s_* = \inf O$ and $t_* = \sup O$.
As the geodesics are continuous, we have that $s_*, t_* \in O$.
Now we define $\pi_1':[s_1, t_1] \to \R$ and $\pi_2':[s_2, t_2] \to \R$ as follows:
for any $t'\in [s_1,t_1]\setminus [s_*, t_*]$, we let $\pi_1'(t')=\pi_{(x_1,s_1;y_1,t_1)}(t')$, and for any $t'\in [s_2,t_2]\setminus [s_*, t_*]$ we let $\pi_2'(t')=\pi_{(x_2,s_2;y_2,t_2)}(t')$.
For any $t' \in (s_*, t_*)$ we let $\pi_1'(t')=\pi_{(x_2,s_2;y_2,t_2)}(t')$ and $\pi_2'(t')=\pi_{(x_1,s_1;y_1,t_1)}(t')$.
See Figure \ref{fig:211} for an illustration.
Then it is straightforward to check that \[\|\pi_1'\|_\cL + \|\pi_2'\|_\cL = \|\pi_{(x_1,s_1;y_1,t_1)}\|_\cL + \|\pi_{(x_2,s_2;y_2,t_2)}\|_\cL = \cL(x_1,s_1;y_1,t_1) + \cL(x_2,s_2;y_2,t_2).\]
This implies that $\pi_1'$ and $\pi_2'$ are also geodesics. 

In the case where both $\pi_{(x_1,s_1;y_1,t_1)}$ and $\pi_{(x_2,s_2;y_2,t_2)}$ are the leftmost geodesics, for any $t' \in (s_*, t_*)$ we must have that $\pi_1'(t')\ge \pi_{(x_1,s_1;y_1,t_1)}(t')$, which implies that $\pi_{(x_2,s_2;y_2,t_2)}(t')\ge \pi_{(x_1,s_1;y_1,t_1)}(t')$. By a symmetric argument we have that $\pi_{(x_2,s_2;y_2,t_2)}(t')\le \pi_{(x_1,s_1;y_1,t_1)}(t')$, so we must have that $t' \in O$, and $O=[s_*, t_*]$.

We next consider the case where $\pi_{(x_1,s_1;y_1,t_1)}$ is the unique geodesic.
Then by uniqueness we must have $\pi_1' = \pi_{(x_1,s_1;y_1,t_1)}$, so we conclude that $\pi_{(x_2,s_2;y_2,t_2)}(t')= \pi_{(x_1,s_1;y_1,t_1)}(t')$ for any $t' \in (s_*, t_*)$, and $O=[s_*, t_*]$.
By symmetry the case where $\pi_{(x_2,s_2;y_2,t_2)}$ is the unique geodesic also follows.
\end{proof}
\subsubsection{Disjointness of geodesics}

Recall that for any $x_1< x_2$, $y_1<y_2$, and $s<t$,  $\Dist_{(x_1, x_2)\to (y_1, y_2)}^{[s, t]}$ is the event where $\pi_{(x_1, s;y_1, t)}$ and $\pi_{(x_2, s;y_2, t)}$ are unique and disjoint.
This disjointness would be equivalent to strict inequality of the quadrangle inequality (Lemma \ref{lem:DL-quad}). 
\begin{lemma}  \label{lem:dis-equiv}
For any $x_1< x_2$, $y_1<y_2$, and $s<t$, the event $\Dist_{(x_1, x_2)\to (y_1, y_2)}^{[s, t]}$ is equivalent to that the geodesics $\pi_{(x_1, s; y_1, t)}$ and $\pi_{(x_2, s; y_2, t)}$ are both unique, and
\begin{equation}  \label{eq:dis-equiv}
\cL(x_1,s;y_1,t) + \cL(x_2,s;y_2,t) > \cL(x_1,s;y_2,t) + \cL(x_2,s;y_1,t).    
\end{equation}
\end{lemma}
Statements similar to this have also appeared in the proof of \cite[Proposition 4.1]{BGHhau}.
\begin{proof}
Whenever $\pi_{(x_1, s; y_1, t)}$ and $\pi_{(x_2, s; y_2, t)}$ are not disjoint, we can take $x_*, t_*$ with \[\pi_{(x_1, s; y_1, t)}(t_*)=\pi_{(x_2, s; y_2, t)}(t_*)=x_*,\]
and construct $\pi_1'$ from $(x_1, s)$ to $(y_2, t)$ and $\pi_2'$ from $(x_2, s)$ to $(y_1, t)$, by switching the paths after $t_*$.
Then $\|\pi_1'\|_\cL+\|\pi_2'\|_\cL=\cL(x_1,s;y_1,t) + \cL(x_2,s;y_2,t)$, and \eqref{eq:dis-equiv} cannot hold.

Now we assume that \eqref{eq:dis-equiv} does not hold; then by the quadrangle inequality (Lemma \ref{lem:DL-quad}), we must have that the left hand side equals the right hand side.
We then consider any geodesic $\pi_{(x_1, s; y_2, t)}$ and $\pi_{(x_2, s; y_1, t)}$.
Take $x_*, t_*$ with \[\pi_{(x_1, s; y_2, t)}(t_*)=\pi_{(x_1, s; y_2, t)}(t_*)=x_*,\]
then we construct $\pi_1'$ from $(x_1, s)$ to $(y_1, t)$ and $\pi_2'$ from $(x_2, s)$ to $(y_2, t)$, by switching the paths after $t_*$.
Now we have 
\[\|\pi_1'\|_\cL+\|\pi_2'\|_\cL=\cL(x_1,s;y_2,t) + \cL(x_1,s;y_2,t)=\cL(x_1,s;y_1,t) + \cL(x_2,s;y_2,t),\]
so $\pi_1'$ and $\pi_2'$ must be geodesics.
By uniqueness we have $\pi_1'=\pi_{(x_1, s; y_1, t)}$ and $\pi_2'=\pi_{(x_2, s; y_2, t)}$, and they are not disjoint.
\end{proof}

We recall the following  estimate on the disjoint probability. This is a limiting version of \cite[Theorem 1.1]{H20Exp}, and is slightly more general than Theorem \ref{thm:disj1}.
\begin{theorem}  \label{thm:disj}
There is a constant $C$ depending on $g, h$, such that for any $g\le t \le h$, and $\epsilon>0$, $l>1$ we have
\[
\PP[ \Dist_{(-l,l)\to (-\epsilon, \epsilon)}^{[0,t]} ] < Cl\epsilon^{1/2} e^{C |\log(\epsilon)|^{5/6}}.
\]
\end{theorem}
\begin{proof}
By Lemma \ref{lem:dis-equiv}, $\Dist_{(-l,l)\to (-\epsilon, \epsilon)}^{[0,t]}$ implies that
\[
\cL(-1,0;-\epsilon,t) - \cL(-1,0;\epsilon,t) > \cL(1,0;-\epsilon,t) - \cL(1,0;\epsilon,t);
\]
thus for some $k\in \llbracket 0, \lceil l \rceil$, we have
\[
\cL(-\lceil l \rceil + 2k,0;-\epsilon,t) - \cL(-\lceil l \rceil + 2k,0;\epsilon,t) > \cL(-\lceil l \rceil + 2k+2,0;-\epsilon,t) - \cL(-\lceil l \rceil + 2k+2,0;\epsilon,t).
\]
Assuming (the probability $1$ event) that these geodesics are unique, this implies
$$\bigcup_{k=0}^{\lceil l\rceil}\Dist_{(-\lceil l \rceil + 2k, -\lceil l \rceil + 2k+2)\to (-\epsilon, \epsilon)}^{[0,t]},$$ by Lemma \ref{lem:dis-equiv}.
Then by the skew-shift invariance of the directed landscape, we have
\[
\PP[ \Dist_{(-l,l)\to (-\epsilon, \epsilon)}^{[0,t]} ] \le \lceil l \rceil \PP[ \Dist_{(-1,1)\to (-\epsilon, \epsilon)}^{[0,t]} ].
\]
We now claim that
\[
\PP[ \Dist_{(-1,1)\to (-\epsilon, \epsilon)}^{[0,t]} ] < C\epsilon^{1/2} e^{C |\log(\epsilon)|^{5/6}}.
\]
for some constant $C$ depending on $g, h$.
Such a disjointness probability estimate is proved in the setting of Brownian LPP, as
\cite[Theorem 1.1]{H20Exp}; and we pass it to the directed landscape limit using \cite[Theorem 1.8]{DOV}, which states that (under an appropriate coupling) any unique geodesic in the directed landscape is the limit (under uniform convergence) of the corresponding ones in Brownian LPP.
The above two inequalities imply the conclusion.
\end{proof}

\subsection{Airy line ensemble and properties}  \label{ss:ale}

The parabolic Airy line ensemble $\cA=\{\cA_n\}_{n=1}^\infty$ is an ordered family of random processes on $\R$, constructed in \cite{CH}.
We summarize some crucial properties of it that we will make use of.
Define $\ocA_n(x):=\cA_n(x)+x^2$ for any $x\in \R$, $n\in \N$.
\begin{enumerate}
    \item Stationarity: $\ocA=\{\ocA_n\}_{n=1}^\infty$ is stationary with respect to horizontal shift.
    \item Ergodicity: $\ocA=\{\ocA_n\}_{n=1}^\infty$ is ergodic with respect to horizontal shift.
    \item Brownian Gibbs property: for any $m\le n\in \N$ and $x\le y$, let $S = \llbracket m,n\rrbracket \times [x, y]$. We can think of $\cA$ as a (random) function on $\N\times \R$, then the process $\cA|_S$ given $\cA|_{S^c}$ is just $n-m + 1$ Brownian bridges connecting up the points $\cA_i(x)$ and $\cA_i(y)$, for $m\le i \le n$, conditioned on that $\cA_i(z)>\cA_{i+1}(z)$ for any $1\vee (m-1) \le i \le n$ and $z\in [x,y]$.
\end{enumerate}
The first property is from its construction, and the third property is proved in \cite{CH}. The second property is proved in \cite{CS}.
We note that here and throughout this paper, all Brownian motions and Brownian bridges have diffusive coefficients equal to $2$.

The first line $\cA_1$ (of $\cA$) is the parabolic Airy$_2$ process on $\R$, and it has the same distribution as the function $x\mapsto \cL(\boo;x,1)$.
This process also has the following Brownian property.
For $K\in \R, d>0$, let $\cB^{[K,K+d]}$ denote a Brownian motion on $[K,K+d]$, taking value $0$ at $K$. Let $\cA_1^{[K,K+d]}$ denote the random function on $[K,K+d]$ defined by 
$$\cA_1^{[K,K+d]}(x):=\cA_1(x)-\cA_1(K), ~\forall x\in [K,K+d].$$

\begin{theorem}[\protect{\cite[Theorem 1.1]{HHJ}}] \label{thm:airytail}
There exists an universal constant $G>0$ such that the following holds.
For any fixed $M>0,$ there exists $a_0=a(M)$ such that for all intervals $[K,K+d]\subset [-M,M]$ and for all measurable $A \subset C^0([K,K+d])$ with $0<\PP[\cB^{[K,K+d]}\in A]=a\le a_0,$
\[\PP\left(\cA_1^{[K,K+d]} \in A\right) \leq a e^{ G {M} ( \log a^{-1} )^{5/6} }.\]
Here $C^0 ([K,K+d])$ is the space of all continuous functions on $[K,K+d]$.
\end{theorem}

We will also need the Brownian absolute continuity for the evolution induced by the directed landscape, with some general initial conditions.
Namely, take any $t>0$ and any function $f_0:\R\to\R\cup\{-\infty\}$, such that $f(x)\neq -\infty$ for some $x$, and $\frac{f_0(x)-x^2/t}{|x|} \to -\infty$ as $|x|\to \infty$. Define
\[
f_t(y) := \sup_{x\in\R} f_0(x) + \cL(x,0;y,t).
\]
Such evolution of functions is also called the KPZ fixed point (see \cite{MQR}).
\begin{theorem}[\protect{\cite[Theorem 1.2]{SV21}}]  \label{thm:abs-bm-gen}
For any fixed $t>0$ and $f_0$ satisfying the above conditions, and any $y_1<y_2$, the law of $y\mapsto f_t(y+y_1)-f(y_1)$ for $y \in [0, y_2-y_1]$ is absolutely continuous with respect to the law of a Brownian motion on $[0, y_2-y_1]$.
\end{theorem}

\section{Busemann function and semi-infinite geodesics}  \label{sec:buse-semiinf-dl}

In this section we construct semi-infinite geodesics in the directed landscape. As stated in Section \ref{s:iop}, we shall first construct the Busemann function, and use it to define the semi-infinite geodesics. We will further deduce various useful properties of these constructed objects.

\subsection{Construction of Busemann function}  \label{ssec:constbuse}
For each $l>0,$ let
\begin{equation} \label{eq:defnUl}
U_l = [-2l^{2/3}, 2l^{2/3}]\times [l, 2l].
\end{equation}
For any compact set $B\subset \R^2$, we denote as $\cC_{B,l}$  the following coalescence event: there exists $(x_*, t_*)$ such that for any $p\in B$ and $q \in U_l$, any geodesic $\pi_{(p;q)}$ contains $t_*$ in its range, and $\pi_{(p;q)}(t_*)=x_*$. The following shows that this is rather likely. 
\begin{lemma}  \label{lem:asy-coal}
For any compact $B \subset \R^2$, there is a constant $C>0$, such that $\PP[\cC_{B,l}]> 1-Cl^{-1/3}e^{C|\log(l)|^{5/6}}$ for any $l>0$.
\end{lemma}

\begin{figure}[hbt!]
    \centering
\begin{tikzpicture}[line cap=round,line join=round,>=triangle 45,x=0.6cm,y=0.4cm]
\clip(-4,-3) rectangle (14,10);

\fill[line width=0.pt,color=green,fill=green,fill opacity=0.2]
(1,4) -- (9,4) -- (9,10) -- (1,10) -- cycle;

\fill[line width=0.pt,color=blue,fill=blue,fill opacity=0.2]
(4,-2) -- (6,-2) -- (6,-1) -- (4,-1) -- cycle;

\draw [dashed] (-15,4) -- (25,4);
\draw [dashed] (-15,-1) -- (25,-1);

\draw [line width=1pt] (2.5,-1) -- (5,0) -- (5.5,1) -- (4.5,2) -- (-2,3) -- (-2.5,4);
\draw [line width=1pt] (7.5,-1) -- (5,0) -- (5.5,1) -- (4.5,2) -- (6,3) -- (12.5,4);

\draw [line width=0.5pt,color=blue] (4.3,-1.8) -- (5.6,-1) -- (5,0) -- (5.5,1) -- (4.5,2) -- (5.3,3) -- (3.3,4) -- (4.8,5) -- (4.1,6) -- (5.5,7) -- (5.6,7.5) -- (4.8,8);

\draw [fill=uuuuuu] (2.5,-1) circle (1.5pt);
\draw [fill=uuuuuu] (7.5,-1) circle (1.5pt);

\draw [fill=uuuuuu] (-2.5,4) circle (1.5pt);
\draw [fill=uuuuuu] (12.5,4) circle (1.5pt);

\draw [blue] [fill=blue] (4.3,-1.8) circle (1.5pt);
\draw [blue] [fill=blue] (4.8,8) circle (1.5pt);

\begin{scriptsize}
\draw (2.5,-1) node[anchor=south east]{$(-\log(l)^5l^{-2/3}, l^{-1})$};
\draw (7.5,-1) node[anchor=south west]{$(\log(l)^5l^{-2/3}, l^{-1})$};

\draw (6,-1.5) node[anchor=north west]{$[-l^{-2/3}, l^{-2/3}] \times [-l^{-1}, l^{-1}]$};
\draw (9,7) node[anchor=north west]{$U_1=[-2,2]\times [1,2]$};

\draw (-2.5,4) node[anchor=south]{$(-\log(l)^5,1)$};
\draw (12.5,4) node[anchor=south]{$(\log(l)^5,1)$};

\draw [blue] (4.3,-1.8) node[anchor=east]{$(x,t)$};
\draw [blue] (4.8,8) node[anchor=east]{$q$};

\end{scriptsize}

\end{tikzpicture}
\caption{An illustration of the proof of Lemma \ref{lem:asy-coal}.}
\label{fig:301}
\end{figure}
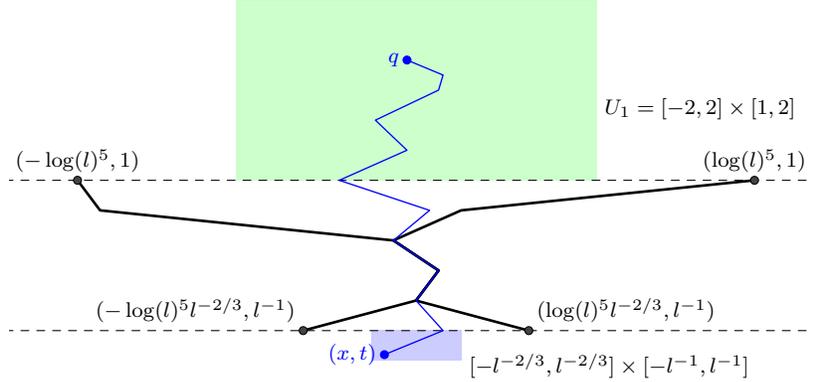

The above is a general version of statements that have already appeared in the literature before (see e.g., \cite{BSS19,BGZ,Z20} which works in the exponential LPP setting). We follow a similar proof relying on the bound given by Theorem \ref{thm:disj} (which accounts for the extra factor of $e^{C|\log(l)|^{5/6}}$), combined with transversal fluctuation estimates and the ordering of geodesics.
\begin{proof}
By scale invariance, we can assume that $B=[-1,1] \times [-1,1]$ without loss of generality. We also assume that $l$ is large enough (since otherwise the conclusion holds obviously by taking $C$ large enough). Let the event $\cC_l'$ be the scaled version of $\cC_{B,l},$  i.e., there exists $(x_*,t_*)$ such that $l^{-1}\le t_* \le 1$, and $\pi_{(x,t;q)}(x_*)=t_*$ for any $|x|\le l^{-2/3}$, $|t|\le l^{-1}$, and $q\in U_1$, and any geodesic $\pi_{(x,t;q)}$. Hence  $\PP[\cC_{B,l}]=\PP[\cC_l']$,

Below we let $C, c>0$ denote large and small constants, whose values may change from line to line.
Let $R$ be the random varialbe given by Lemma \ref{lem:dl-trans-f}.
We now consider the following events.
\begin{enumerate}
    \item Let $\cE_1$ be the event that $R>\log(l)$. 
    \item Let $\cE_2$ be the event $\Dist_{(-\log(l)^5l^{-2/3},\log(l)^5l^{-2/3})\to (-\log(l)^5, \log(l)^5)}^{[l^{-1},1]}$; recall that this is the event where the geodesics $\pi_{(-\log(l)^5l^{-2/3}, l^{-1}; -\log(l)^5, 1)}$ and $\pi_{(\log(l)^5l^{-2/3}, l^{-1}; \log(l)^5, 1)}$ are unique and disjoint (while the coordinates are somewhat  technical looking, we hope Figure \ref{fig:301} would convey the basic setup).
\end{enumerate}
We next show that (up to a probability zero event) $\cE_1^c\cap \cE_2^c \subset \cC_l'$ (see Figure \ref{fig:301} for an illustration).
Using Lemma \ref{lem:dl-trans-f}, for any $|x|\le l^{-2/3}$, $|t|\le l^{-1}$, and $q\in U_1$, we have $|\pi_{(x,t;q)}(1)|\le \log(l)^5$ and $|\pi_{(x,t;q)}(l^{-1})|\le \log(l)^5l^{-2/3}$ under $\cE_1^c$.
Now on $\cE_2^c$ and hence on $\cE_1^c\cap \cE_2^c$, we can find $(x_*, t_*)$ such that
\[\pi_{(-\log(l)^5l^{-2/3}, l^{-1}; -\log(l)^5, 1)}(t_*)=\pi_{(\log(l)^5l^{-2/3}, l^{-1}; \log(l)^5, 1)}(t_*)=x_*.\]
We assume (the probability $1$ event) that $\pi_{(-\log(l)^5l^{-2/3}, l^{-1}; -\log(l)^5, 1)}$ and $\pi_{(\log(l)^5l^{-2/3}, l^{-1}; \log(l)^5, 1)}$ are unique. By the ordering of finite geodesics (Lemma \ref{lem:ordering-geo}), for any $|x|\le l^{-2/3}$, $|t|\le l^{-1}$, and $q\in U_1$, one also has $\pi_{(x,t;q)}(t_*)=x_*$.

By Lemma \ref{lem:dl-trans-f}, we have $\PP[\cE_1]<Ce^{-c(\log(l))^{9/4}\log(\log(l))^{-4}}$,
by Theorem \ref{thm:disj} we have 
\[\PP[\cE_2]<C\log(l)^{15/2}l^{-1/3}e^{C|\log(l)|^{5/6}}.\]
Adding these bounds together yields the conclusion.
\end{proof}

We then immediately get the following convergence of difference of passage times.
\begin{cor}  \label{cor:buse-prep}
For any compact $B \subset \R^2$ and any $r\in \R$, almost surely there is a (random) number $h_0$ such that the following is true.
For any $p, q \in B$, $h\ge h_0$, and $z \in \R$ with $|z-hr|<h^{2/3}$, we have that
$\cL(p;z,h)-\cL(q;z,h) = \cL(p;h_0r,h_0)-\cL(q;h_0r,h_0)$.
\end{cor}
\begin{proof}
By skew-shift invariance of the directed landscape (Lemma \ref{lem:DL-symmetry}), we can assume that $r=0$ without loss of generality.
By Lemma \ref{lem:asy-coal} and the Borel-Cantelli lemma, there is a random $k_0\in \N$, such that $\cC_{B,2^k}$ holds for any $k\ge k_0$.
Then for any $p,q \in B$, $k\ge k_0$ and $p' \in U_{2^k}$, we have $\cL(p;0,2^k)-\cL(q;0,2^k) = \cL(p;p')-\cL(q;p')$; in particular $\cL(p;0,2^k)-\cL(q;0,2^k) = \cL(p;0,2^{k+1})-\cL(q;0,2^{k+1})$.
This implies that for any $h\ge 2^{k_0}$ and $|z|<h^{2/3}$, we have
$\cL(p;0,2^{k_0})-\cL(q;0,2^{k_0}) = \cL(p;z,h)-\cL(q;z,h)$, since $(z,h) \in \cup_{k\ge k_0}U_{2^k}$.
Thus the conclusion follows by taking $h_0=2^{k_0}$.
\end{proof}
As in the Exponential LPP setting, the above coalescence results allow us to define the Busemann function in the directed landscape.
\begin{defn}  \label{defn:DL-Buse}
We define the Busemann function as
\[
\cM^r(x,s;y,t) := \lim_{n\to\infty}\cL(x,s;z_n,h_n) - \cL(y,t;z_n,h_n)
\]
for any $(x,s;y,t) \in \R^4$ and $r\in \R$, where $\{z_n\}_{n=1}^\infty$ and $\{h_n\}_{n=1}^\infty$ are any sequences such that $z_n\to \infty$ as $n\to\infty$, and $|z_n-h_nr|<h_n^{2/3}$ for each $n\in \N$.
By Corollary \ref{cor:buse-prep} we see that for any fixed $r\in\R$, almost surely this limit exists for any $(x,s;y,t) \in \R^4$, and is independent of the choice of the sequences. Also, $\cM^r$ is almost surely continuous as a function on $\R^4$.
For simplicity of notations, we would often write $\cM$ for $\cM^0$.
\end{defn}

The next lemma asserts Brownianity of the Busemann functions.

\begin{lemma}  \label{lem:buse-bm}
For any $r, t\in\R$, the process of
$x\mapsto \cM^r(x,t;0,t) - 2xr$ is a Brownian motion.
\end{lemma}
\begin{proof}
By skew-shift invariance (Lemma \ref{lem:DL-symmetry}) it suffices to assume $r=t=0$.
For any $h>0$, by scaling invariance, the process
$x\mapsto \cL(x,0;0,h)-\cL(\boo;0,h)$ has the same law as $x\mapsto h^{-1/3}\cL(xh^{-2/3},0;0,1)-\cL(\boo;0,1)$, which is the same as $x\mapsto h^{-1/3}(\cA_1(h^{-2/3}x)-\cA_1(0))$, where $\cA_1$ is the parabolic Airy$_2$ process.
By \cite[Theorem 4.14]{MQR}, as $h\to \infty$ it converges to Brownian motion, in finite dimensional distributions.
Also by Corollary \ref{cor:buse-prep}, 
$x\mapsto \cL(x,0;0,h)-\cL(\boo;0,h)$ converges to $x\mapsto \cM(x,t;0,t)$ uniformly in compact sets,
so $x\mapsto \cM(x,0;\boo)$ has the same finite dimensional distributions as a Brownian motion.
Since $\cM$ is continuous, the conclusion follows.
\end{proof}

\begin{lemma}  \label{lem:quadra}
For any $r_1 < r_2$, the following is true almost surely.
For any $x_1 < x_2$ and $t\in \R$, there is
\begin{equation}  \label{eq:quadra}
 \cM^{r_1}(x_1,t;x_2,t) \ge \cM^{r_2}(x_1,t;x_2,t).   
\end{equation}
\end{lemma}
\begin{proof}
It suffices to prove this for fixed (thus all rational) $x_1<x_2$ and $t\in \R$ since then by continuity of $\cM^{r_1}$ and $\cM^{r_2}$ we can get \eqref{eq:quadra} to hold for all $x_1<x_2$ and $t\in \R$ simultaneously.

From the construction, we can find some rational $h_0>t$, such that
$\cM^{r_1}(x_1,t;x_2,t) = \cL(x_1,t;h_0r_1,h_0) - \cL(x_2,t;h_0r_1,h_0)$, and 
$\cM^{r_2}(x_1,t;x_2,t) = \cL(x_1,t;h_0r_2,h_0) - \cL(x_2,t;h_0r_2,h_0)$.
Then by the quadrangle inequality (Lemma \ref{lem:DL-quad}), the conclusion follows.
\end{proof}

We now move on to the formal construction of semi-infinite geodesics. 

\subsection{Construction of semi-infinite geodesics}

\begin{lemma}  \label{lem:semi-inf-exist}
For any $r\in \R$, almost surely the following is true.
For any $p=(x,s)\in\R^2$, there exists a continuous path $\pi_p^r:[s,\infty)\to \R$
satisfying the following condition.

Consider any compact $B\subset \R^2$, and any sequences $\{z_n\}_{n=1}^\infty$ and $\{h_n\}_{n=1}^\infty$ satisfying  $z_n\to \infty$ as $n\to\infty$ and $|z_n-h_nr|<h_n^{2/3}$ for each $n\in \N$, and  $t\in \R$.
Then there is some $n_0\in\N$ such that $\pi_{(p;z_n,h_n)}(t')=\pi_p^r(t')$ for any $n\ge n_0$, any $p=(x,s)\in B$ with $s\le t$, and any $s\le t' \le t$. Here $\pi_{(p;z_n,h_n)}$ is taken as the leftmost geodesic.
\end{lemma}

\begin{proof}
The finite geodesics in this proof are taken to be the leftmost one.
We also just prove for the case where $r=0$, since the argument for general $r$ follows from the skew-shift invariance.

By Lemma \ref{lem:asy-coal} and using the Borel-Cantelli lemma, almost surely, for each $m\in \N$ there is a random $k_m\in \N$, such that the coalescence event $\cC_{[-m,m]\times[-m,m],2^k}$ (recall the definition from the beginning of Section \ref{ssec:constbuse}) holds for any $k\ge k_m$.
Then for any $p=(x,s)\in [-m,m]\times[-m,m]$, $k\ge k_m$, $s\le t\le m$, and $p' \in U_{2^k}$ (defined in \eqref{eq:defnUl}), we have $\pi_{(p;0,2^k)}(t)= \pi_{(p;p')}(t)$ by Lemma \ref{lem:ordering-geo};
in particular $\pi_{(p;0,2^k)}(m)= \pi_{(p;0,2^{k+1})}(m)$.
This implies that for any $h\ge 2^{k_m}$ and $|z|<h^{2/3}$, and $s\le t\le m$, we have
$\pi_{(p;0,2^{k_m})}(t)= \pi_{(p;z,h)}(t)$, since $(z,h) \in \cup_{k\ge k_m}U_{2^k}$.

For any $p=(x,s)\in [-m,m]\times[-m,m],$ we just define $\pi_p^0(t) = \pi_{(p;0,2^{k_m})}(t)$ for any $s\le t \le m$.
Since $\pi_{(p;0,2^{k_m})}(t) = \pi_{(p;0,2^{k_{m+1}})}(t)$ (for any $s\le t \le m$) we would have that this definition is consistent for different $m$.
Now for a general compact set $B\subset \R^2$ and $t\in \R$, the lemma follows by taking $m$ large enough so that $B\subset [-m,m]\times[-m,m]$ and $m>t$. 
\end{proof}

As already indicated in the introduction, given $r\in \R$, unless otherwise noted, for any $p \in \R^2,$ we will use $\pi_p^r$ to denote the path constructed in Lemma \ref{lem:semi-inf-exist} and $\pi_p=\pi_p^0$ for simplicity of notations.

Just like the composition law \eqref{eq:DL-compo} satisfied by finite geodesics, the number $\pi_p^r(t)$ would be the maximizer of the sum of two profiles: the passage time from $p$, and the Busemann function in the $r$-direction.
\begin{lemma}  \label{lem:semi-compo}
For any fixed $r, t\in \R$, almost surely the following is true. 
For any $p=(x,s)\in\R^2$ with $s<t$, we have
\[
\cL(p;\pi_p^r(t),t) + \cM^r(\pi_p^r(t),t;0,t) = \max_x \cL(p;x,t) + \cM^r(x,t;0,t).
\]
\end{lemma}

\begin{proof}
By Corollary \ref{cor:buse-prep} and Lemma \ref{lem:semi-inf-exist},
for any $n\in \N$, we can find $h_n$ large enough such that 
$\cM^r(x,t;0,t) = \cL(x,t;h_nr,h_n)-\cL(0,t;h_nr,h_n)$ for any $|x| \le n$, and $\pi_p^r(t) = \pi_{(p;h_nr,h_n)}(t)$ for any $p\in [-n,n]\times [-n, t)$.

Then for any $p=(x,s)\in\R^2$ with $s<t$, and any $n\in \N$ with $n>2|\pi_p^r(t)|$ and $p\in [-n,n]\times [-n, n]$, by the composition law \eqref{eq:DL-compo} we have that \[\cL(p;\pi_p^r(t),t) + \cM^r(\pi_p^r(t),t;0,t) = \max_{|x|\le n} \cL(p;x,t) + \cM^r(x,t;0,t).\]
By sending $n\to\infty$ the conclusion follows.
\end{proof}

We next verify that the constructed paths $\pi_p^r$ are semi-infinite geodesics.

It is straight forward to see that the paths are locally geodesics. 
Fix any $r\in\R$.
By Lemma \ref{lem:semi-inf-exist}, for each $p=(x,s)\in\R^2$, and any $s\le t_1 < t_2$, the restriction of $\pi_p^r$ to $[t_1, t_2]$ is a subpath of $\pi_{(p;hr,h)}$ for some large enough $h$; so it is also a geodesic. 

It remains to show that $\pi_p^r$ have the desired asymptotic directions.
For that, we need the following transversal fluctuation estimate on semi-infinite geodesics. It will also be used in several other places in the rest of this paper, and can be viewed as a limiting version of Lemmas \ref{lem:trans-dis-inf} and \ref{lem:trans-semi-inf}. While we do not prove this, up to the constants, this is also the best upper bound one can hope for, as a matching lower bound can be obtained by analyzing the location of the $\argmax$ of the parabolic Airy$_2$ process plus a Brownian motion, and using the Brownian Gibbs property of the Airy line ensemble (the latter is an invariance property which is crucially used in Appendix \ref{sec:appa}). 
\begin{lemma}  \label{lem:dl-semi-trans}
There exist constants $C, c>0$ such that the following is true.
For any $r, t\in \R$ and $p=(x,s)\in \R^2$ satisfying $s<t$, we have $\PP[|\pi_p^r(t) - (x+r(t-s))| > w(t-s)^{2/3}] < Ce^{-cw^3}$.
\end{lemma}
\begin{proof}
In this proof, we let $C, c>0$ be universal constants, whose value can change from line to line.
By translation and skew-shift invariance of the directed landscape (thus of semi-infinite geodesics), we just need to prove this for $p=\boo$, $r=0$, and $t=1$.
By Lemma \ref{lem:semi-compo}, the event $|\pi_\boo(1)| > w$ implies that
\[
\max_{|y|>w} \cL(\boo;y,1) + \cM(y,1;0,1) \ge \cL(\boo;0,1).
\]
Denote this event as $\cE$. For any $n\in\Z$ we let $\cE_n$ be the event where
\[
\max_{y\in [n, n+1]} \cL(\boo;y,1) + \cM^r(y,1;0,1) \ge \cL(\boo;0,1).
\]
Then we have $\PP[\cE]\le \sum_{n\ge w} \PP[\cE_n] + \sum_{n\le w-1} \PP[\cE_n]$. 

Next, note that $\cE_n$ implies one of the following three events:
\begin{enumerate}
\item $\max_{y\in [n, n+1]} \cL(\boo;y,1) > -2n^2/3$.
\item $\cL(\boo;0,1) < -n^2/3$.
\item $\max_{y\in [n, n+1]} \cM^r(y,1;0,1) > n^2/3$.
\end{enumerate}
For the first event, by skew-shift invariance of the directed landscape, it has the same probability as the follow event:
\[
\max_{y\in [0, 1]} \cL(\boo;y,1) - (y+n)^2 + y^2 > -2n^2/3.
\]
By Lemma \ref{lem:DLbound}, its probability is bounded by $Ce^{-cn^3}$. Also, the probabilities of the second and third event above are also bounded by $Ce^{-cn^3}$, using Lemma \ref{lem:DLbound} and estimates on Brownian motions, respectively.
Thus we have $\PP[\cE_n] < Ce^{-cn^3}$.
By summing over $n\ge w$ and $n\le w-1$, the conclusion follows.
\end{proof}

We also need the following result which is on the overlaps between the constructed paths, and essentially says that the paths are ``ordered''. It is a consequence of the construction and Lemma \ref{lem:ordering-geo}.
It will be improved later (Lemma \ref{lem:semi-inf-tree-fixed}) for the case of paths in the same direction.
\begin{lemma}  \label{lem:semi-inf-tree}
For any $r_1, r_2\in \R$, almost surely the following is true.
For any $p_2=(x_1, s_1), p_2=(x_2, s_2) \in \R^2$, let $O=\{t\ge s_1\vee s_2: \pi_{p_1}^{r_1}(t)=\pi_{p_2}^{r_2}(t)\}$ be the set of ``overlap times''.
If $r_1\neq r_2$, then $O$ is either empty or a closed interval.
If $r_1=r_2$, then $O$ is either empty, or $O=[t_*,\infty)$ for some $t_*\ge s_1\vee s_2$.
\end{lemma}
\begin{proof}
We first consider the case where $r_1\neq r_2$.
In this case we must have that $O$ is contained in a compact set, since the asymptotic directions of $\pi_{p_1}^{r_1}$ and $\pi_{p_2}^{r_2}$ are different.
For any $t>s_1\vee s_2$, by Lemma \ref{lem:semi-inf-exist} there is some (random) $h$ large enough (depending on $p_1$, $p_2$, $r_1$, $r_2$ and $t$) such that
$\pi_{p_1}^{r_1}(t')=\pi_{(p_1;hr_1,h)}(t')$ for any $s_1\le t' \le t$, and $\pi_{p_2}^{r_2}(t')=\pi_{(p_2;hr_2,h)}(t')$ for any $s_2\le t' \le t$.
Thus by Lemma \ref{lem:ordering-geo}, the set $O\cap (-\infty, t]$ is either empty, or a closed interval. By sending $t\to\infty$ we get the conclusion.

We next consider the case where $r_1=r_2$.
Again we take any $t>s_1\vee s_2$, and by Lemma \ref{lem:semi-inf-exist} there is some (random) $h$ large enough (depending on $p_1$, $p_2$, $r_1=r_2$ and $t$) such that
$\pi_{p_1}^{r_1}(t')=\pi_{(p_1;hr_1,h)}(t')$ for any $s_1\le t' \le t$, and $\pi_{p_2}^{r_2}(t')=\pi_{(p_2;hr_1,h)}(t')$ for any $s_2\le t' \le t$.
By Lemma \ref{lem:ordering-geo}, the 
set $\{s_1\vee s_2 \le t' \le h: \pi_{(p_1;hr_1,h)}(t') = \pi_{(p_2;hr_1,h)}(t') \}$ is a closed interval with right endpoint being $h$;
so the set $O\cap (-\infty, t]$ is either empty, or a closed interval with right endpoint being $t$. Again by sending $t\to\infty$ we get the conclusion.
\end{proof}

Now we finish proving that the constructed paths are semi-infinite geodesics.

\begin{lemma}  \label{lem:is-geo}
For any fixed $r\in \R$, almost surely the constructed path $\pi_p^r$ is a semi-infinite geodesic, for each $p\in \R^2$.
\end{lemma}

\begin{proof}
Using skew-shift invariance, it suffice to prove for the case where $r=0$.

We start with a fixed $p\in \R^2$.
Above we have stated that any segment of $\pi_p$ is a geodesic, and it remains to show that it has the desired asymptotic direction. 
By Lemma \ref{lem:dl-semi-trans}, and the Borel-Cantelli lemma, almost surely there exists some random $k_0 > 0$, such that for any $k>k_0$, $k\in\N$, we have $|\pi_p(2^k)|<2^{2k/3+0.1k}$.
For any $k>k_0$, since $\pi_p$ between $(\pi_p(2^k), 2^k)$ and $(\pi_p(2^{k+1}), 2^{k+1})$ is a geodesic, using Lemma \ref{lem:dl-trans-f} we can bound $\sup_{t\in [2^k, 2^{k+1}]}|\pi_p(t)|$. Thus we conclude that $\lim_{t\to\infty} t^{-1}\pi_p(t) = 0$.

Now we have that almost surely, for any $p\in \Q^2$, we have $\lim_{t\to\infty} t^{-1}\pi_p(t) = 0$. Assuming this event, and consider some $p=(x,s)\not\in \Q^2$.
Then we can find some $t, x_-, x_+ \in \Q$, such that $t>s$ and $x_- < \pi_p(s) < x_+$. By Lemma \ref{lem:semi-inf-tree}, we must have that $\pi_{(x_-,t)}(t') \le \pi_p(t') \le \pi_{(x_+,t)}(t')$ for any $t'\ge t$.
Since $\lim_{t\to\infty} t^{-1}\pi_{(x_-,t)}(t) = \lim_{t\to\infty} t^{-1}\pi_{(x_+,t)}(t) = 0$, we must have $\lim_{t\to\infty} t^{-1}\pi_p(t) = 0$.
\end{proof}

\subsection{Properties of semi-infinite geodesics}

The first property is about the overlap of semi-infinite geodesics for starting points around a fixed point.
For any continuous path $\pi_n:[s_n,t_n]\to \R$, $n\in\N$, and $\pi:[s,t]\to \R$, we say that $\pi_n\to \pi$ \emph{in the overlap topology}, if for all large enough $n$, $O_n=\{t'\in [s,t]\cap [s_n,t_n]: \pi_n(t')=\pi(t')\}$ is an interval, and the endpoints of $O_n$ converge to $s$ and $t$.
We quote the following result on convergence in the overlap topology for finite geodesics.
\begin{lemma}[\protect{\cite[Lemma 3.1, Lemma 3.3]{DSV}}]\label{lem:overlap-finite}
Almost surely the following is true.
Let $(p_n; q_n) \to (p; q) \in \R^4_\uparrow$, and let $\pi_n$ be any sequence of geodesics from $p_n$ to $q_n$.
Suppose that there is a unique geodesic $\pi$ from $p$ to $q$.
Then $\pi_n \to \pi$ in the overlap topology.
\end{lemma}
From this, we can get overlap of semi-infinite geodesics.
\begin{lemma}  \label{lem:semi-overlap}
For any $r\in\R$, $p=(x,s) \in \R$ and $t>s$, almost surely there is an open neighborhood $\hat{B}$ of $p$, such that $\pi_p^r(t')=\pi_{p'}^r(t')$ for any $p' \in \hat{B}$ and $t'\ge t$.
\end{lemma}
\begin{proof}
Take a compact set $B$ containing an open neighborhood of $p$.
We assume (the probability $1$ event) that for any $h\in\N$, $h>s$, there is a unique geodesic $\pi_{(p;hr,h)}$ from $p$ to $(hr, h)$.

{By Lemma \ref{lem:semi-inf-exist}, we can find }some large enough (random) $h \in \N$ with $h>t$, such that $\pi_{(p';hr,h)}(t)=\pi_{p'}^r(t)$ for any $p'=(x',s') \in B$ with $s'\le t$.
By Lemma \ref{lem:overlap-finite}, there is an open neighborhood $\hat{B}$ of $p$, such that $\hat{B}\subset B$, and for any $p' \in \hat{B}$ there is $\pi_{(p';hr,h)}(t)=\pi_{(p;hr,h)}(t)$.
We then have that $\pi_p^r(t)=\pi_{p'}^r(t)$ for any $p' \in \hat{B}$.
By Lemma \ref{lem:semi-inf-tree} we have $\pi_p^r(t')=\pi_{p'}^r(t')$ for any $p' \in \hat{B}$ and $t'\ge t$.
\end{proof}

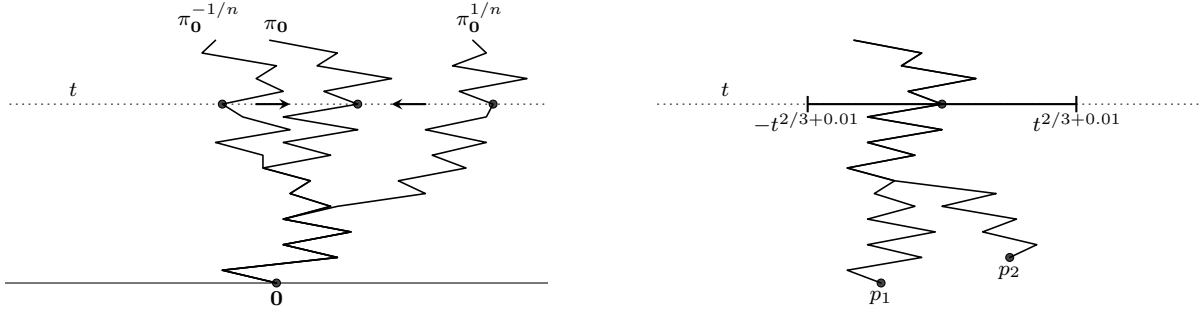
\begin{figure}
     \centering
     \begin{subfigure}[t]{0.45\textwidth}
         \centering
\begin{tikzpicture}
[line cap=round,line join=round,>=triangle 45,x=0.9cm,y=0.17cm]
\clip(1,-2) rectangle (9,22);

\draw (-1,0) -- (11,0);

\draw [dotted] (-1,14) -- (11,14);

\draw [fill=uuuuuu] (5,0) circle (1.5pt);

\draw [fill=uuuuuu] (6.2,14) circle (1.5pt);
\draw [fill=uuuuuu] (4.2,14) circle (1.5pt);
\draw [fill=uuuuuu] (8.2,14) circle (1.5pt);

\draw [line width=0.6pt] (5,0) -- (4.2,1) -- (5.9,2) -- (5.1,3) -- (6.1,4) -- (5.1,5) -- (5.8,6) -- (5.2,7) -- (5.5,8) -- (4.8,9) -- (5.8,10) -- (5.1,11) -- (6.2,12) -- (5.1,13) -- (6.2,14) -- (5.7,15) -- (6.7,16) -- (5.6,17) -- (5.9,18) -- (4.9,19);
\draw [line width=0.6pt] (5,0) -- (4.2,1) -- (5.9,2) -- (5.1,3) -- (6.1,4) -- (5.1,5) -- (5.8,6) -- (5.2,7) -- (5.5,8) -- (4.8,9) -- (4.8,10) -- (4.1,11) -- (5.2,12) -- (4.5,13) -- (4.2,14) -- (5.1,15) -- (4.7,16) -- (5,17) -- (3.9,18) -- (4.1,19);
\draw [line width=0.6pt] (5,0) -- (4.2,1) -- (5.9,2) -- (5.1,3) -- (6.1,4) -- (5.1,5) -- (5.9,6) -- (7.2,7) -- (6.8,8) -- (7.8,9) -- (7.3,10) -- (8.1,11) -- (7.2,12) -- (8.1,13) -- (8.2,14) -- (7.7,15) -- (8.7,16) -- (7.6,17) -- (8.1,18) -- (7.9,19);

\begin{scriptsize}
\draw (5,0) node[anchor=north]{$\boo$};
\draw (4,19) node[anchor=south]{$\pi_\boo^{-1/n}$};
\draw (5,19) node[anchor=south]{$\pi_\boo$};
\draw (8,19) node[anchor=south]{$\pi_\boo^{1/n}$};
\draw (2,14) node[anchor=south]{$t$};

\draw [-stealth] [thick] (4.7,14) -- (5.2,14);
\draw [-stealth] [thick] (7.2,14) -- (6.7,14);
\end{scriptsize}

\end{tikzpicture}
         \caption{Lemma \ref{lem:semi-inf-unique}: both $\pi_\boo^{-1/n}(t)$ and $\pi_\boo^{1/n}(t)$ converge to $\pi_\boo(t)$.}
         \label{fig:312}
     \end{subfigure}
     \hfill
     \begin{subfigure}[t]{0.5\textwidth}
         \centering
\begin{tikzpicture}
[line cap=round,line join=round,>=triangle 45,x=0.9cm,y=0.17cm]
\clip(2,-2) rectangle (10,22);

\draw [dotted] (-1,14) -- (15,14);

\draw [fill=uuuuuu] (5.3,0) circle (1.5pt);

\draw [fill=uuuuuu] (7.2,2) circle (1.5pt);

\draw [fill=uuuuuu] (6.2,14) circle (1.5pt);
\draw [|-|] [thick] (4.2,14) -- (8.2,14);

\draw [line width=0.6pt] (5.3,0) -- (4.8,1) -- (5.9,2) -- (5.1,3) -- (6.1,4) -- (5.1,5) -- (5.8,6) -- (5.2,7) -- (5.5,8) -- (4.8,9) -- (5.8,10) -- (5.1,11) -- (6.2,12) -- (5.1,13) -- (6.2,14) -- (5.7,15) -- (6.7,16) -- (5.6,17) -- (5.9,18) -- (4.9,19);
\draw [line width=0.6pt] (7.2,2) -- (7.6,3) -- (6.8,4) -- (7.3,5) -- (6.2,6) -- (7,7) -- (5.5,8) -- (4.8,9) -- (5.8,10) -- (5.1,11) -- (6.2,12) -- (5.1,13) -- (6.2,14) -- (5.7,15) -- (6.7,16) -- (5.6,17) -- (5.9,18) -- (4.9,19);

\begin{scriptsize}
\draw (3,14) node[anchor=south]{$t$};
\draw (5.3,0) node[anchor=north]{$p_1$};
\draw (7.2,2) node[anchor=north]{$p_2$};
\draw (4.2,14) node[anchor=north]{$-t^{2/3+0.01}$};
\draw (8.2,14) node[anchor=north]{$t^{2/3+0.01}$};

\end{scriptsize}

\end{tikzpicture}        \caption{Lemma \ref{lem:semi-inf-tree-fixed}: under $\cE_1^c$ where $|\pi_{p_1}(t)| \vee |\pi_{p_2}(t)| \le t^{2/3+0.01}$ and $\cE_2^c$ where $x\mapsto \cL(p_1;x,t)-\cL(p_2;x,t)$ is constant in $[-t^{2/3+0.01}, t^{2/3+0.01}]$, $\pi_{p_1}$ and $\pi_{p_2}$ would coalesce before $t$.}
         \label{fig:313}
     \end{subfigure}
        \caption{Illustrations for the proofs of Lemma \ref{lem:semi-inf-unique} and \ref{lem:semi-inf-tree-fixed}.}
        \label{fig:3e}
\end{figure}

If we fix both the starting point and the direction, then almost surely the semi-infinite geodesic is unique.
\begin{lemma}  \label{lem:semi-inf-unique}
For fixed $r\in \R$ and $p\in\R^2$, almost surely there is a unique semi-infinite geodesic from $p$ in the $r$ direction.
\end{lemma}
The proof uses the ordering of semi-infinite geodesics to sandwich any semi-infinite geodesic between two other semi-infinite geodesics with the same starting point but slightly different directions. We will also use the maximizer formulation of semi-infinite geodesics (Lemma \ref{lem:semi-compo}), and monotonicity of Busemann functions (Lemma \ref{lem:quadra}).
\begin{proof}
It suffices to prove uniqueness. 
By skew-shift invariance of the directed landscape (Lemma \ref{lem:DL-symmetry}), without loss of generality we can assume that $r=0$ and $p=\boo$.
We will work on the almost sure event that for any $m, n\in\N$, there is a unique geodesic from $\boo$ to $(m/n, m)$, and a unique geodesic from $\boo$ to $(-m/n, m)$.

Let $\pi:\R_{\ge 0}\to \R$ be a semi-infinite geodesic in the $0$ direction.
Take any fixed $t>0$, and $n, n'\in\N$ with $n<n'$. By the asymptotic directions and Lemma \ref{lem:semi-inf-exist}, we can find some $m\in \N$ such that $\pi_\boo^r(t) = \pi_{(\boo;mr,m)}(t)$, for each $r\in\{-1/n, 1/n, -1/n', 1/n'\}$,
and $|\pi(m)|<m/n'$.
Thus by Lemma \ref{lem:ordering-geo} we have
\[
\pi_\boo^{-1/n}(t)\le \pi_\boo^{-1/n'}(t)\le \pi(t)\le
\pi_\boo^{1/n'}(t) \le \pi_\boo^{1/n}(t).
\]
Then the limits $x_-:=\lim_{n\to\infty}\pi_\boo^{-1/n}(t)$ and $x_+:=\lim_{n\to\infty}\pi_\boo^{1/n}(t)$ exist, and $x_-\le \pi(t) \le x_+$. We would next show that $x_-=x_+$, and then $\pi(t)$ must equal $x_-=x_+$. See Figure \ref{fig:312} for an illustration.

Before proceeding further, we explain how to prove that $x_-=x_+$. Recall (from Lemma \ref{lem:semi-compo}) that the function $\cL(\boo;\cdot,t) + \cM^{-1/n}(\cdot,t;0,t)$ takes its maximum value at $\pi_\boo^{-1/n}(t)$. Then we prove uniform convergence of $\cM^{-1/n}(\cdot,t;0,t)$ to $\cM(\cdot,t;0,t)$, in a compact interval. 
By passing to the $n\to\infty$ limit, we get that $x_-$ is the $\argmax$ of the function $\cL(\boo;\cdot,t) + \cM(\cdot,t;0,t)$; and similarly this is also true for $x_+$. By uniqueness of the $\argmax$ we conclude that $x_-=x_+$.

By Lemma \ref{lem:semi-compo}, for any $y\in\R$ we have
\begin{equation}  \label{eq:pre-lim-qua}
\cL(\boo;\pi_\boo^{-1/n}(t),t) + \cM^{-1/n}(\pi_\boo^{-1/n}(t),t;x_-,t) \ge \cL(\boo;y,t) + \cM^{-1/n}(y,t;x_-,t).
\end{equation}
We would next send $n\to\infty$ for this equation.
By continuity of $\cL$ and $\lim_{n\to\infty}\pi_\boo^{-1/n}(t) = x_-$, we have that
\[
\lim_{n\to\infty} \cL(\boo;\pi_\boo^{-1/n}(t),t) = \cL(\boo;x_-,t).
\]
By Lemma \ref{lem:quadra}, we have that 
\[
\cM(\pi_\boo^{-1/n}(t),t;x_-,t)
\le
\cM^{-1/n}(\pi_\boo^{-1/n}(t),t;x_-,t)
\le
\cM^{-1}(\pi_\boo^{-1/n}(t),t;x_-,t).
\]
By continuity of $\cM$ and $\cM^{-1}$, we have $\cM(\pi_\boo^{-1/n}(t),t;x_-,t) \to 0$ and $\cM^{-1}(\pi_\boo^{-1/n}(t),t;x_-,t)\to 0$ as $n\to\infty$; so \[\lim_{n\to\infty}\cM^{-1/n}(\pi_\boo^{-1/n}(t),t;x_-,t) = 0.\]
For $\cM^{-1/n}(y,t;x_-,t)$, by Lemma \ref{lem:quadra} it is non-decreasing as $n\to\infty$, and $\cM^{-1/n}(y,t;x_-,t)\le \cM(y,t;x_-,t)$.
We claim that almost surely, for any $y\in\R$,
\begin{equation} \label{eq:lim-qua-n}
\lim_{n\to\infty}\cM^{-1/n}(y,t;x_-,t)= \cM(y,t;x_-,t).
\end{equation}
It suffices to show that for any $N\in\N$, almost surely this convergence holds for all $y\in [-N, N]$, conditional on that $|x_-|<N$.
We note that when $|x_-|, |y|<N$, by Lemma \ref{lem:quadra} there is
\begin{align*}
&\cM^{-1/n}(-N,t;N,t) - \cM^{-1/n}(y,t;x_-,t) \\
=&
\cM^{-1/n}(-N,t;y,t) + \cM^{-1/n}(x_-,t;N,t) \\
\ge &
\cM(-N,t;y,t) + \cM(x_-,t;N,t) \\
= & \cM(-N,t;N,t) - \cM(y,t;x_-,t),
\end{align*}
i.e., 
\[
\cM^{-1/n}(y,t;x_-,t)- \cM(y,t;x_-,t) \le \cM^{-1/n}(-N,t;N,t)- \cM(-N,t;N,t).
\]
The right hand side is almost surely positive, while the expectation equals $4N/n$; it is also non-increasing as $n\to\infty$. Thus by sending $n\to\infty$ we have that $\cM^{-1/n}(-N,t;N,t)- \cM(-N,t;N,t) \to 0$ almost surely.
Hence, conditional on $|x_-|<N$, we must have that \[\limsup_{n\to\infty}\cM^{-1/n}(y,t;x_-,t)\le \cM(y,t;x_-,t).\]
Similarly we have $\liminf_{n\to\infty}\cM^{-1/n}(y,t;x_-,t)\ge \cM(y,t;x_-,t)$. So conditional on $|x_-|<N$,  \eqref{eq:lim-qua-n} holds for all $y\in[-N, N]$. Then \eqref{eq:lim-qua-n} almost surely holds for all $y\in \R$.

Then by sending $n\to\infty$ for \eqref{eq:pre-lim-qua}, we conclude that almost surely, $\cL(\boo;x_-,t) \ge \cL(\boo;y,t) + \cM(y,t;x_-,t)$ for any $y\in\R$; by symmetry we also have
$\cL(\boo;x_+,t) \ge \cL(\boo;y,t) + \cM(y,t;x_+,t)$ for any $y\in\R$.
Thus we have that the function $x\mapsto \cL(\boo;x,t) + \cM(x,t;0,t)$ achieves its maximum at both $x_-$ and $x_+$.
However, $x\mapsto \cL(\boo;x,t)$ is a scaled Airy$_2$ process, so it is absolutely continuous with respect to Brownian motion (by \cite{CH} or Theorem \ref{thm:airytail}) in any compact interval; also $x\mapsto\cM(x,t;0,t)$ is a Brownian motion by Lemma \ref{lem:buse-bm}.
These imply that $x\mapsto \cL(\boo;x,t) + \cM(x,t;0,t)$ is also absolutely continuous with respect to Brownian motion in any compact interval.
As Brownian motion almost surely has a unique maximum in any compact interval, we conclude that $x_-=x_+$ almost surely.

Thus the above implies that for a given $t>0$, almost surely $\pi(t)$ is uniquely determined, which then by considering all rational $t>0$ along with continuity of $\pi$, implies the uniqueness of $\pi$.
\end{proof}

Finally, we state an improvement of Lemma \ref{lem:semi-inf-tree} for semi-infinite geodesics in a fixed direction; namely, for fixed two points, the semi-infinite geodesics to this direction coalesce almost surely.
\begin{lemma}  \label{lem:semi-inf-tree-fixed}
For any fixed $r\in \R$ and $p_1=(x_1, s_1), p_2=(x_2, s_2) \in \R^2$, almost surely the following is true.
Let $O=\{t\ge s_1\vee s_2: \pi_{p_1}^{r}(t)=\pi_{p_2}^{r}(t)\}$, then $O=[t_*,\infty)$ for some $t_*\ge s_1\vee s_2$.
\end{lemma}
\begin{proof}
By Lemma \ref{lem:DL-symmetry}, without loss of generality we will only consider the case $r=0$. It suffices to show that $\PP[\pi_{p_1}(t)=\pi_{p_2}(t)] \to 1$ as $t\to \infty$; then we have that almost surely $O\neq \emptyset$, and by Lemma \ref{lem:semi-inf-tree} the conclusion follows.

For fixed $t> s_1\vee s_2$, we consider the following events.
\begin{enumerate}
    \item Let $\cE_1$ be the event where $|\pi_{p_1}(t)| \vee |\pi_{p_2}(t)| > t^{2/3+0.01}$.
    \item Let $\cE_2$ be the event where the function $x\mapsto \cL(p_1;x,t)-\cL(p_2;x,t)$ is not constant for $x\in [-t^{2/3+0.01}, t^{2/3+0.01}]$.
\end{enumerate}
By Lemma \ref{lem:semi-compo}, the functions $x\mapsto \cL(p_1;x,t)+\cM(x,t;0,t)$ and $x\mapsto \cL(p_2;x,t)+\cM(x,t;0,t)$ take their maximum values at $\pi_{p_1}(t)$ and $\pi_{p_2}(t)$, respectively.
Under $\cE_2^c$, these two functions differ by a constant in the interval $[-t^{2/3+0.01}, t^{2/3+0.01}]$.
As argued in the proof of Lemma \ref{lem:semi-inf-unique}, these two functions are absolutely continuous with respect to a Brownian motion in any compact interval, so almost surely they have a unique $\argmax$ in the interval $[-t^{2/3+0.01}, t^{2/3+0.01}]$.
Thus (up to a probability zero event)
$\cE_1^c\cap \cE_2^c$ implies that $\pi_{p_1}(t)=\pi_{p_2}(t)$. See Figure \ref{fig:313} for an illustration.

As $t\to\infty$, by Lemma \ref{lem:dl-semi-trans} we have that $\PP[\cE_1]\to 0$. We net show that $\PP[\cE_2]\to 0$.
For each $i\in \Z$, $|i|<\lceil t^{0.01}\rceil$, let $\cE_{2,i}$ be the event where $x\mapsto \cL(p_1;x,t)-\cL(p_2;x,t)$ is not constant for $x\in [it^{2/3}, (i+1)t^{2/3}]$.
Then $\cE_2\subset \cup_{|i|<\lceil t^{0.01}\rceil} \cE_{2,i}$, and we need to bound each $\PP[\cE_{2,i}]$.
Using scaling and skew-shift invariance of the directed landscape, we get that $\PP[\cE_{2,i}]=\PP[\cE_{2,i}']$, with $\cE_{2,i}'$ being the event where
\[
x\mapsto \cL(t^{-1/6}x_1-it^{-1/2}s_1, t^{-1/4}s_1;x,t^{3/4})-\cL(t^{-1/6}x_2-it^{-1/2}s_2, t^{-1/4}s_2;x,t^{3/4})
\]
is not constant for $x\in [0, t^{1/2}]$.
Let $B=[-1,1]\times [-1,1]$.
Then for $t$ large enough, we would have $(t^{-1/6}x_1-it^{-1/2}s_1, t^{-1/4}s_1)$, $(t^{-1/6}x_2-it^{-1/2}s_2, t^{-1/4}s_2) \in B$ for any $|i|<\lceil t^{0.01}\rceil$.
So by Lemma \ref{lem:asy-coal}, we conclude that there is a constant $C>0$ such that 
\[
\PP[\cE_{2,i}']<Ct^{-1/4}e^{C|\log(t)|^{5/6}},
\]
for each $|i|<\lceil t^{0.01}\rceil$.
Thus we have $\PP[\cE_2] \le \sum_{|i|<\lceil t^{0.01}\rceil} \PP[\cE_{2,i}'] < Ct^{-1/4}e^{C|\log(t)|^{5/6}}(2\lceil t^{0.01}\rceil+1)$ when $t$ is large enough.
We then have that $\PP[\cE_2] \to 0$ as $t\to \infty$, and hence $\PP[\pi_{p_1}(t)=\pi_{p_2}(t)] \to 1$ as $t\to \infty$, and the conclusion follows.
\end{proof}

\begin{rem}
We choose to state Lemma \ref{lem:semi-inf-tree-fixed} in the current form as it will be convenient for a later application. In fact, via the ordering of geodesics and bounds on transversal fluctuations, one can prove that almost surely, the statement holds simultaneously for any starting points in $\R^2$.
\end{rem}

\section{Duality in the limit}\label{dualsec}
In this section we prepare the necessary framework of competition interfaces in the directed landscape. We will start from Lemma \ref{lem:dual-prelim}, i.e., the equality of law of the competition interface starting from the stationary initial condition and the semi-infinite geodesic in the prelimiting model of exponential LPP, and use Theorem \ref{thm:exp-to-dl} from \cite{DV21} to pass to the limit and conclude the same for the directed landscape.

Let $\cB$ be a standard two sided Brownian motion.
Recall from Section \ref{s:iop} that
\begin{align}
\cL^L(x,t)&=\sup_{y\le 0} \cL(y,0;x,t) + \cB(y),\\
\cL^R(x,t)&=\sup_{y\ge 0} \cL(y,0;x,t) + \cB(y),
\end{align}
for any $x\in\R, t>0$.
Using the quadrangle inequality (Lemma \ref{lem:DL-quad}), we can show that the function $x\mapsto \cL^L(x,t)-\cL^R(x,t)$ is non-increasing in $x$.
Indeed, for any $x_1<x_2$ and $y_1\le 0 \le y_2$, by Lemma \ref{lem:DL-quad} we have
\[
(\cL(y_1,0;x_1,t) + \cB(y_1)) + (\cL(y_2,0;x_2,t) + \cB(y_2)) \ge (\cL(y_1,0;x_2,t) + \cB(y_1)) + (\cL(y_2,0;x_1,t) + \cB(y_2)).
\]
By taking the supremum over $y_1\le 0 \le y_2$ in the right hand side, we have
\[
\cL^L(x_1,t) + \cL^R(x_2,t) \ge \cL^L(x_2,t) + \cL^R(x_1,t).
\]
We then let
\[\LV^\cB(t)=\inf\{x\in \R: \cL^L(x,t)-\cL^R(x,t) \le 0\},\]
as given in \eqref{leftrightprofile}.

The goal of this section is to establish Propositon \ref{prop:equal-dist-geo-comp}, i.e., $\LV^\cB$ and $\pi_{\boo}$ are equal in distribution.
Relying on the same statement for exponential LPP (Lemma \ref{lem:dual-prelim}), the Proposition \ref{prop:equal-dist-geo-comp} follows immediately from the following convergence in distribution. Recall the maps $\fG_n$ and $\Graph$ from Definition \ref{rangegraph}.

\begin{prop}  \label{prop:cov-semi-geo}
For any $0<h$, $\fG_n(\Gamma_\boo) \cap (\R\times [0,h]) \to \Graph(\pi_\boo) \cap (\R\times [0,h])$ in distribution in the Hausdorff topology, as $n\to\infty$.
\end{prop}

\begin{prop}  \label{prop:cov-competi}
For any $0<g<h$, $\fG_n(\Delta) \cap (\R\times [g,h]) \to \Graph(\LV^\cB) \cap (\R\times [g,h])$ in distribution in the Hausdorff topology, as $n\to\infty$.
\end{prop}
In Proposition \ref{prop:cov-competi} we need $g>0$ due to using Theorem \ref{thm:exp-to-dl}, where a compact set needs to be taken for uniform convergence.

It remains to prove Proposition \ref{prop:cov-semi-geo} and \ref{prop:cov-competi}, in the following two subsections respectively.

\subsection{Weak convergence of semi-infinite geodesics}

Let us reiterate from Section \ref{s:iop}, the general strategy to prove Proposition \ref{prop:cov-semi-geo}.
For simplicity and by scaling invariance, it suffices to prove the result for $h=1$.

Our proof would rely on the convergence of the exponential LPP to the directed landscape.
However, a technical issue arises since the convergence stated in Theorem \ref{thm:exp-to-dl} is only for compact sets. Thus to use the latter,  we will define a ``truncated'' version of the geodesic $\Gamma_\boo$, parameterized by some $H>0$. Theorem \ref{thm:exp-to-dl} then allows us to deduce a scaling convergence to a ``truncated'' version of $\pi_\boo$.
The final step would be to show that the truncated versions agree with the original versions of the geodesics with probabilities increasing to $1$ as $H\to \infty$. For the pre-limiting semi-infinite geodesics, this needs to be done uniformly in $n$. This last step can also be understood as proving tightness for the pre-limiting semi-infinite geodesics, i.e., to give a uniform bound on the transversal fluctuations.

We first truncate the semi-infinite geodesic $\pi_\boo$ in $\cL$, and the idea is to force $\pi_\boo(1)$ to stay in a compact interval $[-H, H]$.
Recall the Busemann function $\cM$ from Definition \ref{defn:DL-Buse}.
For each $H>0$, we define the \emph{$H-$geodesic in $\cL$} as follows.
We take $x_* = \argmax_{x\in[-H,H]} \cL(\boo;x,1) + \cM(x,1;0,1)$, and let $\pi_{\boo,H}= \pi_{(\boo;x_*,1)}$.

\begin{lemma}  \label{lem:H-geo-unique}
Almost surely the $H$-geodesic $\pi_{\boo,H}$ is unique.
\end{lemma}
This is by the uniqueness of the $\argmax$ of the sum of two weight profiles, at each rational time. Similar arguments have been used to prove uniqueness of finite geodesics in \cite{DOV}.
\begin{proof}
By Theorem \ref{thm:abs-bm-gen} and Lemma \ref{lem:buse-bm}, for any $0< t < 1$ the function $$y\mapsto \sup_{x\in [-H,H]} \cL(y, t; x, 1) + \cM(x,1;0,1)$$ is absolutely continuous with respect to a Brownian motion, on any compact set.
Moreover, on account of being a scaled parabolic Airy$_2$ process, the function $y \mapsto \cL(\boo; y, t)$ is also absolutely continuous with respect to a Brownian motion on any compact set.
Thus, almost surely, the function
\[
y \mapsto \cL(\boo; y, t) + \sup_{x\in [-H,H]} \cL(y, t; x, 1) + \cM(x,1;0,1)
\]
has a unique maximum in any compact interval. 
Since this function also has parabolic decay as $|y|\to\infty$ (by Lemma \ref{lem:DLbound}), it has a unique global maximum, which we denote as $y_{t,*}$. Then we must have $\pi_{\boo,H}(t)=y_{t,*}$ and hence almost surely the uniqueness of $\pi_{\boo,H}(t)$ holds for all rational $0<t<1$. By continuity of $\pi_{\boo,H}$ the conclusion follows.
\end{proof}
We next carry out the counterpart truncation for the pre-limiting semi-infinite geodesics.
Similarly as above, we define the \emph{exponential LPP $(n, H)$-geodesic} $\Gamma_{\boo,H,n}$ as follows.
Recalling the Busemann function $G(\cdot,\cdot)$ from Definition \ref{expbusemann}, let 
\[m_*[n] = \argmax_{m\in\Z, |m-n|\le {2^{5/3}H n^{2/3}}} T_{\boo,(m,n)} + G((m,n+1)),\]
which is almost surely unique (due to continuity of the distribution functions).
Let $\Gamma_{\boo,H,n}=\Gamma_{\boo,(m_*[n],n)}$.
Noting that the constant $2^{5/3}$ is a result of the scaling function in Section \ref{ssec:lpptodl}, we have the following convergence result.
\begin{lemma}   \label{lem:cov-semi-geo-H}
For any $H>0$, we have $\fG_n(\Gamma_{\boo,H,n}) \to \Graph(\pi_{\boo,H})$ in distribution in the Hausdorff topology.
\end{lemma}
\begin{proof}
We take the coupling from Theorem \ref{thm:exp-to-dl} of all $\cK_n$ and $\cL$, on the set $\{(x,s;y,t)\in\R^4: s<t\le 1\}\subset \R^4_\uparrow$.
We also consider the function $\cM_n:\R\to\R$, where
\[\cM_n(x) = 2^{-4/3}n^{-1/3}(G((\lfloor n+2^{5/3}n^{2/3}x \rfloor, n+1)) - G((n,n+1)) + 2(\lfloor n+2^{5/3}n^{2/3}x \rfloor -n)).\]
This is in the same scale as and converges to the function $x\mapsto\cM(x,1;0,1)$.
As stated in Section \ref{ssec:dualityexp}, along any down-right path, $G$ is a random walk with i.i.d. steps. Thus $\cM_n$ is a (rescaled) random walk with i.i.d. steps, and (using Skorokhod's representation theorem) we couple $\cM_n$ with $x\mapsto \cM(x,1;0,1)$ such that $\cM_n$ converges to the latter almost surely uniformly in each compact set.
We show that $\fG_n(\Gamma_{\boo,H,n}) \to \Graph(\pi_{\boo,H})$ almost surely under this coupling.

We work under the probability one event where $\pi_{\boo,H}$ is unique.
In particular, this implies that $\pi_{\boo,H}(1)$ is the unique maximum of the function $x\mapsto \cL(\boo;x,1) + \cM(x,1;0,1)$ in $[-H, H]$, and $\pi_{\boo,H}$ is the unique geodesic from $\boo$ to $(\pi_{\boo,H}(1), 1)$.

Using previous notations, $\fG_n((m_*[n],n))$ is a maximum of the function $x\mapsto \cK_n(\boo;x,1) + \cM_n(x)$ $\fG_n((m_*[n],n))$ in $[-H, H]$.
Since $x\mapsto \cK_n(\boo;x,1) + \cM_n(x)$ uniformly converges to $x\mapsto \cL(\boo;x,1) + \cM(x,1;0,1)$ in $[-H, H]$, we have that $\fG_n((m_*[n],n)) \to \pi_{\boo,H}(1)$.
Since $\Gamma_{\boo,H,n}$ is the geodesic from $\boo$ to $(m_*[n],n)$, $\fG_n(\Gamma_{\boo,H,n})$ is a maximal geodesic set of $\cK_n$ from $\boo$ to $\fG_n((m_*[n],n))$.
Thus we conclude that any subsequential limit of $\fG_n(\Gamma_{\boo,H,n})$ would be $\Graph(\pi)$ for some $\cL$ geodesic $\pi$ from $\boo$ to $(\pi_{\boo,H}(1), 1)$ (by Theorem \ref{thm:exp-to-dl}), thus we must have $\pi=\pi_{\boo,H}$, and the conclusion follows.
\end{proof}

We next prove tightness, which says that in the pre-limiting setting, the truncated semi-infinite geodesics are the same as the original semi-infinite geodesics with high probability.
This is achieved via estimates on the transversal fluctuations of semi-infinite geodesics, given by Lemma \ref{lem:trans-semi-inf}.

\begin{lemma}  \label{lem:semi-geo-tight}
There are constants $c,C>0$, such that for any $H>0$ and any $n$ large enough (depending on $H$), we have $\PP[\fG_n(\Gamma_{\boo,H,n}) \neq \fG_n(\Gamma_\boo) \cap (\R\times [0,1])] < Ce^{-cH}$.
\end{lemma}
\begin{proof}
Let $m_* \in \Z$ such that $(m_*,n), (m_*,n+1) \in \Gamma_\boo$.
Then for any $m\in\Z$, $m\neq m_*$, we must have
\[
T_{\boo,(m_*,n)} + G((m_*,n+1)) > T_{\boo,(m,n)} + G((m,n+1))
\]
since otherwise we would have $T_{\boo,(m_*,n)} + T_{(m_*,n+1),\bc} \le T_{\boo,(m,n)} + T_{(m,n+1),\bc}$ where $\bc$ is any point in $\Gamma_{(m_*,n+1)}\cap\Gamma_{(m,n+1)}$.
This implies that the path obtained by concatenating $\Gamma_{\boo,(m,n)}$ and $\Gamma_{(m,n+1),\bc}$ has a total weight no less than the total weight of the part of $\Gamma_\boo$ up to $\bc$.
However, from the definition of $\Gamma_\boo$ we have that the part of $\Gamma_\boo$ up to $\bc$ is the geodesic $\Gamma_{\boo,\bc}$.
Thus we get a contradiction with the assumption that there is a unique geodesic between $\boo$ and $\bc$. 

Then the event where $\fG_n(\Gamma_{\boo,H,n}) \neq \fG_n(\Gamma_\boo) \cap (\R\times [0,1])$ is equivalent to the event $|m_*-n|> 2^{5/3}Hn^{2/3}$.
By Lemma \ref{lem:trans-semi-inf} its probability is upper bounded by $Ce^{-cH}$ for some constants $c,C>0$, so the conclusion follows.
\end{proof}

The relation between $\pi_\boo$ and $\pi_{\boo,H}$ is more obvious, since we have that $\pi_{\boo,H}$ is the same as $\pi_\boo$ on $[0,1]$, when $H>|\pi_\boo(1)|$.
Now we can finish the proof of Proposition \ref{prop:cov-semi-geo}.
\begin{proof}[Proof of Proposition \ref{prop:cov-semi-geo}]
We just prove for the case where $h=1$, and the case of general $h>0$ follows similar arguments.
For any continuous function $f$ on the space of all compact subsets of $\R^2$ with $0\le f\le 1$, by Lemma \ref{lem:cov-semi-geo-H} we have
$\E[f(\fG_n(\Gamma_{\boo,H,n}))] \to \E[f(\Graph(\pi_{\boo,H}))]$ as $n\to\infty$ for any $H>0$.
Then we have 
\begin{multline*}
\limsup_{n\to\infty} \big|\E[f(\fG_n(\Gamma_\boo)\cap (\R\times[0,1]))] - \E[f(\Graph(\pi_\boo)\cap (\R\times [0,1]))]\big| \\ 
\le \limsup_{n\to\infty} \PP[\fG_n(\Gamma_{\boo,H,n}) \neq \fG_n(\Gamma_\boo) \cap (\R\times [0,1])] + \PP[\Graph(\pi_{\boo,H}) \neq \Graph(\pi_\boo) \cap (\R\times [0,1])]
\\
\le Ce^{-cH} + \PP[\Graph(\pi_{\boo,H}) \neq \Graph(\pi_\boo) \cap (\R\times [0,1])],
\end{multline*}
where the last inequality is by Lemma \ref{lem:semi-geo-tight}.
As $H\to\infty$, we have $Ce^{-cH}\to 0$ and $\PP[\Graph(\pi_{\boo,H}) \neq \Graph(\pi_\boo) \cap (\R\times [0,1])] \le \PP[|\pi_\boo(1)|<H] \to 0$, by Lemma \ref{lem:dl-semi-trans}.
So as $H\to\infty$ the last line decays to zero. Thus we have $\E[f(\fG_n(\Gamma_\boo)\cap (\R\times[0,1]))] \to \E[f(\Graph(\pi_\boo)\cap (\R\times [0,1]))]$, and the conclusion follows.
\end{proof}

We next carry out the counterpart convergence arguments for the competition interface.

\subsection{Weak convergence of the competition interface}
For Proposition \ref{prop:cov-competi}, we again just prove it for $h=1$, without loss of generality.
The proof follows a similar general framework as the proof of Proposition \ref{prop:cov-semi-geo}, since one encounters the same technical difficulty that the convergence stated in Theorem \ref{thm:exp-to-dl} only applies to compact sets. 
We similarly define a ``truncated'' version of  $\Delta$, and use Theorem \ref{thm:exp-to-dl} to deduce that its scaling limit is a ``truncated'' version of $\LV^\cB$.
As before, the last step entails showing that the truncated versions of $\Delta$ are the same as the original version with high probability. 
However, compared to the geodesics setting in the last subsection, now it is not immediate the truncated versions of $\LV^\cB$ are close (in total variation distance) to $\LV^\cB$, and this requires some further argument.

We first truncate $\LV^\cB$. The idea is to only use the boundary condition $\cB$ only in a compact interval.
For each $H>0$, we define the \emph{$H$-competition interface in $\cL$} as follows. 
For any $x\in\R, t>0$, we let
\[
\cL^L_H(x,t)=\sup_{-H \le y\le 0} \cL(y,0;x,t) + \cB(y),
\]
and
\[
\cL^R_H(x,t)=\sup_{0\le y \le H} \cL(y,0;x,t) + \cB(y).
\]
By arguments similar to the proof of the quadrangle inequality (Lemma \ref{lem:DL-quad}), or the arguments before Propositoin \ref{prop:cov-semi-geo}, we have that the function $x\mapsto \cL^L_H(x,t)-\cL^R_H(x,t)$ is non-increasing in $x$. So we can let 
\[\LV^\cB_H(t)=\inf\{x\in \R: \cL^L_H(x,t)-\cL^R_H(x,t) \le 0\}.\]
We can show that actually, almost surely, for each $t> 0$ there is exactly one $x\in\R$ with $\cL^L_H(x,t)=\cL^R_H(x,t)$.
\begin{lemma}  \label{lem:lim-compi-uni}
For any $H>0$, almost surely the following is true.
For any $t> 0$, $|\{x\in\R: \cL^L_H(x,t)-\cL^R_H(x,t)=0 \}|=1$.
Thus via continuity of $\cL$ and $\cB$, $\LV^\cB_H$ must be continuous.
\end{lemma}
The proof proceeds as follows: first we show that $\cL^L_H(x,t)\neq \cL^R_H(x,t)$, for any $x$ and $t>0$, when they are both rational.
Then if there is some $t>0$ and $x_1<x_2$ such that $\cL^L_H(x_1,t)-\cL^R_H(x_1,t)= \cL^L_H(x_2,t)-\cL^R_H(x_2,t)=0$, the difference $\cL^L_H-\cL^R_H$ would be zero in the segment connecting $(x_1, t)$ and $(x_2, t)$. We then show that  one can find a rational point $(x_*,t_*)$ where $\cL^L_H-\cL^R_H$ equals zero, which is already ruled out by the previous line.
Finally, using continuity of $\cL^L_H- \cL^R_H$, and that (for any $t> 0$) $\cL^L_H(x,t)- \cL^R_H(x,t)$ tends to $\infty$ (resp. $-\infty$) as $x\to-\infty$ (resp. $x\to\infty$), we conclude that there must exist some $x$, and hence exactly one such, with $\cL^L_H(x,t)=\cL^R_H(x,t)$.
\begin{proof}
We first show that for any fixed $x\in \R, t>0$, almost surely we have $\cL^L_H(x,t)\neq \cL^R_H(x,t)$.
Indeed, this would imply that the function $y\mapsto \cL(y,0;x,t)+\cB(y)$ has the same maximum in $[-H, 0]$ and $[0, H]$, and almost surely this cannot happen since $y\mapsto \cL(y,0;x_*,t)+\cB(y)$ in $[-H, H]$ is absolutely continuous with respect to a Brownian motion, by \cite{CH} or Theorem \ref{thm:airytail}. 
Thus, almost surely $\cL^L_H(x,t)\neq \cL^R_H(x,t)$, for any rational $x$ and $t>0$.

We next show that $|\{x\in\R: \cL^L_H(x,t)-\cL^R_H(x,t)=0 \}|\le 1$ for any (not necessarily rational) $t>0$.
For this, we use continuity of $\cL$ and geodesics, and basic inequalities of the directed landscape.

\begin{figure}
         \centering
\begin{tikzpicture}[line cap=round,line join=round,>=triangle 45,x=0.6cm,y=0.3cm]
\clip(-2,-1) rectangle (12,10);

\fill[line width=0.pt,color=cyan,fill=cyan,fill opacity=0.2]
(7,8) -- (3,8) -- (3.8,7) -- (3.6,6) -- (6.6,7) -- cycle;
\draw (-10,0) -- (30,0);
\draw [dashed] (-10,8) -- (30,8);
\draw [dashed] (-10,7) -- (30,7);

\draw [line width=1pt] (0.4,0) -- (2.5,1) -- (2.1,2) -- (4.2,3) -- (2.3,4) -- (4.8,5) -- (3.6,6) -- (3.8,7) -- (3,8);
\draw [line width=1pt] (9,0) -- (8.5,1) -- (8.7,2) -- (5.7,3) -- (5,4) -- (4.8,5) -- (3.6,6) -- (6.6,7) -- (7,8);

\draw [line width=1pt] (5.1,6.5) -- (4.8,7);

\draw [fill=uuuuuu] (0.4,0) circle (1.5pt);
\draw [fill=uuuuuu] (5,0) circle (1.5pt);
\draw [fill=uuuuuu] (9,0) circle (1.5pt);

\draw [fill=uuuuuu] (3,8) circle (1.5pt);
\draw [fill=uuuuuu] (7,8) circle (1.5pt);

\draw [fill=uuuuuu] (3.8,7) circle (1.5pt);
\draw [fill=uuuuuu] (4.8,7) circle (1.5pt);
\draw [fill=uuuuuu] (6.7,7) circle (1.5pt);
\begin{scriptsize}
\draw (0.4,0) node[anchor=north]{$(y_2,0)$};
\draw (9,0) node[anchor=north]{$(y_1,0)$};
\draw (5,0) node[anchor=north]{$\boo$};

\draw (3,8) node[anchor=south]{$(x_1, t)$};
\draw (7,8) node[anchor=south]{$(x_2, t)$};

\draw (4.8,7) node[anchor=south]{$(x_*, t_*)$};

\draw (3.8,7) node[anchor=north east]{$(x_{1,*}, t_*)$};
\draw (6.7,7) node[anchor=north west]{$(x_{2,*}, t_*)$};

\end{scriptsize}

\end{tikzpicture}
 \caption{An illustration of the arguments for $|\{x\in\R: \cL^L_H(x,t)-\cL^R_H(x,t)=0 \}|\le 1$ in the proof of Lemma \ref{lem:lim-compi-uni}: suppose that $\cL^L_H(x_1,t)-\cL^R_H(x_1,t)=\cL^L_H(x_2,t)-\cL^R_H(x_2,t)=0$, the difference $\cL^L_H-\cL^R_H$ would be zero in the whole blue region, which must include a rational point $(x_*,t_*)$.}
\label{fig:crox}
\end{figure}
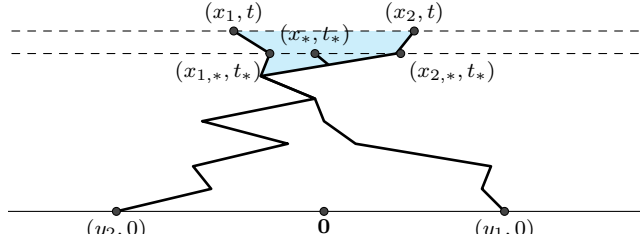

We argue by contradiction (see Figure \ref{fig:crox} for an illustration).
Suppose that there is some $t>0$ and $x_1<x_2$ such that $\cL^L_H(x_1,t)-\cL^R_H(x_1,t)= \cL^L_H(x_2,t)-\cL^R_H(x_2,t)=0$.
Take any $y_1=\argmax_{y\in[0,H]}\cL(y,0;x_1,t)+\cB(y)$.
Take any geodesic $\pi_{(y_1,0;x_1,t)}$.
Similarly, take any $y_2=\argmax_{y\in[-H,0]}\cL(y,0;x_2,t)+\cB(y)$, and any geodesic $\pi_{(y_2,0;x_2,t)}$.
Then we can find a rational $0<t_*<t$ such that $x_{1,*}<x_{2,*}$, where $x_{1,*}=\pi_{(y_1,0;x_1,t)}(t_*)$ and $x_{2,*}=\pi_{(y_2,0;x_2,t)}(t_*)$.
We then have that
\begin{align*}
\cL^L_H(x_{1,*},t_*) &=
\sup_{-H\le y \le 0} \cL(y,0;x_{1,*},t_*) + \cB(y) \\ &\le
\sup_{-H\le y \le 0} \cL(y,0;x_1,t_*) - \cL(x_{1,*},t_*;x_1,t) + \cB(y) 
\\ & = \cL^L_H(x_1,t) - \cL(x_{1,*},t_*;x_1,t) \\ &= \cL^R_H(x_1,t) - \cL(x_{1,*},t_*;x_1,t) \\ &= \cL^R_H(x_{1,*},t_*),
\end{align*}
where the inequality is by the reverse triangle inequality of the directed landscape, or the composition law \eqref{eq:DL-compo}.
Thus we have $\cL^L_H(x_{1,*},t_*) - \cL^R_H(x_{1,*},t_*) \le 0$.
Similarly we have $\cL^L_H(x_{2,*},t_*) - \cL^R_H(x_{2,*},t_*) \ge 0$.

As above mentioned, by arguments similar to Lemma \ref{lem:DL-quad}, or those at the beginning of this section, we have that the function $x\mapsto \cL^L_H(x,t_*) - \cL^R_H(x,t_*)$ is non-increasing.
Now we take any rational $x_*$ such that $x_{1,*}\le x_* \le x_{2,*}$.
Then we have
\[
0
\le \cL^L_H(x_{2,*},t_*) - \cL^R_H(x_{2,*},t_*)
\le \cL^L_H(x_{*},t_*) - \cL^R_H(x_{*},t_*)
\le \cL^L_H(x_{1,*},t_*) - \cL^R_H(x_{1,*},t_*) \le 0.
\]
This implies that $\cL^L_H(x_*,t_*) = \cL^R_H(x_*,t_*)$.
As we've shown that this almost surely cannot happen, we get a contradiction.
We thus conclude that $|\{x\in\R: \cL^L_H(x,t)-\cL^R_H(x,t)=0 \}|\le 1$ for any $t>0$.

Finally, by Lemma \ref{lem:DLbound} and growth estimates of $\cB$, we have that for any $t>0$, $\cL^L_H(x,t)$ has parabolic decay as $x\to\infty$.
Also, for $x\le 0$ we have
\[
\cL^L_H(x,t) \ge \cB(x) + \cL(x,0;x,t),
\]
so by Lemma \ref{lem:DLbound} and growth estimates of $\cB$ we have that $\cL^L_H(x,t) \ge -|x|^{1/2}\log^3(|x|)$ for all $x$ small enough.
Similalry, $\cL^R_H(x,t)$ has parabolic decay as $x\to-\infty$, and $\cL^R_H(x,t)\ge -|x|^{1/2}\log^3(|x|)$ for all $x$ large enough.
Thus we conclude that $\lim_{x\to\infty}\cL^L_H(x,t)- \cL^R_H(x,t)=-\infty$ and $\lim_{x\to-\infty}\cL^L_H(x,t)- \cL^R_H(x,t)=\infty$. With continuity of $\cL^L_H- \cL^R_H$ we get that $|\{x\in\R: \cL^L_H(x,t)-\cL^R_H(x,t)=0 \}|\ge 1$, so the conclusion follows.
\end{proof}

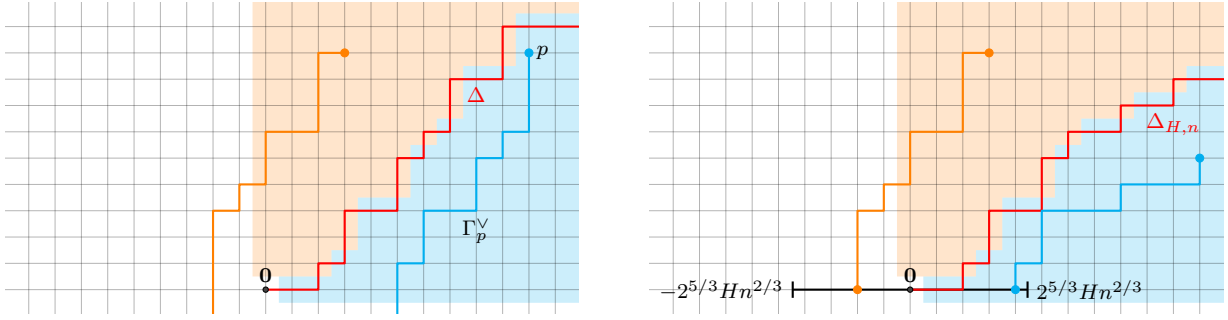
\begin{figure}
     \centering
     \begin{subfigure}[t]{0.45\textwidth}
         \centering
\begin{tikzpicture}[line cap=round,line join=round,>=triangle 45,x=3.5cm,y=3.5cm]
\clip(-1.39,-0.19) rectangle (0.79,.99);

\fill[line width=0.pt,color=orange,fill=orange,fill opacity=0.2]
(-0.45,-0.05) -- (-0.15,-0.05) -- (-0.15,0.05)-- (-0.05,0.05) --(-0.05,0.25) --(0.15,0.25) --(0.15,0.45) -- (0.25,0.45) --(0.25,0.55) -- (0.35,0.55) -- (0.35,0.75) -- (0.55,0.75) -- (0.55,0.95) -- (0.85,0.95) -- (0.85,1) -- (-0.45,1) -- cycle;
\fill[line width=0.pt,color=cyan,fill=cyan,fill opacity=0.2]
(-0.35,-0.15) -- (-0.35,-0.05) -- (-0.15,-0.05) -- (-0.15,0.05)-- (-0.05,0.05) --(-0.05,0.25) --(0.15,0.25) --(0.15,0.45) -- (0.25,0.45) --(0.25,0.55) -- (0.35,0.55) -- (0.35,0.75) -- (0.55,0.75) -- (0.55,0.95) -- (0.85,0.95) -- (0.85,-0.15) -- cycle;

\foreach \i in {-20,...,30}
{
\draw [line width=.1pt, opacity=0.3] (-1.5,\i/10) -- (2.6,\i/10);
\draw [line width=.1pt, opacity=0.3] (\i/10,-0.2) -- (\i/10,2.2);
}

\draw [thick] [red] plot coordinates {(-0.4,-0.1) (-0.2,-0.1) (-0.2,0.) (-0.1,0.) (-0.1,0.2) (0.,0.2) (0.1,0.2) (0.1,0.3) (0.1,0.4) (0.2,0.4) (0.2,0.5) (0.3,0.5) (0.3,0.7) (0.5,0.7) (0.5,0.9) (0.9,0.9)};

\draw [orange] [fill=orange] (-0.1,0.8) circle (1.5pt);
\draw [thick] [orange] plot coordinates {(-0.6,-0.2) (-0.6,0.2) (-0.5,0.2) (-0.5,0.3) (-0.4,0.3) (-0.4,0.5) (-0.2,0.5) (-0.2,0.8) (-0.1,0.8)};

\draw [cyan] [fill=cyan] (0.6,0.8) circle (1.5pt);
\draw [thick] [cyan] plot coordinates {(0.1,-0.2) (0.1,0.) (0.2,0.) (0.2,0.2) (0.4,0.2) (0.4,0.4) (0.5,0.4) (0.5,0.5) (0.6,0.5) (0.6,0.8)};
\draw [fill=uuuuuu] (-0.4,-0.1) circle (1.0pt);

\begin{scriptsize}
\draw (-0.4,-0.1) node[anchor=south]{$\boo$};
\draw [red] (0.4,0.7) node[anchor=north]{$\Delta$};
\draw (0.6,0.8) node[anchor=west]{$p$};
\draw (0.4,0.2) node[anchor=north]{$\Gamma_p^\vee$};
\end{scriptsize}

\end{tikzpicture}
     \end{subfigure}
     \hfill
     \begin{subfigure}[t]{0.5\textwidth}
         \centering
\begin{tikzpicture}[line cap=round,line join=round,>=triangle 45,x=3.5cm,y=3.5cm]
\clip(-1.39,-0.19) rectangle (0.79,.99);

\fill[line width=0.pt,color=orange,fill=orange,fill opacity=0.2]
(-0.45,-0.05) -- (-0.15,-0.05) -- (-0.15,0.05)-- (-0.05,0.05) --(-0.05,0.25) --(0.05,0.25) --(0.15,0.25) --(0.15,0.45) --(0.25,0.45) --(0.25,0.55) -- (0.45,0.55) -- (0.45,0.65) -- (0.65,0.65) -- (0.65,0.75) -- (0.85,0.75) -- (0.85,1) -- (-0.45,1) -- cycle;
\fill[line width=0.pt,color=cyan,fill=cyan,fill opacity=0.2]
(-0.35,-0.15) -- (-0.35,-0.05) -- (-0.15,-0.05) -- (-0.15,0.05)-- (-0.05,0.05) --(-0.05,0.25) --(0.05,0.25) --(0.15,0.25) --(0.15,0.45) --(0.25,0.45) --(0.25,0.55) -- (0.45,0.55) -- (0.45,0.65) -- (0.65,0.65) -- (0.65,0.75) -- (0.85,0.75) -- (0.85,-0.15) -- cycle;

\foreach \i in {-20,...,30}
{
\draw [line width=.1pt, opacity=0.3] (-1.5,\i/10) -- (2.6,\i/10);
\draw [line width=.1pt, opacity=0.3] (\i/10,-0.2) -- (\i/10,2.2);
}

\draw [|-|] [thick] (-0.85,-0.1) -- (0.05,-0.1);

\draw [thick] [red] plot coordinates {(-0.4,-0.1) (-0.2,-0.1) (-0.2,0.) (-0.1,0.) (-0.1,0.2) (0.,0.2) (0.1,0.2) (0.1,0.3) (0.1,0.4) (0.2,0.4) (0.2,0.5) (0.3,0.5) (0.4,0.5) (0.4,0.6) (0.6,0.6) (0.6,0.7) (0.9,0.7) (0.9,0.9) (1.,0.9) (1.2,0.9) (1.2,1.) (1.3,1.) (1.3,1.3) (1.4,1.3) (1.4,1.5) (1.5,1.5) (1.5,1.7) (1.9,1.7) (1.9,1.8) (2.1,1.8) (2.1,2.) (2.3,2.) (2.3,2.1)  (2.4,2.1)};

\draw [fill=uuuuuu] (-0.4,-0.1) circle (1.0pt);
\draw [orange] [fill=orange] (-0.1,0.8) circle (1.5pt);
\draw [orange] [fill=orange] (-0.6,-0.1) circle (1.5pt);
\draw [thick] [orange] plot coordinates {(-0.6,-0.1) (-0.6,0.2) (-0.5,0.2) (-0.5,0.3) (-0.4,0.3) (-0.4,0.5) (-0.2,0.5) (-0.2,0.8) (-0.1,0.8)};

\draw [cyan] [fill=cyan] (0.7,0.4) circle (1.5pt);
\draw [cyan] [fill=cyan] (0.,-0.1) circle (1.5pt);
\draw [thick] [cyan] plot coordinates {(0.,-0.1) (0.,0.) (0.1,0.) (0.1,0.2) (0.4,0.2) (0.4,0.3) (0.7,0.3) (0.7,0.4)};

\begin{scriptsize}
\draw (-0.4,-0.1) node[anchor=south]{$\boo$};
\draw [red] (0.6,0.6) node[anchor=north]{$\Delta_{H,n}$};

\draw (0.05,-0.1) node[anchor=west]{$2^{5/3}Hn^{2/3}$};
\draw (-0.85,-0.1) node[anchor=east]{$-2^{5/3}Hn^{2/3}$};
\end{scriptsize}
\end{tikzpicture} 
     \end{subfigure}
        \caption{A comparison between the competition interfaces $\Delta$ and $\Delta_{H,n}$: for $\Delta$ we consider semi-infinite geodesics in the $(-1,-1)$-direction; while for $\Delta_{H,n}$, we consider the optimal paths in the $(-1,-1)$-direction passing through $[-2^{5/3}Hn^{2/3}, 2^{5/3}Hn^{2/3}] \times \{0\}$.}
        \label{fig:49}
\end{figure}

We also define the \emph{exponential LPP $(n, H)$-competition interface} $\Delta_{H,n}$ as follows.
Recall that $G^\vee$ is the $(-1,-1)$ direction Busemann function in exponential LPP, defined in Section \ref{ssec:dualityexp}. Analogous to the un-truncated definitions, we split $\Z_{\ge 0}^2\setminus\{\boo\}$ into two parts $P_{H,n,1}$ and $P_{H,n,2}$, as follows.
Let $\N\times \{0\} \subset P_{H,n,1}$; and for any $p\in \Z_{\ge 0}\times \N$, if we let $$m_*=\argmax_{m\in\Z, |m|\le 2^{5/3}H n^{2/3}} G^\vee(m,0)+T_{(m,1),p}$$ (which is almost surely unique), then $p\in P_{H,n,1}$ if and only if $m_*>0$.
We again note that the constant $2^{5/3}$ is a result of the scaling in Section \ref{ssec:lpptodl}.
Let $\Delta_{H,n}$ be the boundary of $P_{H,n,1}$ and $P_{H,n,2}$: for any $p\in\Z_{\ge 0}^2$, we let $p\in \Delta_{H,n}$ if and only if $p+(1,0)\in P_{H,n,1}$, and $p+(0,1)\in P_{H,n,2}$.
See the right panel in Figure \ref{fig:49} for an illustration.

We next prove convergence of these $(n, H)$-competition interfaces to the truncated competition interface in the directed landscape.
\begin{lemma}   \label{lem:cov-comp-H}
For any $H>0$ and $0<g<1$, we have $\fG_n(\Delta_{H,n})\cap ([-H, H]\times [g,1]) \to \Graph(\LV^\cB_H)\cap ([-H, H]\times [g,1])$ in distribution in the Hausdorff topology.
\end{lemma}
We note that this where the requirement of $g>0$ in Proposition \ref{prop:cov-competi} comes from.

In the proof we work under the coupling from Theorem \ref{thm:exp-to-dl}, and show almost sure convergence.
Note that (compared to Lemma \ref{lem:cov-semi-geo-H}), here we only consider the sets $\fG_n(\Delta_{H,n})$ and $\Graph(\LV^\cB_H)$ in $[-H, H]\times [0,1]$ rather than $\R\times [0,1]$. This can be viewed as a further truncation, and is also due to the compactness limitation of Theorem \ref{thm:exp-to-dl}.
\begin{proof}
We take the coupling from Theorem \ref{thm:exp-to-dl} of all $\cK_n$ and $\cL$ on the set $\{(x,s;y,t)\in\R^4: 0< s<t\}\subset \R^4_\uparrow$.
We also consider the function $\cB_n:\R\to\R$, where
\[\cB_n(x) = 2^{-4/3}n^{-1/3}(G^\vee((\lceil 2^{5/3}n^{2/3}x \rceil, 0)) - 2\lceil 2^{5/3}n^{2/3}x \rceil).\]
This is a discrete version of the Brownian motion $\cB$.
Using the properties of $G^\vee$ stated in Section \ref{ssec:dualityexp}, $\cB_n$ is a rescaled random walk with i.i.d. steps, and is independent of $\cK_n$ on $\{(x,s;y,t)\in\R^4: 0< s<t\}$.
We can couple $\cB_n$ with $\cB$ such that $\cB_n\to\cB$ almost surely, uniformly on each compact set, by Skorokhod's representation theorem.

We show that under this coupling, almost surely, we have  $\fG_n(\Delta_{H,n})\cap ([-H,H]\times [g,1]) \to \Graph(\LV^\cB_H)\cap ([-H,H]\times [g,1])$.

For each $x\in\R$ and $t>0$, we let
\[
\cL^L_{H,n}(x,t)=\sup_{-H \le y\le 0} \lim_{s\searrow 0}\cK_n(y,s;x,t) + \cB_n(y),
\]
and
\[
\cL^R_{H,n}(x,t)=\sup_{0\le y \le H} \lim_{s\searrow 0}\cK_n(y,s;x,t) + \cB_n(y).
\]
We use the limit $\lim_{s\searrow 0}\cK_n(y,s;x,t)$ (instead of $\cK_n(y,0;x,t)$) due to the rounding setup in Section \ref{ssec:lpptodl}:
under the rouding setup, the limit depends on the values of $\cK_n$ on $\{(x,s;y,t)\in\R^4: 0< s<t\}$, thus is independent of $\cB_n$.

From the convergence of $\cK_n$ to $\cL$ (Theorem \ref{thm:exp-to-dl}) and the convergence of $\cB_n$ to $\cB$, we have $\cL^L_{H,n} \to \cL^L_H$ and $\cL^R_{H,n} \to \cL^R_H$, both uniformly in the set $[-2H, 2H] \times [g/2, 2]$.
From the above constructions and the definitons of $P_{H,n,1}$ and $P_{H,n,1}$, we have that
\[
P_{H,n,1}=\{p\in \Z_{\ge 0}^2\setminus\{\boo\}: \cL^L_{H,n}(\fG_n(p)) - \cL^R_{H,n}(\fG_n(p)) < 0\},
\]
and
\[
P_{H,n,2}=\{p\in \Z_{\ge 0}^2\setminus\{\boo\}: \cL^L_{H,n}(\fG_n(p)) - \cL^R_{H,n}(\fG_n(p)) > 0\}.
\]
Further, from the definition of $\Graph$, for any $(x,t) \in (-H, H) \times (g, 1)$, we have that $\LV^\cB_H(t)=x$ if $(x,t) \in \Graph(\LV^\cB_H)$. Take any $\delta >0$ small enough (depending on $g$, $H$, and $(x,t)$). By Lemma \ref{lem:lim-compi-uni} we have $\cL^L_H(x-\delta/2,t) > \cL^R_H(x-\delta/2,t)$ and $\cL^L_H(x+\delta/2,t) < \cL^R_H(x+\delta/2,t)$.
By continuity of $\cL$, thus of $\cL^L_H$ and $\cL^R_H$, and uniform convergence, we have that for any $n$ large enough, we can find $(x_{L,n}, t_n)\in \fG_n(P_{H,n,2})$ and $(x_{R,n}, t_n)\in \fG_n(P_{H,n,1})$, in the $\delta/3$ neighborhood of $(x-\delta/2, t)$ and $(x+\delta/2, t)$ respectively.
Now consider the $\delta$ neighborhood of $(x,t)$: its pre-image under $\fG_n$ must contain a lattice path from $\fG_n^{-1}((x_{L,n}, t_n))$ and $\fG_n^{-1}((x_{L,n}, t_n))$.
Since $\Delta_{H,n}$ is the boundary of $P_{H,n,1}$ and $P_{H,n,2}$, the lattice path must intersect $\Delta_{H,n}$, so there must be a vertex of $\fG_n(\Delta_{H,n})$ in the $\delta$ neighborhood of $(x,t)$.
Thus we conclude that $(x,t)$ is in any subsequential limit of $\fG_n(\Delta_{H,n})\cap ([-H, H]\times [g,1])$.

If $(x,t) \not\in \Graph(\LV^\cB_H)$, by Lemma \ref{lem:lim-compi-uni}, we have $\cL^L_H(x,t) \neq \cL^R_H(x,t)$ and without loss of generality we assume that $\cL^L_H(x,t) < \cL^R_H(x,t)$.
Then by continuity of $\cL^L_H$ and $\cL^R_H$, and uniform convergence, there is a neighborhood of $(x,t)$, that is disjoint from $\fG_n(P_{H,n,1})$ for any $n$ large enough.
Thus $(x,t)$ cannot be in any subsequential limit of $\fG_n(\Delta_{H,n})\cap ([-H, H]\times [g,1])$.
This completes the proof. \end{proof}

We next prove that the competition interfaces are with high probability the same as the truncated versions.
We start with the pre-limiting ones, and the arguements are via estimates on transversal fluctuations of semi-infinite geodesics in the $(-1, -1)$ direction.

\begin{figure}
         \centering
\begin{tikzpicture}[line cap=round,line join=round,>=triangle 45,x=0.4cm,y=0.4cm]
\clip(-4,-1) rectangle (16,10);

\fill[line width=0.pt,color=cyan,fill=cyan,fill opacity=0.2]
(-2,0) -- (6,8) -- (10,8) -- (2,0) -- cycle;

\fill[line width=0.pt,color=blue,fill=blue,fill opacity=0.2]
(6,0) -- (14,8) -- (10,8) -- (2,0) -- cycle;

\draw (-10,0) -- (30,0);
\draw (-10,8) -- (30,8);

\draw [line width=1pt] (0.4,0) -- (2.5,1) -- (2.7,2) -- (3.1,3) -- (4.9,4) -- (5.3,5) -- (6.6,6) -- (7,7) -- (8.5,8);
\draw [line width=1pt] (5,0) -- (5.1,1) -- (7.6,2) -- (7.9,3) -- (8.5,4) -- (8.7,5) -- (9.1,6) -- (10.6,7) -- (11.5,8);

\draw [line width=1pt,color=red] (2,0) -- (2.9,1) -- (5.9,2) -- (7.1,3) -- (7.2,4) -- (8.1,5) -- (8.8,6) -- (9.3,7) -- (10.5,8) -- (11.1,9) -- (12.4,10);

\draw [fill=uuuuuu] (8.5,8) circle (1.5pt);
\draw [fill=uuuuuu] (11.5,8) circle (1.5pt);
\draw [fill=uuuuuu] (2,0) circle (1.5pt);

\begin{scriptsize}
\draw (2,0) node[anchor=north]{$\boo$};
\draw (8.5,8) node[anchor=south]{$p_1$};
\draw (11.5,8) node[anchor=south]{$p_2$};

\draw (2,4) node[anchor=west]{$U_1$};
\draw (10,4) node[anchor=east]{$U_2$};

\draw (7,7) node[anchor=east]{$\Gamma^\vee_{p_1}$};
\draw (10.6,7) node[anchor=west]{$\Gamma^\vee_{p_2}$};

\draw (7.2,4) node[anchor=east,color=red]{$\Delta$};

\end{scriptsize}

\end{tikzpicture}
 \caption{An illustration of the proof of Lemma \ref{lem:comp-tight}: if $\Gamma^\vee_{p_1}$ and $\Gamma^\vee_{p_2}$ in $\Z\times \llbracket 0,n\rrbracket$ are contained in the parallelograms $U_1$ and $U_2$, respectively, the competition interface $\Delta$ (in $\Z\times \llbracket 0,n\rrbracket$) should be between them and equal $\Delta_{H,n}$.}
\label{fig:comptight}
\end{figure}
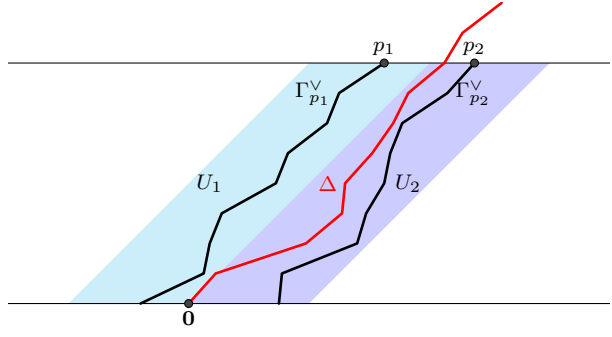

\begin{lemma}  \label{lem:comp-tight}
There are constants $c,C>0$, such that for any $H>0$ and any $n$ large enough (depending on $H$), we have $\PP[\fG_n(\Delta_{H,n})\cap ([-H,H]\times [0,1]) \neq \fG_n(\Delta) \cap (\R\times [0,1])] < Ce^{-cH}$.
\end{lemma}
\begin{proof}
Let $p_1 = (n-\lfloor n^{2/3}H \rfloor , n)$ and $p_2 = (n+\lfloor n^{2/3}H \rfloor , n)$,
and $U_1$ be the parallelogram whose vertices are $(-2^{5/3}n^{2/3}H, 0)$, $\boo$, and $(n-2^{5/3}n^{2/3}H, n)$, $(n, n)$,
and $U_2 = U_1+(2^{5/3}n^{2/3}H, 0)$. See Figure \ref{fig:comptight} for an illustration of these objects.

We now consider the events where $\Gamma^\vee_{p_1} \cap (\Z\times \llbracket 0,n\rrbracket) \subset U_1$ and $\Gamma^\vee_{p_2} \cap (\Z\times \llbracket 0,n\rrbracket) \subset U_2$.
By Lemma \ref{lem:trans-semi-inf} these happen with probability at least $1-Ce^{-cH}$, for some constants $c,C>0$.

It now suffices to show that these events imply that $\fG_n(\Delta_{H,n})\cap ([-H,H]\times [0,1]) = \fG_n(\Delta) \cap (\R\times [0,1])$.
Assuming that the events occur,
let $m_{1,*}, m_{2,*}\in \Z$ such that $(m_{1,*},0), (m_{1,*},1) \in \Gamma^\vee_{p_1}$, and $(m_{2,*},0), (m_{2,*},1) \in \Gamma^\vee_{p_2}$.
Then we have
\[
m_{1,*} = \argmax_{m\in\Z} G^\vee((m,0))+T_{(m,1),p_1},
\quad
m_{2,*} = \argmax_{m\in\Z} G^\vee((m,0))+T_{(m,1),p_2}.
\]
Then we would have that $-2^{5/3}n^{2/3}H \le m_{1,*}\le 0$ and $0<m_{2,*}\le 2^{5/3}n^{2/3}H$.
We now take $p=(x,y) \in \Z_{\ge 0}\times \llbracket 0,n\rrbracket \setminus \{\boo\}$ that is between $\Gamma^\vee_{p_1}$ and $\Gamma^\vee_{p_2}$;
i.e. there are $x_1\le x \le x_2$ such that $(x_1,y) \in \Gamma^\vee_{p_1}$ and $(x_2,y) \in \Gamma^\vee_{p_2}$.

Let $m_*\in \Z$ such that $(m_*,0), (m_*,1) \in \Gamma^\vee_p$. Due to the tree structure of the $(-1, -1)$-direction semi-infinite geodesics, 
 we must have that $m_{1,*} \le m_*\le m_{2,*}$, and thus $m_*\in [-2^{5/3}H n^{2/3}, 2^{5/3}H n^{2/3}]$.
This implies that $p\in P_{H,n,1}$ if and only if $p \in P_1$, and that $p\in P_{H,n,2}$ if and only if $p \in P_2$.
Since $\Gamma^\vee_{p_1} \cap \Z_{\ge 0}^2\setminus \{\boo\} \subset P_1 \cap U_1$ and $\Gamma^\vee_{p_2} \cap \Z_{\ge 0}^2\setminus \{\boo\} \subset P_2 \cap U_2$, we conclude that $\Delta_{H,n} = \Delta \cap (\Z\times \llbracket 0,n\rrbracket)\subset U_1\cup U_2$, which implies that $\fG_n(\Delta_{H,n})\cap ([-H,H]\times [0,1]) = \fG_n(\Delta) \cap (\R\times [0,1])$.
\end{proof}

As we have pointed out, unlike for semi-infinite geodesics, it is not immediate that $\LV^\cB_H$ becomes $\LV^\cB$ when $H$ is large enough. We prove this result now.
\begin{lemma}  \label{lem:comp-tight-lim}
Almost surely there is a (random) $H_0>0$, such that $\Graph(\LV^\cB_H) \cap ([-H,H]\times [0,1])= \Graph(\LV^\cB) \cap (\R\times [0,1])$ for any $H>H_0$.
\end{lemma}
\begin{proof}
We let $C$ denote a large constant in this proof, whose value may change from line to line.
Let $R_1$ be given by Lemma \ref{lem:DLbound}, and let $R_2 = \sup_x |\cB(x)| - |x|$. Let $R = R_1 \vee R_2 \vee C$.
For any $x, y \in \R$ and $0< t \le 1$, suppose that $|y|>10|x| \vee 10$, then we have
\begin{align*}
\cL(0,0;x,t) - \cL(y,0;x,t) - \cB(y) >& 
-\frac{x^2}{t} + \frac{(x-y)^2}{t} - |y| -R\\
&- 2Rt^{1/3} \log^{4/3}(2(|x|+|y|+3)^{3/2}/t) \log^{2/3}(|x|+|y|+3)
\\
\ge &
\frac{y^2-2xy}{t} - |y| - R - 2R\log^2(2(|x|+|y|+3)^{3/2}) \\
> &
\frac{y^2}{2} - 5R \log^2(10(|y|+1)).
\end{align*}
Since the last line is positive if $y$ is large enough, 
this implies that we can find some (random) $H_1>0$, such that for any $H>10H_1$, $|x|<H_1$ and $0<t\le 1$, we have $\cL^L_H(x,t) = \cL^L(x,t)$ and $\cL^R_H(x,t) = \cL^R(x,t)$.
We also have that for any $x, y\in \R$ and $0<t\le 1$, 
\begin{align*}
\cL(x,0;x,t) + \cB(x) - \cL(y,0;x,t) - \cB(y) >& 
\frac{(x-y)^2}{t} - |x| - |y| -2R\\
&- 2Rt^{1/3} \log^{4/3}(2(2|x|+|y|+3)^{3/2}/t) \log^{2/3}(2|x|+|y|+3)
\\
\ge &
\frac{(x-y)^2}{t} - |x| - |y| - 2R - 2R\log^2(2(2|x|+|y|+3)^{3/2}) \\
> &
\frac{(x-y)^2}{t} - |x| - |y| - 5R \log^2(10(|x|+|y|+1)).
\end{align*}
This implies that we can find some (random) $H_2>0$, such that $\cL^L(x,t) < \cL^R(x,t)$ for any $x>H_2$, $0<t\le 1$; and $\cL^L(x,t) > \cL^R(x,t)$ for any $x<-H_2$, $0<t\le 1$.
By taking $H_0=10H_1\vee H_2$ the conclusion follows.
\end{proof}

Now we can deduce Proposition \ref{prop:cov-competi}, and its proof is similar to the proof of Proposition \ref{prop:cov-semi-geo}.
\begin{proof}[Proof of Proposition \ref{prop:cov-competi}]
Again we just prove for the case where $h=1$, and the case of general $h>0$ follows similar arguments.
For any continuous function $f$ on the space of all compact subsets of $\R^2$ with $0\le f\le 1$, by Lemma \ref{lem:cov-comp-H} we have
$$\E[f(\fG_n(\Delta_{H,n})\cap ([-H,H]\times [g,1]))] \to \E[f(\Graph(\LV^\cB_H)\cap ([-H,H]\times [g,1]))]$$ as $n\to\infty$ for any $H>0$.
Then by Lemma \ref{lem:comp-tight} we have 
\begin{multline*}
\limsup_{n\to\infty} |\E[f(\fG_n(\Delta)\cap (\R\times [g,1]))] - \E[f(\Graph(\LV^\cB)\cap (\R\times [g,1]))]| \\
\le Ce^{-cH} + \PP[\Graph(\LV^\cB_H)\cap ([-H,H]\times [g,1])\neq \Graph(\LV^\cB) \cap (\R\times [g,1])].
\end{multline*}
As $H\to\infty$ the right hand side decays to zero, by Lemma \ref{lem:comp-tight-lim}. Thus we have $\E[f(\fG_n(\Delta)\cap (\R\times [g,1]))] \to \E[f(\Graph(\LV^\cB)\cap (\R\times [g,1]))]$, and the conclusion follows.
\end{proof}

With the above, somewhat long, sequence of preparatory results, we can finally delve in to the proof of Theorem \ref{thm:2dhaus}.

\section{Hausdorff dimension of $\NC$: Proof of Theorem \ref{thm:2dhaus}}  \label{sec:2ddphau}

Our first observation is that for the difference profile $\cD$, being constant in a 2D neighborhood is equivalent to being constant in the space direction. Recall the sets $\hNC_t$ and $\NC$ from the statement of Theorem \ref{thm:2dhaus}.
\begin{lemma}  \label{lem:equivD}
Almost surely, we have $\NC\cap(\R\times\{t\})=\hNC_t$ for any $t>0$.
\end{lemma}

\begin{figure}[hbt!]
    \centering
\begin{tikzpicture}[line cap=round,line join=round,>=triangle 45,x=0.6cm,y=0.45cm]
\clip(-2,-1) rectangle (12,10);

\fill[line width=0.pt,color=blue,fill=blue,fill opacity=0.2]
(6,8) -- (7,8) -- (6.5,7) -- (3.6,6) -- (4.5,7) -- (3,8) -- (5,8) -- (5,8.5) -- (6,8.5) -- cycle;

\draw (-10,0) -- (30,0);
\draw [dashed] (-10,8) -- (30,8);

\draw [line width=1pt] (0,0) -- (2.5,1) -- (1.7,2) -- (4.2,3) -- (2.3,4) -- (4.8,5) -- (3.6,6) -- (4.5,7) -- (3,8);
\draw [line width=1pt] (10,0) -- (8.5,1) -- (8.7,2) -- (5.7,3) -- (5,4) -- (4.8,5) -- (3.6,6) -- (6.5,7) -- (7,8);

\draw [line width=0.5pt,color=blue] (4.5,7) -- (4.8,7.5) -- (4.35,8) -- (5.15,8.2);
\draw [line width=0.5pt,color=blue] (5.05,6.5) -- (5.1,7) -- (6.3,7.5);

\draw [fill=uuuuuu] (0,0) circle (1.5pt);
\draw [fill=uuuuuu] (10,0) circle (1.5pt);

\draw [fill=uuuuuu] (3,8) circle (1.5pt);
\draw [fill=uuuuuu] (5.5,8) circle (1.5pt);
\draw [fill=uuuuuu] (7,8) circle (1.5pt);
\draw [blue] [fill=blue] (5.15,8.2) circle (1.5pt);
\draw [blue] [fill=blue] (6.3,7.5) circle (1.5pt);

\begin{scriptsize}
\draw (0,0) node[anchor=north]{$(-1,0)$};
\draw (10,0) node[anchor=north]{$(1,0)$};

\draw (3,8) node[anchor=north east]{$(x_-, t)$};
\draw (7,8) node[anchor=north west]{$(x_+, t)$};
\draw (5.5,8) node[anchor=north]{$(x,t)$};

\draw (1.7,2) node[anchor=east]{$\pi_-$};
\draw (8.7,2) node[anchor=west]{$\pi_+$};

\draw (6,8.5) node[anchor=west]{$(x-\delta,x+\delta)\times(t,t+\delta)$};
\end{scriptsize}

\end{tikzpicture}
 \caption{An illustration of the proof of Lemma \ref{lem:equivD}: we prove that $\cD$ is flat in the blue region.}
\label{fig:501}
\end{figure}
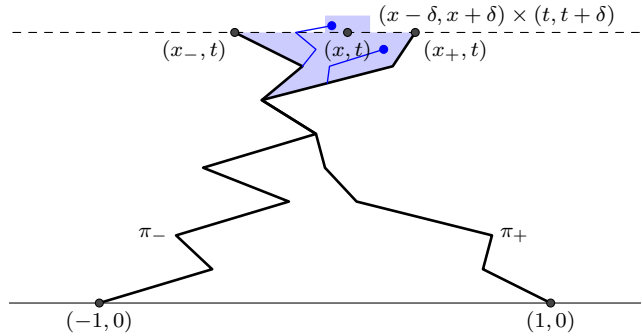

The idea is to show that, if $\cD$ is constant on a segment $(x_-,x_+)\times \{t\}$, it would also be constant in a 2D open neighborhood of it, using monotonicity of $\cD$ and transversal fluctuation estimates of geodesics.
\begin{proof}
The inclusion of $\hNC_t\subset \NC\cap(\R\times\{t\})$ is obvious, and we prove the other inclusion now.
Take any $(x,t) \in \HH$ such that $(x,t)\not\in \hNC_t$, i.e., $\cD(\cdot,t)$ is constant in a neighborhood of $x$.
Denote $X=\cD(x,t)$.
We can find $x_-<x<x_+$ such that
\[
\cD(x_-,t) = \cD(x_+,t) = X.
\]
We will show that the difference profile is constant in a 2D neighborhood of $(x, t)$. 
Although the proof below mainly consists of equations, there is a geometric heuristic behind this, as illustrated by Figure \ref{fig:501}. From the equality, we would expect that the geodesics $\pi_{(-1,0;x_-,t)}$ and $\pi_{(1,0;x_+,t)}$ coalesce, assuming their uniqueness (by Lemma \ref{lem:dis-equiv}). This then implies that for all $p$ and $p'$ in a neighborhood of $(x,t)$, $\pi_{(-1,0;p)}$ and $\pi_{(1,0;p')}$ coalesce, due to transversal fluctuation estimates and the ordering of the geodesics, and this would imply that the difference profile is constant in the neighborhood of $(x,t)$, whose generic point is, say, $(z,s).$

We argue separately for $(z,s)$ with $s<t$ (case (1)) or $s>t$ (case (2)).

\noindent
\textbf{Case 1.} We let $\pi_-$ and $\pi_+$ be any geodesics from $(1,0)$ to $(x_-,t)$, and $(-1,0)$ to $(x_+,t)$ respectively.
Take any $0<s< t$ with $\pi_-(s) < \pi_+(s)$.
We then have \[
\cD(\pi_-(s),s) = \cL(-1,0;\pi_-(s),s) - \cL(1,0;x_-,t) + \cL(\pi_-(s),s;x_-,t)
\le 
\cD(x_-,t) =: X    
\]
and
\[
\cD(\pi_+(s),s) = \cL(-1,0;x_+,t) - \cL(1,0;\pi_+(s),s) - \cL(\pi_+(s),s;x_+,t)
\ge 
\cD(x_+,t) = X.    
\]
By monotonicity of $\cD(\cdot,s)$, for any $z\in (\pi_-(s), \pi_+(s))$ we have
\[
\cD(x_-(s),s) \ge \cD(z,s) \ge \cD(x_+(s),s),
\]
so $\cD(z,s) = X$.
Since the geodesics $\pi_-$ and $\pi_+$ are continuous, we conclude that there exists some $\delta > 0$, such that for any $t-\delta<s\le t$ and $|z-x|<\delta$, we have $\cD(z,s) = X$.\\

\noindent
\textbf{Case 2:}
 Take any $s\ge t$ and $z \in \R$.
Consider any geodesics $\pi_-'$ from $(-1,0)$ to $(z,s)$, and $\pi_+'$ from $(1,0)$ to $(z,s)$.
We claim that if $\pi_-'(t), \pi_+'(t)\in [x_-, x_+]$, we must have  $\cD(z,s)= X$.
Indeed, assuming that $\pi_-'(t), \pi_+'(t)\in [x_-, x_+]$, we have
\[
\cD(z,s) \le \cD(\pi_-'(t),t) \le \cD(x_-,t) = X,
\]
where the first inequality is by $\cL(-1,0;z,s) = \cL(-1,0;\pi_-'(t),t) + \cL(\pi_-'(t),t;z,s)$ and $\cL(1,0;z,s) \ge \cL(1,0;\pi_-'(t),t) + \cL(\pi_-'(t),t;z,s)$, and the second inequality is by Lemma \ref{lem:DL-quad}. Similarly we have
\[
\cD(z,s) \ge \cD(\pi_+'(t),t) \ge \cD(x_+,t) = X,
\]
so the claim follows.

By Lemma \ref{lem:dl-trans-f}, (almost surely) we can find some $\delta>0$, such that for any $t\le s\le t+\delta$ and $|z-x|<\delta$, and any geodesic $\pi'$ from $(-1,0)$ to $(z,s)$ or from $(1,0)$ to $(z,s)$, there is $\pi'(t)\in [x_-, x_+]$.
Thus by the above claim we have $\cD(z,s) = X$ (for each $(z,s) \in [t,t+\delta] \times (x-\delta, x+\delta)$) and hence the proof is complete.
\end{proof}

We now are in a position to embark on proving Theorem \ref{thm:2dhaus}. 
The upper bound is straightforward, and the arguments are similar to those in \cite{BGH}.
The lower bound is significantly more delicate, and we now briefly reiterate our strategy from Section \ref{s:iop}.

We define a measure $\zeta$ on $\HH$, such that for any $a < b$ and $0\le g<h$, 
\begin{equation} \label{eq:def-zeta}
\zeta([a,b]\times (g,h)) = \int_g^h [\cD(a,t)-\cD(b,t)] dt.    
\end{equation}
We first observe that the support of $\zeta$ is contained in $\NC$, since for any $(x,t)\in \HH \setminus \NC$, there is an open neighborhood of $(x,t)$ with zero measure under $\zeta$.
We further prove that $\zeta$ is not degenerate, i.e., $\zeta(\HH)>0$.
In fact, by Lemma \ref{lem:DLbound}, as $a\to \infty$ we have
\[
\cD(-a,t)-\cD(a,t) \to \infty
\]
uniformly for $t \in [g, h]$, for any $0<g<h$.
Thus we conclude that 
\begin{equation}\label{nondeg1}
\zeta(\HH) \ge \lim_{a\to\infty} \zeta([-a,a]\times [g,h]) = \infty.
\end{equation}
Given the above no-degeneracy of the measure, 
to show that $\NC$ has Hausdorff dimension at least $5/3$,  it now suffices to prove the following ``2D H\"{o}lder'' bound.
\begin{prop}\label{prop:2d-53-H\"{o}lder}
For any $0<g<h$ and $N, \epsilon >0$, almost surely
\[\sup_{A\subset [-N,N]\times [g,h]} \frac{\zeta(A)}{\diam(A)^{5/3-\epsilon}} < \infty.\]
\end{prop}
As indicated in Section \ref{s:iop}, to prove such bound on the measure $\zeta$, we will obtain a ``level set decomposition'' of it. 

\subsection{Level set decomposition}
We consider ``level functions'' of the difference profile, as given in \eqref{levelcurve1}; i.e., for each $\ell \in \R$ and $t\in\R_+$, we take
\[\LV_\ell(t)=\inf\{x\in \R: \cD(x,t) \le \ell\}.\]
Then by continuity of $\cL$ we have $\cD(\LV_\ell(t),t) = \ell$. Actually, for all but countably many $\ell$, the infimum in the above definition is not needed.
\begin{lemma} \label{lem:uni-kappa}
For each $\ell \in \R,$ except for countably many of them, there is a unique $x\in\R$ with $\cD(x,t)=\ell$ for every $t\in\R_+$.  
\end{lemma}
\begin{proof}
Suppose that for some $\ell\in \R$ there are $t\in \R_+$ and $x_1<x_2$, such that $\cD(x_1,t)=\cD(x_2,t)=\ell$. 
By monotonicity of $\cD(\cdot, t)$, we have $\cD(x,t)=\ell$ for any $x_1\le x \le x_2$.
Then by Lemma \ref{lem:equivD}, one can find some rational $x'$ and $t'>0$ with $\cD(x',t')=\ell$. Thus the conclusion follows.
\end{proof}
Another basic property of $\LV_\ell$ is its continuity.
\begin{lemma} \label{lem:conti-kappa}
For each $\ell \in \R$, $\LV_\ell$ is a continuous function on $\R_+$.
\end{lemma}
As a standard exercise in measure theory, using Lemma \ref{lem:conti-kappa} and the fact that $\ell\mapsto\LV_\ell(t)$ is right continuous for each $t$, one can show that $(\ell, t)\mapsto \LV_\ell(t)$ is a measurable function on $\HH$.
Now, for each $\ell \in \R$, we define the occupation measure of $\LV_\ell(\cdot)$ to be denoted $\zeta_\ell$ on $\HH$, such that for any $0\le g<h$ and $a<b$, we have
\begin{equation} \label{eq:def-zeta-ell}
\zeta_\ell([a,b]\times (g,h)) = \int_g^h \don[\LV_\ell(t)\in [a,b]] dt.
\end{equation}
Recalling \eqref{eq:def-zeta}, observing that $\cD(a,t)-\cD(b,t)=\int \don[\LV_\ell(t)\in [a,b]] d\ell,$ using Fubini's theorem, we have that 
\[
\zeta([a,b]\times (g,h)) = \int_g^h \int \don[\LV_\ell(t)\in [a,b]] d\ell dt = \int \int_g^h \don[\LV_\ell(t)\in [a,b]]  dt d\ell = \int \zeta_\ell([a,b]\times (g,h)) d\ell.
\]
Thus (by standard measurable theory arguments)
\begin{equation} \label{eq:zeta-ell}
\zeta(A)=\int \zeta_\ell(A) d\ell    
\end{equation}
for any measurable $A\subset \HH$.
Then to bound $\zeta(A)$, we shall bound an average of $\zeta_\ell(A)$ over $\ell$.

Before proceeding further we prove the continuity of each $\LV_\ell$.
\begin{proof}[Proof of Lemma \ref{lem:conti-kappa}]
We argue by contradiction.
Suppose that $\LV_\ell$ is not continuous at some $t>0$. Denote $x=\LV_\ell(t)$.
Then there is a sequence $\{t_k\}_{k=1}^\infty$, such that $t_k\to t$ while $\LV_\ell(t_k)$ does not converge to $x$, as $k\to\infty$.
By Lemma \ref{lem:DLbound}, $|\LV_\ell(t_k)|$ are bounded as $k\to\infty$, so by taking a subsequence we can assume that $\lim_{k\to\infty}\LV_\ell(t_k)=y$ for some $y\neq x$.
Since $\cD(\LV_\ell(t_k),t_k) = \ell$ for each $k$, by continuity of $\cL$ we have $\cD(y,t) = \ell$. Thus we have $y>x$.
We take some $x<z<y$.
By Lemma \ref{lem:equivD}, there is some $\delta > 0$, such that for any $s\in (t-\delta, t+\delta)$, $\cD(z,s) = \ell$.
In particular, we can find some large enough $k$, such that $t_k\in (t-\delta, t+\delta)$, and $z<\LV_\ell(t_k)$, whereas $\cD(z,t_k) = \ell$.
This contradicts the definition of $\LV_\ell(t_k)$.
\end{proof}

\subsection{Comparison to competition interface}

To analyze the measures $\zeta_\ell$, we study the level functions $\LV_\ell$, by connecting them with the competition interface starting from the stationary initial condition.

We start by recalling certain objects. Namely, recall that we let $\cB$ be a standard two sided Brownian motion, $\cL^L(x,t)$ and $\cL^R(x,t)$ be the left and right weight profiles from $\cB$, and that $\LV^\cB$ denotes the competition interface  (from \eqref{leftrightprofile}).
Recall also the equality in law statement in Proposition \ref{prop:equal-dist-geo-comp}.
We define $\zeta^\cB$ as the occupation measure for  $\LV^\cB$ on $\HH,$  i.e., for any $a<b$ and $0\le g<h$, 
\begin{equation}\label{compinterfaceoccmeasure}
\zeta^\cB([a,b]\times (g,h)) = \int_g^h \don[\LV^\cB(t)\in [a,b]] dt.
\end{equation}
As indicated in Section \ref{s:iop}, the main idea is to show that, on an event with positive probability, $\LV^\cB$ equals $\LV_\ell$ for some random $\ell$; thus bounds on $\zeta^\cB$ would imply bounds on $\zeta_\ell$ averaged over $\ell$.
The former would be obtained by comparison with the semi-infinite geodesic $\pi_\boo$.\\

The following is a crucial estimate on the occupation measure $\zeta^\cB$ (and is equivalent to Lemma \ref{lem:bound-zb-s}).
\begin{lemma}  \label{lem:bound-zb}
There are universal constants $C,c>0$ such that the following is true.
For any closed interval $I\subset \R$ and $J\subset \R_{\ge 0}$, we have
\[
\PP[\zeta^\cB(I\times J) > M\LE(I)\LE(J)^{1/3}] < Ce^{-cM}
\]
for any $M>0$.
\end{lemma}
However, it turns out that for later applications, we will need a stronger version (see Lemma \ref{lem:geo-ocp-tail}) of the above. As already mentioned, the proofs are adaptations of arguments appearing previously in \cite{SSZ}. For completeness we include the proofs in the Appendix \ref{sec:appb}.

In the remainder of this section we finish the proof of Theorem \ref{thm:2dhaus}.
In Section \ref{ssec:spikebe} we set up the events where $\LV^\cB$ equals $\LV_\ell$ for some random $\ell$.
We prove Proposition \ref{prop:2d-53-H\"{o}lder} in Section \ref{ssec:53-lower}, using Lemma \ref{lem:bound-zb}, and the relations between $\LV^\cB$ and the level functions.
In Section \ref{ssec:53-upper} we prove the upper bound on the Hausdorff dimension.

Before proceeding to establish the comparison estimates between $\LV_\ell$ and the competition interface $\LV^\cB$ in the next section, it might be instructive to point out that such results only hold in a compact set (i.e. away from zero and infinity).
In fact, it is not too difficult to pin down the behaviors near zero (short times) and infinity (large times) respectively. It turns out that for $t$ near zero, $\LV_\ell(t)$ is roughly $-\ell t/4$ plus an order $t^{4/3}$ fluctuation, while as $t\to\infty$, $\LV_\ell(t)$ is linear in $t$ with a random slope, as evident from Figure \ref{f.simulations}.
See Remark \ref{rem:lc-ls-b} below for more details regarding this.

\subsection{Spiky boundary behavior}  \label{ssec:spikebe}

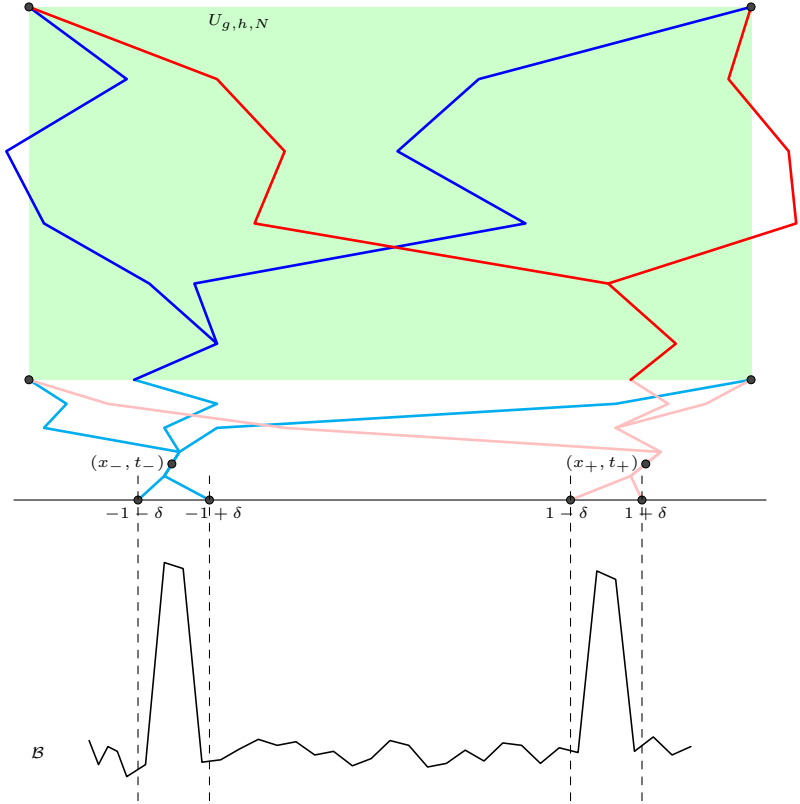
\begin{figure}[hbt!]
    \centering
\begin{tikzpicture}[line cap=round,line join=round,>=triangle 45,x=1cm,y=.16cm]
\clip(-0.5,-0.8) rectangle (10.5,68);

\fill[line width=0.pt,color=green,fill=green,fill opacity=0.2]
(0.2,35)-- (9.8,35)-- (9.8,66)-- (0.2,66) -- cycle;

\draw [line width=1pt,color=cyan] (1.65,25) -- (2,27) -- (2.2,29) -- (0.4,31) -- (0.7,33) -- (0.2,35);
\draw [line width=1pt,color=cyan] (2.6,25) -- (2,27) -- (2.2,29) -- (2.7,31) -- (8,33) -- (9.8,35);

\draw [line width=1pt,color=cyan] (2.2,29) -- (2,31) -- (2.7,33) -- (1.6,35);

\draw [line width=1pt,color=blue] (1.6,35) -- (2.7,38) -- (1.8,43) -- (0.4,48) -- (-0.1,54) -- (1.5,60) -- (0.2,66);

\draw [line width=1pt,color=blue] (2.7,38) -- (2.4,43) -- (6.8,48) -- (5.1,54) -- (6.18,60) -- (9.8,66);

\draw [line width=1pt,color=pink] (8.35,25) -- (8.2,27) -- (8.6,29) -- (8,31) -- (9.2,33) -- (9.8,35);
\draw [line width=1pt,color=pink] (7.4,25) -- (8.2,27) -- (8.6,29) -- (3.6,31) -- (1.25,33) -- (0.2,35);
\draw [line width=1pt,color=pink] (8,31) -- (8.7,33) -- (8.2,35);

\draw [line width=1pt,color=red] (8.2,35) -- (8.8,38) -- (7.9,43) -- (10.4,48) -- (10.3,54) -- (9.5,60) -- (9.8,66);

\draw [line width=1pt,color=red] (7.9,43) -- (3.2,48) -- (3.6,54) -- (2.7,60) -- (0.2,66);

\draw (0,25) -- (10,25);

\draw [fill=uuuuuu] (0.2,35) circle (1.5pt);
\draw [fill=uuuuuu] (9.8,35) circle (1.5pt);

\draw [fill=uuuuuu] (1.65,25) circle (1.5pt);
\draw [fill=uuuuuu] (2.6,25) circle (1.5pt);

\draw [fill=uuuuuu] (2.1,28) circle (1.5pt);
\draw [fill=uuuuuu] (8.4,28) circle (1.5pt);

\draw [fill=uuuuuu] (7.4,25) circle (1.5pt);
\draw [fill=uuuuuu] (8.35,25) circle (1.5pt);
\draw [fill=uuuuuu] (0.2,66) circle (1.5pt);
\draw [fill=uuuuuu] (9.8,66) circle (1.5pt);

\draw [line width=0.6pt] (1,5) -- (1.125,3.) -- (1.25,4.5) -- (1.375,4.1) -- (1.5,2) -- (1.75,3) -- (2,19.8) -- (2.25,19.3) -- (2.5,3.2) -- (2.75,3.4) -- (3,4.3) -- (3.25,5.1) -- (3.5,4.6) -- (3.75,4.9) -- (4,3.8) -- (4.25,4.1) -- (4.5,2.9) -- (4.75,3.5) -- (5,5) -- (5.25,4.6) -- (5.5,2.8) -- (5.75,3.1) -- (6,4.2) -- (6.25,3.3) -- (6.5,4.8) -- (6.75,4.6) -- (7,3.1) -- (7.25,4.4) -- (7.5,4) -- (7.75,19.1) -- (8,18.4) -- (8.25,4.1) -- (8.5,5.3) -- (8.75,3.8) -- (9,4.5);

\foreach \i in {1.65, 2.6, 7.4, 8.35}
{
\draw [dashed] (\i,0) -- (\i,27);
}

\begin{tiny}
\draw (1.6,25) node[anchor=north]{$-1-\delta$};
\draw (2.65,25) node[anchor=north]{$-1+\delta$};
\draw (7.35,25) node[anchor=north]{$1-\delta$};
\draw (8.4,25) node[anchor=north]{$1+\delta$};

\draw (3,66) node[anchor=north]{$U_{g,h,N}$};

\draw (2.1,28) node[anchor=east]{$(x_-,t_-)$};
\draw (8.4,28) node[anchor=east]{$(x_+,t_+)$};

\draw (0.5,4) node[anchor=east]{$\cB$};
\end{tiny}

\end{tikzpicture}
\caption{The top half of this figure illustrates the event $\Good_{g,h,N,\delta}$, and the bottom half of this figure illustrates the event where $\cB$ has spikes around $-1$ and $1$. These two events imply that for any $p$ in the box $U_{g,h,N}$, the optimal paths would land near $-1$ or $1$.}
\label{fig:1}
\end{figure}

We will first establish how one can relate $\LV^\cB$, and $\LV_\ell$ for different $\ell$.
For each $0<g<h$, $N>0$, and $\delta>0$, let $U_{g,h,N}$ denote the box $[-N,N] \times [g/2,2h]$. We consider the event $\Good_{g,h,N,\delta}$, where the following natural  conditions (with somewhat technical looking quantifiers which the reader could without hesitation pay no attention to on first read) hold.
\begin{enumerate}
    \item Coalescence of geodesics to $U_{g,h,N}$: there exist some $t_-, t_+>0$, $x_-, x_+ \in \R$, such that for any $x\in (-1-\delta, -1+\delta)$ and $q \in U_{g,h,N}$, there is $\pi_{(x,0;q)}(t_-)=x_-$; and for any $x\in (1-\delta, 1+\delta)$ and $q \in U_{g,h,N}$, there is $\pi_{(x,0;q)}(t_+)=x_+$. See the top half of Figure \ref{fig:1} for an illustration of the event.
    \item Regularity of the directed landscape: $R<\delta^{-0.01}$, where $R$ is the maximum of the random variables given by Lemma \ref{lem:DLbound}, Lemma \ref{lem:modcont} for $K=([-4N, 4N]\times [-1, g/6])\times ([-4N, 4N]\times [g/3, 2h])$, and Lemma \ref{lem:dl-trans-f}.
\end{enumerate}
The next result says that the above event holds with high probability as $\delta$ tends to $0.$
\begin{lemma}  \label{lem:low-bound-G}
For any $0<g<h$, $N>0$, we have $\lim_{\delta\searrow 0}\PP[\Good_{g,h,N,\delta}]=1$. 
\end{lemma}
\begin{proof}
We denote the two events as $\cE_1, \cE_2$ respectively.

By Lemmas \ref{lem:DLbound}, \ref{lem:modcont}, and \ref{lem:dl-trans-f}, we have $\PP[\cE_2]\to 1$ as $\delta \searrow 0$.
Thus it remains to analyze $\cE_1$.
Assuming that $\cE_2$ holds, by Lemma \ref{lem:dl-trans-f}, there exists some $Q\in \N$, determined by $g$, $h$, $\delta$, $N$ (which are the parameters appearing in the definition of the good event), such that 
\[
\pi_{(-1-\delta,0;Q,2h)}(t) > N, \;
\pi_{(1+\delta,0;-Q,2h)}(t) < -N
\]
for any $g/2\le t \le 2h$.
By Theorem \ref{thm:disj}, almost surely there is some random $\delta'>0$, such that $\pi_{(-1-\delta',0;-Q,2h)}$ intersects $\pi_{(-1+\delta',0;Q,2h)}$, and $\pi_{(1-\delta',0;-Q,2h)}$ intersects $\pi_{(1+\delta',0;Q,2h)}$, and all of these geodesics are unique.
Then we just take $t_-, x_-$ such that $\pi_{(-1-\delta',0;-Q,2h)}(t_-)=\pi_{(-1+\delta',0;Q,2h)}(t_-)=x_-$, and $t_+, x_+$ such that $\pi_{(1-\delta',0;-Q,2h)}(t_+)=\pi_{(1+\delta',0;Q,2h)}(t_+)=x_+$.
Thus by the ordering of geodesics (Lemma \ref{lem:ordering-geo}), we have that for any $x\in (-1-\delta', -1+\delta')$ and $q \in U_{g,h,N}$, there is $\pi_{(x,0;q)}(t_-)=x_-$; and for any $x\in (1-\delta', 1+\delta')$ and $q \in U_{g,h,N}$, there is $\pi_{(x,0;q)}(t_+)=x_+$.
This implies that $\PP[\cE_2\setminus\cE_1]\to 0$ as $\delta \searrow 0$. Thus the conclusion follows.
\end{proof}

We next consider the key event $\ConsDiff_{g,h,N}$ which allows us to relate the difference profile rooted at $(\pm1, 0)$ to the one corresponding to Brownian initial data. This is defined as follows: there exists a (random) number $\alpha$, such that for any $p\in U_{g,h,N}$, there is
\[
\alpha = \cD(p) - (\cL^L(p)-\cL^R(p)),
\]
i.e., the two difference profiles are a shift of each other when restricted to $U_{g,h,N}.$
As an useful consequence, we have that under $\ConsDiff_{g,h,N}$,
\[
\{(\LV_\alpha(t),t): t\ge 0\} \cap U_{g,h,N} = \{(\LV^\cB(t),t): t\ge 0\} \cap U_{g,h,N}.
\]
We next come to the crucial result which allows us to go between the level sets of $\cD$ and $\LV^\cB$.
\begin{lemma}  \label{lem:low-bound-D}
For any $0<g<h$, $N>0$, and any $\delta >0$, $H>0$, there is some $\delta'>0$, such that for any interval $I\subset [-H,H]$, and any sampling of $\cL$ such that $\Good_{g,h,N,\delta}$ holds, we have
\[
\PP[ \ConsDiff_{g,h,N}, \alpha \in I \mid \cL] > \delta' \LE(I).
\]
where $\LE(\cdot)$ denotes the Lebesgue measure.
\end{lemma}
Re-iterating from Section \ref{s:iop}, the proof of this lemma is to force $\cB$ to have spikes around $-1$ and $1$. Thus for any $p$ in a box $U_{g,h,N}$, the optimal path to the left half boundary (i.e. the path corresponding to $\cL^L(p)$) would land near $-1$; and the optimal path to the right half boundary (i.e. the path corresponding to $\cL^R(p)$) would land near $1$, respectively. See Figure \ref{fig:1} for an illustration. Further, we will show that on the coalescence event $\Good_{g,h,N,\delta}$, the event $\ConsDiff_{g,h,N}$ would hold.
Finally, by a simple argument involving resampling  $\cB,$ we will also establish that the random difference $\alpha$ is independent of the events $\Good_{g,h,N,\delta}$ and $\cB$ being spiky and has a smooth distribution. 

\begin{figure}[hbt!]
    \centering
\begin{tikzpicture}[line cap=round,line join=round,>=triangle 45,x=.8cm,y=.16cm]
\clip(-2.5,-10.8) rectangle (12.5,42);

\draw (-5,4) -- (15,4);

\draw [fill=uuuuuu] (5,4) circle (1.5pt);

\draw [line width=0.6pt] (2.75,3.4-5) -- (3,4.3-4) -- (3.25,5.1-3.5) -- (3.5,4.6-3) -- (3.75,4.9-2.5) -- (4,3.8-2) -- (4.25,4.1-1.5) -- (4.5,2.9-1) -- (4.75,3.5-0.5) -- (5,4) -- (5.25,4.6+0.5) -- (5.5,2.8+1) -- (5.75,3.1+2) -- (6,4.2+2.5) -- (6.25,3.3+3) -- (6.5,4.8+3.5) -- (6.75,4.6+4) -- (7,3.1+4.5) -- (7.25,4.4+5);

\draw [blue] [line width=0.9pt] (-1,5.3-5) -- (-0.75,3.2-5) -- (-0.5,6.1-5) -- (-0.25,4.5-5) -- (0,5.9-5) -- (0.25,4.1-5) -- (0.5,5.3-5) -- (0.75,4.2-5) -- (1,5-5) -- (1.25,4.5-5) -- (1.5,2-5) -- (1.75,3-5) -- (2,29.8-5) -- (2.25,29.3-5) -- (2.5,3.2-5) -- (2.75,3.4-5);

\draw [blue] [line width=0.9pt] (7.25,4.4+5) -- (7.5,4+5) -- (7.75,29.1+5) -- (8,28.4+5) -- (8.25,4.1+5) -- (8.5,5.3+5) -- (8.75,3.8+5) -- (9,4.5+5) -- (9.25,5.6+5) -- (9.5,3.8+5) -- (9.75, 4.1+5) -- (10,4.4+5) -- (10.25,3.6+5) -- (10.5,4.9+5) -- (10.75,4.7+5) -- (11,3.6+5);

\foreach \i in {1.5, 2.75, 7.25, 8.5}
{
\draw [dotted] (\i,-15) -- (\i,55);
}

\draw [blue] [dashed] (2.75,27-5) -- (1.5,27-5);
\draw [blue] [dashed] (2.75,31-5) -- (1.5,31-5);

\draw [blue] [dashed] (7.25,28+5) -- (8.5,28+5);
\draw [blue] [dashed] (7.25,32+5) -- (8.5,32+5);

\draw [blue] [dashed] (1.5,8-5) -- (-6,10-5);
\draw [blue] [dashed] (8.5,9+5) -- (16,11+5);

\draw [dashed] (7.25,-7) -- (2.75,-7);
\draw [dashed] (7.25,15) -- (2.75,15);

\begin{tiny}
\draw (1.5,-9) node[anchor=north]{$-1-\delta$};
\draw (2.75,-9) node[anchor=north]{$-1+\delta$};
\draw (7.25,-9) node[anchor=north]{$1-\delta$};
\draw (8.5,-9) node[anchor=north]{$1+\delta$};

\draw (2.75,15) node[anchor=south west]{$A/2$};
\draw (2.75,-7) node[anchor=south west]{$-A/2$};

\draw (-1,3-5) node[anchor=east]{$\cB$};

\draw [blue] (1.5,27-5) node[anchor=east]{$\cB(-1+\delta)+A$};
\draw [blue] (1.5,31-5) node[anchor=east]{$\cB(-1+\delta)+A+1$};

\draw [blue] (7.25,28+5) node[anchor=east]{$\cB(1-\delta)+A$};
\draw [blue] (7.25,32+5) node[anchor=east]{$\cB(1-\delta)+A+1$};

\draw [blue] (-1,9-5) node[anchor=south]{$\cB(-1+\delta)+1-x/A$};
\draw [blue] (11,10+5) node[anchor=south]{$\cB(1-\delta)+1+x/A$};

\draw (5,3.8) node[anchor=north]{$0$};
\end{tiny}

\end{tikzpicture}
\caption{An illustration of the events in the proof of Lemma \ref{lem:low-bound-D}: the blue curves are the processes $\cB(1-\delta+x)-\cB(1-\delta)$ and  $\cB(-1+\delta-x)-\cB(-1+\delta)$ for $x>0$, which determine the event $\cE$. Conditional on these two processes, $\cB$ within $[-1-\delta, -1+\delta]$ is still a two-sided Brownian motion, and determines $\cE'$.}
\label{fig:11}
\end{figure}
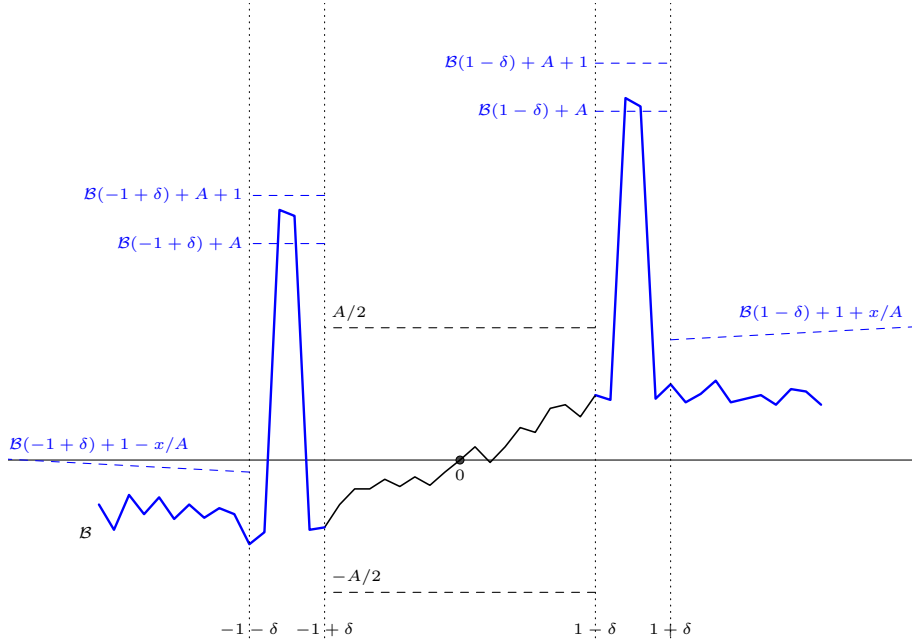
\begin{proof}
Since $\cL$ is conditioned on, the only randomness is from the Brownian motion $\cB$. The general idea is to make $\cB$ ``spiky'' in the intervals $[-1-\delta, -1+\delta]$ and $[1-\delta, 1+\delta]$, such that for any $p\in U_{g,h,N}$, we have 
\begin{equation}  \label{eq:low-bound-D-pf2}
\cL^L(p)=\sup_{|y+1|\le \delta} \cL(y,0;p) + \cB(y),\quad \cL^R(p)=\sup_{|y-1|\le \delta} \cL(y,0;p) + \cB(y).
\end{equation}
Assuming $\Good_{g,h,N,\delta}$, we have that \eqref{eq:low-bound-D-pf2} implies $\ConsDiff_{g,h,N}$.
Indeed, by the coalecence condition of the event $\Good_{g,h,N,\delta}$ (as illustrated by Figure \ref{fig:1}, and recall the points $(x_-, t_-)$ and $(x_+, t_+)$ there), we have for any $p_1, p_2 \in U_{g,h,N}$,
\begin{align*}
\cL^L(p_1)-\cL^L(p_2) &= \cL(x_-,t_-;p_1) - \cL(x_-,t_-;p_2) = \cL(-1,0;p_1) - \cL(-1,0;p_2),\\
\cL^R(p_1)-\cL^R(p_2) &= \cL(x_+,t_+;p_1) - \cL(x_+,t_+;p_2) = \cL(1,0;p_1) - \cL(1,0;p_2).
\end{align*}

These imply that $\cD(p_1)-(\cL^L(p_1)-\cL^R(p_1)) = \cD(p_2)-(\cL^L(p_2)-\cL^R(p_2))$, so $\ConsDiff_{g,h,N}$ holds.
The value of $\alpha$ (from the definition of the event $\ConsDiff_{g,h,N}$), can be altered by ``adjusting the difference in the heights of the spikes'', i.e., by adding different constants to $\cB$ in $[-1-\delta, -1+\delta]$ and $[1-\delta, 1+\delta]$. This is done by decomposing the process $\cB$ into three independent parts: the processes $\cB(1-\delta+x)-\cB(1-\delta)$ and  $\cB(-1+\delta-x)-\cB(-1+\delta)$ for $x>0$, and $\cB$ in $[-1+\delta, 1-\delta]$, and modifying the last part while fixing the first two.
See Figure \ref{fig:11} for an illustration of some of the events we define next.

Take $A$ to be a large enough number depending on $g$, $h$, $\delta$, $N$ (which are parameters in the definition of the good event) and $H$.
We consider the event $\cE$, given by the following conditions:
\begin{itemize}
    \item $\sup_{x<-1-\delta} \cB(x)-\cB(-1+\delta)+x/A < 1$,
   \item   $\sup_{x>1+\delta} \cB(x)-\cB(1-\delta)- x/A < 1$,
    \item $\sup_{x\in(-1-\delta,-1+\delta)} \cB(x)-\cB(-1+\delta) \in [A, A+1]$, 
    \item
    $\sup_{x\in(1-\delta,1+\delta)} \cB(x)-\cB(1-\delta) \in [A, A+1]$.
\end{itemize}
In words, the first two conditions stipulate a bound for  the processes $\cB(1-\delta+x)-\cB(1-\delta)$ and  $\cB(-1+\delta-x)-\cB(-1+\delta)$for $x> 2\delta$ that force them to be essentially bounded by $1$ for an interval of size $A$ (a large number) while beyond that, not growing faster than the linear function $|x|/A.$ Note that the last part can be ensured in fact by any function that grows faster than $\sqrt{x}$.

The third and fourth conditions say that the processes $\cB(1-\delta+x)-\cB(1-\delta)$ and  $\cB(-1+\delta-x)-\cB(-1+\delta)$ attain the value $A$ in $[0,2\delta]$.
Thus together, the above ensures the spiky behavior we desire.

As indicated above, we also want to adjust the difference in the spike heights and to be able to do that we define  
 for any $x\in\R, t>0$, the following proxies for $\cL^L(x,t)$ and $\cL^R(x,t)$, by restricting the range of $y$ and shifting by the value $\cB$ at the endpoint. Namely,
\[
\hcL^L(x,t)=\sup_{y\le -1+\delta} \cL(y,0;x,t) + \cB(y)-\cB(-1+\delta),
\]
and
\[
\hcL^R(x,t)=\sup_{y\ge 1-\delta} \cL(y,0;x,t) + \cB(y)-\cB(1-\delta).
\]
Observe that $\hcL^L$ and $\hcL^R$ are independent of $\cB$ in $[-1+\delta, 1-\delta]$ and we will use the randomness in the latter to alter the heights of the spikes. 
Recalling the parameter $g$ appearing in the definition of $\Good_{g,h,N,\delta},$ let $p_0=(0,g)$, and define $\alpha'$ via
\[
\alpha'= \cD(p_0) - (\hcL^L(p_0)-\hcL^R(p_0)).
\]
Recalling that $I$ is the interval in the statement of Lemma \ref{lem:low-bound-D}, we next consider the event $\cE'$, where
\begin{itemize}
    \item $\cB(-1+\delta)-\cB(1-\delta) \in \alpha' - I$;
    \item $\sup_{|x|\le 1-\delta} |\cB(x)| < A/2$.
\end{itemize}
Note that for a given $\alpha'$, the event $\cE'$ is independent of $\cE$.
Using basic properties of Brownian motions, we can find some $\delta'>0$, depending on all parameters, such that $\PP[\cE] > \sqrt{\delta'}$; and conditional on $\cE$, and any $\cL$ such that $\Good_{g,h,N,\delta}$ holds, we have $\PP[\cE'\mid \cE, \cL] > \sqrt{\delta'}\LE(I)$.
Here for the second part, we use that $\PP[\cE'\mid \cE, \cL, \alpha'] = \PP[\cE'\mid \cL, \alpha'] = \PP[\cE'\mid \Good_{g,h,N,\delta}, \alpha']$, and that under the event $\Good_{g,h,N,\delta}$, $\alpha'-I\subset [-A/2, A/2]$ by taking $A$ large enough.

Thus we have $$\PP[\cE\cap\cE'\mid \cL] = \PP[\cE'\mid \cE, \cL]\PP[\cE] > \delta'\LE(I),$$ for any $\cL$ such that $\Good_{g,h,N,\delta}$ holds. 

It remains to show that $\cE\cap\cE'\cap\Good_{g,h,N,\delta}\subset\ConsDiff_{g,h,N}\cap\{\alpha \in I\}$.
Under the event $\cE\cap\cE'\cap\Good_{g,h,N,\delta}$, since $A$ is taken large enough, we have that \eqref{eq:low-bound-D-pf2} holds for any $p \in U_{g,h,N}$. Indeed, for this, we just use the event $\cE\cap\cE'$ to lower bound 
$\sup_{|y+1|\le \delta}\cB(y)$ and $\sup_{|y-1|\le \delta}\cB(y)$, and upper bound $\cB(y)$ for $y \not\in [-1-\delta,-1+\delta]\cup [1-\delta, 1+\delta]$; and we use $R<\delta^{-0.01}$ (where $R$ is the random variable given by Lemma \ref{lem:DLbound}) from the event $\Good_{g,h,N,\delta}$ to bound $\cL(y,0;p)$ for any $y\in \R$.
Then \eqref{eq:low-bound-D-pf2} implies $\ConsDiff_{g,h,N}$, as argued at the beginning of this proof. Under $\ConsDiff_{g,h,N}$, the event $\cE'$ implies that $\alpha'-\alpha = \cB(-1+\delta)-\cB(1-\delta)$, so $\alpha \in I$. Thus we conclude that $\cE\cap\cE'\cap\Good_{g,h,N,\delta}\subset\ConsDiff_{g,h,N}\cap\{\alpha \in I\}$.
\end{proof}

\subsection{Proof of the lower bound}
\label{ssec:53-lower}

Proposition \ref{prop:2d-53-H\"{o}lder} is an immediate consequence of the following result.
\begin{lemma}\label{lem:2d-53-H\"{o}lder}
For any $0<g<h$ and $N, \epsilon >0$,
there is a random number $Q$ such that the following is true.
For any $k\in \N$ and $i\in \llbracket -2^k, 2^k-1\rrbracket$, $j\in \llbracket 0, 2^k-1\rrbracket$ we have
\[
\zeta([i2^{-k}N, (i+1)2^{-k}N] \times [j2^{-k}(h-g)+g, (j+1)2^{-k}(h-g)+g] ) < Q 2^{-(5/3-\epsilon)k}.
\]
\end{lemma}
Let us iterate from Section \ref{s:iop} the source of the exponent $5/3$ in the right hand side. 
By the regularity condition of the directed landscape under the event $\Good_{g,h,N,\delta}$, the measure $\zeta_\ell$ on a set of size $2^{-k}\times 2^{-k}$ is nonzero only for $\ell$ in an interval of length (in the order of) $2^{-k/3}$.
Thus to bound $\zeta$ on a set of size $2^{-k}\times 2^{-k}$, we just bound an integral of $\zeta_\ell$ (on the set) for $\ell$ in an interval of length $2^{-k/3}$, and take the maximum over all such intervals.
By Lemma \ref{lem:low-bound-D}, (under $\Good_{g,h,N,\delta}$) such a quantity is bounded by $2^{-k/3}$ times $\zeta^\cB$ (on the set of size $2^{-k}\times 2^{-k}$).
Using Lemma \ref{lem:bound-zb}, we can bound $\zeta^\cB$ (on a set of size $2^{-k}\times 2^{-k}$) by the order of $2^{-4k/3}$, with exponential tails.
These combined would give the power of $5/3$.
\begin{proof}
We let $C,c>0$ denote large and small constants depending on $g, h, N, \epsilon$, and their values may change from line to line.
All other constants and parameters in this proof would also depend on $g, h, N, \epsilon$.

For any $\delta>0$, it suffices to prove the result under the event $\Good_{g,h,N,\delta}$.
To control $\zeta$ on a dyadic box, for any $\ell,$ we let
\[
Q_\ell = \sup_{k\in\N, i\in \llbracket -2^k, 2^k-1\rrbracket, j\in \llbracket 0, 2^k-1\rrbracket} 
2^{(4/3-\epsilon/2)k} \zeta_\ell([i2^{-k}N, (i+1)2^{-k}N] \times [j2^{-k}(h-g)+g, (j+1)2^{-k}(h-g)+g] ).
\]
We note that, as we expect that $\zeta_\ell$ would assign mass of the order of $2^{-k/3}\times 2^{-k}=2^{-4k/3}$ to any dyadic box of size $2^{-k}\times 2^{-k}$, we multiply by the factor of $2^{(4/3-\epsilon/2)k}$ to make $Q_\ell$ to be of an order slightly smaller than $1$.
Now for any $H>0$ large enough (depending on $\delta$),
for any $\varepsilon>0$, we denote
\[
J_\varepsilon = \sup_{A\subset [-H,H], \LE(A)<\delta^{-1} \varepsilon^{1/3}\log^{2/3}(1+\varepsilon^{-1}) } \int_A Q_\ell d\ell,
\]
and
\begin{equation}  \label{eq:defJ}
J= \sup_{\varepsilon\in (0,1)} \frac{J_\varepsilon}{\varepsilon^{1/3}\log(1+\varepsilon^{-1})} .
\end{equation}
The definition of $J_\varepsilon$ is motivated by the following (above mentioned) heuristic: by continuity estimates of $\cL$ (ensured by the second condition of $\Good_{g,h,N,\delta}$), for any $k\in\N, i\in \llbracket -2^k, 2^k-1\rrbracket, j\in \llbracket 0, 2^k-1\rrbracket$, the measure $\zeta_\ell([i2^{-k}N, (i+1)2^{-k}N] \times [j2^{-k}(h-g)+g, (j+1)2^{-k}(h-g)+g] )$ is non-zero only for $\ell$ in an interval of length roughly in the order of $2^{-k/3}$.

This number $J$ is essentially (up to a constant factor) the desired $Q$ (in the statement of this lemma).
This will be shown in two steps: (1) we bound the upper tail of $\don[\Good_{g,h,N,\delta}]J$; (2) we bound $\zeta$ on small rectangles, using $J$ under the event $\Good_{g,h,N,\delta}$.\\

\noindent
\textbf{Step 1.}
We carry out the first step next relying on the following input about the competition interface with stationary initial condition.  Let the counterpart of $Q_\ell$ for the stationary initial condition be defined as 
\[
Q^\cB = \sup_{k\in\N, i\in \llbracket -2^k, 2^k-1\rrbracket, j\in \llbracket 0, 2^k-1\rrbracket} 
2^{(4/3-\epsilon/2)k} \zeta^\cB([i2^{-k}N, (i+1)2^{-k}N] \times [j2^{-k}(h-g)+g, (j+1)2^{-k}(h-g)+g] ).
\]
By applying Lemma \ref{lem:bound-zb} to the boxes $[i2^{-k}N, (i+1)2^{-k}N] \times [j2^{-k}(h-g)+g, (j+1)2^{-k}(h-g)+g]$, and a union bound over all $k, i, j$, we have
\[
\PP[Q^{\cB}>M]<C_1e^{-c_1M},
\]
for any $M>0$, for some constants $c_1,C_1>0$. 
Thus by Lemma \ref{lem:low-bound-D}, we have 
\begin{equation}\label{eq:tail-Q-ell}
\E\left[\don[\Good_{g,h,N,\delta}]\int_{-H}^H \don[Q_\ell > M] d\ell \right] \le {\delta'}^{-1} \PP[Q^{\cB}>M]  < {\delta'}^{-1} C_1e^{-c_1M}.    
\end{equation}
where $\delta'$ appears in the statement of  Lemma \ref{lem:low-bound-D}.
We also assume that $\delta'\le \delta$, since otherwise we can repalce $\delta'$ by $\delta'\wedge \delta$.

The next few lines are technical and can be a bit hard to parse, but are used to spell out the algebra needed to translate this expectation bound \eqref{eq:tail-Q-ell} into bounds on the tail of $\don[\Good_{g,h,N,\delta}]J$.
For any $\varepsilon, M>0,$ let
\[
J_{\varepsilon,M} = \sup_{A\subset [-H,H], \LE(A)<\delta^{-1} \varepsilon^{1/3}\log^{2/3}(1+\varepsilon^{-1}) } \int_A \don[Q_\ell>M] d\ell .
\]
By simply bounding the integrand by $1$, we have
\begin{equation}  \label{eq:E-varM-r}
\don[\Good_{g,h,N,\delta}]J_{\varepsilon,M} \le 
\delta^{-1} \varepsilon^{1/3}\log^{2/3}(1+\varepsilon^{-1}) .
\end{equation}
On the other hand, since $J_{\varepsilon,M} \le \int_{-H}^H \don[Q_\ell > M] d\ell$, by \eqref{eq:tail-Q-ell} we have
\begin{equation}  \label{eq:E-varM}
\E[\don[\Good_{g,h,N,\delta}]J_{\varepsilon,M}] \le {\delta'}^{-1} C_1e^{-c_1M}.
\end{equation}
We take
\[
E_\varepsilon = 
\begin{cases}
\sup_M c_1^{-1}\log(C_1^{-1}\delta' J_{\varepsilon,M}) + M/2 & \text{if } \Good_{g,h,N,\delta} \text{ holds,} \\
0 & \text{otherwise.}
\end{cases}
\]
which then rewritten yields
\begin{equation}  \label{eq:E-varep}
\don[\Good_{g,h,N,\delta}]J_{\varepsilon,M} \le {\delta'}^{-1} C_1e^{-c_1M/2+c_1E_\varepsilon},\;\; \forall M>0.    
\end{equation}
We next bound the upper tail of $E_\varepsilon$ (which in particular implies that almost surely it is finite).
Assume that $\Good_{g,h,N,\delta}$ holds, 
then from \eqref{eq:E-varM-r} and the fact that $\delta'<\delta$, we have $c_1^{-1}\log(C_1^{-1}\delta' J_{\varepsilon,M}) + M/2<C$ for any $M<1$.
Also, by \eqref{eq:E-varM}, and Markov's inequality, for any $i\in\N$ and $M>0$, we have
\[
\PP[\{c_1^{-1}\log(C_1^{-1}\delta' J_{\varepsilon,i}) + (i+1)/2 > M\} \cap \Good_{g,h,N,\delta}] < e^{-c_1((i-1)/2+M)}.
\]
Thus by a union bound over $i\in\N$, and noting that $J_{\varepsilon,M}$ is non-increasing in $M$, we have $\PP[E_\varepsilon > M] < C e^{-c_1M}$ for any $M>0$.

Denote by $M_\varepsilon$  the smallest non-negative number such that $${\delta'}^{-1} C_1e^{-c_1M_\varepsilon/2+c_1E_\varepsilon} \le \delta^{-1} \varepsilon^{1/3}\log^{2/3}(1+\varepsilon^{-1}).$$
Then we have
\begin{align*}
\don[\Good_{g,h,N,\delta}]J_\varepsilon &\le \int_0^\infty \don[\Good_{g,h,N,\delta}]J_{\varepsilon, M} dM
\\
&\le \delta^{-1} \varepsilon^{1/3}\log^{2/3}(1+\varepsilon^{-1})M_\varepsilon + \int_{M_\varepsilon}^\infty \don[\Good_{g,h,N,\delta}]J_{\varepsilon, M} dM
\\
& \le
\delta^{-1} \varepsilon^{1/3}\log^{2/3}(1+\varepsilon^{-1})M_\varepsilon + \int_{M_\varepsilon}^\infty {\delta'}^{-1} C_1e^{-c_1M/2+c_1E_\varepsilon} dM
\\
& =
\delta^{-1} \varepsilon^{1/3}\log^{2/3}(1+\varepsilon^{-1})M_\varepsilon + 2c_1^{-1}{\delta'}^{-1} C_1e^{-c_1M_\varepsilon/2+c_1E_\varepsilon}
\\
& \le
(\delta^{-1}M_\varepsilon+ 2c_1^{-1}\delta^{-1}) \varepsilon^{1/3}\log^{2/3}(1+\varepsilon^{-1}).
\end{align*}
Here the first inequality is by the definition of $J_\varepsilon$ and $J_{\varepsilon,M}$, the second inequality is by splitting the integral at $M_\varepsilon$ and \eqref{eq:E-varM-r}, the third inequality is by \eqref{eq:E-varep}, and the last inequality is by the definition of $M_\varepsilon$.
From the tail of $E_\varepsilon$ and the relation between $E_\varepsilon$ and $M_\varepsilon$, we conclude that for any $M>0$,
\[
\PP[\don[\Good_{g,h,N,\delta}]J_\varepsilon > M\varepsilon^{1/3} \log^{2/3}(1+\varepsilon^{-1})] < C_\delta' e^{-c_\delta' M},
\]
where $C_\delta', c_\delta'$ are constants depending on $\delta$.
Recall $J$ from \eqref{eq:defJ}. By taking the union bound over an exponentially decreasing sequence of $\varepsilon$, and noting that $J_{\varepsilon}$ is monotone in $\varepsilon$ (for $\varepsilon$ small enough), 
we conclude that $\PP[\don[\Good_{g,h,N,\delta}]J > M] < C_\delta e^{-c_\delta M}$, for any $M>0$, where $C_\delta, c_\delta$ are constants depending on $\delta$.

\noindent 
\textbf{Step 2.}
We now implement the second step.
Under the event $\Good_{g,h,N,\delta}$, using its second condition and Lemma \ref{lem:DLbound}, we can bound $|\cD|$ in $[-N,N]\times[g,h]$.
Then as $H$ is taken to be large enough depending on $\delta$, we have $\zeta_\ell([-N,N]\times[g,h]) = 0$ for any $|\ell|>H$.
Thus for any measurable $A \subset[-N,N]\times[g,h]$, we have $\zeta(A) = \int_{-H}^H \zeta_\ell(A) d\ell$.
Hence, for any $k\in \N$ large enough (depending on $\delta$) and $i\in \llbracket -2^k, 2^k-1\rrbracket$, $j\in \llbracket 0, 2^k-1\rrbracket$ we have
\begin{align*}
&\zeta([i2^{-k}N, (i+1)2^{-k}N] \times [j2^{-k}(h-g)+g, (j+1)2^{-k}(h-g)+g] )\\
=&\int_{-H}^H \zeta_\ell([i2^{-k}N, (i+1)2^{-k}N] \times [j2^{-k}(h-g)+g, (j+1)2^{-k}(h-g)+g] ) d\ell
\\
<& 2^{-(4/3-\epsilon/2)k} 2^{-k/3}\log(1+2^k) J\\
<& CJ 2^{-(5/3-\epsilon)k} .
\end{align*}
Here the first inequality is by the definition of $J$, and the following reasoning: under the second condition of $\Good_{g,h,N,\delta}$ and using Lemma \ref{lem:modcont}, for any $\ell_1, \ell_2$, if $\zeta_{\ell_1}([i2^{-k}N, (i+1)2^{-k}N] \times [j2^{-k}(h-g)+g, (j+1)2^{-k}(h-g)+g] )$ and $\zeta_{\ell_2}([i2^{-k}N, (i+1)2^{-k}N] \times [j2^{-k}(h-g)+g, (j+1)2^{-k}(h-g)+g] )$ are both nonzero, we must have
\[
|\ell_1-\ell_2| < C\delta^{-0.01}2^{-k/3}\log^{2/3}(1+2^k) < \delta^{-1}2^{-k/3}\log^{2/3}(1+2^k),
\]
where the second inequality is by taking $\delta$ small enough.
Finally, by taking $Q=CJ$ the conclusion follows.
\end{proof}

This lemma immediately implies Proposition \ref{prop:2d-53-H\"{o}lder}. Recall from Section \ref{s:iop} the general strategy to prove lower bound for Hausdorff dimension, Proposition \ref{prop:2d-53-H\"{o}lder} and \eqref{nondeg1} imply that almost surely the Hausdorff dimension of $\NC$ is at least $5/3$.

\subsection{Proof of the upper bound}
\label{ssec:53-upper}

We next upper bound the Hausdorff dimension of the set $\NC$. This is a more straightforward argument involving  splitting $\HH$ into boxes, and computing the expected number of boxes intersecting $\NC$.

\begin{prop} \label{prop:upper-2D}
For any $0<g<h$, $N>0$, and any $d>5/3$, almost surely the $d$-dimensional Hausdorff measure of the set $\NC\cap ([-N,N]\times [g, h])$ is zero.
\end{prop}

The following is the key estimate relying on non-coalescence of geodesics similar to the argument in \cite{BGH}.
\begin{prop}   \label{prop:empty-prob}
For any $0<g<h$ and $N>0$, there is a constant $C$, such that for any $g\le t< t+\epsilon \le h$ and $-N \le x < x+\epsilon \le N$, we have
\[
\PP[ \NC \cap ((x, x+\epsilon) \times (t, t+\epsilon)) \neq \emptyset ] < C\epsilon^{1/3} e^{C |\log(\epsilon)|^{5/6}}.
\]
\end{prop}

Before jumping into the proof we take a moment to explain the form of the above bound.
For any geodesic from $(-1, 0)$ or $(1, 0)$ to $(x, x+\epsilon) \times (t, t+\epsilon)$, the fluctuation of its intersection with the line $\R\times \{t\}$ would be of the order of $\epsilon^{2/3}$.
Thus the event $\NC\cap ((x, x+\epsilon) \times (t, t+\epsilon))\neq \emptyset$ would roughly imply disjointness of geodesics from $(-1, 0)$ and $(1, 0)$ to $(x-\epsilon^{2/3},t)$ and $(x+\epsilon^{2/3},t)$, respectively.
Theorem \ref{thm:disj} implies that the probability is of the order of $\epsilon^{1/3}$, with the extra sub-polynomial factor.

Assuming Proposition \ref{prop:empty-prob}, the proof of Proposition \ref{prop:upper-2D} is quick.
\begin{proof}[Proof of Proposition \ref{prop:upper-2D}]
Take any large enough $n\in \N$ and let
\[
\sS_n = \{i, j\in \Z: ([i/n, (i+1)/n] \times [j/n, (j+1)/n]) \cap \NC \cap ([-N,N]\times [g, h]) \neq \emptyset\}.
\]
By Proposition \ref{prop:empty-prob}, we have that $\E[|\sS_n|] < C n^{5/3}e^{C(\log(n))^{5/6}}$, where $C$ is a constant depending only on $g, h, N$.
Taking $0<\alpha<d-5/3$, we then have that
\[
\PP[n^{-d}|\sS_n|>n^{-\alpha}] < C n^{5/3-d+\alpha}e^{C(\log(n))^{5/6}}.
\]
This implies by the Borel-Cantelli lemma, that almost surely, we can find a random $k_0$, such that for any integer $k>k_0$,
\[
2^{-dk}|\sS_{2^k}| \le 2^{-\alpha k}.
\]
Since the right hand side converges to $0$ as $k\to\infty$, and that $\NC\cap ([-N,N]\times [g, h])\subset \cup_{i\in\sS_n} [i/n, (i+1)/n]\times [j/n, (j+1)/n]$, the conclusion follows. 
\end{proof}

We now finish the proof of Proposition \ref{prop:empty-prob}. 
Recall that for any $x_1< x_2$, $y_1<y_2$, and $s<t$, $\Dist_{(x_1, x_2)\to (y_1, y_2)}^{[s, t]}$ denotes the event where $\pi_{(x_1, s), (y_1, t)}$ and $\pi_{(x_2, s), (y_2, t)}$ are unique and disjoint.

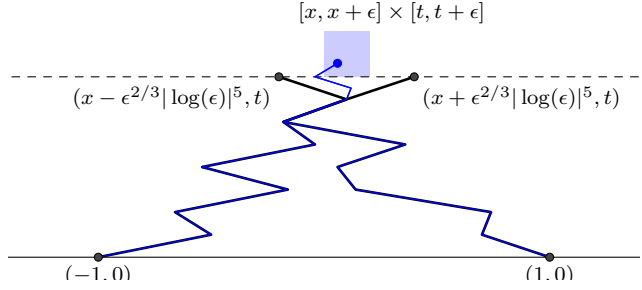
\begin{figure}[hbt!]
    \centering
\begin{tikzpicture}[line cap=round,line join=round,>=triangle 45,x=0.6cm,y=0.3cm]
\clip(-2,-1) rectangle (12,11.5);

\fill[line width=0.pt,color=blue,fill=blue,fill opacity=0.2]
(5,8) -- (6,8) -- (6,10) -- (5,10) -- cycle;

\draw (-10,0) -- (30,0);
\draw [dashed] (-10,8) -- (30,8);

\draw [line width=1pt] (0,0) -- (2.5,1) -- (1.7,2) -- (4.2,3) -- (2.3,4) -- (4.8,5) -- (4.1,6) -- (5.5,7) -- (4,8);
\draw [line width=1pt] (10,0) -- (8.5,1) -- (8.7,2) -- (5.7,3) -- (5.3,4) -- (6.8,5) -- (4.1,6) -- (5.5,7) -- (7,8);

\draw [line width=0.5pt,color=blue] (0,0) -- (2.5,1) -- (1.7,2) -- (4.2,3) -- (2.3,4) -- (4.8,5) -- (4.1,6) -- (5.5,7) -- (5.6,7.5) -- (4.8,8) -- (5.3,8.6);
\draw [line width=0.5pt,color=blue] (10,0) -- (8.5,1) -- (8.7,2) -- (5.7,3) -- (5.3,4) -- (6.8,5) -- (4.1,6) -- (5.5,7) -- (5.6,7.5) -- (4.8,8) -- (5.3,8.6);

\draw [fill=uuuuuu] (0,0) circle (1.5pt);
\draw [fill=uuuuuu] (10,0) circle (1.5pt);

\draw [fill=uuuuuu] (4,8) circle (1.5pt);
\draw [fill=uuuuuu] (7,8) circle (1.5pt);
\draw [blue] [fill=blue] (5.3,8.6) circle (1.5pt);

\begin{scriptsize}
\draw (0,0) node[anchor=north]{$(-1,0)$};
\draw (10,0) node[anchor=north]{$(1,0)$};

\draw (4,8) node[anchor=north east]{$(x-\epsilon^{2/3}|\log(\epsilon)|^{5}, t)$};
\draw (7,8) node[anchor=north west]{$(x+\epsilon^{2/3}|\log(\epsilon)|^{5}, t)$};
\draw (6.5,10) node[anchor=south]{$[x,x+\epsilon]\times [t, t+\epsilon]$};
\end{scriptsize}

\end{tikzpicture}
\caption{An illustration of the proof of Proposition \ref{prop:empty-prob}: the events $\cE_1^c$ and $\cE_2^c$.}
\label{fig:510}
\end{figure}

\begin{proof}[Proof of Proposition \ref{prop:empty-prob}]
We can assume that $\epsilon$ is small enough since otherwise the inequality is obvious by taking the constant large enough.

We consider the following two events.
\begin{enumerate}
    \item Let $\cE_1$ be the event $\Dist_{(-1,1)\to (x-\epsilon^{2/3}|\log(\epsilon)|^{5}, x+\epsilon^{2/3}|\log(\epsilon)|^{5})}^{[0,t]}$.
    \item Let $\cE_2$ be the event that there is some $y\in[x, x+\epsilon]$ and $s\in [t, t+\epsilon]$, such that
    \[|\pi_{(-1,0;y,s)}(t)-x| \vee |\pi_{(1,0;y,s)}(t)-x|>\epsilon^{2/3}|\log(\epsilon)|^{5}.\]
\end{enumerate}
Under $\cE_1^c\cap\cE_2^c$, we have that (almost surely) the function $(y,s)\mapsto \cL(-1,0;y,s)-\cL(1,0;y,s)$ is constant on $[x, x+\epsilon]\times [t, t+\epsilon]$ (see Figure \ref{fig:510} for an illustration).
Indeed, by $\cE_1^c$ we can find some $t_* \in (0, t)$ and $x_*$, such that \[\pi_{(-1,0;x-\epsilon^{2/3}|\log(\epsilon)|^{5},t)}(t_*)=\pi_{(1,0;x+\epsilon^{2/3}|\log(\epsilon)|^{5},t)}(t_*)=x_*.\]
We also assume (the probability $1$ event) that the geodesics $\pi_{(-1,0;x-\epsilon^{2/3}|\log(\epsilon)|^{5},t)}$ and $\pi_{(1,0;x+\epsilon^{2/3}|\log(\epsilon)|^{5},t)}$ are unique.
Then by $\cE_2^c$ and the ordering of geodesics (Lemma \ref{lem:ordering-geo}), we must have that
\[\pi_{(1,0;y,s)}(t_*)=\pi_{(-1,0;y,s)}(t_*)=x_*,\]
for any $(y,s) \in [x, x+\epsilon]\times [t, t+\epsilon]$.
We then have that 
\[
\cL(-1,0;y,s)-\cL(1,0;y,s) = \cL(-1,0;x_*,t_*)-\cL(1,0;x_*,t_*),
\]
which implies that the left hand side is a constant (in $[x, x+\epsilon]\times [t, t+\epsilon]$).
Thus we have that the event $\NC\cap ((x, x+\epsilon)\times (t, t+\epsilon))\neq \emptyset$ (almost surely) implies $\cE_1\cup\cE_2$.

We next bound the probabilities of these events.
Below we let $C, c>0$ denote large and small constants depending on $g, h, N$  whose values may change from line to line.
By Theorem \ref{thm:disj} we have that $\PP[\cE_1] < C\epsilon^{1/3} e^{C |\log(\epsilon)|^{5/6}}$.
For $\PP[\cE_2]$, we take $R$ to be the random varialbe given by Lemma \ref{lem:dl-trans-f}.
Then $\PP[R>|\log(\epsilon)|] < Ce^{-c|\log(\epsilon)|^{9/4}/\log(|\log(\epsilon)|)^4}$, and $R\le |\log(\epsilon)|$ implies $\cE_2^c$.
So we have $\PP[\cE_2] < Ce^{-c|\log(\epsilon)|^{9/4}/\log(|\log(\epsilon)|)^4}$.
Thus $\PP[\cE_1\cup \cE_2] < C\epsilon^{1/3} e^{C |\log(\epsilon)|^{5/6}}$, and the conclusion follows.
\end{proof}

Having completed our analysis of the two dimensional process $\cD,$ we now shift our focus on one dimensional objects for which, as outlined in Section \ref{s:iop}, the arguments are more delicate. We start with GLT and the geodesic zero set.

\section{Constructing GLT: Proof of Theorem \ref{thm:geohaus}}  \label{sec:geo-loc}

Recall from the discussion following the statement of the theorem that we construct a non-decreasing function $L:\R_{\ge 0} \to \R$, which we call GLT, in analogue with the well known BLT.
This $L$ induces a (random) measure $\mu$ on $\R_{\ge 0}$, supported on the set $Z$ defined in \eqref{zeroset12}.
We show that  almost surely, $L$ is $1/3-\epsilon$ H\"{o}lder for any $\epsilon >0$, and that it is almost surely non-degenerate, i.e., $L>0$ on $\R_+$. This implies that the set $Z$ has Hausdorff dimension at least $1/3$.
The upper bound of the Hausdorff dimension is easier, by splitting the line $\{(0,t):t\ge 0\}$ into segments, and bounding the expected number of segments that intersect with $\pi_\boo$.

In Section \ref{ssec:geo-ocp-mea} we construct $L$.
The H\"{o}lder property of $L$ is proved in Section \ref{ssec:geo-ocp-H\"{o}lder}.
In Section \ref{ssec:L-non-zero} we show that $L$ is almost surely non-zero.
Finally we prove the Hausdorff dimension upper bound in Section \ref{ssec:geo-haus-upper}.

\subsection{Construction: occupation measure and the limit}  \label{ssec:geo-ocp-mea}
We recall from \eqref{occmeasure1} the occupation measure $\nu_h$ for the geodesic $\pi_\boo$ at height $h$, where for any measurable $A\subset \R$ 
\begin{equation}\label{occmeasure2}
\nu_h(A) = \int_0^h \don[\pi_\boo(t)\in A] dt.
\end{equation}
The function $L$ will be constructed as in the classical case of Brownian motion by showing the existence of the density of $\nu_h$ at $0$.
For this, we prove Proposition \ref{prop:om-limit}, which asserts that almost surely the limit $\lim_{w\searrow 0} (2w)^{-1}\nu_h([-w, w])$ exists and is finite. In fact, in the forthcoming work \cite{GZ22}, we will show that properly normalized counting measure on the zero set of the geodesic in the discrete prelimiting model of Exponential LPP converges to $L,$ similar to how it is classically known that the random walk local time converges to BLT. 

Reiterating the outline from Section \ref{s:iop},
note that by Lebesgue's theorem for the differentiability of monotone functions, for any (locally finite) measure on $\R,$ its density exists almost everywhere. However, for our purposes we would need it to exist at $0.$ To ensure this we will construct the occupation measure $\nu_h'$ of semi-infinite geodesics starting from random points. The latter has the advantage of being translation invariant which is enough to imply that almost surely its density exists at a fixed point.
We will then use coalescence of geodesics to compare $\nu_h$ and $\nu_h'$, and deduce the same for $\nu_h$.

We start by recording the following bound on the number of intersections of semi-infinite geodesics with a line segment in the space direction.
\begin{lemma}  \label{lem:semi-cross}
For any $t>0$, and interval $A\subset \R$, we have \[\E[|\{\pi_{(x,0)}(t): x\in \R, \pi_{(x,0)} \text{ is unique}\} \cap A|]<C\LE(A)t^{-2/3},\] where $C$ is a universal constant.
\end{lemma}
We start by briefly explaining the occurrence of the factor of $t^{-2/3}$: by the KPZ scaling (of the directed landscape $\cL$), we expect $\pi_{(x,0)}$ for $x$ varying in an interval of length (of the order of) $t^{2/3}$ to coalesce before height $t$.
Thus there would be constant many points in $\{\pi_{(x,0)}(t): x\in \R, \pi_{(x,0)} \text{ is unique}\}$, in an interval of length $t^{2/3}$.

Next we prove Proposition \ref{prop:om-limit} assuming Lemma \ref{lem:semi-cross}.

\begin{figure}
     \centering
     \begin{subfigure}[t]{0.45\textwidth}
         \centering
\begin{tikzpicture}
[line cap=round,line join=round,>=triangle 45,x=0.5cm,y=0.17cm]
\clip(-1,-2) rectangle (11,19.5);
\draw (-1,0) -- (11,0);

\draw [fill=uuuuuu] (0.5,0) circle (1.5pt);
\draw [fill=uuuuuu] (4,0) circle (1.5pt);
\draw [fill=uuuuuu] (9,0) circle (1.5pt);

\draw [line width=0.6pt] (0.5,0) -- (0.2,1) -- (1.9,2) -- (0.7,3) -- (-.1,4) -- (.7,5) -- (-.8,6) -- (1.2,7) -- (-.5,8) -- (1.8,9) -- (.8,10) -- (2.1,11) -- (.2,12) -- (1.1,13) -- (-.2,14) -- (2.7,15) -- (.7,16) -- (1.6,17) -- (.9,18) -- (2.5,19);
\draw [line width=0.6pt] (4,0) -- (4.2,1) -- (2.9,2) -- (3.7,3) -- (2.1,4) -- (3.7,5) -- (1.8,6) -- (2.2,7) -- (1.2,8) -- (1.8,9);
\draw [line width=0.6pt] (9,0) -- (8.2,1) -- (9.9,2) -- (8.7,3) -- (8.1,4) -- (10.7,5) -- (9.8,6) -- (7.2,7) -- (9.5,8) -- (8.8,9) -- (9.8,10) -- (9.1,11) -- (10.2,12) -- (8.1,13) -- (9.2,14) -- (7.7,15) -- (8.7,16) -- (7.6,17) -- (7.9,18) -- (6.5,19);

\begin{scriptsize}
\draw (2,0) node[anchor=north]{$\xi$};
\end{scriptsize}

\end{tikzpicture}
         \caption{Semi-infinite geodesics starting from $\xi \times \{0\}$.}
         \label{fig:62p}
     \end{subfigure}
     \hfill
     \begin{subfigure}[t]{0.5\textwidth}
         \centering
\begin{tikzpicture}
[line cap=round,line join=round,>=triangle 45,x=0.5cm,y=0.17cm]
\clip(-1,-2) rectangle (11,19.5);

\draw (-1,0) -- (11,0);

\draw [|-|] [blue] [thick] (4,0) -- (6,0);
\draw [|-|] [blue] [thick] (-10,0) -- (0,0);
\draw [|-|] [blue] [thick] (10,0) -- (19,0);

\draw [fill=blue] (5,0) circle (1.5pt);
\draw [fill=blue] (-.7,0) circle (1.5pt);
\draw [fill=blue] (10.8,0) circle (1.5pt);

\draw [fill=uuuuuu] (5.5,0) circle (1.5pt);

\draw [blue] [line width=1pt] (5,0) -- (5.2,1) -- (4.9,2) -- (5.7,3) -- (5.1,4) -- (6.7,5) -- (4.8,6) -- (6.2,7) -- (5.5,8) -- (6.8,9) -- (4.8,10) -- (5.1,11) -- (4.2,12) -- (5.1,13) -- (4.2,14) -- (4.7,15) -- (5.1,16) -- (4.6,17) -- (5.9,18) -- (4.5,19);
\draw [line width=0.6pt] (5.5,0) -- (5.4,1) -- (6.3,2) -- (5.7,3) -- (5.1,4) -- (6.7,5) -- (4.8,6) -- (6.2,7) -- (5.5,8) -- (6.8,9) -- (4.8,10) -- (5.1,11) -- (4.2,12) -- (5.1,13) -- (4.2,14) -- (4.7,15) -- (5.1,16) -- (4.6,17) -- (5.9,18) -- (4.5,19);

\draw [line width=0.6pt] (-0.7,0) -- (0.2,1) -- (-.8,2) -- (0.7,3) -- (-.1,4) -- (-.7,5) -- (.2,6) -- (-1,7);
\draw [line width=0.6pt] (10.8,0) -- (10.2,1) -- (10.9,2) -- (9.7,3) -- (10.1,4) -- (10.7,5) -- (11.8,6) -- (10.2,7) -- (10.5,8) -- (9.8,9) -- (10.8,10) -- (11.1,11) -- (10.2,12) -- (11.1,13) -- (10.2,14) -- (11.7,15) ;

\begin{scriptsize}
\draw (2,0) node[anchor=north]{$\xi$};
\draw [blue] (4,0) node[anchor=north]{$-a$};
\draw [blue] (5,0) node[anchor=north]{$\boo$};
\draw [blue] (6,0) node[anchor=north]{$a$};
\draw [blue] (0,0) node[anchor=north]{$-a^{-1}$};
\draw [blue] (10,0) node[anchor=north]{$a^{-1}$};
\end{scriptsize}

\end{tikzpicture}        \caption{$\pi_\boo$ and semi-infinite geodesics starting from $\xi \times \{0\}$, conditional on $\cE_a$.}
         \label{fig:62f}
     \end{subfigure}
        \caption{An illustration of the proof of Proposition \ref{prop:om-limit}.}
        \label{fig:62}
\end{figure}
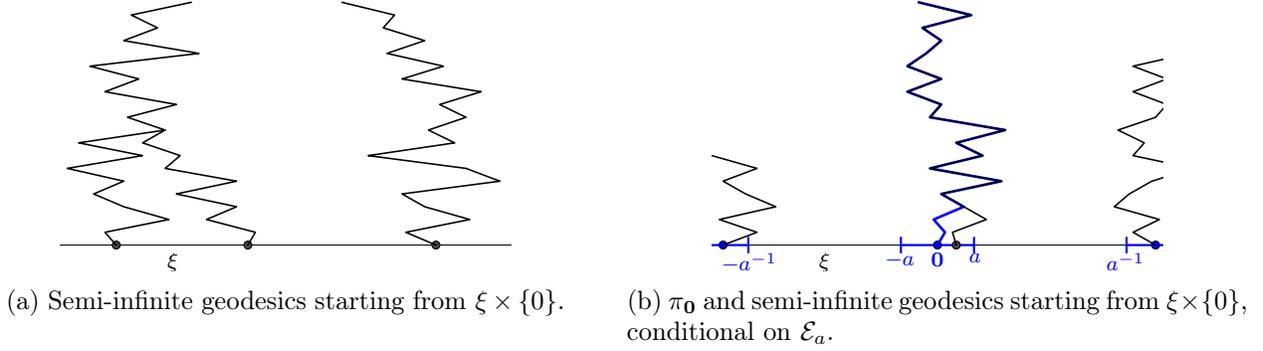

\begin{proof}[Proof of Proposition \ref{prop:om-limit}]
The proof goes in two steps: in the first step we prove that the density of $\nu_h-\nu_g$ exists at $0$ for any $0<g<h$, using a translation invariant argument; in the second step we bound $\limsup_{w\searrow 0}(2w)^{-1}\nu_g([-w, w])$, and show that it decays to zero as $g\searrow 0$.\\

\noindent
\textbf{Step 1.}
We start with the first step, which is illustrated by Figure \ref{fig:62}.
Let $\xi \subset \R$ be a Poisson point process with rate $1$. We consider the measure $\nu_h'$ such that for any measurable $A\subset \R$,
\[
\nu_h'(A) = \int_0^h \sum_{x\in\xi}\don[\pi_{(x,0)}(t)\in A] dt.
\]
Recall that for any (locally finite) measure, the density exists and is finite Lebesgue almost everywhere, by Lebesgue's theorem for the differentiability of monotone functions.
We next show that $\nu_h'$ is locally finite.
Note that under the probability one event where $\pi_{(x,0)}$ is unique for any $x\in\xi$, we have $\nu_h'(A) \le \int_0^h |\{\pi_{(x,0)}(t): x\in \R, \pi_{(x,0)} \text{ is unique}\} \cap A| dt$. So by Lemma \ref{lem:semi-cross} we have that almost surely $\nu_h'(A)<\infty$ for any compact $A$.
Then if we let $E$ be the set where such limit $\lim_{w\searrow 0} (2w)^{-1}\nu_h'([x-w, x+w])$ does not exist or is infinite, we would have that $E$ is Lebesgue measurable\footnote{The fact that $E$ is measurable can be argued as follows: for each $w>0$, the function $x\mapsto (2w)^{-1}\nu_h'([x-w, x+w])$ is continuous except for countably many points, thus it is measurable. Then $E$ is measurable, since it is the set where a sequence of measurable functions does not converge or converges to infinity.} with $\LE(E)=0$.
By translation invariance of the law of $\nu_h'$ and Fubini's theorem, we have $\PP[0\in E] = \E[\LE(E\cap [0,1])] = 0$.
Thus we have that almost surely the limit $\lim_{w\searrow 0} (2w)^{-1}\nu_h'([-w, w])$ exists and is finite.

We next compare $\nu'$ and $\nu.$ Take $0<g<h$.
For $0<a<1$, consider the following event $$\cE_a:=\{|\xi \cap [-a, a]| = |\xi \cap [-a^{-1}, a^{-1}]| = 1\}.$$
The idea being that on $\cE_a,$ the presence of only a single point near the origin would allow us to carry out the desired comparison.  
Since $\PP[\cE_a]>0$, 
conditional on $\cE_a$ we still have that almost surely the limit $\lim_{w\searrow 0} (2w)^{-1}(\nu_h'([-w, w])- \nu_g'([-w, w]))$ exists and is finite.

We note that, almost surely, as $x\to 0$, the geodesic $\pi_{(x,0)}\to\pi_\boo$ in the overlap topology (by Lemma \ref{lem:semi-overlap}).
Also, by Lemma \ref{lem:dl-semi-trans} and the Borel-Cantelli lemma, we have that almost surely $|\pi_{(x,0)}(h)| \to \infty$ as $|x|\to \infty$.
Then by Lemma \ref{lem:dl-trans-f}, and ordering of finite geodesics, we further have that almost surely $\min_{t\in[0,h]}|\pi_{(x,0)}(t)| \to \infty$ as $|x|\to \infty$.
Hence, given $0<g<h$, almost surely for all small enough $a$, on $\cE_a$,  $$\nu_h([-w, w])- \nu_g([-w, w]=\nu'_h([-w, w])- \nu'_g([-w, w].$$
Thus we have that (almost surely for the directed landscape) the limit $\lim_{w\searrow 0} (2w)^{-1}(\nu_h([-w, w])- \nu_g([-w, w]))$ exists and is finite.\\

\noindent
\textbf{Step 2.}
We next carry out the second step, i.e., to bound $\limsup_{w\searrow 0}(2w)^{-1}\nu_g([-w, w])$. 
The idea is to bound $\nu_g$ by the occupation measure of the union of $\pi_{(x,0)}$ \textit{for all} $x\in\R$, and to apply Lemma \ref{lem:semi-cross}.
We begin by denoting by $\nu_g''$ the measure such that for any Borel measurable set $A,$
\[
\nu_g''(A) = \int_0^g
|\{\pi_{(x,0)}(t): x\in \R, \pi_{(x,0)} \text{ is unique}\} \cap A| dt.
\]
From this definition, almost surely $\nu_g(A)\le \nu''_g(A)$ for any compact $A$. Also, for compact $A$, by Lemma \ref{lem:semi-cross} we have $\E[\nu_g''(A)] < \infty$, so almost surely we have $\nu_g''(A)<\infty$.
Since the density of $\nu_g''$ exists Lebesgue almost everywhere, by translation invariance of the law of $\nu_g''$, we have that almost surely $\lim_{w\searrow 0}(2w)^{-1}\nu_g''([-w, w])$ exists.
Assuming this limit exists, we have
\[\limsup_{w\searrow 0}(2w)^{-1}\nu_g([-w, w])\le \lim_{w\searrow 0}(2w)^{-1}\nu_g''([-w, w]).\]
We next bound the upper tail of this limit.
Let $E_z\subset \R$ consist of all $x$, where $\lim_{w\searrow 0}(2w)^{-1}\nu_g''([x-w, x+w])$ exists and is at least $z$.
Then $E_z$ is measurable, by the same arguments as proving that $E$ is measurable (see the footnote above).
Since the law of $\nu_g''$ is translation invariant, we have $\PP[0\in E_z] = (2l)^{-1}\E[\LE(E_z\cap [-l, l])]$ for any $l>0$.
We also have\footnote{The first inequality is by the following argument: for any $w>0$, we have 
$\nu_g''([-l-w,l+w]) \ge \int_{-l}^l (2w)^{-1}\nu_g''([x-w, x+w]) dx$.
Then the inequality follows by sending $w\to 0$ and using Fatou's lemma for the functions $x\mapsto (2w)^{-1}\nu_g''([x-w, x+w])$.}
\[
\nu_g''([-l,l]) \ge \int_{E_z\cap[-l,l]}\lim_{w\searrow 0}(2w)^{-1}\nu_g''([x-w, x+w]) dx \ge z\LE(E_z\cap [-l, l]),
\]
so
\[
\PP[0\in E_z] \le (2l)^{-1}z^{-1} \E[\nu_g''([-l,l])] < Cg^{1/3}z^{-1}
\]
for some constant $C>0$.
Thus we conclude that
\[
\PP[\limsup_{w\searrow 0}(2w)^{-1}\nu_g([-w, w]) > z] < Cg^{1/3}z^{-1}.
\]
By the Borel-Cantelli lemma, this implies that as $g\searrow 0$, $\limsup_{w\searrow 0}(2w)^{-1}\nu_g([-w, w]) \to 0$ almost surely. Thus we get Proposition \ref{prop:om-limit}.
\end{proof}
It remains to prove Lemma \ref{lem:semi-cross}.

\begin{proof}[Proof of Lemma \ref{lem:semi-cross}]
Recall that to justify the factor of $t^{-2/3}$, we would need to prove coalescence of semi-infinite geodesics, and for that we rely on Lemma \ref{lem:number-of-inter}.

By scaling invariance it suffices to prove this result for $t=1$; and by translation invariance it suffices to consider a fixed interval $A=[-1,1]$;
i.e. we just need to show that $\E[|\{\pi_{(x,0)}(1): x\in \R, \pi_{(x,0)} \text{ is unique}\} \cap [-1,1]|]$ is finite.
As throughout the article, below we use $C,c>0$ to denote large and small constants whose values may change from line to line.

Take $m\in\N$ and we consider $\pi_{(-m,0)}$ and $\pi_{(m,0)}$.
By passing Lemma \ref{lem:trans-semi-inf} to the limit using Proposition \ref{prop:cov-semi-geo}, we have that
\begin{equation}  \label{eq:semi-c-bd}
\PP[\sup_{s\in[0,2]} |\pi_{(-m,0)}(s)+m| > m/2 ]=\PP[\sup_{s\in[0,2]} |\pi_{(m,0)}(s)-m| > m/2 ] < Ce^{-cm^3}.
\end{equation}
By the ordering of semi-infinite geodesics (Lemma \ref{lem:semi-inf-tree}), when $\sup_{s\in[0,2]} |\pi_{(-m,0)}(s)+m|\le m/2$ and $\sup_{s\in[0,2]} |\pi_{(m,0)}(s)-m| \le m/2$, we must have that $\sup_{s\in[0,2], x\in[-m,m]} |\pi_{(x,0)}(s)| \le 3m/2$, and hence 
\begin{multline*}
|\{\pi_{(x,0)}(1): x\in \R, \pi_{(x,0)} \text{ is unique}\} \cap [-1,1]| \\ \le |\{\pi_{(x,0),(y,2)}(1): |x|,|y|\le 3m/2,  \pi_{(x,0),(y,2)} \text{ is unique}\} \cap [-1,1]|.   
\end{multline*}
Then by \eqref{eq:semi-c-bd} and Lemma \ref{lem:number-of-inter}, we have that $|\{\pi_{(x,0)}(1): x\in \R, \pi_{(x,0)} \text{ is unique}\} \cap [-1,1]|$ has a stretched exponential tail, so the conclusion holds.
\end{proof}
We record the following useful immediate corollary.
\begin{cor}  \label{cor:cross-int}
For any $t>0$ and interval $A\subset \R$, we have
$\PP[\pi_\boo(t)\in A] < C\LE(A)t^{-2/3}$, where $C$ is a universal constant.
\end{cor}
Such an estimate in the pre-limit has appeared as \cite[Theorem 2]{BB21} where the proof also uses translation invariance of the models.

We finish this subsection by recording the following absolute continuity property of the occupation measure $\nu_h$, which will be used later.
\begin{prop}  \label{prop:abs-cont-leb}
Almost surely the measure $\nu_h$ is absolutely continuous to the Lebesgue measure for any $h>0$.
\end{prop}
The proof is similar to the one showing absolute continuity for the Brownian occupation measure, in studying BLT (see e.g., \cite[Theorem 3.25]{morters2010brownian}).
\begin{proof}
Since $\nu_h$ is monotone in $h$, it suffices to prove this for fixed $h>0$.
By an equivalent definition of absolute continuity, we need to show that almost surely, for $\nu_h$-almost every $x\in \R$, we have
$$\liminf_{w\searrow 0}\frac{\nu_h([x-w,x+w])}{w} < \infty.$$
We take arbitrary $H>0$, and by Fatou's lemma we have
\begin{align*}
&\E \int_{-H}^H \liminf_{w\searrow 0}\frac{\nu_h([x-w,x+w])}{w}  d\nu_h(x) 
\\
\le & \liminf_{w\searrow 0}w^{-1} \E \int_{-H}^H \nu_h([x-w,x+w]) d\nu_h(x) \\
= & \liminf_{w\searrow 0}w^{-1} \int_0^h \int_0^h \PP[|\pi_\boo(t)|<H, |\pi_\boo(t)-\pi_\boo(s)|<w] dt ds.
\end{align*}
For each $i\in \Z$, we have
\[
\int_0^h \int_0^h \PP[\pi_\boo(t), \pi_\boo(s)\in [(i-1)w, (i+1)w]] dtds
=\E[\nu_h([(i-1)w, (i+1)w])^2] \le C w^2,
\]
where $C>0$ is a constant depending on $h, H$, and the inequality is by Lemma \ref{lem:bound-zb}.
Then we have
\begin{multline*}
\int_0^h \int_0^h \PP[|\pi_\boo(t)|<H, |\pi_\boo(t)-\pi_\boo(s)|<w] dt ds
\\ \le \sum_{i\in\Z: w(|i|-1)\le H} \int_0^h \int_0^h \PP[\pi_\boo(t), \pi_\boo(s)\in [(i-1)w, (i+1)w]] dtds
\le (2\lfloor H/w\rfloor +3)C w^2 <  3CHw
\end{multline*}
when $w$ is small enough. 
Thus we conclude that $\nu_h$ is absolutely continuous with respect to Lebesgue measure.
\end{proof}

\subsection{$1/3$ H\"{o}lder regularity of GLT}  \label{ssec:geo-ocp-H\"{o}lder}
By Proposition \ref{prop:om-limit}, almost surely we can define $L(h)=\lim_{w\searrow 0} (2w)^{-1}\nu_h([-w, w])$ for any non-negative rational $h$.
The following H\"{o}lder estimate of $L$ in particular implies continuity and will allow us to extend $L$ to $\R_{\ge 0}$, and conclude that $L$ is almost surely $\frac{1}{3}-\epsilon$ H\"{o}lder for any $\epsilon>0$.
\begin{prop} \label{prop:mea-L-tail}
There exist constants $c,C>0$ such that the following is true.
For any $g,h \in \Q$, $0\le g < h$, and $M>0$, we have
\[
\PP[L(h)-L(g)>M(h-g)^{1/3}] < Ce^{-cM}.
\]
\end{prop}

In fact, we will prove a more general estimate.
Namely, for some intervals $I, J$, we will obtain a uniform bound on the occupation measure in the set $I\times J$, for any geodesic whose endpoints are at the boundary of $I\times J$.
It would imply both Proposition \ref{prop:mea-L-tail} and Lemma \ref{lem:bound-zb}, and have many other applications in the remainder of this article.

For any $u=(p;q)=(x,s;y,t) \in \R^4_\uparrow$, recall that $\pi_u$ denotes any geodesic from $p=(x,s)$ to $q=(y,t)$. For any interval $I$ we define 
\[W[u, I] = \int_s^t \don[\pi_u(r) \in I, \pi_u \text{ is unique}] dr.\]
In words, $W[u, I]$ is the total amount of time $\pi_u$ spends in $I$, if $\pi_u$ is the unique geodesic; and $W[u, I]=0$ otherwise.

For any intervals $I, J$, we denote
\begin{equation}  \label{eq:defn-phi-ij}
\Phi_{I, J} = \{(p;q) \in \R^4_\uparrow: p, q \in \partial(I\times J)\},
\end{equation}
where $\partial(I\times J)$ denotes the boundary of $I\times J$, i.e., $\partial(I\times J)=(\{\inf I, \sup I\}\times J) \cup (I \times \{\inf J, \sup J\})$,
and
\[
W_{I,J}=\sup_{u \in \Phi_{I, J}} W[u, I].
\]
We note that $W_{I,J}$ is also an upper bound for the occupation measure in $I\times J$, for any unique geodesic $\pi_u$, where $u=(p;q)=(x,s;y,t) \in \R^4_\uparrow$ with $p, q \not\in I\times J$. Here this occupation measure is defined as
\[
\int_{[s,t]\cap J} \don[\pi_u(r) \in I] dr.
\]
Indeed, if the graph of $\pi_u$ intersects $I\times J$, the above expression equals $W[(p';q'),I]$, for $p',q'$ being the first and last intersection of the graph of $\pi_u$ with $I\times J$. Since $\pi_u$ is unique, $\pi_{(p';q')}$ must also be unique, so $W[(p';q'),I]\le W_{I,J}$.
Similarly, for any $p\not\in I\times J$, if $\pi_p$ is unique, its occupation measure in $I\times J$ is also upper bounded by $W_{I,J}$. We now state the following uniform estimate.
\begin{lemma}  \label{lem:geo-ocp-tail}
 There are universal constants $C,c>0$, such that for any $w, M>0$ we have
 \[
 \PP[W_{[-w,w], [0,1]} > Mw] < Ce^{-cM}.
 \]
\end{lemma}
Its proof follows induction arguments that are essentially the same as those in \cite[Section 3]{SSZ}, where a version for exponential LPP is proved.
We provide the proof adapted to the directed landscape setting in Appendix \ref{sec:appb}.

\subsection{Non-degeneracy of GLT}  \label{ssec:L-non-zero}
The goal of this section is to show that $L$ is strictly positive on $\R_+$.
\begin{prop}  \label{prop:L-not-dege}
Almost surely, $L(h)>0$ for any $h>0$.
\end{prop}
Before proving the above, we quickly prove the lower bound of the Hausdorff dimension.
\begin{proof}[Proof of Theorem \ref{thm:geohaus}: lower bound]
By Proposition \ref{prop:mea-L-tail} $L$ is continuous, and we can define a measure $\mu$ on $\R_{\ge 0}$, via $\mu([g,h]) = L(h)-L(g)$ for any $0\le g < h$.
We next claim that $\mu$ is supported on $Z$ (although we believe that the support is exactly equal to $Z,$ the proof of this unfortunately still eludes us).
Indeed, as $\pi_\boo$ is continuous, for any $t\in \R_+\setminus Z$ we can find some $\delta>0$ with $|\pi_\boo(s)|>\delta$ for any $|s-t|\le \delta$.
Thus $\nu_{t+\delta}([-\delta, \delta]) - \nu_{t-\delta}([-\delta, \delta]) = 0$, and hence $L(t+\delta)=L(t-\delta)$, leading to $\mu([t-\delta, t+\delta]) = 0$.

Take any $h, \epsilon>0$.
By Proposition \ref{prop:L-not-dege} we have (almost surely) $\mu([0,h])>0$.
By Proposition \ref{prop:mea-L-tail}, (almost surely) there is a constant $C>0$ such that $L([a,b])<C(b-a)^{1/3-\epsilon}$ for any $0\le a < b \le h$.
Hence by the general strategy from Section \ref{s:iop}, the Hausdorff dimension of $Z\cap [0,h]$ is at least $1/3-\epsilon$.
Since $h, \epsilon$ are arbitrary, the Hausdorff dimension of $Z$ is at least $1/3$.
\end{proof}
The proof of Proposition \ref{prop:L-not-dege} is via the monotonicity and scaling invariance of the function $L$, and the following 0-1 law of the tail sigma algebra of the directed landscape.

For any $r \in \N,$ we let $\cF_{\ge r}$ be the sigma algebra generated by $\cL(x,s;y,t)$ for $r\le s<t$, and let $\cF_\infty = \cap_{r>0} \cF_{\ge r}$. 
\begin{lemma}  \label{lem:F-inf-tri}
Any $\cF_\infty$-measurable event has probability $0$ or $1$.
\end{lemma}
\begin{proof}
This follows from the definition of $\cF_\infty$. For any $(x,s;y,t)\in\R^4_\uparrow$, and $r>t$, $\cL(x,s;y,t)$ is independent of $\cF_{\ge r}$. 
This implies that $\cL$ is independent of $\cF_\infty$. So any $\cF_\infty$-measurable event is independent of itself, thus has probability $0$ or $1$.
\end{proof}

\begin{proof}[Proof of Proposition \ref{prop:L-not-dege}]
By scaling invariance and monotonicity of the function $L$, it suffices to show that $\PP[L(1)>0]=1$.
We first show that $\PP[L(1)>0]>0$.
For this it suffices to show that $\E[L(1)]>0$.
We note that by Lemma \ref{lem:geo-ocp-tail}, the random variables $(2w)^{-1}\nu_1([-w, w])$ for all small enough $w$ are uniformly integrable.
This implies that $\E[(2w)^{-1}\nu_1([-w, w])] \to \E[L(1)]$ as $w\searrow 0$; so it suffices to lower bound $\E[(2w)^{-1}\nu_1([-w, w])]$. 

By Lemma \ref{lem:lbound-geo-intv} below and scaling invariance, we have that $\PP[\pi_\boo(t) \in [-w, w]] > cw$ for any $0\le t \le 1$; and by integrating over $t$ we get $\E[(2w)^{-1}\nu_1([-w, w])]>c/2$ for any $0<w<1$. We thus conclude that $\E[L(1)]>0$ and hence $\PP[L(1)>0]>0$.

Now suppose that $\PP[L(1)>0]=\alpha$ for some $0<\alpha<1$.
We then claim that conditional on $L(1)>0$, we have $\lim_{h\to\infty} L(h)-L(g) > 0$ for any $g>0$.
Indeed, we take
\[
g_* = \inf\{g>0: L(h) = L(g), \forall h>g \} \cup \{\infty \}.
\]
By scaling invariance of $L$, the probability $\PP[g_*>g]$ is same for any $g>0$; thus we have $\PP[g<g_*\le h]=0$ for any $0<g<h$, so $\PP[0<g_*<\infty] = 0$. Then conditional on $L(1)>0$ we have that $g_*=\infty$ almost surely. 
Also by scaling invariance we have that $\PP[L(h)>0]=\alpha$ for any $h>0$, so by monotonicity of $L$ we have that conditional on $L(1)>0$, almost surely  $L(h)>0$ for any $h>0$, and $g_*=\infty$.
Then we conclude that the events $L(1)>0$ and $g_*=\infty$ are almost surely equivalent; i.e., the probability of exactly one of them happening equals zero.

We will now show $g_*=\infty$ happens almost surely.
Towards this for any $r\in \N$ and $h>0$, we denote $\nu_h^r$ as the occupation measure for $\pi_{(0,r)}$ at height $h+r$; i.e., for any measurable $A\subset \R$, we let
\[
\nu_h^r(A) = \int_0^h \don[\pi_{(0,r)}(r+t)\in A] dt,
\]
and we let $L^r(h)=\lim_{w\searrow 0} (2w)^{-1}\nu_h^r([-w, w])$. Almost surely this limit exists and is finite, by Proposition \ref{prop:om-limit} and translation invariance.
Then we can define $g_*^r = \inf\{g>0: L^r(h) = L^r(g), \forall h>g \} \cup \{\infty \}$.
By the coalescence of semi-infinite geodesics (Lemma \ref{lem:semi-inf-tree-fixed}), we have that $g_*=\infty$ is almost surely equivalent to $g_*^r=\infty$, thus to $L^r(1)>0$, for any $r\in \N$.
Thus we can consider the $\cF_\infty$ measurable event $|\{r\in \N: L^r(1)>0\}|=\infty$.
By the arguments above, this is almost surely equivalent to $L(1)>0$.
These imply that $\PP[|\{r\in \N: L^r(1)>0\}|=\infty] = \alpha$.
By Lemma \ref{lem:F-inf-tri}, this contradicts with the fact that $0<\alpha<1$, so we must have that $\PP[L(1)>0]=1$.
\end{proof}
It remains to prove the following technical result. 
\begin{lemma} \label{lem:lbound-geo-intv}
For any $0<w<1$ we have $\PP[\pi_\boo(1) \in [-w, w]] > cw$, where $c>0$ is a universal constant.
\end{lemma}
The matching upper bound for this probability is already given by Corollary \ref{cor:cross-int}. A similar lower bound in the pre-limit is given by \cite[Theorem 2]{BB21} and our arguments are similar.

Unlike the upper bound, such a lower bound does not directly follow from translation invariance. In addition, we need to use the ordering of geodesics, and some additional geometric arguments.

\begin{figure}[t]
         \centering
\begin{tikzpicture}
[line cap=round,line join=round,>=triangle 45,x=0.9cm,y=0.2cm]
\clip(0,-2) rectangle (10,20);

\draw (-1,0) -- (11,0);

\draw [dotted] (-1,14) -- (11,14);

\draw [line width=0.6pt] (2,0) -- (1.2,1) -- (1.9,2) -- (1.1,3) -- (2.1,4) -- (1.1,5) -- (2.8,6) -- (2.2,7) -- (2.5,8) -- (1.8,9) -- (2.8,10) -- (1.6,11) -- (3.2,12) -- (3.1,13) -- (3.4,14) -- (2.7,15) -- (4.1,16) -- (3.8,17) -- (4.5,18) -- (3.9,19);
\draw [line width=0.6pt] (3,0) -- (2.6,1) -- (2.9,2) -- (2.3,3) -- (3.1,4) -- (2.5,5) -- (2.8,6) -- (2.2,7) -- (2.5,8) -- (1.8,9) -- (2.8,10) -- (1.6,11) -- (3.2,12) -- (3.1,13) -- (3.4,14) -- (2.7,15) -- (4.1,16) -- (3.8,17) -- (4.5,18) -- (3.9,19);

\draw [line width=0.6pt] (5,0) -- (4.8,1) -- (5.9,2) -- (5.1,3) -- (6.1,4) -- (5.1,5) -- (5.8,6) -- (5.2,7) -- (5.5,8) -- (4.8,9) -- (4.8,10) -- (4.1,11) -- (5.2,12) -- (4.1,13) -- (3.9,14) -- (4.7,15) -- (4.1,16) -- (3.8,17) -- (4.5,18) -- (3.9,19);

\draw [line width=0.6pt] (6,0) -- (6.8,1) -- (6.4,2) -- (7.1,3) -- (6.7,4) -- (7.3,5) -- (7.2,6) -- (6.2,7) -- (6.5,8) -- (5.3,9) -- (5.8,10) -- (5.1,11) -- (6.1,12) -- (5.9,13) -- (4.7,14) -- (5.7,15) -- (5.4,16) -- (5.8,17) -- (5.1,18) -- (5.6,19);
\draw [line width=0.6pt] (7,0) -- (7.4,1) -- (7.2,2) -- (7.8,3) -- (6.7,4) -- (7.3,5) -- (7.2,6) -- (6.2,7) -- (6.5,8) -- (5.3,9) -- (5.8,10) -- (5.1,11) -- (6.1,12) -- (5.9,13) -- (4.7,14) -- (5.7,15) -- (5.4,16) -- (5.8,17) -- (5.1,18) -- (5.6,19);
\draw [line width=0.6pt] (8,0) -- (9.4,1) -- (8.2,2) -- (8.8,3) -- (7.7,4) -- (8.6,5) -- (8.1,6) -- (7.9,7) -- (8.5,8) -- (7.3,9) -- (7.8,10) -- (6.6,11) -- (7.1,12) -- (6.3,13) -- (6.7,14) -- (6.6,15) -- (5.6,16) -- (5.8,17) -- (5.1,18) -- (5.6,19);

\draw [line width=1pt, color=blue] (4,0) -- (3.7,1) -- (4.3,2) -- (3.6,3) -- (5.1,4) -- (4.2,5) -- (4.8,6) -- (3.7,7) -- (4.5,8) -- (2.8,9) -- (3.3,10) -- (2.6,11) -- (3.2,12) -- (3.1,13) -- (3.4,14) -- (2.7,15) -- (4.1,16) -- (3.8,17) -- (4.5,18) -- (3.9,19);

\draw [fill=uuuuuu] (2,0) circle (1.5pt);
\draw [fill=uuuuuu] (3,0) circle (1.5pt);
\draw [fill=uuuuuu] (4,0) circle (1.5pt);
\draw [fill=uuuuuu] (5,0) circle (1.5pt);
\draw [fill=uuuuuu] (6,0) circle (1.5pt);
\draw [fill=uuuuuu] (7,0) circle (1.5pt);
\draw [fill=uuuuuu] (8,0) circle (1.5pt);

\draw [fill=uuuuuu] (3.4,14) circle (1.5pt);
\draw [fill=uuuuuu] (3.9,14) circle (1.5pt);
\draw [fill=uuuuuu] (4.7,14) circle (1.5pt);
\draw [fill=uuuuuu] (6.7,14) circle (1.5pt);

\begin{scriptsize}
\draw (5,0) node[anchor=north]{$\boo$};
\draw (2,0) node[anchor=north]{$(-a,0)$};
\draw (8,0) node[anchor=north]{$(a,0)$};
\draw (3.2,13.8) node[anchor=east]{$(\pi_{(-a,0)}(1),1)$};
\draw (6.7,14) node[anchor=west]{$(\pi_{(a,0)}(1),1)$};
\end{scriptsize}

\end{tikzpicture}
         \caption{An illustration for the proof of Lemma \ref{lem:lbound-geo-intv}: consider $\pi_{(ai/n,0)}(1)$ for each $i\in\llbracket -n, n\rrbracket$; assuming $\pi_{(-a,0)}(1)>-a$ and $\pi_{(a,0)}(1)<a$, there must exist an $i$ with $|\pi_{(ai/n,0)}(1)-ai/n|\le a/n$ (such semi-infinite geodesic $\pi_{(ai/n,0)}$ is colored blue). Then $\PP[\pi_{(-a,0)}(1)>-a, \pi_{(a,0)}(1)<a] > 0$ implies that $\PP[\pi_\boo(1) \in [-a/n, a/n]]$ is at least in the order of $1/n$.}
         \label{fig:612}
\end{figure}
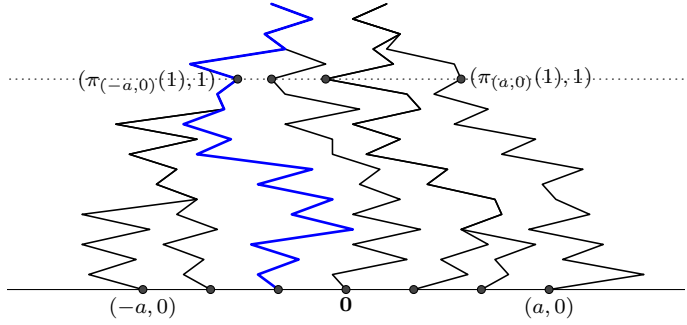
\begin{proof}
We claim that there is some $a\in [1, 2]$ such that
\begin{equation}  \label{eq:le-ri-pos}
\PP[\pi_{(-a,0)}(1)>-a, \pi_{(a,0)}(1)<a] > 0.
\end{equation}
We note that here $[1, 2]$ is chosen arbitrarily and can be replaced by any closed interval in $\R_+$.
We show this by contradiction. 
Using Corollary \ref{cor:cross-int}, and translation and scaling invariance of $\cL$, we have that almost surely $\pi_{(x,0)}(t)\neq x$ for any $x, t \in \Q$.
If the claim is not true, then by translation and scaling invariance, almost surely the following is true: for any $x, y, t\in\Q$, with $t>0$ and $y-x\in [2,4]$, we have $(\pi_{(x,0)}(t)-x)(\pi_{(y,0)}(t)-y)> 0$.
By using this repeatedly, we have that almost surely, $\pi_\boo(t)(\pi_{(x,0)}(t)-x)> 0$ for any $x, t\in \Q$, $t>0$.
Consider $\pi_\boo$ and $\pi_{(1,0)}$: there is some (random) $t_0>0$ such that $\pi_\boo(t)=\pi_{(1,0)}(t)$ for any $t>t_0$.
Then either $\liminf_{t\to\infty}\pi_\boo(t)> 1$ or $\limsup_{t\to\infty}\pi_{(1,0)}(t)< 0$.
By flip invariance, these two events have equal probability, so we have
\begin{equation}  \label{eq:right-event-prob}
0<\PP[\liminf_{t\to\infty}\pi_\boo(t)> 1]<1.    
\end{equation}
However, by coalescing of semi-infinite geodesics, the event $\liminf_{t\to\infty}\pi_\boo(t)> 1$ is (up to a probability zero event) equivalent to the event $\liminf_{t\to\infty}\pi_{(0,r)}(t)> 1$, for any $r\in \R$; thus it is (up to a probability zero event) equivalent to an event that is $\cF_\infty$ measurable. So \eqref{eq:right-event-prob} contradicts Lemma \ref{lem:F-inf-tri}, and we can find an $a\in [1, 2]$ such that \eqref{eq:le-ri-pos} holds.

Now we fix such an $a\in [1, 2]$.
It suffices to prove this lemma for $w=\frac{a}{n}$, $n\in\N$.
We claim that under $\pi_{(-a,0)}(1)>-a, \pi_{(a,0)}(1)<a$, there exists some $i \in \llbracket -n, n \rrbracket$, such that $|\pi_{(ai/n,0)}(1)-ai/n|\le a/n$ (see Figure \ref{fig:612}).
Indeed, otherwise we denote
\[
I_- = \{i \in \llbracket -n, n \rrbracket: \pi_{(ai/n,0)}(1)<a(i-1)/n\}
\]
and
\[
I_+ = \{i \in \llbracket -n, n \rrbracket: \pi_{(ai/n,0)}(1)>a(i+1)/n\}.
\]
We have $-n \in I_+$ and $n\in I_-$, and $\llbracket -n, n \rrbracket \subset I_-\cup I_+$.
Then we can find some $i\in \llbracket -n, n-1 \rrbracket$ such that $i\in I_+$ and $i+1\in I_-$.
These would lead to a contradiction with the ordering of semi-infinite geodesics (Lemma \ref{lem:semi-inf-tree}).
Then we conclude that $|\pi_{(ai/n,0)}(1)-ai/n|\le a/n$ for some $i \in \llbracket -n, n \rrbracket$; thus by translation invariance we have
\[
\PP[\pi_\boo(1) \in [-a/n, a/n]] > (2n+1)^{-1} \PP[\pi_{(-a,0)}(1)>-a, \pi_{(a,0)}(1)<a],
\]
so our conclusion follows.
\end{proof}

\subsection{Hausdorff dimension of the zero set: upper bound}  \label{ssec:geo-haus-upper}
To finish the proof of Theorem \ref{thm:geohaus}, it remains to prove the upper bound for the Hausdorff dimension of the set $Z=\{t\ge 0: \pi_\boo(t)=0\}$. This is immediately implied by the following result. 

\begin{prop} \label{prop:upper-geo}
For any $0<g<h$, and any $d>1/3$, almost surely the $d$-dimensional Hausdorff measure of the set $Z\cap [g, h]$ is zero.
\end{prop}
Recall the formal definition of Hausdorff dimension in Definition~\ref{d.hd}.
As indicated in Section \ref{s:iop}, we would split $[g,h]$ into segments, and bound the expected number of segments that intersect with $Z$.
The key estimate is as follows.
\begin{prop}   \label{prop:empty-prob-geo}
There is a constant $C$ depending on $g, h$, such that for any $g\le t< t+\epsilon \le h$, we have
\[
\PP[ Z\cap (t, t+\epsilon) \neq \emptyset ] < C\epsilon^{2/3} |\log(\epsilon)|^5.
\]
\end{prop}
There are two mechanisms for $Z\cap (t, t+\epsilon) \neq \emptyset$: either $\pi_\boo(t)$ is in a neighborhood of $0$ of size in the order of $\epsilon^{2/3}$, or $\pi_\boo$ has a transversal fluctuation in $(t, t+\epsilon)$ larger than its typical order (of $\epsilon^{2/3}$). We shall bound the probabilities of both events.
The proof is similar in spirit to the proof of Proposition \ref{prop:empty-prob}.
\begin{proof}
We can assume that $\epsilon$ is small enough since otherwise the inequality is obvious by taking the constant large enough.

We consider the following two events, corresponding to the above two mechanisms.
\begin{itemize}
    \item Let $\cE_1$ be the event where $|\pi_\boo(t)|\le \epsilon^{2/3}|\log(\epsilon)|^{5}$.
    \item Let $\cE_2$ be the event where there is some $s\in [t, t+\epsilon]$, such that
    \[|\pi_{(\boo;0,s)}(t)|>\epsilon^{2/3}|\log(\epsilon)|^{5}.\]
\end{itemize}
Then we have that the event $Z\cap (t, t+\epsilon) \neq \emptyset$ implies $\cE_1\cup\cE_2$.
Indeed, suppose that there is some $s\in Z\cap (t, t+\epsilon)$, then $\pi_{(\boo;0,s)}$ is a part of $\pi_\boo$; i.e., for any $s'\in [0, s]$ we have $\pi_{(\boo;0,s)}(s')=\pi_\boo(s')$.
In particular, we would have either $|\pi_\boo(t)|\le \epsilon^{2/3}|\log(\epsilon)|^{5}$, or $\pi_{(\boo;0,s)}(t)|>\epsilon^{2/3}|\log(\epsilon)|^{5}$, implying that either $\cE_1$ or $\cE_2$ occurs.

It now suffices to bound the probabilities of these events by the order of $\epsilon^{2/3} |\log(\epsilon)|^5$.
The bound on $\PP[\cE_1]$ is by Corollary \ref{cor:cross-int}.
For $\PP[\cE_2]$, we us Lemma \ref{lem:dl-trans-f}, as bounding $\PP[\cE_2]$ in the proof of Proposition \ref{prop:empty-prob}. We omit the details here. 
\end{proof}

\begin{proof}[Proof of Proposition \ref{prop:upper-geo}]
For any large enough $n\in \N$, let
\[
\sZ_n = \{i\in \Z: [i/n, (i+1)/n] \cap Z \cap [g, h] \neq \emptyset\}.
\]
By Proposition \ref{prop:empty-prob-geo}, $\E[|\sZ_n|] < C n^{2/3}(\log(n))^5$, where $C$ is a constant depending only on $g, h$.
For $0<\alpha<d-1/3$, we then have 
\[
\PP[n^{-d}|\sZ_n|>n^{-\alpha}] < C n^{1/3-d+\alpha}(\log(n))^5.
\]
This implies that almost surely, we can find a random $k_0$, such that for any $k>k_0$,
\[
2^{-dk}|\sZ_{2^k}| \le 2^{-\alpha k}.
\]
Since the right hand side $\to 0$ as $k\to\infty$, and $Z\cap [g, h]\subset \cup_{i\in\sZ_n} [i/n, (i+1)/n]$, the conclusion follows.
\end{proof}

We have now spelled out all the key ingredients and ideas for us to be able to plunge in to the proof of Theorem \ref{thm:1dhaus}.

\section{Hausdorff dimension of $\vNC$: Proof of Theorem \ref{thm:1dhaus}}  \label{sec:1D-diff-prof}

Our general strategy would be similar to that implemented in Section \ref{sec:2ddphau}. However while the upper bound straight forwardly follows from Proposition \ref{prop:empty-prob}, there are several key differences in the lower bounds, which we explain below. 

\subsection{Proof of the upper bound}
The proof is similar to that of the proof of Proposition \ref{prop:upper-geo} with the application of Proposition \ref{prop:empty-prob-geo} replaced by an application of  Proposition \ref{prop:empty-prob}. We provide the details below for completeness at the risk of being slightly repetitive. 
For any $0<g<h$, and any $d>2/3$, we show that almost surely, the $d$-dimensional Hausdorff measure of the set $\vNC\cap [g, h]$ is zero.
Take any large enough $n\in \N$, and let
\[
\sS_n = \{i\in \Z: [i/n, (i+1)/n] \cap S \cap [g, h] \neq \emptyset\}.
\]
By Proposition \ref{prop:empty-prob}, we have that $\E[|\sS_n|] < C n^{2/3}e^{C(\log(n))^{5/6}}$, where $C$ is a constant depending only on $g, h$.
Taking $0<\alpha<d-2/3$, we then have 
\[
\PP[n^{-d}|\sS_n|>n^{-\alpha}] < C n^{2/3-d+\alpha}e^{C(\log(n))^{5/6}}.
\]
This implies that almost surely, we can find a random $k_0$, such that for any $k>k_0$,
\[
2^{-dk}|\sS_{2^k}| \le 2^{-\alpha k}.
\]
The proof can then be finished as for Proposition \ref{prop:upper-geo}.
\qed\\\\

We next move on to the proof of the lower bound.
For this, we need to construct a measure $\KA$ on $\R_+$, which is the density of the measure $\zeta$ from Section \ref{sec:2ddphau}, as described in Section \ref{s:iop} (see \eqref{onedimmeasure} and the ensuing discussion). We shall call it \emph{the measure induced by the difference profile density}, and show that it is supported on $\vNC$.
We then prove a H\"{o}lder property of $\KA$, and that $\KA$ is strictly positive.
Recall the level functions $\LV_\ell$ defined as
\[\LV_\ell(t)=\inf\{x\in \R: \cD(x,t) \le \ell\}.\]
The H\"{o}lder property of $\KA$ is proved by writing $\KA$ as local times of $\LV_\ell$, averaged over $\ell \in \R$. Via results in Section \ref{ssec:spikebe}, this is reduced to analyzing the local time of the competition interface $\LV^\cB$, which has the same law as the semi-infinite geodesic $\pi_\boo$ (Proposition \ref{prop:equal-dist-geo-comp}).

There are two main additional difficulties in the setting of $\vNC$ (compared to that of $\NC$ in Section \ref{sec:2ddphau}).
First, the construction of $\KA$ is not straightforward, and we need arguments similar to those in Section \ref{ssec:geo-ocp-mea}.
In addition, there is some work needed to show that $\KA$ is an average of occupation measures of $\LV_\ell$.
We resolve these difficulties in Section \ref{ssec:constr-avg-ocp}.
Second, the fact that $\KA$ is non-degenerate is not obvious, and requires some work using estimates about the Airy line ensemble. We conduct these analyses in Section \ref{ssec:non-deg-eta}.
The desired H\"{o}lder regularity of $\KA$ is proved in Section \ref{ssec:H\"{o}lder-eta}.

Recall the notations $U_{g,h,N}$ and $\Good_{g,h,N,\delta}$ from Section \ref{ssec:spikebe}.  
For this section, we shall usually take $N=1/2$, and use the abbreviations $U_{g,h}=U_{g,h,1/2}$ and $\Good_{g,h,\delta}=\Good_{g,h,1/2,\delta}$, for simplicity of notations.

\subsection{Construction and average of occupation measures}  \label{ssec:constr-avg-ocp}

Recall $\zeta^\cB$, $\zeta_\ell$ for each $\ell \in \R$, and $\zeta$ on $\HH$, which are defined in \eqref{compinterfaceoccmeasure}, \eqref{eq:def-zeta-ell}, and \eqref{eq:def-zeta} respectively, and are given by
\[
\zeta^\cB([a,b]\times (g,h)) = \int_g^h \don[\LV^\cB(t)\in [a,b]] dt,
\]
and
\[
\zeta_\ell([a,b]\times (g,h)) = \int_g^h \don[\LV_\ell(t)\in [a,b]] dt,
\]
\[
\zeta([a,b]\times (g,h)) = \int_g^h [\cD(a,t)-\cD(b,t)] dt,
\]
for any $a < b$ and $0\le g<h$.
In this section, we will mostly consider these measures on rectangles. To ease notations, for each $h\ge 0$, we denote by $\lambda^\cB_h$ and $\lambda_{\ell,h}$ the measures on $\R$, such that
\[
\lambda^\cB_h(A)=\zeta^\cB(A\times (0,h])=\int_0^h \don[\LV^\cB(t) \in A] dt,
\]
and
\[
\lambda_{\ell,h}(A) = \zeta_\ell(A\times (0,h]) =\int_0^h \don[\LV_\ell(t) \in A] dt,
\]
for any measurable $A\subset \R$. Thus $\lambda^\cB_h$ and $\lambda_{\ell,h}$ can be viewed as the occupation measure for $\LV^\cB$ and $\LV_\ell$, respectively, for the time interval $[0, h]$. Also for each $0\le g < h$, we denote by $\lambda_{(g, h]}$ the measures on $\R$, such that
\begin{equation}  \label{eq:def-la-h}
\lambda_{(g, h]}(A)=\zeta(A\times (g,h])
\end{equation}
for any measurable $A\subset \R$.
\begin{rem}  \label{rem:la-gh}
For $\lambda_{(g, h]}$, we shall usually take $g>0$, because $\lambda_{(0, h]}$ is infinite. Indeed, by symmetry of the directed landscape, for any $a<b$ and $t>0$ we have $\E[\cD(a,t)-\cD(b,t)]=4(b-a)/t$, and the fluctuation of $\cD(a,t)-\cD(b,t)$ is $O(t^{1/3})$ (by Lemma \ref{lem:DLbound}).
Thus when $g=0$, we have $\lambda_{(0, h]}([a,b]) = \int_0^h [\cD(a,t)-\cD(b,t)] dt = \infty$. On the other hand, when $g>0$, we have $\E[\lambda_{(g, h]}([a,b])]<\infty$, so almost surely we have $\lambda_{(g, h]}([a,b]) < \infty$.
\end{rem}

The relation \eqref{eq:zeta-ell} for the set $A\times (g,h]$ can be written as
\begin{equation} \label{eq:mean-lambda-h}
   \lambda_{(g, h]}(A) = \int [\lambda_{\ell,h}(A) - \lambda_{\ell,g}(A)] d\ell .
\end{equation}
For the remainder of this section, we will use  $\lambda_{(g, h]}$, $\lambda^\cB_h$, $\lambda_{\ell,h}$, instead of $\zeta$, $\zeta^\cB$, $\zeta_\ell$, at the risk of introducing new notation. This is because we will mostly consider sets of the form of $A\times (g,h]$, and $\lambda^\cB_{h}$ has the same law as $\nu_h$ (defined in Section \ref{ssec:geo-ocp-mea}) by Proposition \ref{prop:equal-dist-geo-comp}.

We next consider the densities of these measures at zero, by defining, for any $0\le g<h$, 
\begin{equation} \label{eq:defn-kgh}
K_{g,h}:=\lim_{w\searrow 0} (2w)^{-1} \lambda_{(g, h]}([-w,w]),
\end{equation}
\[
K^\cB_{g,h}:=\lim_{w\searrow 0} (2w)^{-1} (\lambda^\cB_{h}([-w,w]) - \lambda^\cB_{g}([-w,w])),
\]
and
\[
K_{\ell, g,h}:=\lim_{w\searrow 0} (2w)^{-1} (\lambda_{\ell,h}([-w,w]) - \lambda_{\ell,g}([-w,w])).
\] 
We now show that these limits are (in some sense) well-defined and finite.

For the stationary initial setting, as $\lambda^\cB_{h}$ has the same law as $\nu_h$, almost surely $K^\cB_{g,h}$ is well defined and finite by Proposition \ref{prop:om-limit}, and $\lambda^\cB_h$ is absolutely continuous with respect to Lebesgue measure by Proposition \ref{prop:abs-cont-leb}.

We next study the limits in defining $K_{g,h}$ and $K_{\ell,g,h}$.
\begin{lemma}  \label{lem:ci-sum-limit}
For any $0<g<h$, almost surely $K_{g,h}$ is well defined, i.e., the limit exists and is finite.
\end{lemma}
\begin{proof}
This proof is similar to the proof of Proposition \ref{prop:om-limit}. We construct a translation invariant version of the difference profile.
Let $\xi \subset \R$ be a Poisson point process with rate $1$. For each $0\le g < h,$ we consider a measure $\lambda_{(g, h]}'$ such that for any $a<b$, we have
\begin{multline*}
\lambda_{(g, h]}'([a,b]) = \int_g^h \sum_{x\in\xi}(\cL(x-1,0;a\vee (x-1),t)-\cL(x+1,0;a\vee (x-1),t)) \\ - (\cL(x-1,0;b\wedge(x+1),t)-\cL(x+1,0;b\wedge(x+1),t)) dt.    
\end{multline*}
Basically this is just the sum of $\lambda_{(g, h]}$ centered at each $x\in \xi$ and truncated by $[x-1, x+1]$.
The truncation is taken so that at any point $y$, the measure only depends on passage times of finitely many starting points, and each is on the segment $[y-2,y+2]\times \{0\}$. Crucially however, even with this truncation the law of the measure $\lambda_{(g, h]}'$ is  translation invariant.

Via the translation invariance, we have that when $g>0$, $\lim_{w\searrow 0} (2w)^{-1}\lambda_{(g, h]}'([-w, w])$ exists and is finite almost surely, using the same arguments as in the proof of Proposition \ref{prop:om-limit} and the fact that $\lambda_{(g, h]}$ is finite in any compact interval (as mentioned in Remark \ref{rem:la-gh}).

Next, for small enough $\delta$, we consider the event $\cE_\delta$, where $|\xi \cap [-\delta, \delta]| = |\xi \cap [-3, 3]| = 1$. 
Then, given $0<g<h$ and recalling the event $\Good_{g,h,\delta}$ from Section \ref{ssec:spikebe}, we notice that under the latter, we have 
\begin{multline*}
(\cL(x_1-1,0;p_1)-\cL(x_1+1,0;p_1)) - (\cL(x_1-1,0;p_2)-\cL(x_1+1,0;p_2)) \\= (\cL(x_2-1,0;p_1)-\cL(x_2+1,0;p_1)) - (\cL(x_2-1,0;p_2)-\cL(x_2+1,0;p_2))    
\end{multline*}
for any $|x_1|, |x_2|<\delta$ and $p_1, p_2\in U_{g,h}=[-1/2, 1/2]\times [g/2, 2h]$.
In words, for any $|x_1|, |x_2|<\delta$, the difference profiles centered at $x_1$ and $x_2$ are the same inside the box $U_{g,h}$, up to adding a (random) constant.

From the above arguments, we have that $$(2w)^{-1}(\lambda_{(g, h]}([-w, w]) = (2w)^{-1}(\lambda_{(g, h]}'([-w, w])$$ for any $0<w<1/2$, assuming the event $\cE_\delta\cap \Good_{g,h,\delta}$.
By Lemma \ref{lem:low-bound-G} we have that $\PP[\Good_{g,h,\delta}]\to 1$ as $\delta\searrow 0$.
So by the independence of $\cE_\delta$ and $\Good_{g,h,\delta}$, we have $\PP[\cE_\delta\cap \Good_{g,h,\delta}] = \PP[\cE_\delta]\PP[\Good_{g,h,\delta}] > 0$ when $\delta$ is small enough.
These imply that conditional on the event $\cE_\delta\cap \Good_{g,h,\delta}$, we still have the (almost sure) existence and finiteness of the limit $\lim_{w\searrow 0} (2w)^{-1}(\lambda_{(g, h]}'([-w, w]))$, thus we have the (almost sure) existence and finiteness of the limit $\lim_{w\searrow 0} (2w)^{-1}(\lambda_{(g, h]}([-w, w]))$.
Since $\PP[\Good_{g,h,\delta}]\to 1$ as $\delta\searrow 0$, and $\cE_\delta$ is independent of $\lambda_{(g, h]}$, the existence and finiteness of the limit hold almost surely.
\end{proof}

\begin{lemma}  \label{lem:ci-ae-limit}
For any $0<g<h$, almost surely, for almost every $\ell \in \R$, $K_{\ell,g,h}$ is well defined (i.e., the limit exists and is finite), and $\lambda_{\ell,h}-\lambda_{\ell,g}$ is absolutely continuous with respect to Lebesgue measure.
By \eqref{eq:mean-lambda-h} we then conclude that (almost surely) $\lambda_{(g, h]}$ is also absolutely continuous with respect to Lebesgue measure.
\end{lemma}
\begin{proof}
Take any $\delta, H>0$.
By Proposition \ref{prop:equal-dist-geo-comp} and Proposition \ref{prop:om-limit}, there is a probability one event, under which $K_{g,h}^\cB$ is well-defined and finite; and by Proposition \ref{prop:abs-cont-leb}, there is a probability one event, under which $\lambda^\cB_{h}-\lambda^\cB_{g}$ is absolutely continuous with respect to Lebesgue measure.
On the other hand, by Lemma \ref{lem:low-bound-D}, under the event $\Good_{g,h,\delta}$ (which has positive probability by Lemma \ref{lem:low-bound-G} and taking $\delta$ small) almost surely we have $\PP[K_{g,h}^\cB = K_{\alpha,g,h}, \lambda^\cB_{h}-\lambda^\cB_{g}=\lambda_{\alpha,h}-\lambda_{\alpha,g} \mid \cL] > 0$, with $\alpha$ being random and has (conditional on $\cL$) probability density $\ge \delta'$ on $[-H, H]$.
These imply that under the event $\Good_{g,h,\delta}$, almost surely there is a full measure subset of $[-H, H]$, which is random and determined by $\cL$, such that for each $\ell$ in this set, $K_{\ell,g,h}$ is well-defined and finite, and $\lambda_{\ell,h}-\lambda_{\ell,g}$ is absolutely continuous with respect to Lebesgue measure.
Then the conclusion follows by taking $\delta\searrow 0$ and $H\to\infty$, and using Lemma \ref{lem:low-bound-G}.
\end{proof}
\begin{rem}  \label{rem:lc-ls-b}
We explain the reason why the above two lemmas require that $g>0$  (unlike in the case of semi-infinite geodesics).
For Lemma \ref{lem:ci-sum-limit}, it is simply because $\lambda_{(g, h]}$ is finite only when $g>0$, as discussed in Remark \ref{rem:la-gh}. 
For Lemma \ref{lem:ci-ae-limit}, one can actually show convergence to $K_{\ell, g, h}$ when $g=0$; but since that is not used in the remainder of this paper we omit the details.
In fact, the occupation measures for short times can be (roughly) described as follows. 
For any $(x, t)\in\HH$ and $\ell \in\R$ with $|x+\ell t/4|\gg t^{4/3}$, it is unlikely to have $\cD(x,t)=\ell$, by Lemma \ref{lem:DLbound}.
This implies that $|\LV_\ell(t)+\ell t/4|$ is at most in the order of $t^{4/3}$.
From these we expect that $\lambda_{\ell, 1}([-w,w])$ is in the order of $|\ell|^{-1}w$, if $\ell \neq 0$; and it is in the order of $w^{3/4}$ if $\ell=0$. On the other hand, for large times, $\LV_\ell$ grows linearly, as depicted in Figure \ref{f.simulations}. More precisely, it is proved in \cite[Proposition 4.5]{RV21} that $\LV_\ell(t)/t$ converges almost surely to a Gaussian random variable (while \cite[Proposition 4.5]{RV21} proves this for the $\ell=0$ case, the arguments for general $\ell$ are essentially verbatim). 
A discrete (or pre-limiting) analogue of this large $t$ behavior is the motion of a second class particle in a step initial TASEP or equivalently a single particle in the TASEP speed process (see \cite{MG05,AAV11}).
\end{rem}
By Lemma \ref{lem:ci-ae-limit}, we now have that almost surely, $K_{g,h}$ is well-defined and is positive and finite for each rational $0<g<h$.
We also have that almost surely for $\cL$, we can find a random full measure set $\Theta \subset \R$, which can be viewed as a function of $\cL$, such that for any $\ell \in \Theta$, and each rational $0<g<h$, $K_{\ell,g,h}$ is well-defined and finite.

Recall the measures $\KA$ and $\KA_\ell$ described in Section \ref{s:iop}. We are now ready to define them.
For any rational $0<h_1<h_2<h_3$, we have $K_{h_1,h_2}+K_{h_2,h_3}=K_{h_1,h_3}$.
Hence we can define the measure $\KA$ on $\R_+$, by letting
\begin{equation}  \label{eq:def-KA}
\KA((g,h)) = \sup_{g<g'<h'<h,\; g', h' \in \Q} K_{g',h'},
\end{equation}
for any $0<g<h$. We take the supremum from interior since we need $\KA((g,h))$ be right continuous in $g$ and left continuous in $h$.
Using standard measure theory arguments, this suffices to define $\KA$ as a measure on $\R_+$.

We now move on to the construction of the decomposition of $\KA$ indicated in \eqref{onedimmeasure}, which is the one dimensional analogue of \eqref{eq:mean-lambda-h}.
For each $\ell\in \Theta$ and $0<g<h$, we let
\begin{equation}\label{eq:defkl}
\KA_\ell((g,h)) = \sup_{g<g'<h'<h,\; g', h' \in \Q} K_{\ell,g',h'},
\end{equation}
and we extend this to define a measure $\KA_\ell$ on $\R_+$.
\begin{rem}
Actually, almost surely $K_{\ell,g,h}$ is continuous in $g$ and $h$ for $\ell$ in a random set with full measure.
This is by Proposition \ref{prop:mea-L-tail} and the equality in distribution between $\pi_\boo$ and $\LV^\cB$ (Lemma \ref{prop:equal-dist-geo-comp}), and taking $\KA^\cB$ as $\KA_\ell$ for a random $\ell$ using Lemma \ref{lem:low-bound-D} (as done in the proof of Lemma \ref{lem:ci-ae-limit}).
Thus $\KA_\ell((g,h)) = K_{\ell,g,h}$ for almost every $\ell$ and any rational $0<g<h$. 
\end{rem}
We next show that $\KA$ is the integral of $\KA_\ell$ over $\ell$. For the convenience of notations, we let $\KA_\ell$ be an arbitrary measure on $\R_+$ for each $\ell \in \R\setminus \Theta$. Since $\R\setminus \Theta$ has zero Lebesgue measure, the choice of such $\KA_\ell$ does not affect the integral.
\begin{lemma}  \label{lem:decomp-eta-to-ell}
Almost surely, $\KA(A)=\int \KA_\ell(A) d\ell$, for any measurable set $A$ that is contained in a compact subset of $\R_+$.
\end{lemma}
To prove this, we first prove the same relation for $K_{g,h}$ and $K_{\ell,g,h}$, for fixed $g, h$.
\begin{lemma}  \label{lem:connect}
For any $0<g<h$, almost surely we have that $K_{g,h}=\int K_{\ell, g, h} d\ell$.
\end{lemma}
We note that this is a ``density'' and hence a limiting version of \eqref{eq:mean-lambda-h}.
The proof of this lemma is via a uniform integrability result, when passing from $\lambda_{(g, h]}$ and $\lambda_{\ell, h}-\lambda_{\ell,g}$ to $K_{g,h}$ and $K_{\ell, g, h}$.
We would again rely on the connections with $\lambda^\cB_h-\lambda^\cB_g$, and with semi-infinite geodesics, and the estimate in Lemma \ref{lem:geo-ocp-tail}.
\begin{proof}
Fixing $0<g< h,$ for brevity, let 
\[
f_w(\ell) = (2w)^{-1} (\lambda_{\ell,h}([-w,w]) - \lambda_{\ell,g}([-w,w])),
\]
for $w>0$, and
\[
f(\ell) = K_{\ell,g,h}.
\]
From Lemma \ref{lem:ci-sum-limit} and Lemma \ref{lem:ci-ae-limit}, we have that almost surely, $f_w\to f$ almost everywhere, and also $\|f_w\|_{L_1}$ converges as $w\searrow 0$. Now we want to show (almost surely) the limit of $\lim_{w\searrow 0}\|f_w\|_{L_1}=\|f\|_{L_1}$.
For this we show for any $H>0$, there is some (random) positive sequence $w_1, w_2, \cdots \to 0$, such that $f_{w_n}\to f$ in $L_1([-H, H])$. 

Note that almost surely, all $f_w$ with $0<w<1$ and $f$ are supported on some compact interval. This is because, by Lemma \ref{lem:DLbound}, for the level function $\LV_\ell$, we have $|\LV_\ell(t)|>1$ for $t\in[g,h]$ when $|\ell|$ is large enough.
Thus it suffices to show that $\{f_{w_n}\}_{n=1}^\infty$ are uniformly integrable in $[-H, H]$, for any $H>0$.

Analogous to $f_w(\cdot)$, let 
\[
f^\cB_w=(2w)^{-1} (\lambda^\cB_{h}([-w,w]) - \lambda^\cB_{g}([-w,w])).
\]
We then have (by Lemma \ref{lem:geo-ocp-tail}), 
\[
\PP[f^\cB_w > M] < Ce^{-cM}
\]
for some $c,C>0$ that are independent of $w$ (but depend on $g, h$).
To show uniform integrablity, we let $\gamma_{w,m}=\int_{-H}^H \don[f_w(\ell)>m]f_w(\ell) d\ell$, for any $m\in \N$.
Now, for any $\delta>0$, let $\delta'$ be given by Lemma \ref{lem:low-bound-D}. Then,
\[
\E[\don[\Good_{g,h,\delta}]\liminf_{w\searrow 0} \gamma_{w,m}] \le \liminf_{w\searrow 0} \E[\don[\Good_{g,h,\delta}] \gamma_{w,m}] \le
1/\delta' \liminf_{w\searrow 0} \E[\don[f^\cB_w>m]f^\cB_w] < Ce^{-cm}/\delta',
\]
where the first inequality is by Fatou's lemma, and $C, c>0$ are constants independent of $w, m$.
Recall (from Lemma \ref{lem:low-bound-G}) that $\lim_{\delta\searrow 0}\PP[\Good_{g,h,\delta}]=1$, and note that $\Good_{g,h,\delta}$ does not depend on $w$ or $m$. These imply that $\lim_{m\to\infty}\liminf_{w\searrow 0} \gamma_{w,m} = 0$ almost surely.

Via a diagonal argument, almost surely we can take a (random) positive subsequence $w_1, w_2, \ldots \to 0$ such that $\lim_{n\to \infty} \gamma_{w_n,m}=\liminf_{w\searrow 0} \gamma_{w,m}$ for each $m\in \N$.
Hence
\[
\lim_{m\to \infty}\lim_{n\to \infty} \gamma_{w_n,m} = \lim_{m\to \infty}\liminf_{w\searrow 0} \gamma_{w,m} = 0.
\]
That is, $\{f_{w_n}\}_{n=1}^\infty$ are uniformly integrable in $[-H, H]$.
For $H$ large enough, we have all $f_w$ and $f$ are supported on $[-H, H]$ (and the largeness is random).
Thus for the subsequence $w_1, \ldots,$ we have $f_{w_n}\to f$ in $L_1$; this implies that $\|f_{w_n}\|_{L_1}\to\|f\|_{L_1}$, and hence $\|f_{w}\|_{L_1}\to\|f\|_{L_1}$.
\end{proof}

Using the above, the proof of Lemma \ref{lem:decomp-eta-to-ell} is straightforward. 

\begin{proof}[Proof of Lemma \ref{lem:decomp-eta-to-ell}]
For the arguments below, we assume the probability one event that Lemma \ref{lem:connect} holds for any rational $0<g<h$.

We first consider the case where $A$ is an interval $(g, h)$ for some rational $0<g<h$.
Take a rational sequence $g_n\searrow g$ and a rational sequence $h_n \nearrow h$.
Then $\KA(A)=\lim_{n\to\infty} K_{g_n,h_n}$ and $\KA_\ell(A)=\lim_{n\to\infty} K_{\ell, g_n,h_n}$ for any $\ell \in \Theta$.
Since $K_{\ell, g_n,h_n} \le K_{\ell, g,h}$ and $K_{\ell, g,h}$ is integrable as a function of $\ell$ (by Lemma \ref{lem:connect}), we can apply dominated convergence theorem to conclude that \[\int \KA_\ell(A)d\ell = \lim_{n\to\infty}  \int K_{\ell, g_n,h_n}d\ell.\]
By Lemma \ref{lem:connect}, the right hand side equals $\lim_{n\to\infty} K_{g_n,h_n} = \KA(A)$, so the conclusion holds in the case where $A$ is a rational interval.

Then we can extend this to any measurable set $A$ contained in a compact set in $\R_+$, using $\pi-\lambda$ theorem. For this, we need to show that for any rational $0< g < h$, the collection $\{A\subset (g,h): \KA(A)=\int\KA_\ell(A)d\ell \}$ is a Dynkin system: closed under complements is obvious, and closed under countable disjoint union is by dominated convergence theorem (since $\KA_\ell(A)\le \KA_\ell((g, h))$ for any $A\subset (g,h)$ and any $\ell\in \Theta$).
\end{proof}

As was indicated in Section \ref{s:iop}, the following statement highlights the usefulness of $\KA.$
\begin{lemma}  \label{lem:eta-sup-S}
Almost surely, the measure $\KA$ is supported inside $\vNC$.
\end{lemma}
\begin{proof}
In this proof we assume the event where we can find a measurable set $\Theta'\subset \Theta \subset \R$ with $\LE(\R\setminus \Theta')=0$, such that each $\ell \in \Theta'$ satisfies the following: (1) for each $t\in\R_+$, there is a unique $x\in\R$ with $\cD(x,t)=\ell$; (2) $K_{\ell,g,h}$ is well-defined and finite for each rational $0<g<h$.
By Lemma \ref{lem:uni-kappa} and Lemma \ref{lem:ci-ae-limit}, this event happens almost surely.

We now fix some $0<g<h$, and show that $\KA$ restricted on $(g,h)$ is supported inside $\vNC\cap (g, h)$.

Denote $S$ to be the closure of $\cup_{\ell \in \Theta'}\{t: g<t<h, \LV_\ell(t)=0\}$. Here the closure is taken to ensure the measurability of $S$.
We claim that $S \subset \vNC$.
Indeed, for any $t\not\in \vNC$, we have $\cD(0,s)=\ell_*$ for each $s\in A$, where $A$ is a neighborhood of $t$ and $\ell_*\in \R$ is a constant.
If $t \in S$, by the definition of $S$ we can find some $s \in A$ and some $\ell' \in \Theta'$, such that $\LV_{\ell'}(s)=0$. Thus we must have that $\ell'=\ell_*$, so $\ell_* \in \Theta'$.
This implies that $\LV_{\ell_*}(s)=0$ for each $s\in A$ and $K_{\ell_*,g,h}$ is well-defined and finite.
However, from its definition we have $K_{\ell_*,g,h}=\infty$ when $\LV_{\ell_*}(s)=0$ for each $s\in A$.
Thus we must have $t\not\in S$.

We shall next show that $\KA$ restricted to $(g,h)$ is supported inside $S$.
Take any $\ell \in \Theta'$. By Lemma \ref{lem:conti-kappa}, the set $\{t: g\le t\le h, \LV_\ell(t)=0\}$ is closed, and for any $g<t<h$ with $\LV_\ell(t)\neq 0$, we can find an open neighborhood around $(0,t)$ that is disjoint from the trajectory of $\LV_\ell$.
The later implies that for any $g<t<h$ with $\LV_\ell(t)\neq 0$, one can find some $\delta>0$ such that $\KA_\ell((t-\delta, t+\delta))=0$.
Thus $\KA_\ell(\{t:g<t<h, \LV_\ell(t)\neq 0\}) = 0$. Since $(g,h)\setminus S \subset \{t:g<t<h, \LV_\ell(t)\neq 0\}$ we have $\KA_\ell((g,h)\setminus S)=0$.
By Lemma \ref{lem:decomp-eta-to-ell}, we have that $\KA((g,h)\setminus S) =\int \KA_\ell((g,h)\setminus S) d\ell = 0$.
This implies that $\KA$ restricted to $(g,h)$ is supported inside $S$, thus supported inside $\vNC$.
\end{proof}

To show that the Hausdorff dimension of $\vNC$ is at least $2/3$,
it remains to show that $\KA$ is $2/3-\epsilon$ H\"{o}lder in any interval $[g,h]$ for any $\epsilon>0$, and $\KA([g,h])>0$ for small enough $g>0$.
These two steps are carried out in the next two subsections.

\subsection{$2/3$ H\"{o}lder of the difference profile density} \label{ssec:H\"{o}lder-eta}

\begin{lemma}  \label{lem:vertical-23-H\"{o}lder}
For any $0<g<h$ and $\epsilon>0$, almost surely there is some random $Q$, such that $\KA([a,b])<Q(b-a)^{2/3-\epsilon}$ for any $g\le a<b\le h$.
\end{lemma}
The proof follows a similar strategy as the proof of Lemma \ref{lem:2d-53-H\"{o}lder}.
Specifically, it suffices to prove $K_{a,b}<Q(b-a)^{2/3-\epsilon}$ for any $a,b \in \Q$ satisfying $g\le a<b\le h$ and $b-a<1$.
By the regularity condition of the directed landscape under the event $\Good_{g,h,\delta}$, $K_{a,b}^\ell$ is nonzero only for $\ell$ in an interval of length (in the order of) $(b-a)^{1/3}$.
To bound $K_{a,b}$, we just bound an integral of $K_{a,b}^\ell$ over $\ell$ in an interval of length (in the order of) $(b-a)^{1/3}$, and take the supremum over all such intervals.
Using Lemma \ref{lem:low-bound-D}, (under $\Good_{g,h,\delta}$) such a supremum of integrals is bound by $(b-a)^{1/3}$ times $K_{a,b}^\cB$; and by Proposition \ref{prop:mea-L-tail}, $K_{a,b}^\cB$ is in the order of $(b-a)^{1/3}$ with exponential tails.
\begin{proof}
Recall the almost sure event where we have a set $\Theta\subset \R$ with full measure, such that for every $\ell\in \Theta$, and any $0<a <b$, $a, b \in \Q$, the numbers $K_{a,b}$, $K_{\ell, a,b}$, and $K_{a,b}^\cB$ are well defined and finite. We assume this event for the rest of this proof.

Fix any small enough $\delta>0$. It suffices to bound $K_{a,b}$ for any rational $a, b$ satisfying $g\le a<b\le h$ and $b-a<1$, under the event $\Good_{g,h,\delta}$.
For each $\ell$, let
\[
Q_\ell = \sup_{g\le a < b\le h, a,b \in \Q} \frac{K_{\ell,a,b}}{(b-a)^{1/3}(1+|\log(b-a)|^2)}.
\]
Here the extra factor of $1+|\log(b-a)|^2$ is due to the fact that we take the supremum over $a, b$, and the power of $2$ can be replaced by any number $>1$.
To bound $K_{a,b}$, we would bound the integral of $Q_\ell$ over $\ell$.
We take $H>0$, large enough and depending on $g, h, \delta$.
For any $\varepsilon>0$ we denote 
\[
J_\varepsilon = \sup_{A\subset [-H,H], \LE(A)<\delta^{-1} \varepsilon^{1/3}\log^{2/3}(1+\varepsilon^{-1}) } \int_A Q_\ell d\ell ,
\]
and let
\[
J= \sup_{\varepsilon\in (0,1)} \frac{J_\varepsilon}{\varepsilon^{1/3}\log^3(1+\varepsilon^{-1})} .
\]
The motivation behind these definitions is that, by regularity estimates of $\cL$, $K_{\ell,a,b}$ is non-zero only for $\ell$ in an interval of length of the order of $(b-a)^{1/3}$.

As in the proof of Lemma \ref{lem:2d-53-H\"{o}lder}, this random variable $J$ is essentially the desired $Q$ in the statement of the lemma, and the remainder of this proof mainly consists of two steps: (1) we bound the upper tail of $\don[\Good_{g,h,\delta}]J$; (2) we bound $K_{a,b}$ for any $g\le a< b\le h$, $a,b\in \Q$ using $J$ under the event $\Good_{g,h,\delta}$. Unavoidably, the following sequence of arguments, like in the proof of Lemma \ref{lem:2d-53-H\"{o}lder}, will be somewhat technical.\\

\noindent
\textbf{Step 1.}
We first bound $\don[\Good_{g,h,\delta}]J$.
As in the proof of Lemma \ref{lem:2d-53-H\"{o}lder}, we take inputs from estimates on the competition interface with stationary initial condition.
Define
\[
Q^\cB = \sup_{g\le a < b\le h, a,b \in \Q} \frac{K^\cB_{a,b}}{(b-a)^{1/3}(1+|\log(b-a)|^2)}.
\]
By Lemma \ref{lem:low-bound-D} and taking $\delta'$ there, for any $M>0$ we have
\begin{equation} \label{eq:good-J-pre}
\E\left[\don[\Good_{g,h,\delta}]\int_{-H}^H \don[Q_\ell > M] d\ell \right] < {\delta'}^{-1}\PP[Q^\cB > M]  < {\delta'}^{-1}C_1e^{-c_1M},    
\end{equation}
where $c_1,C_1>0$ are some constants (depending on $g, h$) and the last inequality is by Proposition \ref{prop:mea-L-tail}.

Recall the arguments in Step 1 in the proof of Lemma \ref{lem:2d-53-H\"{o}lder}, where \eqref{eq:tail-Q-ell} is translated into an upper tail bound of $\don[\Good_{g,h,\delta}]J$ there. 
Using the same algebra, \eqref{eq:good-J-pre} can be translated into the following bound on the upper tail of $\don[\Good_{g,h,\delta}]J$:
for any $M>0$, we have that $$\PP[\don[\Good_{g,h,\delta}]J > M] < C_\delta e^{-c_\delta M},$$ for $C_\delta, c_\delta$ being constants depending on $g, h, \delta$.\\

\noindent
\textbf{Step 2.} 
Under $\Good_{g,h,\delta}$, using its second condition and Lemma \ref{lem:DLbound}, we can bound $|\cD|$ in $[-1/2,1/2]\times[g,h]$.
Consequently, we have that $K_{\ell,g,h} = 0$ for any $|\ell|>H$, as $H$ is large enough depending on $g, h, \delta$.
Then by Lemma \ref{lem:connect}, for any rational $a, b$ with $g\le a<b\le h$, we have $K_{a,b} = \int_{-H}^H K_{\ell,a,b}d\ell$.
We note that by the second condition of $\Good_{g,h,\delta}$ and using Lemma \ref{lem:modcont}, for any $\ell_1, \ell_2$, if both  $K_{\ell_1,a,b}$ and $K_{\ell_2,a,b}$ are nonzero, we must have
\[
|\ell_1-\ell_2|\le C\delta^{-0.01} (b-a)^{1/3}\log^{2/3}(1+(b-a)^{-1}) < \delta^{-1} (b-a)^{1/3}\log^{2/3}(1+(b-a)^{-1}),
\]
where $C>0$ is a constant depending on $g, h$, and the last inequality is by taking $\delta$ small enough depending on $g, h$.
Thus by the definitions of $J_{b-a}$ and $J$, we have (for $b-a<1$)
\[
K_{a,b} = \int_{-H}^H K_{\ell,a,b}d\ell \le (b-a)^{1/3}(1+|\log(b-a)|^2) J_{b-a} < (b-a)^{2/3}(1+|\log(b-a)|^2) \log^3(1+(b-a)^{-1}) J.
\]
By taking $Q$ to be $J$ multiplied by a constant depending on $g, h, \epsilon$, we get the desired bound $K_{a,b}<Q(b-a)^{2/3-\epsilon}$ under $\Good_{g,h,\delta}$.
Finally, by taking $\delta\searrow 0$, we get that the desired bound holds almost surely.
\end{proof}

\subsection{Non-degeneracy of the difference profile density}  \label{ssec:non-deg-eta}
In this subsection we implement the last step in our strategy to lower bound the Hausdorff dimension of $\vNC$, namely showing that $\KA$ is strictly positive This along with Lemma \ref{lem:eta-sup-S} and Lemma \ref{lem:vertical-23-H\"{o}lder} concludes the proof of the lower bound of Theorem \ref{thm:1dhaus}, as indicated by the general strategy in Section \ref{s:iop}.
\begin{prop}  \label{prop:eta-non-deg}
Almost surely, there exist $0<g<h$ such that $\KA([g,h])>0$.
\end{prop}
Note that for fixed $0<g<h,$ with positive probability,
owing to coalescence of geodesics, the difference profile $\cD$ can be flat in the box $[-1,1]\times [g,h]$, thus $\KA([g,h])=0$, with positive probability.

The proof of Proposition \ref{prop:eta-non-deg} follows a similar strategy as the proof of Proposition \ref{prop:L-not-dege}: we start by showing that for a fixed interval, $\KA$ is positive with positive probability; then we use ergodicity to show that $\KA$ is almost surely strictly positive.
Recall the notation of $K_{g,h}$ for any $0<g<h$, from \eqref{eq:defn-kgh}.
\begin{lemma}  \label{lem:K-pos-posprob}
We have $\PP[K_{1,2}>0]>0$.
\end{lemma}
To prove this we need estimates using the Airy line ensemble. We defer the Airy line ensemble arguments to Appendix \ref{sec:appa}, and state the rest of the arguments here.

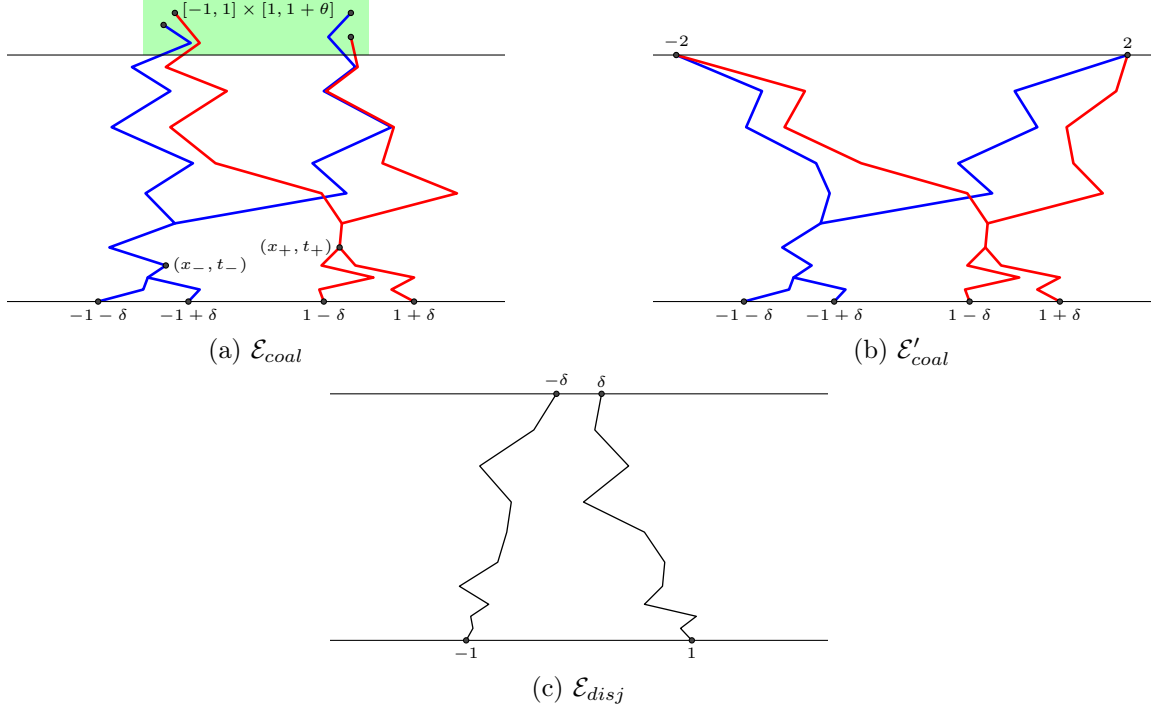
\begin{figure}[hbt!]
    \centering
    
\begin{subfigure}{0.48\textwidth}
    \centering
\begin{tikzpicture}[line cap=round,line join=round,>=triangle 45,x=0.3cm,y=.08cm]
\clip(-6,22) rectangle (16,75);

\fill[line width=0.pt,color=green,fill=green,fill opacity=0.3]
(0,66) -- (10,66) -- (10,75) -- (0,75) -- cycle;

\draw [line width=1pt,color=blue] (2,25) -- (2.5,27) -- (0.2,29);
\draw [line width=1pt,color=blue] (-2,25) -- (0,27) -- (0.2,29) -- (1,31) -- (-1.5,34) -- (1.4,38) -- (0.1,43) -- (2.2,48) -- (-1.4,54) -- (1.2,60) -- (-0.5,64) -- (2.1,68) -- (0.9,71);

\draw [line width=1pt,color=blue] (1.4,38) -- (9,43) -- (7.5,48) -- (11,54) -- (8,60) -- (9.4,64) -- (8.2,69) -- (9.2,73);

\draw [line width=1pt,color=red] (8,25) -- (7.8,27) -- (10.2,29) -- (7.9,31) -- (8.7,34);

\draw [line width=1pt,color=red] (12,25) -- (11,27) -- (12,29) -- (9.4,31) -- (8.7,34) -- (8.8,38) -- (13.9,43) -- (10.6,48) -- (11.1,54) -- (8.1,60) -- (9.5,64) -- (9.2,69);

\draw [line width=1pt,color=red] (8.8,38) -- (7.9,43) -- (3.2,48) -- (1.2,54) -- (3.7,60) -- (1,64) -- (2.5,68) -- (1.4,73);

\draw (-10,25) -- (20,25);
\draw (-10,66) -- (20,66);

\draw [fill=uuuuuu] (0.9,71) circle (1pt);
\draw [fill=uuuuuu] (9.2,73) circle (1pt);
\draw [fill=uuuuuu] (1.4,73) circle (1pt);
\draw [fill=uuuuuu] (9.2,69) circle (1pt);

\draw [fill=uuuuuu] (1,31) circle (1pt);
\draw [fill=uuuuuu] (8.7,34) circle (1pt);

\draw [fill=uuuuuu] (-2,25) circle (1pt);
\draw [fill=uuuuuu] (2,25) circle (1pt);
\draw [fill=uuuuuu] (8,25) circle (1pt);
\draw [fill=uuuuuu] (12,25) circle (1pt);

\begin{tiny}
\draw (-2,25) node[anchor=north]{$-1-\delta$};
\draw (2,25) node[anchor=north]{$-1+\delta$};
\draw (8,25) node[anchor=north]{$1-\delta$};
\draw (12,25) node[anchor=north]{$1+\delta$};

\draw (1,31) node[anchor=west]{$(x_-,t_-)$};
\draw (8.7,34) node[anchor=east]{$(x_+,t_+)$};

\draw (5.1,76) node[anchor=north]{$[-1,1]\times [1,1+\theta]$};

\end{tiny}
\end{tikzpicture}
\subcaption{$\cE_{coal}$}
\end{subfigure}
\hfill
\begin{subfigure}{0.48\textwidth}
    \centering
\begin{tikzpicture}[line cap=round,line join=round,>=triangle 45,x=0.3cm,y=.08cm]
\clip(-6,22) rectangle (16,75);

\draw [line width=1pt,color=blue] (2,25) -- (2.5,27) -- (0.2,29);
\draw [line width=1pt,color=blue] (-2,25) -- (0,27) -- (0.2,29) -- (1,31) -- (-0.3,34) -- (1.4,38) -- (1.8,43) -- (1.2,48) -- (-1.9,54) -- (-1.2,60) -- (-5,66);

\draw [line width=1pt,color=blue] (1.4,38) -- (9,43) -- (7.5,48) -- (11,54) -- (10,60) -- (15,66);

\draw [line width=1pt,color=red] (8,25) -- (7.8,27) -- (10.2,29) -- (7.9,31) -- (8.7,34);

\draw [line width=1pt,color=red] (12,25) -- (11,27) -- (12,29) -- (9.4,31) -- (8.7,34) -- (8.8,38) -- (13.9,43) -- (12.6,48) -- (12.3,54) -- (14.5,60) -- (15,66);

\draw [line width=1pt,color=red] (8.8,38) -- (7.9,43) -- (3.2,48) -- (-0.2,54) -- (0.7,60) -- (-5,66);

\draw (-10,25) -- (20,25);
\draw (-10,66) -- (20,66);

\draw [fill=uuuuuu] (-5,66) circle (1pt);
\draw [fill=uuuuuu] (15,66) circle (1pt);

\draw [fill=uuuuuu] (-2,25) circle (1pt);
\draw [fill=uuuuuu] (2,25) circle (1pt);
\draw [fill=uuuuuu] (8,25) circle (1pt);
\draw [fill=uuuuuu] (12,25) circle (1pt);

\begin{tiny}
\draw (-2,25) node[anchor=north]{$-1-\delta$};
\draw (2,25) node[anchor=north]{$-1+\delta$};
\draw (8,25) node[anchor=north]{$1-\delta$};
\draw (12,25) node[anchor=north]{$1+\delta$};

\draw (-5,66) node[anchor=south]{$-2$};
\draw (15,66) node[anchor=south]{$2$};
\end{tiny}

\end{tikzpicture}
\subcaption{$\cE_{coal}'$}
\end{subfigure}
\hfill
\begin{subfigure}{0.48\textwidth}
    \centering
\begin{tikzpicture}[line cap=round,line join=round,>=triangle 45,x=0.3cm,y=.08cm]
\clip(-6,22) rectangle (16,70);

\draw [line width=0.5pt] (10,25) -- (9.5,27) -- (10.2,29) -- (7.9,31) -- (8.7,34) -- (8.8,38) -- (7.9,43) -- (5.2,48) -- (7.2,54) -- (5.7,60) -- (6,66);

\draw [line width=0.5pt] (0,25) -- (0.3,27) -- (0.2,29) -- (1,31) -- (-0.3,34) -- (1.4,38) -- (1.8,43) -- (2,48) -- (0.6,54) -- (3,60) -- (4,66);

\draw (-10,25) -- (20,25);
\draw (-10,66) -- (20,66);

\draw [fill=uuuuuu] (4,66) circle (1pt);
\draw [fill=uuuuuu] (6,66) circle (1pt);

\draw [fill=uuuuuu] (0,25) circle (1pt);
\draw [fill=uuuuuu] (10,25) circle (1pt);

\begin{tiny}
\draw (0,25) node[anchor=north]{$-1$};
\draw (10,25) node[anchor=north]{$1$};

\draw (4,66) node[anchor=south]{$-\delta$};
\draw (6,66) node[anchor=south]{$\delta$};

\end{tiny}

\end{tikzpicture}
\subcaption{$\cE_{disj}$}
\end{subfigure}

\caption{Illustrations of the events $\cE_{coal}$, $\cE_{coal}'$, and $\cE_{disj}$ in the proof of Lemma \ref{lem:K-pos-posprob}:
assuming $\cE_{disj}$, we have $\int_{-\delta}^\delta \oK^s_{1,1+\theta} ds>0$ for some random $\theta>0$, where $\oK^s_{1,1+\theta}$ is $K_{1,1+\theta}$ translated by $(s,0)$ in its definition.
Then by translation invariance (and assuming $\PP[\cE_{coal}\cap\cE_{disj}]>0$), we have $\E[K_{1,1+\theta}]>0$ for small enough $\theta$. To show $\PP[\cE_{coal}\cap\cE_{disj}]>0$ we reduce $\cE_{coal}$ to $\cE_{coal}'$, and use tools from the Airy line ensemble to show $\PP[\cE_{coal}'\cap\cE_{disj}]>0$ (in Appendix \ref{sec:appa}).}
\label{fig:3}
\end{figure}

\begin{proof}
Take some small $\theta>0$ to be determined.
It suffices to prove that $\E[K_{1,1+\theta}]>0$.
The general idea is to use translation invariance and coalescence of geodesics, to reduce this to proving that, for some small $\delta>0$, $\lambda_{(1, 1+\theta]}([-\delta, \delta])>0$ with positive probability (recall \eqref{eq:def-la-h} for the definition of $\lambda_{(1, 1+\theta]}$); and this event is further implied by some coalescence and disjointness of geodesics. Such coalescence and disjointness of geodesics are studied using the Airy line ensemble, and we leave the arguments to Appendix \ref{sec:appa}.

We first define translated versions of $K_{1,1+\theta}$. For any $s\in\R$ we let
\begin{multline*}
\oK^s_{1,1+\theta} = \lim_{w\searrow 0}(2w)^{-1} \int_1^{1+\theta} [(\cL(-1+s,0;-w+s,t)-\cL(1+s,0;-w+s,t)) \\
-(\cL(-1+s,0;w+s,t)-\cL(1+s,0;w+s,t))] dt.
\end{multline*}
By Lemma \ref{lem:ci-sum-limit} and translation invariance of the directed landscape, for fixed $s, \theta$, almost surely $\oK^s_{1,1+\theta}$ is well-defined and finite, and has the same distribution as $K_{1,1+\theta}$.
We also define
\begin{align*}
\hK_{1,1+\theta}^s=&\lim_{w\searrow 0} (2w)^{-1} (\lambda_{(1, 1+\theta]}([s-w,s+w]))
\\
=&\lim_{w\searrow 0}(2w)^{-1} \int_1^{1+\theta} \big[(\cL(-1,0;-w+s,t)-\cL(1,0;-w+s,t)) \\
&\quad\quad\quad\quad\quad\quad-(\cL(-1,0;w+s,t)-\cL(1,0;w+s,t))\big] dt.
\end{align*}
We note that (by Lebesgue's theorem for the differentiability of monotone functions as before), almost surely this is well-defined for almost every $s$.
The difference between these two definitions is that, for $\hK_{1,1+\theta}^s$ only the upper endpoints of the considered passage times are shifted by $(s,0)$, while everything is shifted by $(s,0)$ for $\oK_{1,1+\theta}^s$.

We now explain the reason behind defining $\oK_{1,1+\theta}^s$ and $\hK_{1,1+\theta}^s$. We will write $\E[K_{1,1+\theta}]$ as an integral of $\E[\oK^s_{1,1+\theta}]$ over $s$.
We will show that $\oK^s_{1,1+\theta}=\hK^s_{1,1+\theta}$ under certain coalescence events, so it suffices to consider the expectation of some integral of $\hK^s_{1,1+\theta}$ over $s$, which can be written as $\lambda_{(1, 1+\theta]}([-\delta, \delta])$.

We first define the coalescence event under which $\oK^s_{1,1+\theta}=\hK^s_{1,1+\theta}$ holds.
Let $\delta>0$ be a small number. Let $\cE_{coal}$ be the following event (see Figure \ref{fig:3} for an illustration): there exist some $t_-, t_+>0$, $x_-, x_+ \in \R$, such that for any $x\in (-1-\delta, -1+\delta)$ and $p \in [-1,1]\times [1, 1+\theta]$, there is $\pi_{(x,0;p)}(t_-)=x_-$; and for any $x\in (1-\delta, 1+\delta)$ and $p \in [-1,1]\times [1, 1+\theta]$, there is $\pi_{(x,0;p)}(t_+)=x_+$. 
Under $\cE_{coal}$, for any $p_1, p_2 \in [-1,1]\times [1, 1+\theta]$ and $|s|<\delta$, we have
\[
\cL(-1,0;p_1)-\cL(-1,0;p_2)=\cL(-1+s,0;p_1)-\cL(-1+s,0;p_2).
\]
This is because, from the coalescence of geodesics as illustrated by Figure \ref{fig:3}, we have that both sides would equal $\cL(x_-,t_-;p_1)-\cL(x_-,t_-;p_2)$.
Similarly,
\[
\cL(1,0;p_1)-\cL(1,0;p_2)=\cL(1+s,0;p_1)-\cL(1+s,0;p_2).
\]
By taking the difference of the above two equations, letting $p_1=(-w+s,t)$, $p_2=(w+s,t)$, integrating over $t$ from $1$ to $1+\theta$, dividing both hand sides by $2w$, and sending $w\searrow 0$, we conclude that $\oK^s_{1,1+\theta}=\hK^s_{1,1+\theta}$ under $\cE_{coal}$.
We then have
\begin{equation}  \label{eq:two-ev-key}
2\delta \E[K_{1,1+\theta}] = \int_{-\delta}^\delta \E[\oK^s_{1,1+\theta}] ds
\ge \int_{-\delta}^\delta \E[\don[\cE_{coal}]\hK^s_{1,1+\theta}] ds
=
\E[\don[\cE_{coal}]\int_{-\delta}^\delta \hK^s_{1,1+\theta} ds].   
\end{equation}
It remains to show that (for appropriate choice of $\delta, \theta$) the last term is positive.
Note that as $\lambda_{(1, 1+\theta]}$ is absolutely continuous with respect to Lebesgue measure by Lemma \ref{lem:ci-ae-limit} (almost surely for fixed $\theta$), we can write
\begin{align*}
\int_{-\delta}^\delta \hK^s_{1,1+\theta} ds
&=
\lambda_{(1, 1+\theta]}([-\delta, \delta])
\\
&=
\int_1^{1+\theta} [(\cL(-1,0;-\delta,t)-\cL(1,0;-\delta,t))
-(\cL(-1,0;\delta,t)-\cL(1,0;\delta,t))]  dt.  
\end{align*}
We now find an event where the right hand side is positive.
Recall that for any $x_1< x_2$, $y_1<y_2$, and $s<t$, $\Dist_{(x_1, x_2)\to (y_1, y_2)}^{[s, t]}$ is the event where $\pi_{(x_1, s;y_1, t)}$ and $\pi_{(x_2, s;y_2, t)}$ are unique and disjoint.
We then let $\cE_{disj}=\Dist_{(-1, 1)\to (-\delta, \delta)}^{[0, 1]}$.
Then by Lemma \ref{lem:dis-equiv} and spatial monotonicity of $\cD$, and the fact that almost surely $\cL$ is continuous, we have that $\cE_{disj}$ implies $\int_{-\delta}^\delta \hK^s_{1,1+\theta} ds>0$.

It remains to show that $\PP[\cE_{coal}\cap \cE_{disj}]>0$. 
The event $\cE_{coal}$ involves geodesics whose upper endpoints vary in the box $[-1,1]\times [1,1+\theta]$, and is not easy to analyze using the Airy line ensemble.
For this, we consider \[\cE_{coal}'=(\Dist_{(-1-\delta, -1+\delta)\to (-2, 2)}^{[0, 1]} \cup \Dist_{(1-\delta, 1+\delta)\to (-2, 2)}^{[0, 1]})^c\] (see Figure \ref{fig:3} for an illustration).
The event $\cE_{coal}'$ only involves geodesics from time $0$ to time $1$, and can be written as an event of the Airy line ensemble. From this we can show that $\PP[\cE_{coal}'\cap \cE_{disj}]>0$. We leave the details to Lemma \ref{lem:sur-ALE}.

By Lemma \ref{lem:dl-trans-f}, given $\cE_{coal}'$ almost surely we have that $\cE_{coal}$ holds for some random $\theta>0$. 
Indeed, by Lemma \ref{lem:dl-trans-f} we can find some $\theta>0$, such that for any $x\in (-1-\delta, -1+\delta) \cup (1-\delta, 1+\delta)$ and $p \in [-1,1]\times [1, 1+\theta]$, we have $|\pi_{(x,0;p)}(1)|\le 2$.
Then by ordering of geodesics, for any $t_-, t_+$ and $x_-, x_+ \in \R$ such that $\pi_{(-1-\delta,0; -2,1)}(t_-)=\pi_{(-1+\delta,0; 2,1)}(t_-)=x_-$ and $\pi_{(1-\delta,0; -2,1)}(t_+)=\pi_{(1+\delta,0; 2,1)}(t_+)=x_+$, they satisfy the properties in the definition of $\cE_{coal}$ above, so $\cE_{coal}$ holds.
Thus $\PP[\cE_{coal}'\cap \cE_{disj}]>0$ implies that for some $\theta$ small enough (depending on $\delta$), we have $\PP[\cE_{coal}\cap \cE_{disj}]>0$. 
We then have that $\E[\don[\cE_{coal}]\int_{-\delta}^\delta \hK^s_{1,1+\theta} ds]>0$, which with \eqref{eq:two-ev-key} implies the conclusion.
\end{proof}

We need one more ingredient to finish the proof of Proposition \ref{prop:eta-non-deg}, namely a 0-1 law.
For any $r > 0,$ we let $\cF_r$ be the sigma algebra generated by $\cL(x,s;y,t)$, for $0\le s<t\le r$.
Let $\cF_{0+} = \cap_{r>0} \cF_r$ be the corresponding germ sigma-algebra. 
\begin{lemma}  \label{lem:F-tri}
Any $\cF_{0+}$-measurable event has probability $0$ or $1$.
\end{lemma}

\begin{proof}
We show that $\cL(x,s;y,t)$ is independent of $\cF_{0+}$ for any $(x,s;y,t)\in\R^4_\uparrow$.
If $s>0$ or $t\le 0$, this is obvious.
If $s=0$, we have $\cL(x,s;y,t) = \lim_{\delta\searrow 0} \cL(x,s+\delta;y,t)$, which also implies the statement, since each $\cL(x,s+\delta;y,t)$ is independent of $\cF_{0+}$. 
Now for $s<0<t$, we have the composition law
\[
\cL(x,s;y,t) = \sup_z \cL(x,s;z,0) + \cL(z,0;y,t).
\]
Since $\cL(x,s;z,0)$ is independent of $\cF_{0+}$, and that we have just shown that $\cL(z,0;y,t)$ is independent of $\cF_{0+}$, we must have that the left hand side is also independent of $\cF_{0+}$.
Thus the conclusion follows.
\end{proof}

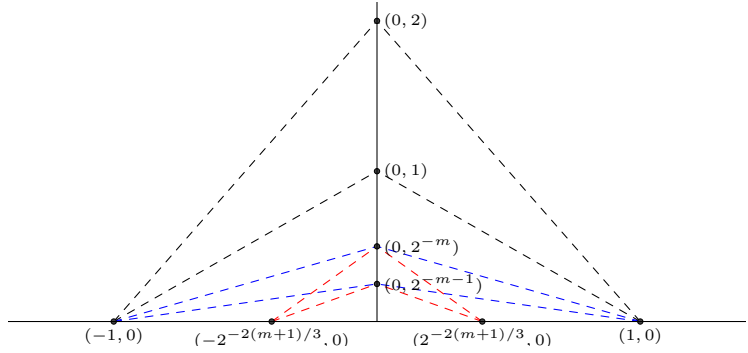
\begin{figure}[hbt!]
    \centering
\begin{tikzpicture}[line cap=round,line join=round,>=triangle 45,x=0.7cm,y=.08cm]
\clip(-2,-4) rectangle (12,53);

\draw [dashed] (0,0) -- (5,50);
\draw [dashed] (10,0) -- (5,50);
\draw [dashed] (0,0) -- (5,25);
\draw [dashed] (10,0) -- (5,25);

\draw [dashed] [blue] (0,0) -- (5,12.5);
\draw [dashed] [blue] (10,0) -- (5,12.5);
\draw [dashed] [blue] (0,0) -- (5,6.25);
\draw [dashed] [blue] (10,0) -- (5,6.25);

\draw [dashed] [red] (3,0) -- (5,12.5);
\draw [dashed] [red] (7,0) -- (5,12.5);
\draw [dashed] [red] (3,0) -- (5,6.25);
\draw [dashed] [red] (7,0) -- (5,6.25);

\draw [fill=uuuuuu] (0,0) circle (1 pt);
\draw [fill=uuuuuu] (10,0) circle (1 pt);
\draw [fill=uuuuuu] (3,0) circle (1 pt);
\draw [fill=uuuuuu] (7,0) circle (1 pt);
\draw [fill=uuuuuu] (5,25) circle (1 pt);
\draw [fill=uuuuuu] (5,50) circle (1 pt);
\draw [fill=uuuuuu] (5,12.5) circle (1 pt);
\draw [fill=uuuuuu] (5,6.25) circle (1 pt);

\draw (-10,0) -- (20,0);
\draw (5,0) -- (5,500);

\begin{tiny}
\draw (0,0) node[anchor=north]{$(-1,0)$};
\draw (3,0) node[anchor=north]{$(-2^{-2(m+1)/3},0)$};
\draw (7,0) node[anchor=north]{$(2^{-2(m+1)/3},0)$};
\draw (10,0) node[anchor=north]{$(1,0)$};
\draw (5,6.25) node[anchor=west]{$(0,2^{-m-1})$};
\draw (5,12.5) node[anchor=west]{$(0,2^{-m})$};
\draw (5,25) node[anchor=west]{$(0,1)$};
\draw (5,50) node[anchor=west]{$(0,2)$};

\end{tiny}

\end{tikzpicture}
\caption{An illustration of the proof of Proposition \ref{prop:eta-non-deg}: the black dashed lines indicate the passage times involved in $K_{1,2}$, and the blue dashed lines indicate the passage times involved in $K_{2^{-m-1},2^{-m}}$. The red dashed lines are rescaled versions of the black ones, which gives a rescaled version of $K_{1,2}$ that lower bounds $K_{2^{-m-1},2^{-m}}$, implying that $\PP[K_{2^{-m-1},2^{-m}}>0] \ge \PP[K_{1,2}>0]$. By the 0-1 law for the germ sigma algebra $\cF_{0+}$, almost surely one can find some $m\in\N$ with $K_{2^{-m-1},2^{-m}}>0$.}
\label{fig:710}
\end{figure}

We can finally put things together to finish the proof of Proposition \ref{prop:eta-non-deg}.
\begin{proof} [Proof of Proposition \ref{prop:eta-non-deg}]
For each $m\in \N,$ we consider $K_{2^{-m-1},2^{-m}}$ (recall its definition from \eqref{eq:defn-kgh}), which we assume to be well-defined and finite.
By the quadrangle inequality (Lemma \ref{lem:DL-quad}), we have
\begin{multline*}
K_{2^{-m-1},2^{-m}} = \lim_{w\searrow 0} (2w)^{-1} \int_{2^{-m-1}}^{2^{-m}}
\cD(-w,t)-\cD(w,t) dt
\\ \ge \limsup_{w\searrow 0} (2w)^{-1} \int_{2^{-m-1}}^{2^{-m}} (\cL(-2^{-2(m+1)/3},0;-w,t)-\cL(2^{-2(m+1)/3},0;-w,t))\\-(\cL(-2^{-2(m+1)/3},0;w,t)-\cL(2^{-2(m+1)/3},0;w,t)) dt,
\end{multline*}
where (by scaling invariance) the right hand side has the same distribution as $2^{-2(m+1)/3}K_{1,2}$ (see Figure \ref{fig:710}).
So we must have $\PP[K_{2^{-m-1},2^{-m}}>0] \ge \PP[K_{1,2}>0]$, which implies that $$\PP[\limsup_{m\to\infty} \don[K_{2^{-m-1},2^{-m}}>0]=1]>0$$ by Lemma \ref{lem:K-pos-posprob}.
Since $K_{2^{-m-1},2^{-m}}$ is $\cF_{2^{-m}}$ measurable, $\limsup_{m\to\infty} \don[K_{2^{-m-1},2^{-m}}>0]$ is $\cF_{0+}$ measurable. Then by Lemma \ref{lem:F-tri} we conclude that 
\[\PP[\limsup_{m\to\infty} \don[K_{2^{-m-1},2^{-m}}>0]=1]=1,\]
which implies the conclusion.
\end{proof}

\bibliographystyle{halpha}
\bibliography{bibliography}

\appendix

\section{Transversal fluctuation in the directed landscape}  \label{sec:apptrans}
In this appendix we prove Lemma \ref{lem:dl-trans-f}.
We first prove the following weaker version, which can be directly deduced from Lemma \ref{lem:DLbound}.
\begin{lemma}  \label{lem:dl-trans1}
For any $u=(x,s;y,t) \in \R_\uparrow^4$, any geodesic $\pi_u$, and $(s+t)/2\le r<t$, we have
\begin{multline}  \label{eq:dl-trans-bd-1}
\left|\pi_u(r) - \frac{x(t-r)+y(r-s)}{t-s}\right| \\ < C(R+1)^{1/2}\log(R+2)^2 (t-s)^{1/6}(t-r)^{1/2} \log^{2/3}(2(2\|u\|+2)^{3/2}/(t-r)) \log^{1/3}(2\|u\|+2).
\end{multline}
where $R$ is the random number from Lemma \ref{lem:DLbound}, and $C>0$ is a universal constant.
Similar result holds when $s<r<(s+t)/2$ by symmetry.
\end{lemma}

\begin{proof}
Denote the right hand side of \eqref{eq:dl-trans-bd-1} by $D$.
Let $z=\pi_u(r)$, and $z'=z - \frac{x(t-r)+y(r-s)}{t-s}$.
We argue by contradiction: we show that if $|z'|\ge D$, we must have
\begin{equation}  \label{eq:dl-trans-bd-pf-1}
\cL(x,s;y,t) - (\cL(x,s;z,r) + \cL(z,r;y,t)) > 0,
\end{equation}
which contradicts with the fact that $\pi_u$ is a geodesic.
Indeed, by Lemma \ref{lem:DLbound}, we have
\[
\cL(x,s;y,t) \ge - \frac{(x-y)^2}{t-s} - R(t-s)^{1/3} \log^{4/3}(2(\|u\|+2)^{3/2}/(t-s)) \log^{2/3}(\|u\|+2),
\] 
\[
-\cL(x,s;z,r) \ge \frac{(x-z)^2}{r-s} - R(r-s)^{1/3} \log^{4/3}(2(2\|u\|+|z'|+2)^{3/2}/(r-s)) \log^{2/3}(2\|u\|+|z'|+2),
\] 
\begin{equation}  \label{eq:dl-trans-bd-pf-15}
-\cL(z,r;y,t) \ge \frac{(z-y)^2}{t-r} - R(t-r)^{1/3} \log^{4/3}(2(2\|u\|+|z'|+2)^{3/2}/(t-r)) \log^{2/3}(2\|u\|+|z'|+2).
\end{equation}
Here we use that $\|(x,s;z,r)\|, \|(z,r;y,t)\| \le 2\|u\|+|z'|$.
By adding up these three inequalities, we have that the left hand side of \eqref{eq:dl-trans-bd-pf-1} is lower bounded by
\begin{equation}  \label{eq:dl-trans-bd-pf-2}
\frac{|z'|^2(t-s)}{(r-s)(t-r)} - 3R(t-s)^{1/3} \log^{4/3}(2(2\|u\|+|z'|+2)^{3/2}/(t-r)) \log^{2/3}(2\|u\|+|z'|+2).
\end{equation}
One can check that when $C$ is large enough, for any $i\in \Z_{\ge 0}$ we have
\[
\frac{(2^iD)^2(t-s)}{(r-s)(t-r)} > 3R(t-s)^{1/3} \log^{4/3}(2(2\|u\|+2^{i+1}D+2)^{3/2}/(t-r)) \log^{2/3}(2\|u\|+2^{i+1}D+2),
\]
so \eqref{eq:dl-trans-bd-pf-2} $>0$ when $|z'|\in [2^iD, 2^{i+1}D]$. Thus \eqref{eq:dl-trans-bd-pf-2} $>0$ whenever $|z'|\ge D$, and the conclusion follows.
\end{proof}
To upgrade this to Lemma \ref{lem:dl-trans-f}, we need more refined estimates in a compact subset of $\R_\uparrow^4$. These can be obtained using Lemma \ref{lem:modcont}.
\begin{lemma}  \label{lem:dl-trans2}
Let $K$ be a compact subset of $\R_\uparrow^4$, and $R$ be the maximum of the random numbers from Lemma \ref{lem:DLbound} and Lemma \ref{lem:modcont} for $K$.
For any $u=(x,s;y,t) \in \R_\uparrow^4$, any geodesic $\pi_u$, and $(s+t)/2\le r<t$, such that $(x,s;\pi_u(r),r), u \in K$,
we must have
\begin{equation}  \label{eq:dl-trans-bd-2}
\left|\pi_u(r) - \frac{x(t-r)+y(r-s)}{t-s}\right| < C(R+1)^{2/3}\log(R+2)^2(t-r)^{2/3}\log^2(C(1+(t-r)^{-1})),
\end{equation}
where $C>0$ is a constant depending only on $K$.
Similar result holds when $s<r<(s+t)/2$ by symmetry.
\end{lemma}

\begin{proof}
This proof is similar to the proof of the previous corollary.
We now denote the right hand side of \eqref{eq:dl-trans-bd-2} by $D$, and define $z, z'$ as before.
We still show that if $|z'|\ge D$, the inequality \eqref{eq:dl-trans-bd-pf-1} holds.

Note that \eqref{eq:dl-trans-bd-pf-15} still holds.
By Lemma \ref{lem:modcont}, we have
\begin{multline*}
\left(\cL(x,s;y,t) + \frac{(x-y)^2}{t-s} \right)- \left(\cL(x,s;z,r) + \frac{(x-z)^2}{r-s} \right) \\
\ge - R\left((t-r)^{1/3}\log^{2/3}(1+(t-r)^{-1}) + |y-z|^{1/2}\log^{1/2}(1+|y-z|^{-1}) \right)
\end{multline*}
Adding this to \eqref{eq:dl-trans-bd-pf-15}, we get that the left hand side of \eqref{eq:dl-trans-bd-pf-1} is lower bounded by
\begin{multline*}
\frac{|z'|^2(t-s)}{(r-s)(t-r)} - R(t-r)^{1/3} \log^{4/3}(2(2\|u\|+|z'|+2)^{3/2}/(t-r)) \log^{2/3}(2\|u\|+|z'|+2) \\ - R\left((t-r)^{1/3}\log^{2/3}(1+(t-r)^{-1}) + |y-z|^{1/2}\log^{1/2}(1+|y-z|^{-1}) \right).
\end{multline*}
Below we let $C>0$ denote a constant depending on $K$, and its value can change from line to line.
Then we have $\|u\|<C$, and $|y-z|< C|z'| + C|t-r|$, so the above can be lower bounded by
\begin{equation}  \label{eq:dl-trans-bd-25}
\frac{|z'|^2}{t-r} - CR(t-r)^{1/3} \log^2(C(1+|z'|)/(t-r)) - CR(|z'|^{1/2} + (t-r)^{1/3}) \log^{2/3}(C(1+(t-r)^{-1})) .
\end{equation}
When $|z'| > D$, the derivative of \eqref{eq:dl-trans-bd-25} with respect to $|z'|$ is $>0$; and when $|z'| = D$, \eqref{eq:dl-trans-bd-25} is positive. Thus we have that \eqref{eq:dl-trans-bd-25} is positive whenever $|z'| > D$, and the conclusion follows.
\end{proof}

\begin{proof}[Proof of Lemma \ref{lem:dl-trans-f}]
We let $C, c>0$ denote large and small universal constants, and the values may change from line to line.
Let $K = ([-1,1]\times [0,1]) \times ([-1,1]\times [2,3])$, and 
\[K_+ = \{\left(x,s; \frac{x(t-r)+y(r-s)}{t-s}, r\right): |x|, |y| \le 2, s \in [0, 1], t \in [2,3], r \in [(s+t)/2, t] \}.\]
Let $R_1$ be the random number from Lemma \ref{lem:DLbound}, and $R_2$ be the random number from Lemma \ref{lem:modcont} with $K_+$.
For any $M>2$, if $R_1, R_2 <M$, by Lemma \ref{lem:dl-trans1} and Lemma \ref{lem:dl-trans2} for $K_+$, for any $u=(x,s;y,t)\in K$ and $r \in [(s+t)/2, t]$, we have
\[
\left|\pi_u(r) - \frac{x(t-r)+y(r-s)}{t-s}\right| < CM^{1/2}\log(M)^2 (t-r)^{1/2} \log^{2/3}(1+(t-r)^{-1}),
\]
and unless the left hand below side is $>1$, we also have
\[
\left|\pi_u(r) - \frac{x(t-r)+y(r-s)}{t-s}\right| < CM^{2/3}\log(M)^2(t-r)^{2/3}\log^2((1+(t-r)^{-1})) .
\]
Denote the right hand sides of the above two estimates by $D_1$ and $D_2$, respectively. Then if $D_1\le 1$, we can bound $\left|\pi_u(r) - \frac{x(t-r)+y(r-s)}{t-s}\right|$ by $D_2$; otherwise we can bound it by $D_1^{4/3}$.
Thus we always have
\[
\left|\pi_u(r) - \frac{x(t-r)+y(r-s)}{t-s}\right| < CM^{2/3}\log(M)^{8/3}(t-r)^{2/3}\log^2((1+(t-r)^{-1})) .
\]
These imply that, if we denote
\[
R_* = \sup_{u=(x,s;y,t)\in K, r \in [(s+t)/2, t]} \frac{\left|\pi_u(r) - \frac{x(t-r)+y(r-s)}{t-s}\right|}{(t-r)^{2/3}\log^2((1+(t-r)^{-1}))},
\]
then $\PP[R_*> M]< Ce^{-cM^{9/4}\log(M)^{-4}}$ for any $M>0$.

Now for any $k, a, b, d\in \Z$, we consider
\[
K[k,a,b,d]=\{(2^{2k/3}(x+b/10), 2^k(s+a/10); 2^{2k/3}(y+d/10), 2^k(t+a/10) )  : (x,s;y,t)\in K\},
\] 
and
\[
R[k,a,b,d] = \sup_{u=(x,s;y,t)\in K[k,a,b,d], r \in [(s+t)/2, t]} \frac{\left|\pi_u(r) - \frac{x(t-r)+y(r-s)}{t-s}\right|}{(t-r)^{2/3}\log^2((1+2^k(t-r)^{-1}))}.
\]
By skew-shift invariance and scaling invariance of the directed landscape, each $R[k,a,b,d]$ has the same distribution as $R_*$. Since $\cup_{k,a,b,d\in\Z} K[k,a,b,d]=(x,s;y,t) = \R_\uparrow^4$, we can take
\[
R= C\sup_{k,a,b,d\in\Z} \frac{R[k,a,b,d]}{\log(1+|k|+|a|+|b|+|d|)}.
\]
Then \eqref{eq:dl-trans-bd-0} holds, and by a union bound over all $k$, $a$, $b$, $d \in \Z$ we have the desired bound on the upper tail of $R$.
\end{proof}

\section{Occupation time in a tube}
\label{sec:appb}

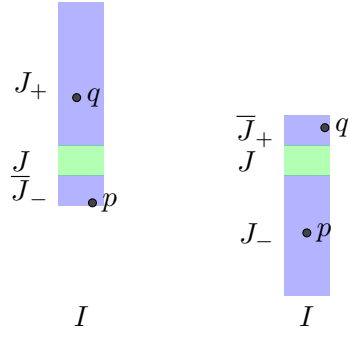
\begin{figure}[hbt!]
    \centering
\begin{tikzpicture}[line cap=round,line join=round,>=triangle 45,x=0.3cm,y=0.2cm]
\clip(-0.5,-6) rectangle (20.5,15.5);

\fill[line width=0.pt,color=blue,fill=blue,fill opacity=0.3]
(6,2) -- (6,4) -- (4,4) -- (4,2) -- cycle;

\fill[line width=0.pt,color=blue,fill=blue,fill opacity=0.3]
(6,22) -- (4,22) -- (4,6) -- (6,6) -- cycle;

\fill[line width=0.pt,color=green,fill=green,fill opacity=0.3]
(6,6) -- (6,4) -- (4,4) -- (4,6) -- cycle;

\draw [fill=uuuuuu] (5.5,2.2) circle (1.5pt);

\draw [fill=uuuuuu] (4.8,9.2) circle (1.5pt);

\draw (5.5,2.2) node[anchor=west]{$p$};
\draw (4.8,9.2) node[anchor=west]{$q$};

\fill[line width=0.pt,color=blue,fill=blue,fill opacity=0.3]
(16,6) -- (16,8) -- (14,8) -- (14,6) -- cycle;

\fill[line width=0.pt,color=blue,fill=blue,fill opacity=0.3]
(16,4) -- (16,-4) -- (14,-4) -- (14,4) -- cycle;

\fill[line width=0.pt,color=green,fill=green,fill opacity=0.3]
(16,6) -- (16,4) -- (14,4) -- (14,6) -- cycle;

\draw [fill=uuuuuu] (15,0.2) circle (1.5pt);

\draw [fill=uuuuuu] (15.8,7.2) circle (1.5pt);

\draw (15,0.2) node[anchor=west]{$p$};
\draw (15.8,7.2) node[anchor=west]{$q$};

\draw (15,-4) node[anchor=north]{$I$};
\draw (5,-4) node[anchor=north]{$I$};

\draw (3.2,5) node[anchor=east]{$J$};
\draw (4,3) node[anchor=east]{$\overline{J}_-$};
\draw (4,10) node[anchor=east]{$J_+$};

\draw (13.2,5) node[anchor=east]{$J$};
\draw (14,0) node[anchor=east]{$J_-$};
\draw (14,7) node[anchor=east]{$\overline{J}_+$};

\end{tikzpicture}
\caption{
An illustration of the set $\Phi_{I,J}'$: for $(p;q)=(x,s;y,t)\in \Phi_{I,J}'$, either $s\in \overline{J}_-$ and $t\in J_+$, as illustrated in the left diagram; or $s\in J_-$ and $t\in \overline{J}_+$, as illustrated in the right diagram. The spatial coordinates $x, y$ are always in $I$.
}
\label{fig:defphi}
\end{figure}
As promised before, in this appendix we provide the proof of Lemma \ref{lem:geo-ocp-tail}. Now, recalling the set $\Phi_{I,J}$ and the random variable $W_{I,J}$ appearing in the discussion preceding the statement of the lemma, it would be simpler to consider the following, similar but different, objects.
Namely, for each interval $J=[s,t]$, let $J_-=(-\infty, s]$, $J_+=[t, \infty)$, and $\overline{J}_-=[2s-t, s]$, $\overline{J}_+=[t, 2t-s]$.
We define
\[
\Phi_{I,J}' = \{ (p;q) \in \R^4_\uparrow: p \in I\times J_-, q \in I\times \overline{J}_+\} \cup \{ (p;q) \in \R^4_\uparrow: p \in I\times \overline{J}_-, q \in I\times J_+\}.
\]
See Figure \ref{fig:defphi} for an illustration.
Then we let
\[
W_{I,J}' =  \int_J \sup_{u\in \Phi_{I,J}'} \don[\pi_u(r)\in I, \pi_u \text{is unique}] dr .
\]
Note that in the definition of $W_{I,J}'$, the supremum over $u\in \Phi_{I,J}'$ is inside the integral over $J$. This upper bounds the occupation measure in $I\times J$ for any unique geodesic $\pi_u$ (with $u\in \Phi_{I,J}'$).

We then have the following crucial estimate.

\begin{lemma}  \label{lem:exp-Wp}
We have $\E[W_{I,J}']<C_1\LE(I)\LE(J)^{1/3}$ for some universal constant $C_1$.
\end{lemma}
An upper bound of the order of $\LE(I)\LE(J)^{1/3}$ can be explained by a simple scaling argument as follows.
First, we expect $\E[W_{I,J}']$ to be linear in $\LE(I)$, at least when $\LE(I)$ is small (compared to $\LE(J)^{2/3}$).
Second, by scaling invariance of $\cL$, when $\LE(J)$ is multiplied by $h>0$ and $\LE(I)$ is multiplied by $h^{2/3}$, $\E[W_{I,J}']$ should be multiplied by $h$.

Compared to $\E[W_{I,J}]$, the expectation $\E[W_{I,J}']$ is easier to bound: since the time coordinates of the endpoints are outside or at the boundary of $J$ (whereas in the definition of $\Phi_{I,J}$ in \eqref{eq:defn-phi-ij}, the endpoints are on the boundary of $I\times J$), we can directly bound $\E[W_{I,J}']$ using translation invariance of $\cL$, and bounds on the number of intersections that geodesics can have with a line (from Lemma \ref{lem:number-of-inter}).

This is then used in an inductive proof, where we work on different dyadic scales. The definition of $\Phi_{I,J}'$ enables one to bound the occupation time for a geodesic in different time intervals.

We also remark that, actually, the same upper bound holds for the expectation of the integral $\int_J \sup_{u\in \Phi_{I,J}', \pi_u \text{ is any geodesic}} \don[\pi_u(r)\in I] dr$; i.e., the uniqueness requirement of the geodesic $\pi_u$ can be dropped and all possible ones can be considered instead. The uniqueness requirement is inherited from Lemma \ref{lem:number-of-inter}, and can actually be removed (see the discussions after Lemma \ref{lem:number-of-inter}).

We now prove Lemma \ref{lem:geo-ocp-tail} assuming Lemma \ref{lem:exp-Wp}, and using an induction argument. A key ingredient is a dyadic decomposition of the time interval $[0,1]$.
Further, the exponential tail (in addition to bounding the expectations) in Lemma \ref{lem:geo-ocp-tail} is proved via a concentration estimate, using independence of the directed landscape across different time increments.

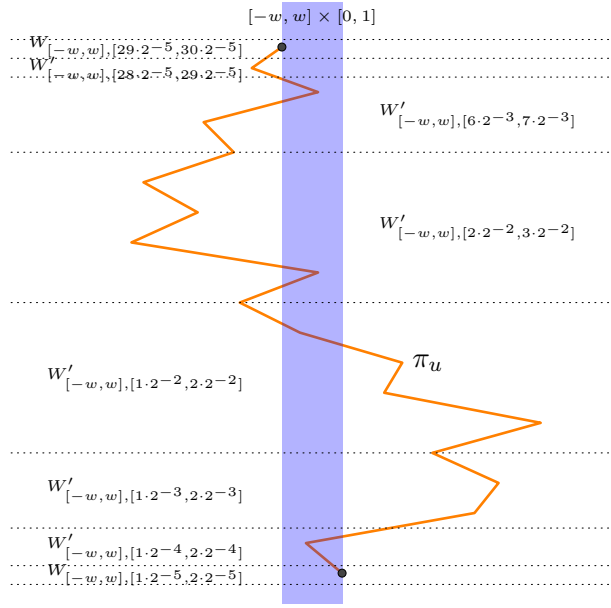
\begin{figure}[hbt!]
    \centering
\begin{tikzpicture}[line cap=round,line join=round,>=triangle 45,x=0.8cm,y=0.4cm]
\clip(-0.5,-0.5) rectangle (10.5,20.5);

\draw [line width=1pt, color=orange] (5.5,1) -- (4.9,2) -- (7.7,3) -- (8.1,4) -- (7,5) -- (8.8,6) -- (6.2,7) -- (6.5,8) -- (4.8,9) -- (3.8,10) -- (5.1,11) -- (2,12) -- (3.1,13) -- (2.2,14) -- (3.7,15) -- (3.2,16) -- (5.1,17) -- (4,17.8) -- (4.5,18.5);

\fill[line width=0.pt,color=blue,fill=blue,fill opacity=0.3]
(4.5,0) -- (5.5,0) -- (5.5,20) -- (4.5,20) -- cycle;

\draw (0,25) -- (10,25);

\draw [fill=uuuuuu] (5.5,1) circle (1.5pt);
\draw [fill=uuuuuu] (4.5,18.5) circle (1.5pt);

\foreach \i in {1.25-5/8,1.25, 2.5, 5, 10, 15, 17.5, 17.5+5/8, 17.5+5/4}
{
\draw [dotted] (0, \i) -- (10, \i);
}

\begin{tiny}
\draw (5,20) node[anchor=north]{$[-w,w]\times [0,1]$};
\draw (4,18.5) node[anchor=east]{$W_{[-w,w],[29\cdot 2^{-5}, 30\cdot 2^{-5}]}$};
\draw (4,17.75) node[anchor=east]{$W'_{[-w,w],[28\cdot 2^{-5}, 29\cdot 2^{-5}]}$};
\draw (4,7.5) node[anchor=east]{$W'_{[-w,w],[1\cdot 2^{-2}, 2\cdot 2^{-2}]}$};
\draw (4,3.75) node[anchor=east]{$W'_{[-w,w],[1\cdot 2^{-3}, 2\cdot 2^{-3}]}$};
\draw (4,1.875) node[anchor=east]{$W'_{[-w,w],[1\cdot 2^{-4}, 2\cdot 2^{-4}]}$};
\draw (4,1.875/2) node[anchor=east]{$W_{[-w,w],[1\cdot 2^{-5}, 2\cdot 2^{-5}]}$};

\draw (6,12.5) node[anchor=west]{$W'_{[-w,w],[2\cdot 2^{-2}, 3\cdot 2^{-2}]}$};
\draw (6,16.25) node[anchor=west]{$W'_{[-w,w],[6\cdot 2^{-3}, 7\cdot 2^{-3}]}$};
\end{tiny}
\draw (6.5,8) node[anchor=west]{$\pi_u$};

\end{tikzpicture}
\caption{
An illustration of \eqref{eq:coverW} in the proof of Lemma \ref{lem:geo-ocp-tail}: for $u\in \Phi_{[-w,w],[0,1]}$, we can cover the domain of $\pi_u$ by intervals $J$ (that are disjoint except for the endpoints), such that for $J$ other than the first and last one, $J\in \Lambda$ and the time of $\pi_u\in [-w,w]$ in $J$ is at most $W'_{[-w,w],J}$; and for $J$ being the first or last one, $J\in \partial\Lambda$ and the time of $\pi_u\in [-w,w]$ in $J$ is at most $W_{[-w,w],J}$.
}
\label{fig:2}
\end{figure}

\begin{proof}[Proof of Lemma \ref{lem:geo-ocp-tail}]
We use induction in $w$, and prove the following estimates: there exist constants $A, C, c>0$, such that for any $w>0$, there is
\begin{equation}  \label{eq:W-exp}
\E[W_{[-w,w], [0,1]}] < A^{2/3+0.01} w,    
\end{equation}
and for any $M>0$,
\begin{equation}  \label{eq:W-tail}
\PP[W_{[-w,w], [0,1]} > AMw] < Ce^{-cM}.
\end{equation}
We assume that $A, C$ are large enough and $c$ is small enough, and the precise relation between them is evident from the proof. 

The case when $w>A^{-2/3}$ is obvious (for \ref{eq:W-tail} just take $C$ large enough), so we assume that $w\le A^{-2/3}$. 
Now to use an induction argument, we assume that the inequalities hold when $w$ is replaced by any number $\ge 2^{2/3}w$.

We first prove \eqref{eq:W-exp}.
Take $g\in\N$ such that $2^g\le A<2^{g+1}$. Let $\Lambda$ consist of all intervals $[a2^{-b}, (a+1)2^{-b}]$ such that $a,b\in \Z$, $0 \le b\le g$, and $0\le a2^{-b} < (a+1)2^{-b} \le 1$.
Let $\partial\Lambda$ consist of all intervals $[a2^{-g}, (a+1)2^{-g}]$ such that $a\in \Z$, $0\le a2^{-g} < (a+1)2^{-g} \le 1$.
We have that
\begin{equation}  \label{eq:coverW}
W_{[-w,w], [0,1]} \le \sum_{J\in \Lambda} W_{[-w,w],J}' + 2\max_{J\in \partial\Lambda}W_{[-w,w],J}.
\end{equation}
This is via a dyadic decomposition of $\pi_u$ for $u\in \Phi_{[-w,w], [0,1]}$, as illustrated by Figure \ref{fig:2}.The precise details of the proof of \eqref{eq:coverW} are left for the moment for later.

By Lemma \ref{lem:exp-Wp}, for each $J=[a2^{-b}, (a+1)2^{-b}]\in \Lambda$, 
we have $\E[W_{[-w,w],J}'] < 2C_1w2^{-b/3}$, thus \[\sum_{J\in \Lambda}\E[W_{[-w,w],J}'] < \sum_{b=0}^g \sum_{a=0}^{2^b-1}2C_1w2^{-b/3} = \sum_{b=0}^g 2C_1w2^{2b/3} < 10C_1w2^{2g/3}.\] So we have
\[
\E[W_{[-w,w], [0,1]}] < 10C_1w 2^{2g/3} + 2\E[\max_{J\in \partial\Lambda}W_{[-w,w],J}].
\]
By induction hypothesis and scale invariance, for each $J\in \partial\Lambda$ we have
\[
\PP[W_{[-w,w],J} > AMw 2^{-g/3}] < Ce^{-cM}.
\]
Thus we have
\[
\E[\max_{J\in \partial\Lambda}W_{[-w,w],J}] < c^{-1}(\log(2^g)+C)A w 2^{-g/3},
\]
and this implies that 
\[
\E[W_{[-w,w], [0,1]}] < (10C_1 + 4c^{-1}(\log(A)+C)) A^{2/3}w < A^{2/3+0.01}w
\]
by requiring $A$ being large enough so that $A^{0.01}>(10C_1 + 4c^{-1}(\log(A)+C))$.

We next prove \eqref{eq:W-tail}. 
Denote $m=\lfloor A^{1/10} \rfloor$.
For each $1\le i \le m$ let $W_i=W_{[-w,w],[(i-1)/m,i/m]}$, and we have $W_{[-w,w],[0,1]}\le \sum_{i=1}^m W_i$.
We shall get the exponential tail of $W_{[-w,w],[0,1]}/(Aw)$ by considering the exponential moments of each $W_i$, using the mutual independence of the latter across $i,$ and applying the induction hypothesis for them.

We write
\[
\sum_{i=1}^m W_i \le \sum_{i=1}^m W_i\don[W_i\le Aw/m] + \sum_{i=1}^m W_i\don[W_i> Aw/m] \le Aw + \sum_{i=1}^m W_i\don[W_i> Aw/m].
\]
When $A$ is large enough, we have $m\ge 2$, so we can apply the induction hypothesis and use scaling invariance to conclude that
$$\E[W_i]\le A^{2/3+0.01}wm^{-1/3} \quad\quad \text{and}\,\,\quad \quad \PP[W_i>MAwm^{-1/3}]<Ce^{-cM},$$ for any $M>0$.
From these bounds we expect that the sum of $W_i$ should have exponential tail at the order of $Aw$, and we analyze it as follows.
Let $X_i = \frac{m^{1/3}}{Aw}W_i\don[W_i> Aw/m]$, then for any $M>0$, by the exponential tail we have
\[
\PP[X_i>M] < Ce^{-cM}.
\]
Using Markov's inequality, for any $M>0$ we also have
\[
 \PP[X_i>M] \le \PP\left[ W_i > Aw/m \right] < \E[W_i]/(Aw/m) \le A^{-1/3+0.01}m^{2/3}.
\]
Take $M_*>0$ such that $Ce^{-cM_*} = A^{-1/3+0.01}m^{2/3}$.
Take $a=c/10$. We then have
\begin{multline*}
\E[e^{aX_i}] = 1 + \int_0^\infty ae^{aM}\PP[X_i>M] dM \\
< 1 + \int_0^{M_*} ae^{aM}A^{-1/3+0.01}m^{2/3} dM + \int_{M_*}^\infty ae^{aM}Ce^{-cM} dM 
< 1 + 2e^{aM_*} A^{-1/3+0.01}m^{2/3},
\end{multline*}
and hence
\[
\E[e^{a\sum_{i=1}^mX_i}] < (1 + 2e^{aM_*} A^{-1/3+0.01}m^{2/3})^m
<
e^{C(A^{-1/3+0.01}m^{2/3})^{1/2}m},
\]
where the last inequality follows by using  
\[
2e^{aM_*} A^{-1/3+0.01}m^{2/3} < C(A^{-1/3+0.01}m^{2/3})^{1/2}m,
\]
which is satisfied by taking $C>10$.
Then we have
\[
\PP[\sum_{i=1}^m W_i > AMw] \le \PP[\sum_{i=1}^mX_i > m^{1/3}(M-1)] < e^{C(A^{-1/3+0.01}m^{2/3})^{1/2}m - am^{1/3}(M-1)} < Ce^{-cM},
\]
where the last inequality is true by taking $A$ large enough (so that $am^{1/3}>2c$ and $e^{Cm^{-1/3}} < C$), for $M>2$.
For $0< M\le 2$ we also have \eqref{eq:W-tail}, by taking $c$ small and $C$ large.
\end{proof}
\begin{proof}[Proof of \eqref{eq:coverW}]
Consider any $u=(x,s;y,t)\in \Phi_{[-w,w], [0,1]}$. The general strategy is to show that $[s, t]$ can be covered the union of (1) two intervals in $\partial \Lambda$ (around the endpoints) and (2) several intervals in $\Lambda$ that are contained in $[s, t]$.
Then we can bound $W[u,[-w,w]]$ by the sum of $W_{[-w,w],J}$ or $W_{[-w,w],J}'$ for $J$ among these intervals.

For each $i\in \N$, we let $a_i^-, a_i^+\in \Z$ be the numbers with $a_i^-2^{-i}\le s < (a_i^-+1)2^{-i}$, and $(a_i^+-1)2^{-i}< t \le a_i^+2^{-i}$.
Consider collections of intervals: 
\[
\partial\Lambda_u=\{ [a_g^-2^{-g}, (a_g^-+1)2^{-g}], [(a_g^+-1)2^{-g}, a_g^+2^{-g}] \},
\]
and
\[
\Lambda_u=\{ [(a_i^-+1)2^{-i}, (a_i^-+2)2^{-i}], [(a_i^+-2)2^{-i}, (a_i^+-1)2^{-i}] : i\in \llbracket 1, g \rrbracket, a_i^+-a_i^- > 2\}.
\]
We note that $\partial\Lambda_u \subset \partial\Lambda$ and $\Lambda_u\subset \Lambda$. We claim that the union of the intervals in $\partial\Lambda_u\cup \Lambda_u$ covers $[s, t]$.
Indeed, if $a_g^+-a_g^-\le 2$, the union of the intervals in $\partial\Lambda_u$ covers $[s, t]$.
Otherwise, let $i_*$ be the smallest integer with $a_{i_*}^+-a_{i_*}^- > 2$, and we must have $1< i_* \le g$.
We next show that the intervals in $\partial\Lambda_u\cup \Lambda_u$ can be written in a sequence:
\begin{multline*}
[a_g^-2^{-g}, (a_g^-+1)2^{-g}], [(a_g^-+1)2^{-g}, (a_g^-+2)2^{-g}] , \cdots, [(a_{i_*}^-+1)^-2^{-i_*}, (a_{i_*}^-+2)2^{-i_*}], 
\\
[(a_{i_*}^+-2)2^{-i_*}, (a_{i_*}^+-1)2^{-i_*}] ,\cdots, [(a_g^+-2)2^{-g}, (a_g^+-1)2^{-g}] , [(a_g^+-1)2^{-g}, a_g^+2^{-g}],
\end{multline*}
such that any two neighboring ones are not disjoint. Then the union of all these intervals is a closed interval.
\begin{enumerate}
  \item The facts that $[a_g^-2^{-g}, (a_g^-+1)2^{-g}]\cap [(a_g^-+1)2^{-g}, (a_g^-+2)2^{-g}] \neq \emptyset$ and $[(a_g^+-2)2^{-g}, (a_g^+-1)2^{-g}] \cap [(a_g^+-1)2^{-g}, a_g^+2^{-g}] \neq \emptyset$ are obvious.
  \item We must have that $a_{i_*}^+-a_{i_*}^- \le 4$. This is because, otherwise we have $a_{i_*}^+2^{-i_*}-a_{i_*}^-2^{-i_*} > 2^{-i_*+2}$. From the definitions (of $a_i^-, a_i^+$) above we have
\[
a_{i_*-1}^-2^{-i_*+1} \le a_{i_*}^-2^{-i_*}, \quad a_{i_*}^+2^{-i_*} \le a_{i_*-1}^+2^{-i_*+1},\]
so $a_{i_*-1}^+2^{-i_*+1}-a_{i_*-1}^-2^{-i_*+1} > 2^{-i_*+2}$, implying that $a_{i_*-1}^+-a_{i_*-1}^- > 2$, which contradicts with the choice of $i_*$.
Thus we have $a_{i_*}^+-a_{i_*}^- =3$ or $4$; and in each case we have that $[(a_{i_*}^-+1)^-2^{-i_*}, (a_{i_*}^-+2)2^{-i_*}]$ and $[(a_{i_*}^+-2)2^{-i_*}, (a_{i_*}^+-1)2^{-i_*}]$ are not disjoint.
  \item From the definition (of $a_i^-, a_i^+$) above, for any $i\in\N$ we have that
\[
a_i^-2^{-i} \le a_{i+1}^-2^{-i-1}, \; (a_{i+1}^-+1)2^{-i-1} \le (a_i^-+1)2^{-i}, \; 
\]
so the intervals $[(a_i^-+1)^-2^{-i}, (a_i^-+2)2^{-i}]$ and $[(a_{i+1}^-+1)^-2^{-i-1}, (a_{i+1}^-+2)2^{-i-1}]$ are not disjoint.
Similarly, the intervals $[(a_i^+-2)2^{-i}, (a_i^+-1)2^{-i}]$ and $[(a_{i+1}^+-2)2^{-i-1}, (a_{i+1}^+-1)2^{-i-1}]$ are not disjoint.
\end{enumerate}
From the above statements on non-disjointness of intervals in $\partial\Lambda_u\cup \Lambda_u$, we conclude that the union of these intervals is a closed interval. The union also contains $s$ and $t$, so it must contain $[s,t]$.

Further, because of the way the definitions were set up we crucially now have the following (see Figure \ref{fig:2}).
For each $J\in \Lambda_u$, we have
\[
\int_J \don[\pi_u(r)\in [-w,w], \pi_u \text{ is unique}] dr \le W_{[-w,w], J}',
\]
and for each $J\in \partial\Lambda_u$, we have
\[
\int_{J\cap[s,t]} \don[\pi_u(r)\in [-w,w], \pi_u \text{ is unique}] dr \le W_{[-w,w], J}.
\]
Thus we have 
\[
W[u,[-w,w]]\le \sum_{J\in \Lambda} W_{[-w,w],J}' + 2\max_{J\in \partial\Lambda}W_{[-w,w],J},
\]
and that \eqref{eq:coverW} holds.
\end{proof}

It remains to prove the bound for $\E[W_{I,J}']$.
We would again do a dyadic decomposition, but now for the locations of the possible starting and ending points of the geodesics.
\begin{proof}[Proof of Lemma \ref{lem:exp-Wp}]
As $\cL$ is translation and scaling invariant we can assume that $J=[0,1]$, and $I=[-w,w]$ for some $w>0$.

For each $0\le t \le 1$, we consider 
\[
A_t = \sup_{u=(p;q):p\in [-w,w]\times [-1,0], q\in [-w,w]\times [1,\infty) } \don[\pi_u(t) \in [-w,w], \pi_u \text{is unique}],
\]
and
\[
A_t' = \sup_{u=(p;q):p\in [-w,w]\times (-\infty,0], q\in [-w,w]\times [1,2] } \don[\pi_u(t) \in [-w,w], \pi_u \text{is unique}].
\]
We then have 
\[
\E[W_{[-w,w],[0,1]}'] \le \int_0^1 \E[A_t] + \E[A_t'] dt = 2\int_0^1 \E[A_t] dt \le 4w^{3/2} + 2\int_{w^{3/2}}^{1-w^{3/2}} \E[A_t] dt,
\]
where the equality is by the fact that $\int_0^1 \E[A_t]dt=\int_0^1 \E[A_t']dt$, due to symmetry.
We next show that, there is a universal constant $C_2$, such that for any $w^{3/2}<t<1/2$, 
\begin{equation}  \label{eq:A-t1}
\E[A_t] < C_2 wt^{-2/3}
\end{equation}
and for $1/2<t<1-w^{3/2}$
\begin{equation}  \label{eq:A-t2}
\E[A_t] < C_2 w(1-t)^{-2/3}.
\end{equation}
Given the above, integrating over $t$ using \eqref{eq:A-t1} and \eqref{eq:A-t2} completes the proof.
Thus it remains to prove these two inequalities.

We prove \eqref{eq:A-t1} first.
The general idea is again a dyadic decomposition, but now for the location of the endpoints. Towards this we define the following sets.
For each $a \in \N$ with $a\le \lceil -\log_2(t) \rceil+1$, we denote
\[
A_{t,a} = \sup_{u=(p;q):p\in [-w,w]\times [-(2^a-1)t,-(2^{a-1}-1)t], q\in [-w,w]\times [1,\infty) } \don[\pi_u(t) \in [-w,w], \pi_u \text{is unique}].
\]
Then we have $\E[A_t] \le \sum_{a=1}^{\lceil -\log_2(t) \rceil+1} \E[A_{t,a}]$.
To bound each $\E[A_{t,a}]$, we use translation invariance.
Denote
\begin{multline*}
B_{t,a} = |[-(2^at)^{2/3}, (2^at)^{2/3}] \cap \\ \{\pi_u(t): u=(p;y,(2^{a-3}+1)t), p\in [-(2^at)^{2/3}, (2^at)^{2/3}] \times [-(2^a-1)t,-(2^{a-1}-1)t], \pi_u \text{is unique} \} |.    
\end{multline*}
In words, for $B_{t,a}$, we allow the lower endpoint to be in a wider box (compared to $[-w,w]\times [-(2^a-1)t,-(2^{a-1}-1)t]$ in $A_{t,a}$), and the upper endpoint to vary in the line $\R \times \{(2^{a-3}+1)t\}$.
Note that since $(2^{a-3}+1)t \le 1$, we can bound $A_{t,a}$ by
\[
\sup_{u=(p;y,(2^{a-3}+1)t):p\in [-w,w]\times [-(2^a-1)t,-(2^{a-1}-1)t] } \don[\pi_u(t) \in [-w,w], \pi_u \text{is unique}].
\]
Also note that $[-(2^at)^{2/3}, (2^at)^{2/3}]$ contains $\lfloor (2^at)^{2/3}/w \lfloor$ non-overlaping translations of $[-w, w]$.
Then {$$\E[A_{t,a}]\le \lfloor (2^at)^{2/3}/w \lfloor^{-1}  \E[B_{t,a}] \le 2w(2^at)^{-2/3} \E[B_{t,a}]$$} by translation invariance.
We then have
\begin{multline*}
\E[A_t] \le \sum_{a=1}^{\lceil -\log_2(t) \rceil+1} \E[A_{t,a}] \le \sum_{a=1}^{\lceil -\log_2(t) \rceil+1} 2w(2^at)^{-2/3}
\E[B_{t,a}] = \sum_{a=1}^{\lceil -\log_2(t) \rceil+1} 2w(2^at)^{-2/3} \E[B] \\< 10wt^{-2/3}\E[B],    
\end{multline*}
where
\[
B = |[-1, 1] \cap \\ \{\pi_u(0): u=(p;y,1/8), p\in [-1, 1] \times [-1,-1/2], \pi_u \text{is unique} \} |,
\]
and by translation and scaling invariance we have each $\E[B_{t,a}] = \E[B]$. Finally, by the following claim we get \eqref{eq:A-t1}.
\begin{cla}  \label{cla:bd-B}
$\E[B]$ is finite.
\end{cla}
\begin{proof}[Proof of the claim]
Take $M>0$ large enough.
Let $\cE_M$ denote the event
\[
|\{\pi_{(x,-1/4;y,1/8)}(0): |x|, |y| \le M^{0.01}, \pi_{(x,-1/4;y,1/8)} \text{is unique} \}| \ge M.
\]
By Lemma \ref{lem:number-of-inter}, and translation and scaling invariance of the directed landscape, we have that $\PP[\cE_M] \le Ce^{-cM^{1/384}}$ for some universal constants $C,c>0$.

We next take $R$ from Lemma \ref{lem:dl-trans-f}. We show that $\{B>M\} \cap \{R\le M^{0.001}\} \subset \cE_M$.
Indeed, by Lemma \ref{lem:dl-trans-f}, if $\pi_{(x,s;y,1/8)}(0)\in[-1,1]$ for some $y\in\R$, $x\in[-1,1]$, $s\in[-1,-1/2]$, we must have that $|y|<C(R+1)\log^3(R+1)$, and $|\pi_{(x,s;y,1/8)}(-1/4)|<C(R+1)\log^3(R+1)$, for some universal constant $C>0$.
Thus if $B>M$, we can find $x_1,\cdots, x_M$ and $y_1,\cdots, y_M$, all in the interval $[-C(R+1)\log^3(R+1), C(R+1)\log^3(R+1)]$, such that for each $1\le i \le M$, $\pi_{(x_i,-1/4;y_i,1/8)}$ is unique; and $\pi_{(x_i,-1/4;y_i,1/8)}(0)$ are mutually different and are all in $[-1,1]$.
If (in addition) $R\le M^{0.001}$, we must have that $C(R+1)\log^3(R+1)< M^{0.01}$, so $\cE_M$ holds.

Thus we have $\{B>M\} \cap \{R\le M^{0.001}\} \subset \cE_M$, which implies that $\PP[B>M] < \PP[R>M^{0.001}] + \PP[\cE_M]$.
Then by $\PP[\cE_M] \le Ce^{-cM^{1/384}}$, and the bound on $\PP[R>M^{0.001}]$ from Lemma \ref{lem:dl-trans-f}, the conclusion follows.
\end{proof}

The proof of \eqref{eq:A-t2} is similar. Denote
\[
\hat{A}_t= \sup_{u=(p;q):p\in [-w,w]\times [-1,0], q\in [-w,w]\times [1,2] } \don[\pi_u(t) \in [-w,w], \pi_u \text{is unique}]
\]
and
\[
\overline{A}_t= \sup_{u=(p;q):p\in [-w,w]\times [-1,0], q\in [-w,w]\times [2,\infty) } \don[\pi_u(t) \in [-w,w], \pi_u \text{is unique}].
\]
We the have $\E[A_t]\le \E[\hat{A}_t] + \E[\overline{A}_t] \le \E[A_{1-t}] + \E[\overline{A}_t]$.
Thus by \eqref{eq:A-t1} it suffices to bound $\E[\overline{A}_t]$.
Via translation invariance, we have $\E[\overline{A}_t]<2w\E[B]$. Thus by Claim \ref{cla:bd-B} we get \eqref{eq:A-t2}. 
\end{proof}

\section{Estimates using the Airy line ensemble}  \label{sec:appa}
In this appendix we include the arguments left from the proof of Lemma \ref{lem:K-pos-posprob}.
Recall the events $\cE_{coal}'$ and $\cE_{disj}$ appearing there and that it remained to prove the following result.
\begin{lemma}  \label{lem:sur-ALE}
We have $\PP[\cE_{coal}'\cap \cE_{disj}]>0$.
\end{lemma}

To prove this we use the connection between the directed landscape and the Airy line ensemble $\cA=\{\cA_n\}_{n=1}^\infty$.
As was the key observation in \cite{DOV} and also later exploited in \cite{DZ}, marginals of the directed landscape can be obtained from last passage percolation across $\cA$, and the disjointness of geodesics in the directed landscape can be reduced to disjointness of geodesics across $\cA$.
What makes it useful to work with $\cA$ is a resampling invariance property introduced and termed as the Brownian Gibbs property in \cite{CH}. Informally, this implies that fixing $\{\cA_n\}_{n=2}^\infty$ and $\cA_1$ outside some interval; the conditional law of the latter is simply given by a Brownian bridge joining the values at the endpoints of the intervals and further conditioned to avoid $\cA_2.$
Given this, we will reduce computing probabilities of events about coalescence properties of geodesics to certain Brownian bridge computations.

We need a bit of preparation to implement the above strategy. First, as $\cA$ has infinitely many lines, the last passage times and geodesics across it are not completely straightforward to define; for this we quote the setup developed in \cite{DZ}, in Section \ref{ssec:p-ale}.
To study the disjointness of geodesics, we will need to understand last passage percolation across $\cA$ whose ending point lives on the second line $\cA_2$; we call such optimal paths the ``second geodesics''  and show their existence using the framework of ``disjoint optimizers'' from \cite{DZ}, in Section \ref{ssec:ale-2nd}.
The proof of Lemma \ref{lem:sur-ALE} is finally given in Section \ref{ssec:sur-ALE}.

\begin{figure}
  \centering
  \begin{tikzpicture}[line cap=round,line join=round,>=triangle 45,x=3cm,y=3.5cm]
  \clip(-0.15,-0.15) rectangle (2.15,1.15);
 
  \draw (0.,0.1) node[anchor=east]{$5$};
  \draw (0.,0.3) node[anchor=east]{$4$};
  \draw (0.,0.5) node[anchor=east]{$3$};
  \draw (0.,0.7) node[anchor=east]{$2$};
  \draw (0.,0.9) node[anchor=east]{$1$};
  
  \draw plot coordinates {(0.1,0.1) (0.2,0.) (0.3,0.11) (0.4,0.16) (0.5,0.04) (0.6,0.08) (0.7,0.06) (0.8,0.01) (0.9,0.02) (1.,0.14) (1.1,0.12) (1.2,0.03) (1.3,0.06) (1.4,0.) (1.5,0.04) (1.6,0.12) (1.7,0.06) (1.8,0.04) (1.9,0.05) (2.,0.07) };
  \draw plot coordinates {(0.1,0.3) (0.2,0.38) (0.3,0.21) (0.4,0.29) (0.5,0.34) (0.6,0.28) (0.7,0.26) (0.8,0.31) (0.9,0.32) (1.,0.23) (1.1,0.22) (1.2,0.33) (1.3,0.26) (1.4,0.2) (1.5,0.36) (1.6,0.33) (1.7,0.26) (1.8,0.28) (1.9,0.37) (2.,0.33) };
  \draw plot coordinates {(0.1,0.5) (0.2,0.54) (0.3,0.59) (0.4,0.46) (0.5,0.44) (0.6,0.48) (0.7,0.46) (0.8,0.51) (0.9,0.58) (1.,0.49) (1.1,0.45) (1.2,0.53) (1.3,0.56) (1.4,0.49) (1.5,0.54) (1.6,0.52) (1.7,0.56) (1.8,0.47) (1.9,0.45) (2.,0.46) };
  \draw plot coordinates {(0.1,0.7) (0.2,0.74) (0.3,0.63) (0.4,0.73) (0.5,0.64) (0.6,0.78) (0.7,0.66) (0.8,0.71) (0.9,0.65) (1.,0.63) (1.1,0.71) (1.2,0.63) (1.3,0.66) (1.4,0.77) (1.5,0.71) (1.6,0.78) (1.7,0.63) (1.8,0.64) (1.9,0.75) (2.,0.72) };
  \draw plot coordinates {(0.1,0.9) (0.2,0.98) (0.3,0.81) (0.4,0.87) (0.5,0.84) (0.6,0.91) (0.7,0.86) (0.8,0.85) (0.9,0.82) (1.,0.93) (1.1,0.82) (1.2,0.96) (1.3,0.96) (1.4,0.91) (1.5,0.88) (1.6,0.91) (1.7,0.87) (1.8,0.94) (1.9,0.95) (2.,0.93) };
  
  \draw [thick] [blue] plot coordinates {(0.3,0.21) (0.4,0.29)};
  \draw [dashed] [blue] plot coordinates {(0.4,0.29) (0.4,0.46)};
  \draw [thick] [blue] plot coordinates {(0.4,0.46) (0.5,0.44) (0.6,0.48) (0.7,0.46) (0.8,0.51) (0.9,0.58)};
  \draw [dashed] [blue] plot coordinates {(0.9,0.58) (0.9,0.65)};
  \draw [thick] [blue] plot coordinates {(0.9,0.65) (1.,0.63) (1.1,0.71) (1.2,0.63) (1.3,0.66) (1.4,0.77) (1.5,0.71) (1.6,0.78)};
  \draw [dashed] [blue] plot coordinates {(1.6,0.78) (1.6,0.91)};
  \draw [thick] [blue] plot coordinates {(1.6,0.91) (1.7,0.87) (1.8,0.94)};

\draw [fill=uuuuuu] (0.3,0.11) circle (1.pt);
\draw [fill=uuuuuu] (1.8,0.94) circle (1.pt);
\draw [blue] [fill=white] (0.4,0.29) circle (1.5pt);
\draw [blue] [fill=white] (0.9,0.58) circle (1.5pt);
\draw [blue] [fill=white] (1.6,0.78) circle (1.5pt);

  \begin{tiny}
\draw (0.3,0.11) node[anchor=south]{$(x,5)$};
\draw (1.8,0.94) node[anchor=south]{$(y,1)$};
\end{tiny}

  \end{tikzpicture}
  \caption{An illstration of a non-increasing c\`adl\`ag function across $\cA$ from $(x,5)$ to $(y,1)$.}   \label{fig:cadlag}
\end{figure}
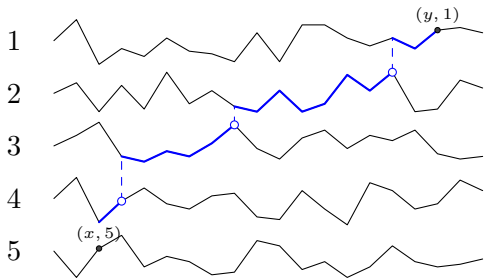

\subsection{Parabolic path across the Airy line ensemble: weight, geodesics, and disjointness}  \label{ssec:p-ale}
For any $x\le y$, and any non-increasing c\`adl\`ag function $\pi:[x,y]\to \N$, we define the \emph{weight of $\pi$ across $\cA$} as
\[
\|\pi\|_\cA = \sum_{i=\pi(y)}^{\pi(x)} \cA_i(y_i)-\cA_i(x_i),
\]
where $x_i=\inf\{z:\pi(z)\le i\}$ and $y_i=\sup\{z:\pi(z)\ge i\}$.
Now, for any $x\le y$ and $n\ge m \in \N$, we define the \emph{last passage value from $(x,n)$ to $(y,m)$} as
\[
\cA[(x,n) \to (y,m)] = \sup_{\pi} \|\pi\|_\cA,
\]
where the supremum is over all non-increasing c\`adl\`ag functions $\pi$ on $[x,y]$, with $\pi(x)\le n$ and $\pi(y)\ge m$ (see Figure \ref{fig:cadlag}). Continuity of $\cA$ and compactness of the set of functions ensures that this supremum is always attained.

We let $\cS:\R^2 \to \R$ be the Airy sheet constructed in \cite{DOV}, which is the fixed time marginal of the directed landscape.
Namely, the function $(x, y)\mapsto \cL(x,0;y,1)$ has the same distribution as $\cS$.
It can be coupled with $\cA$ such that
\begin{itemize}
  \item $\cS(0, y) = \cA_1(y)$ for $y \in \R$.
  \item For $x > 0, z_1<z_2$, we have 
  \begin{equation} \label{eq:ale-diff-asy}
  \cS(x, z_1) - \cS(x, z_2) = \lim_{y \to -\infty} \cA[(y,\pi(y)) \to (z_1, 1)] - \cA[(y,\pi(y)) \to (z_2, 1)],     
  \end{equation}
where $\pi:(-\infty, z_1]\to\N$ is any non-increasing c\`adl\`ag function, such that $\lim_{y \to -\infty} \frac{\pi(y)}{2 y^2} = x$.
\end{itemize}
See \cite[Definition 1.2]{DOV} and \cite[Definition 2.20]{DZ}; also see \cite[Lemma 2.21]{DZ} for the second property.
Below we always assume the coupling between $\cS$ and $\cA$.

A \emph{parabolic path} or \emph{path} across $\cA$ from $x \ge 0$ to $z \in \R$ is a nonincreasing c\`adl\`ag function $\pi:(-\infty, z] \to \N$ such that 
\[
\lim_{y \to -\infty} \frac{\pi(y)}{2 y^2} = x.
\]
We then define the weight of $\pi$ on $\cA$ as
\begin{equation}  \label{eq:def-ale-pw}
\|\pi\|_\cA = \cS(x, z) + \lim_{y \to -\infty} \|\pi|_{[y, z]}\|_\cA - \cA[(y, \pi(y)) \to (z, 1)].    
\end{equation}
We note that we use the same notation $\|\cdot\|_\cA$ to denote weight across $\cA$ for both a non-increasing c\`adl\`ag function in a finite interval, and a parabolic path.
The quantity $\|\pi|_{[y, z]}\|_\cA - \cA[(y, \pi(y)) \to (z, 1)]$ can be understood as the ``loss of weight'' of $\pi|_{[y, z]}$, versus the optimal path from $(y, \pi(y))$ to $(z, 1)$; whereas the optimal parabolic path from $x$ to $z$ should have weight $\cS(x, z)$.
The limit exists since $\|\pi|_{[y, z]}\|_\cA - \cA[(y, \pi(y)) \to (z, 1)]$ is non-decreasing in $y$, by \cite[Lemma 4.1(i)]{DZ}.

We quote the following results on computing the weight of paths.
\begin{lemma}[\protect{\cite[Lemma 4.3]{DZ}}] \label{lem:ale-weight-cal}
The following statement holds almost surely for $\cA$.
Let $\pi_1,\pi_2$ be any two parabolic paths across $\cA$ from $x\ge 0$ to $z_1, z_2$ respectively, such that for some $z_0<z_1\wedge z_2$, we have $\pi_1(y)=\pi_2(y)$ for any $y\le z_0$.
Then
\[
\|\pi_1\|_\cA - \|\pi_1|_{[z_0,z_1]}\|_\cA
=
\|\pi_2\|_\cA - \|\pi_2|_{[z_0,z_2]}\|_\cA
\]
\end{lemma}
The next result states that, for a bounded region in $\R\times\N$, and any parabolic path without a ``jump point'' inside that region, the weight of the path is measurable with respect to $\cA$ outside that region.
\begin{lemma}[\protect{\cite[Lemma 4.4]{DZ}}] \label{lem:measurable-path-weight}
Take any closed interval $I\subset \R$ and $k\in\N$.
Let $\cF_{I,k}$ be the sigma-algebra generated by all null sets, all $\cA_i$ for $i>k$, and $\{\cA_i(x):x\not\in I\}$ for $1\le i \le k$.
Take any $x\ge 0$, and let $\Sigma_x$ be the set of parabolic paths $\pi$ from $x$ to some $z \in \R$ such that either $\pi(y)>k$ for any $y\in(-\infty,z]\cap I$, or else $I \subset (-\infty, z]$ and $\pi$ is constant on $I$. Let $F:\Sigma_x \to \R$ be the random function recording path weight in $\cA$: $F(\pi) = \|\pi\|_\cA$. Then $F$ is $\cF_{I,k}$-measurable.
\end{lemma}

A parabolic path $\pi$ from $x\ge 0$ to $y$ is a \emph{geodesic} if the weight $\|\pi\|_\cA$ is finite, and is maximal among all parabolic paths in $\cA$ from $x$ to $y$; and this is equivalent to $\|\pi\|_\cA = \cS(x,y)$, by \cite[Lemma 4.1(ii)]{DZ}.
It is proved in \cite[Lemma 4.2]{DZ} that, almost surely, there exists a geodesic from any $x\ge 0$ to any $y$; and for any fixed $x\ge 0$ and $y$, almost surely there is a unique geodesic from $x$ to $y$.
We shall use $\pi\{x,y\}$ to denote any geodesic from $x$ to $y$.

For paths across the Airy line ensemble, we use the following notion for disjointness.
For any paths $\pi_1:(-\infty, x_1]\to \N$ and $\pi_2:(-\infty, x_2]\to \N$, we say $\pi_1<\pi_2$, if $x_1<x_2$, and for any $x\le x_1$, we have $\lim_{y\nearrow x} \pi_1(y) < \pi_2(x)$ (we require this instead of $\pi_1(x) < \pi_2(x)$ because $\pi_1$ is defined to be c\`adl\`ag).
We say that $\pi_1$ and $\pi_2$ are \emph{disjoint}, if $\pi_1<\pi_2$ or $\pi_2<\pi_1$; and we say that $\pi_1$ \emph{intersects} $\pi_2$, if they are not disjoint.
Just like in the directed landscape (as stated in Lemma \ref{lem:dis-equiv}), the disjointness of geodesics is also equivalent to an inequality for the Airy sheet.
\begin{lemma}  \label{lem:dis-equiv-AS}
For any $0\le x_1< x_2$, $y_1<y_2$, assuming that $\pi\{x_1,y_1\}$ and $\pi\{x_2,y_2\}$ are unique, $\pi\{x_1,y_1\}$ and $\pi\{x_2,y_2\}$ being disjoint is equivalent to that
\begin{equation}  \label{eq:dis-equiv-AS}
\cS(x_1,y_1) + \cS(x_2,y_2) > \cS(x_1,y_2) + \cS(x_2,y_1).    
\end{equation}
\end{lemma}
\begin{proof}
We follow a similar idea as in the proof of Lemma \ref{lem:dis-equiv}.
Whenever $\pi\{x_1, y_1\}$ intersects $\pi\{x_2, y_2\}$, we can find $x_*\le y_1$ with \[\lim_{y\nearrow x_*}\pi\{x_1, y_1\}(y)\ge \pi\{x_2, y_2\}(x_*).\]
Again note that we use $\lim_{y\nearrow x_*}\pi\{x_1, y_1\}(y)$ instead of $\pi\{x_1, y_1\}(x_*)$ here, due to that $\pi\{x_1, y_1\}$ is c\`adl\`ag.
Using this we construct $\pi_1'$ from $x_1$ to $y_2$ and $\pi_2'$ from $x_2$ to $y_1$, by switching the paths $\pi\{x_1, y_1\}$ on $[x_*, y_1]$ and $\pi\{x_2, y_2\}$ on $[x_*, y_2]$.
Then by Lemma \ref{lem:ale-weight-cal} $$\|\pi_1'\|_\cA+\|\pi_2'\|_\cA=
\|\pi\{x_1, y_1\}\|_\cA+\|\pi\{x_2, y_2\}\|_\cA
=\cS(x_1,y_1) + \cS(x_2,y_2),$$ 
where the second equality is by the definition of $\pi\{x_1, y_1\}$ and $\pi\{x_2, y_2\}$.
Thus \eqref{eq:dis-equiv-AS} cannot hold.

Now we assume that \eqref{eq:dis-equiv-AS} does not hold; then by the quadrangle inequality (Lemma \ref{lem:DL-quad}), we must have that $\cS(x_1,y_1) + \cS(x_2,y_2) = \cS(x_1,y_2) + \cS(x_2,y_1)$.
We then consider any geodesics $\pi\{x_1, y_2\}$ and $\pi\{x_2, y_1\}$.
By planarity they should cross, i.e., we can find some $x_*\le y_1$ with \[\lim_{y\nearrow x_*}\pi\{x_1, y_2\}(y)\le \lim_{y\nearrow x_*}\pi\{x_1, y_2\}(y),\quad \pi\{x_1, y_2\}(x_*) \ge \pi\{x_1, y_2\}(x_*).\]
Then we construct $\pi_1'$ from $x_1$ to $y_1$ and $\pi_2'$ from $x_2$ to $y_2$, by switching the paths $\pi\{x_2, y_1\}$ on $[x_*, y_1]$ and $\pi\{x_1, y_2\}$ on $[x_*, y_2]$.
Now by Lemma \ref{lem:ale-weight-cal} we have 
\[\|\pi_1'\|_\cA+\|\pi_2'\|_\cA=\|\pi\{x_1, y_2\}\|_\cA+\|\pi\{x_1, y_2\}\|_\cA=\cS(x_1,y_2) + \cS(x_2,y_1)
=\cS(x_1,y_1) + \cS(x_2,y_2),\]
so $\pi_1'$ and $\pi_2'$ must be geodesics.
By uniqueness, we have $\pi_1'=\pi\{x_1, y_1\}$ and $\pi_2'=\pi\{x_2, y_2\}$, and they are not disjoint at $x_*$.
\end{proof}
Using this and Lemma \ref{lem:dis-equiv}, disjointness of geodesics in the directed landscape $\cL$ can be reduced to disjointness of geodesics in $\cA$.
\begin{cor}  \label{cor:dis-eq}
Couple $\cL$ with $\cS$ such that $\cS(x,y)=\cL(x,0;y,1)$ for any $x,y\in \R$.
For any fixed $0\le x_1<x_2$, $y_1<y_2$, the event that $\pi\{x_1,y_1\}$ and $\pi\{x_2,y_2\}$ are disjoint and $\Dist^{[0,1]}_{(x_1,x_2)\to(y_1,y_2)}$ are equivalent (up to an event with probability zero). (Recall that $\Dist^{[0,1]}_{(x_1,x_2)\to(y_1,y_2)}$ is the event where $\pi_{(x_1,0;y_1,1))}$ and $\pi_{(x_2,0;y_2,1)}$ are disjoint.)
\end{cor}

\subsection{Disjoint optimizers, ordering, and second geodesics}  \label{ssec:ale-2nd}

We also need the notion of disjoint optimizers from \cite{DZ}.
For any $k\in\N$, we let
\[
\R^k_\le = \{(x_1,\ldots,x_k) : x_1\le \cdots \le x_k\}.
\]
For any $\bx=(x_1,\ldots,x_k)$, $\by=(y_1,\ldots,y_k)$ in $\R^k_\le$ with $0\le x_1$, a \emph{disjoint $k$-path} from $\bx$ to $\by$ is a collection of parabolic paths $(\pi_1, \ldots, \pi_k)$, such that each $\pi_i$ is from $x_i$ to $y_i$, and $\pi_i<\pi_{i+1}$ on $(-\infty, x_i)$, $\pi_i(x_i)\le \pi_{i+1}(x_i)$, for each $1\le i <k$.
Note (due to that the paths are c\`adl\`ag) this is a weaker condition than $\pi_i$ and $\pi_{i+1}$ being disjoint.
A \emph{disjoint optimizer} from $\bx$ to $\by$ is a disjoint $k$-path $(\pi_1, \ldots, \pi_k)$, such that $\sum_{i=1}^k\|\pi_i\|_\cA$ achieves its supremum (over all disjoint $k$-paths from $\bx$ to $\by$).
It has been shown, see \cite[Proposition 5.8]{DZ}, that almost surely for any $\bx=(x_1,\ldots,x_k)$ and $\by=(y_1,\ldots,y_k)$ in $\R^k_\le$ with $0\le x_1$, there is a disjoint optimizer from $\bx$ to $\by$.

For $x\ge 0$ and $y\in \R$, we consider all the parabolic paths $\pi:(-\infty, y] \to \N$ across the Airy line ensemble from $x$, with $\lim_{z\nearrow y}\pi(z)\ge 2$.
Note that at the endpoint we allow $\pi(y)=1$.
We denote any one with the maximum weight as a \emph{second geodesic from $x$ to $y$}.
\begin{lemma}  \label{lem:2nd-geo-exist}
Almost surely, for any $x\ge 0$ and $y$, there is a second geodesic from $x$ to $y$. In addition, if we take $\pi:(-\infty,y]\to \N$ such that $\pi \equiv 1$, then any second geodesic from $x$ to $y$ together with $\pi$ would give a disjoint optimizer from $(0, x)$ to $(y, y)$. 
\end{lemma}
\begin{proof}
For the first part, we just take any disjoint optimizer $(\pi_1, \pi_2)$ from $(0,x)$ to $(y,y)$, and we claim that $\pi_2$ is a second geodesic from $x$ to $y$.
Otherwise, there is a path $\pi_2'$ from $x$ to $y$, with $\lim_{z\nearrow y}\pi_2'(z)\ge 2$, and $\|\pi_2'\|_\cA > \|\pi_2\|_\cA$.
For $\pi:(-\infty, y]\to \N$ with $\pi \equiv 1$, we have $\|\pi\|_\cA = \cA_1(y)=\cS(0,y)\ge \|\pi_1\|_\cA$ (where the inequality is by \eqref{eq:def-ale-pw}), so $\|\pi\|_\cA + \|\pi_2'\|_\cA > \|\pi_1\|_\cA + \|\pi_2\|_\cA$, and this contradicts that $(\pi_1, \pi_2)$ is a disjoint optimizer.

We next prove the second part. Take any second geodesic $\pi_*$ from $x$ to $y$, then $\|\pi_*\|_\cA = \|\pi_2\|_\cA$, implying that $\|\pi\|_\cA + \|\pi_*\|_\cA \ge \|\pi_1\|_\cA + \|\pi_2\|_\cA$. Thus $(\pi, \pi_*)$ must also be a disjoint optimizer from $(0, x)$ to $(y, y)$. 
\end{proof}
From any $x\ge 0,$ to any $y\in \R$, there maybe more than one second geodesic, and below we use $\pi^2\{x,y\}:(-\infty, y] \to \N$ to denote an arbitrary one.

Next we state the following result, which is about the ordering of disjoint optimizers. It can be viewed as a counterpart of Lemma \ref{lem:ordering-geo} or Lemma \ref{lem:semi-inf-tree} in the disjoint optimizer setting.
For any $x_1< x_2$ and $\pi_1:(-\infty, x_1]\to \N$ and $\pi_2:(-\infty, x_2]\to \N$, we define $\pi_1\vee \pi_2:(-\infty, x_2]\to \N$ and $\pi_1\wedge \pi_2:(-\infty, x_1]\to \N$ via
\[
(\pi_1\vee \pi_2)(x) = 
\begin{cases}
\pi_1(x) \vee \pi_2(x),\quad & \forall x \in (-\infty, x_1],\\
\pi_2(x), \quad & \forall x \in (x_1, x_2],
\end{cases}
\]
and
\[
(\pi_1\wedge \pi_2)(x) = \pi_1(x) \wedge \pi_2(x),\quad \forall x \in (-\infty, x_1].
\]
\begin{lemma}  \label{lem:ALE-opt-min-max}
Take any $\bx, \by, \bx', \by' \in \R^k_\le$ with $x_1, x_1'\ge 0$, and $\bx \le \bx', \by \le \by'$ entry-wise. Take any disjoint optimizers $(\pi_1,\ldots, \pi_k)$ from $\bx$ to $\by$, and $(\pi_1',\ldots, \pi_k')$ from $\bx'$ to $\by'$.
Then $(\pi_1\vee \pi_1',\ldots, \pi_k\vee \pi_k')$ is a disjoint optimizer from $\bx'$ to $\by'$, and $(\pi_1\wedge \pi_1',\ldots, \pi_k\wedge \pi_k')$ is a disjoint optimizer from $\bx$ to $\by$.
\end{lemma}
For paths in finite intervals (rather than parabolic paths), such result follows obviously from the following facts: the $\wedge$ and $\vee$ of two disjoint optimizers are still disjoint optimizers (with the respective endpoints), and the sum of their weights is the same as the sum of weights of the disjoint optimizers.
For parabolic paths, there are extra complexities due to infiniteness of the paths.
\begin{proof}
It is straightforward to check that for each $1\le i \le k$, $\pi_i\vee \pi_i'$ is a path from $x_i'$ to $y_i'$, and $\pi_i\wedge \pi_i'$ is a path from $x_i$ to $y_i$.
For $1\le i <k$, we have $\pi_i\wedge \pi_i' < \pi_{i+1}\wedge \pi_{i+1}'$ on $(-\infty, y_i)$, $\pi_i(y_i)\wedge \pi_i'(y_i) \le \pi_{i+1}(y_i)\wedge \pi_{i+1}'(y_i)$; and $\pi_i\vee \pi_i' < \pi_{i+1}\vee \pi_{i+1}'$ on $(-\infty, y_i')$, $\pi_i(y_i')\vee \pi_i'(y_i') \le \pi_{i+1}(y_i')\vee \pi_{i+1}'(y_i')$.
Thus $(\pi_1\vee \pi_1',\ldots, \pi_k\vee \pi_k')$ is a disjoint $k$-path from $\bx'$ to $\by'$, and $(\pi_1\wedge \pi_1',\ldots, \pi_k\wedge \pi_k')$ is a disjoint $k$-path from $\bx$ to $\by$.

We claim that for each $1\le i \le k$,
\begin{equation} \label{eq:ale-switch}
\|\pi_i\vee \pi_i'\|_\cA + \|\pi_i\wedge \pi_i'\|_\cA = \|\pi_i\|_\cA + \|\pi_i'\|_\cA.
\end{equation}
For this, note that for any $y\le y_i$, we have
\[
\|(\pi_i\vee \pi_i')|_{[y, y_i']}\|_\cA + \|(\pi_i\wedge \pi_i')|_{[y, y_i]}\|_\cA -\|\pi_i|_{[y, y_i]}\|_\cA - \|\pi_i'|_{[y, y_i']}\|_\cA=0,
\]
as the involved paths are finite.
By \eqref{eq:def-ale-pw}, to prove \eqref{eq:ale-switch} it now suffices to show that
\begin{multline}  \label{eq:ale-pass-quad}
\lim_{y\to-\infty}
\cA[(y,\pi_i(y)\wedge \pi_i'(y))\to (y_i,1)] + \cA[(y,\pi_i(y)\vee \pi_i'(y))\to (y_i',1)]
\\ - \cA[(y,\pi_i(y))\to (y_i,1)] - \cA[(y,\pi_i'(y))\to (y_i',1)] = 0.    
\end{multline}
If there is some $y_*$ such that $\pi_i(y)\le \pi_i'(y)$ for any $y\le y_*$, \eqref{eq:ale-pass-quad} holds obviously.
Otherwise, for infinitely many $y$ decreasing to $-\infty$, we have $\pi_i(y)> \pi_i'(y)$.
Then since $\lim_{y\to-\infty}\frac{\pi_i(y)}{2y^2}=x_i$, $\lim_{y\to-\infty}\frac{\pi_i'(y)}{2y^2}=x_i'$, and (the assumption that) $x_i\le x_i'$, we must have $x_i=x_i'$.
Thus by \eqref{eq:ale-diff-asy} we have
\begin{multline*}
\lim_{y\to-\infty}
\cA[(y,\pi_i(y))\to (y_i,1)] - \cA[(y,\pi_i(y))\to (y_i',1)]
\\ = \lim_{y\to-\infty} \cA[(y,\pi_i'(y))\to (y_i,1)] - \cA[(y,\pi_i'(y))\to (y_i',1)].
\end{multline*}
Note that (before sending $y\to-\infty$) the expression in \eqref{eq:ale-pass-quad} equals either $0$ or
\begin{multline*}
\cA[(y,\pi_i'(y))\to (y_i,1)] + \cA[(y,\pi_i(y))\to (y_i',1)]
\\ - \cA[(y,\pi_i(y))\to (y_i,1)] - \cA[(y,\pi_i'(y))\to (y_i',1)]. 
\end{multline*}
Thus \eqref{eq:ale-pass-quad} holds, and so does \eqref{eq:ale-switch}. So we have $\sum_{i=1}^k\|\pi_i\vee \pi_i'\|_\cA + \|\pi_i\wedge \pi_i'\|_\cA = \sum_{i=1}^k\|\pi_i\|_\cA + \|\pi_i'\|_\cA$, and the conclusion follows.
\end{proof}

\subsection{Proof of Lemma \ref{lem:sur-ALE}}  \label{ssec:sur-ALE}
We first translate the event $\cE_{coal}'\cap \cE_{disj}$ into an event on disjointness of geodesics in $\cA$.
Recall that 
\[
\cE_{disj}=\Dist_{(-1, 1)\to (-\delta, \delta)}^{[0, 1]} \cap (\Dist_{(-1-\delta, -1+\delta)\to (-2, 2)}^{[0, 1]})^c \cap (\Dist_{(1-\delta, 1+\delta)\to (-2, 2)}^{[0, 1]})^c,
\]
i.e., the geodesics $\pi_{(-1,0;-\delta,1)}$ and $\pi_{(1,0;\delta,1)}$ are disjoint, while, the geodesic pairs given by $$\{\pi_{(-1-\delta,0;-2,1)},\pi_{(-1+\delta,0;2,1)}\}, \{\pi_{(1-\delta,0;-2,1)},\pi_{(1+\delta,0;2,1)}\},$$ coalesce.
For easier analysis and to use Corollary \ref{cor:dis-eq}, we reflect time and then shift the starting points, {since the corollary only allows positive starting points, and the skew-shift invariance property of the directed landscape implies that the disjointness of the geodesics stay invariant as well under such a shift}. More precisely, by Lemma \ref{lem:DL-symmetry}, $\cE_{coal}'\cap \cE_{disj}$ has the same probability as 
\[
\Dist_{(2-\delta, 2+\delta)\to (-1,1)}^{[0, 1]} \cap (\Dist_{(0,4)\to(-1-\delta, -1+\delta)}^{[0, 1]})^c \cap (\Dist_{(0,4)\to(1-\delta, 1+\delta)}^{[0, 1]})^c.
\]
Then by Corollary \ref{cor:dis-eq}, we have $\PP[\cE_{coal}'\cap \cE_{disj}]=\PP[\cE_{ALE}]$, where $\cE_{ALE}$ is the following event:
\begin{enumerate}
    \item $\pi\{0,-1-\delta\}$ intersects $\pi\{4,-1+\delta\}$, and $\pi\{0,1-\delta\}$ intersects $\pi\{4,1+\delta\}$.
    \item $\pi\{2-\delta,-1\}$ and $\pi\{2+\delta,1\}$ are disjoint.
\end{enumerate}
It remains to show that $\cE_{ALE}$ happens with positive probability. 

Below we assume that the geodesic $\pi\{x,y\}$ for any rational $x\ge 0$ and $y$ is unique. This event has probability one by \cite[Lemma 4.2]{DZ}.
An implication of this assumption is that, for any rational $y$, $\pi\{0,y\}$ would be the function staying at $1$ on $(-\infty, y]$. Recall (from the proof of Lemma \ref{lem:K-pos-posprob}) that $\delta>0$ is taken to be an arbitrary small number, and below we assume that it is rational.\\

\noindent\textbf{Step 1: reduce to events on single geodesics.} 
We note that $\pi\{0,-1-\delta\}$ and $\pi\{0,1-\delta\}$ would be constant $1$ on the whole domains $(-\infty, -1-\delta]$ and $(-\infty, 1-\delta]$ respectively, since these geodesics are assumed to be unique.
Then the events $\pi\{0,-1-\delta\}$ intersects $\pi\{4,-1+\delta\}$, and $\pi\{0,1-\delta\}$ intersects $\pi\{4,1+\delta\}$, would be equivalent to $\pi\{4,-1+\delta\}(-1-\delta)=1$ and $\pi\{4,1+\delta\}(1-\delta)=1$; and we denote them as $\cE_{ALE,1}$ and $\cE_{ALE,2}$ respectively.

For the non-coalescence event, we consider any second geodesics $\pi^2\{2-\delta,-1\}$ and $\pi^2\{2+\delta,0\}$.
Then by the asymptotic directions (of $\pi^2\{2-\delta,-1\}$ and $\pi^2\{2+\delta,0\}(y)$), there exists rational $y_0<-1$, such that $\pi^2\{2-\delta,-1\}(y)<\pi^2\{2+\delta,0\}(y)$ for any $y\le y_0$.
Further we claim that we can take $y_0$ small enough, such that 
\begin{equation}  \label{eq:ale4-1}
\|\pi^2\{2+\delta, 0\}\|_\cA - \cA_2(0) > \|\pi^2\{2+\delta, y_0\}\|_\cA - \cA_2(y_0),
\end{equation}
which will ensure that $\pi^2\{2+\delta, 0\}$ stays below the second line up to $y_0$, which will allow us to construct events ensuring disjointness of  $\pi\{2-\delta,-1\}$ and $\pi\{2+\delta,1\}.$ Before describing this, we first show how to establish \eqref{eq:ale4-1} using the following lemma which is essentially the same as  the quadrangle inequality.

\begin{lemma} \label{lem:y0-exist}
Almost surely the following is true.
The function $y\mapsto \|\pi^2\{2+\delta, y\}\|_\cA - \cA_2(y)$ is non-decreasing, and $\limsup_{y\to-\infty}\|\pi^2\{2+\delta, y\}\|_\cA - \cA_2(y) < \|\pi^2\{2+\delta, 0\}\|_\cA - \cA_2(0)$.
\end{lemma}
\begin{proof}
For any $y_1<y_2$, and any second geodesic $\pi^2\{2+\delta, y_1\}$, we can construct a path $\pi$ from $2+\delta$ to $y_2$, which is the same as $\pi^2\{2+\delta, y_1\}$ on $(-\infty, y_1)$, and equals $1$ on $[y_1, y_2]$.
We then have that $\|\pi\|_\cA=\|\pi^2\{2+\delta, y_1\}\|_\cA - \cA_2(y_1) + \cA_2(y_2)$. We also have that  $\|\pi\|_\cA\le \|\pi^2\{2+\delta, y_2\}\|_\cA$, so we have $\|\pi^2\{2+\delta, y_1\}\|_\cA - \cA_2(y_1) \le \|\pi^2\{2+\delta, y_2\}\|_\cA - \cA_2(y_2)$.

For the second part, we assume the contrary; then we have that $\|\pi^2\{2+\delta, y\}\|_\cA - \cA_2(y) = \|\pi^2\{2+\delta, 0\}\|_\cA - \cA_2(0)$ for any $y\le 0$.
By Lemma \ref{lem:ale-weight-cal} we have $\|\pi\{2+\delta,0\}\|_\cA=\sup_{y\le 0}\|\pi^2\{2+\delta,y\}\|_\cA + \cA[(y,2)\to (0,1)]$, which further means
\begin{equation}  \label{eq:pi2del}
\begin{split}
\|\pi\{2+\delta,0\}\|_\cA &= \sup_{y\le y' \le 0}\|\pi^2\{2+\delta,y\}\|_\cA - \cA_2(y) + \cA_2(y') - \cA_1(y') + \cA_1(0) \\
& = \sup_{y' \le 0}\|\pi^2\{2+\delta,0\}\|_\cA - \cA_2(0) + \cA_2(y') - \cA_1(y') + \cA_1(0) \\
& = \|\pi^2\{2+\delta,0\}\|_\cA - \cA_2(0) + \cA_1(0),
\end{split}
\end{equation}
where the last equality is due to the following reason.
By the Brownian Gibbs property of $\cA$, we have $\PP[\cA_1(0)-\cA_2(0)\le \delta]>0$ for any $\delta>0$. 
Then by stationarity and ergodicity of $\cA$ plus a parabola, we have that $\PP[\liminf_{y\to\infty} \cA_1(y)-\cA_2(y) \le \delta]=1$, so almost surely $\liminf_{y\to\infty} \cA_1(y)-\cA_2(y) = 0$, and $\sup_{y' \le 0}\cA_2(y') - \cA_1(y')=0$.

On the other hand, by the existence of the geodesic $\pi\{2+\delta, 0\}$, there exists some $y_*\le 0$ with $\|\pi\{2+\delta,0\}\|_\cA=\|\pi^2\{2+\delta,y_*\}\|_\cA + \cA_1(0)-\cA_1(y_*)$, by Lemma \ref{lem:ale-weight-cal}. Thus with \eqref{eq:pi2del} we have
\[
 \|\pi^2\{2+\delta,0\}\|_\cA - \cA_2(0) = \|\pi^2\{2+\delta,y_*\}\|_\cA + -\cA_1(y_*).
\]
Since we have assumed that $\|\pi^2\{2+\delta,0\}\|_\cA - \cA_2(0) = \|\pi^2\{2+\delta,y_*\}\|_\cA - \cA_2(y_*)$, we have
$\cA_1(y_*)=\cA_2(y_*)$.
However, by the Brownian Gibbs property of $\cA$, we have that $\PP[\cA_1(y)=\cA_2(y),\;\exists y\le 0] = 0$. Thus we conclude that $\PP[\|\pi^2\{2+\delta, y\}\|_\cA - \cA_2(y) = \|\pi^2\{2+\delta, 0\}\|_\cA - \cA_2(0),\; \forall y\le 0] = 0$.
\end{proof}

\begin{figure}[hbt!]
    \centering
\begin{tikzpicture}[line cap=round,line join=round,>=triangle 45,x=1.2cm,y=0.3cm]
\clip(-0.5,-2) rectangle (10.5,10.5);

\draw [line width=2pt, color=orange] (6,10) -- (2.5,10) -- (2.5,8) -- (1.7,8) -- (1.7,6) -- (0.5,6) -- (0.5,4) -- (0,4);
\draw [line width=2pt, color=orange] (8,10) -- (7.7,10) -- (7.7,8) -- (7.4,8) -- (7.4,6) -- (6.6,6) -- (6.6,4) -- (4.6,4) -- (4.6,2) -- (3.2,2) -- (3.2,0) -- (2.5,0);

\draw [line width=1pt, color=blue] (6,8) -- (4.6,8) -- (4.6,6) -- (3.7,6) -- (3.7,4) -- (0.8,4) -- (0.8,2) -- (0,2);

\draw [line width=1pt, color=blue] (7,8) -- (6.2,8) -- (6.2,6) -- (3.7,6) -- (3.7,4) -- (3.4,4) -- (3.4,2) -- (3.2,2) -- (3.2,0) -- (2.5,0);

\draw (0,25) -- (10,25);

\foreach \i in {1,...,6}
{
\draw [dotted] [thick] (0, 12-2*\i) -- (10, 12-2*\i);
\begin{tiny}
\draw (0,12-2*\i) node[anchor=east]{$\cA_{\i}$};
\end{tiny}
}

\draw [dashed] [ultra thin] (3, -1) -- (3, 11);
\draw [dashed] [ultra thin] (6, -1) -- (6, 11);
\draw [dashed] [ultra thin] (7, -1) -- (7, 11);
\draw [dashed] [ultra thin] (8, -1) -- (8, 11);

\begin{tiny}
\draw (3,-1) node[anchor=north]{$y_0$};
\draw (6,-1) node[anchor=north]{$-1$};
\draw (7,-1) node[anchor=north]{$0$};
\draw (8,-1) node[anchor=north]{$1$};

\draw (2.5,9) node[anchor=east]{$\pi\{2-\delta, -1\}$};
\draw (4.6,7) node[anchor=east]{$\pi^2\{2-\delta, -1\}$};
\draw (5.3,6) node[anchor=north]{$\pi^2\{2+\delta, 0\}$};
\draw (5.3,4) node[anchor=north]{$\pi\{2+\delta, 1\}$};
\end{tiny}

\end{tikzpicture}
\caption{
An illustration of the proof of Lemma \ref{lem:geo-sand}: assuming $\pi\{2-\delta,-1\}(y_0)=1$, $\pi\{2-\delta,-1\}$ is to the left of $\pi^2\{2-\delta,-1\}$; and assuming $\pi\{2+\delta,1\}(0)>1$, $\pi\{2+\delta,1\}$ is to the right of $\pi^2\{2+\delta,0\}$.
By the choice of $y_0$, $\pi^2\{2-\delta,-1\}$ and $\pi^2\{2+\delta,0\}$ are disjoint in $(-\infty, y_0]$; thus by the ordering of geodesics, $\pi\{2-\delta,-1\}$ and $\pi\{2+\delta,1\}$ are disjoint.
}
\label{fig:B9}
\end{figure}
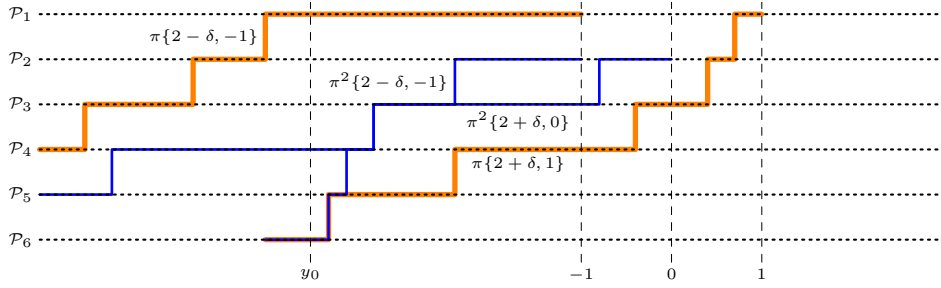

We now construct an event which ensures disjointness of geodesics.

\begin{lemma}  \label{lem:geo-sand}
$\pi\{2-\delta,-1\}(y_0)=1$ and $\pi\{2+\delta,1\}(0)>1$ imply that $\pi\{2-\delta,-1\}$ and $\pi\{2+\delta,1\}$ are disjoint.
\end{lemma}
\begin{proof}
The proof is a consequence of ordering properties of geodesics and second geodesics. In particular, we show that (in $(-\infty, -1]$)
$$\pi\{2-\delta,-1\}\le \pi^2\{2-\delta,-1\}\le  \pi^2\{2+\delta,0\}\le \pi\{2+\delta,1\}$$ and then establish certain strictness of the inequalities (see Figure \ref{fig:B9}). While the inequalities are just a consequence of planarity and uniqueness of the relevant objects, the formal arguments are a bit tedious.

\begin{enumerate}
    \item  
    Take some rational $y_1<-1$ with $\pi\{2-\delta,-1\}(y_1)>1$. Then $(\pi\{0,y_1\},\pi\{2-\delta,-1\})$ would be a disjoint optimizer from $(0,2-\delta)$ to $(y_1,-1)\}$.
    We also know that $(\pi\{0,-1\}, \pi^2\{2-\delta,-1\})$ is a disjoint optimizer from $(0,2-\delta)$ to $(-1,-1)$, by Lemma \ref{lem:2nd-geo-exist}.
    Then by Lemma \ref{lem:ALE-opt-min-max}, $(\pi\{0,y_1\},\pi\{2-\delta,-1\}\wedge \pi^2\{2-\delta,-1\})$ would be another disjoint optimizer from $(0,2-\delta)$ to $(y_1,-1)$.
    This implies that $\|\pi\{2-\delta,-1\}\wedge \pi^2\{2-\delta,-1\}\|_\cA \ge \|\pi\{2-\delta,-1\}\|_\cA = \cS(2-\delta,-1)$, and $\pi\{2-\delta,-1\}\wedge \pi^2\{2-\delta,-1\}$ is also a geodesic from $2-\delta$ to $-1$, 
    So $\pi\{2-\delta,-1\} = \pi\{2-\delta,-1\}\wedge \pi^2\{2-\delta,-1\}$ by the uniqueness.
    Thus we have that $\pi\{2-\delta,-1\}(y)\le \pi^2\{2-\delta,-1\}(y)$ for any $y\in (-\infty, -1]$. 
    \item When $\pi\{2+\delta,1\}(0)>1$, $(\pi\{0,0\},\pi\{2+\delta,1\})$ is a disjoint optimizer from $(0,2+\delta)$ to $(0,1)$. 
    By Lemma \ref{lem:2nd-geo-exist}, $(\pi\{0,0\}, \pi^2\{2+\delta,0\})$ is a disjoint optimizer from $(0,2+\delta)$ to $(0,0)$.
    Then by Lemma \ref{lem:ALE-opt-min-max}, $(\pi\{0,0\},\pi\{2+\delta,1\}\vee \pi^2\{2+\delta,0\})$ is another disjoint optimizer from $(0,2+\delta)$ to $(0,1)$, so $\pi\{2+\delta,1\}\vee \pi^2\{2+\delta,0\}$ would also be a geodesic from $2+\delta$ to $1$.
    So $\pi\{2+\delta,1\} = \pi\{2+\delta,1\}\vee \pi^2\{2+\delta,0\}$ by the uniqueness,
    and $\pi\{2+\delta,1\}(y)\ge\pi^2\{2+\delta,0\}(y)$ for any $y\in(-\infty, 0]$. 
\end{enumerate}
Since $\pi^2\{2-\delta,-1\}$ and $\pi^2\{2+\delta,0\}$ are disjoint on $(-\infty, y_0]$ by definition of $y_0$, we have that $\pi\{2-\delta,-1\}$ and $\pi\{2+\delta,1\}$ are disjoint on $(-\infty, y_0]$. On $[y_0, -1]$ we have $\pi\{2-\delta,-1\}=1$ and $\pi\{2+\delta,1\}>1$ by hypothesis, and so $\pi\{2-\delta,-1\}$ and $\pi\{2+\delta,1\}$ are also disjoint.
\end{proof}
To use the above, we let $\cE_{ALE,3}$ be the event $\pi\{2-\delta,-1\}(y_0)=1$ and $\cE_{ALE,4}$ be the event $\pi\{2+\delta,1\}(0)>1$.

From the above lemma we have that $$\cE_{ALE,1}\cap \cE_{ALE,2} \cap \cE_{ALE,3} \cap \cE_{ALE,4}\subset \cE_{ALE}=\cE_{coal}'\cap \cE_{disj}.$$

We will next construct events $\cE_{ALE,1}'$, $\cE_{ALE,2}'$, $\cE_{ALE,3}'$, and $\cE_{ALE,4}'$, which imply $\cE_{ALE,1}$, $\cE_{ALE,2}$, $\cE_{ALE,3}$, $\cE_{ALE,4}$, respectively.
These events are given in terms of inequalities involving $\cA_1$, $\cA_2$, and weights of certain second geodesics.
We will show that, conditional on $\{\cA_n\}_{n=2}^\infty$, and $\cA_1$ outside some compact interval, the intersection of these events has positive probability.
This is by considering $\cA_1$ inside that compact interval and using the Brownian Gibbs property of $\cA$.

\noindent\textbf{Step 2: reduce to events on $\cA_1$ and weights of second geodesics.} 
Let $\cF$ be the sigma-algebra generated by $\{\cA_n\}_{n=2}^\infty$, and $\cA_1(x)$ for $x\le y_0-1$ and $x\ge 2$.
By Lemma \ref{lem:measurable-path-weight}, for any rational $x\ge 0$, the random variable $\|\pi^2\{x,y\}\|_\cA$ for any $y$ is $\{\cA_n\}_{n=2}^\infty$ measurable; and $y_0$ is $\{\cA_n\}_{n=2}^\infty$ measurable.
Then these random variables are all $\cF$ measurable. By the Brownian Gibbs property of the Airy line ensemble \cite{CH}, conditional on $\cF$, the process $\cA_1$ on $[y_0-1, 2]$ is a Brownian bridge conditional on avoiding $\cA_2$.

It would be useful to remember the following simple relation between geodesics and second geodesics.
By Lemma \ref{lem:ale-weight-cal}, we have that for any rational $x\ge 0$ and $y$, 
\begin{equation} \label{eq:rel-geo-2geo}
\|\pi\{x,y\}\|_\cA = \sup_{z\le y} \|\pi^2\{x,z\}\|_\cA + \cA_1(y) - \cA_1(z),
\end{equation}
and by uniqueness of $\pi\{x,y\}$ this supremum is achieved at a unique $z_*\le y$, such that $\pi\{x,y\}(z)=1$ for $z\in [z_*,y]$, and $\pi\{x,y\}(z)>1$ for $z\in (-\infty,z_*)$.

\begin{figure}[hbt!]
    \centering
\begin{tikzpicture}[line cap=round,line join=round,>=triangle 45,x=.55cm,y=.7cm]
\clip(-3,-2) rectangle (28,8.5);

\draw [line width=1pt,color=darkgreen] (-3,5) -- (-2,3.3) -- (-1,4.2) -- (0,3.1) -- (0.5,4.1) -- (1,3.3);

\draw [line width=1pt,color=darkgreen] (28,5.4) -- (27,3.8) -- (26,5.6) -- (25,3.4) -- (24.5,4.5) -- (24,4);

\draw [line width=1pt] (1,3.3) -- (1.25,2.5) -- (1.5,4) -- (1.75,3) -- (2,2.8) -- (2.25,3.3) -- (2.5,3.2) -- (2.75,3.4) -- (3,4.3) -- (3.25,5.1) -- (3.5,4.6) -- (3.75,4.9) -- (4,3.8) -- (4.25,4.1) -- (4.5,2.9) -- (4.75,3.5) -- (5,5) -- (5.25,4.6) -- (5.5,4.8) -- (5.75,5.1) -- (6,6.2) -- (6.25,5.3) -- (6.5,6.8) -- (6.75,6.6) -- (7,7.1) -- (7.25,7.4) -- (7.5,7) -- (7.75,6.1) -- (8,7.4) -- (8.25,6.1) -- (8.5,8.2) -- (8.75,7.7) -- (9,8.3);

\draw [line width=1pt] (24,4) -- (23.75,4.1) -- (23.5,3.1) -- (23.25,4.5) -- (23,4.8) -- (22.75,3.3) -- (22.5,3.2) -- (22.25,3.4) -- (22,4.3) -- (21.75,5.1) -- (21.5,5.6) -- (21.25,4.2) -- (21,4.8) -- (20.75,4.1) -- (20.5,2.9) -- (20.25,4.9) -- (20,3.8) -- (19.75,4.6) -- (19.5,6.8) -- (19.25,7.1) -- (19,8.2) -- (18.75,7.6) -- (18.5,7.7) -- (18.25,4.6) -- (18,5.1) -- (17.75,5.4) -- (17.5,4) -- (17.25,5.1) -- (17,4.4) -- (16.75,4.6) -- (16.5,3.8) -- (16.25,4.) -- (16,6.2);

\draw [line width=1pt] (9,8.3) -- (9.5,8.2) -- (10,7.1) -- (10.5,6.7) -- (11,7.5) -- (11.5,6.2) -- (12,7.8) -- (12.5,6.3) -- (13,6.5) -- (13.5,7.1) -- (14,6.3) -- (14.5,5.8) -- (15,7.1) -- (15.5,5.7) -- (16,6.2);

\draw [line width=1pt,color=darkgreen] (-3,1.5) -- (-2,0.3) -- (-1,1.8) -- (0,0.7) -- (0.5,1.1) -- (1,0.3) -- (1.25,1.5) -- (1.5,1.2) -- (1.75,1.4) -- (2,0.8) -- (2.25,1.3) -- (2.5,1.2) -- (2.75,1.4) -- (3,2.3) -- (3.25,2.1) -- (3.5,3.1) -- (3.75,2.9) -- (4,1.8) -- (4.25,2.1) -- (4.5,0.9) -- (4.75,1.5) -- (5,2) -- (5.25,2.6) -- (5.5,0.8) -- (5.75,1.1) -- (6,2.2) -- (6.25,1.3) -- (6.5,2.8) -- (6.75,2.6) -- (7,3.1) -- (7.25,3.4) -- (7.5,2) -- (7.75,2.1) -- (8,2.4) -- (8.25,3.1) -- (8.5,3.6) -- (8.75,3.) -- (9,2.5);

\draw [line width=1pt,color=darkgreen] (28,1.4) -- (27,0.8) -- (26,2.2) -- (25,1.4) -- (24.5,2.5) -- (24,2) -- (23.75,2.1) -- (23.5,1.1) -- (23.25,2.5) -- (23,2.8) -- (22.75,2.3) -- (22.5,3.2) -- (22.25,1.4) -- (22,2.3) -- (21.75,2.1) -- (21.5,1.6) -- (21.25,2.2) -- (21,2.8) -- (20.75,2.1) -- (20.5,1.9) -- (20.25,2.9) -- (20,1.8) -- (19.75,2.6) -- (19.5,1.8) -- (19.25,2.1) -- (19,2.2) -- (18.75,1) -- (18.5,3.1) -- (18.25,1.6) -- (18,2.1) -- (17.75,1.4) -- (17.5,2) -- (17.25,1.1) -- (17,2.8) -- (16.75,2.4) -- (16.5,3.6) -- (16.25,2.) -- (16,2.2);

\draw [line width=1pt,color=darkgreen] (9,2.5) -- (9.5,3.5) -- (10,2.1) -- (10.5,3) -- (11,2.5) -- (11.5,3.2) -- (12,1.8) -- (12.5,2.3) -- (13,2.5) -- (13.5,3.1) -- (14,2.3) -- (14.5,3.8) -- (15,2.1) -- (15.5,3.7) -- (16,2.2);

\draw [line width=0.5pt,color=blue] (6,2.7) -- (6.25,1.8) -- (6.5,3.3) -- (6.75,3.1) -- (7,3.6) -- (7.25,3.9) -- (7.5,2.5) -- (7.75,2.6) -- (8,2.9) -- (8.25,3.6) -- (8.5,4.1) -- (8.75,3.5) -- (9,3) -- (9.5,4) -- (10,2.6) -- (10.5,3.5) -- (11,3) -- (11.5,3.7) -- (12,2.3) -- (12.5,2.8) -- (13,3) -- (13.5,3.6) -- (14,2.8);

\draw [line width=0.5pt,color=blue] (6,4.2) -- (6.25,3.3) -- (6.5,4.8) -- (6.75,4.6) -- (7,4.1) -- (7.25,5.4) -- (7.5,4) -- (7.75,4.1) -- (8,4.4) -- (8.25,5.6) -- (8.5,6.1) -- (8.75,5.5) -- (9,5);

\draw [line width=0.5pt,color=blue] (8.5,7.2) -- (8.75,6.7) -- (9,6.4) -- (9.5,7.5);
\draw [line width=0.5pt,color=blue] (18.5,6.3) -- (18.75,5.1) -- (19,5.6) -- (19.25,5.9) -- (19.5,5.7);

\foreach \i in {1,3.5,6,8,8.5,9,9.5,14,16.5,18,18.5,19,19.5,24}
{
\draw [dotted] (\i,0) -- (\i,10);
}

\begin{tiny}
\draw (1,0) node[anchor=north east,rotate=60]{$y_0-1$};
\draw (3.5,0) node[anchor=north east,rotate=60]{$y_0-1/2$};
\draw (6,0) node[anchor=north east,rotate=60]{$y_0$};

\draw (8,0) node[anchor=north east,rotate=60]{$-1-2\delta$};
\draw (8.5,0) node[anchor=north east,rotate=60]{$-1-\delta$};
\draw (9,0) node[anchor=north east,rotate=60]{$-1$};
\draw (9.5,0) node[anchor=north east,rotate=60]{$-1+\delta$};
\draw (14,0) node[anchor=north east,rotate=60]{$0$};
\draw (16.5,0) node[anchor=north east,rotate=60]{$1/2$};
\draw (18,0) node[anchor=north east,rotate=60]{$1-2\delta$};
\draw (18.5,0) node[anchor=north east,rotate=60]{$1-\delta$};
\draw (19,0) node[anchor=north east,rotate=60]{$1$};
\draw (19.5,0) node[anchor=north east,rotate=60]{$1+\delta$};
\draw (24,0) node[anchor=north east,rotate=60]{$2$};

\draw (3.5,5) node[anchor=south, color=red]{$\|\pi^2\{2-\delta, y_0-1/2\}\|_\cA + H$};
\draw (8,7.7) node[anchor=east, color=red]{$\|\pi^2\{4, -1-2\delta\}\|_\cA + H^2$};
\draw (16.5,4) node[anchor=east, color=red]{$\cA_2(1/2)+H^{-1}$};
\draw (18,5.5) node[anchor=east, color=red]{$\|\pi^2\{4,1-2\delta\}\|_\cA+H$};

\draw (9.5,4.1) node[anchor=west, color=blue]{$\cA_2(x)+H^{-1}$};
\draw (9,5) node[anchor=west, color=blue]{$\|\pi^2\{2-\delta, x\}\|_\cA +H$};
\draw (9,6.4) node[anchor=north west, color=blue]{$\|\pi^2\{4,x\}\|_\cA + H^2$};
\draw (19.5,5.7) node[anchor=west, color=blue]{$\|\pi^2\{4,x\}\|_\cA + H^2$};
\end{tiny}

\draw [color=red, fill=red] (3.5,5) circle (1.5pt);
\draw [color=red, fill=red] (8,7.7) circle (1.5pt);
\draw [color=red, fill=red] (16.5,4) circle (1.5pt);
\draw [color=red, fill=red] (18,5.5) circle (1.5pt);

\end{tikzpicture}
\caption{An illustration of the conditions \eqref{eq:A1bounds} to be satisfied by $\cA_1$: the green curves are $\cA_2$ and $\cA_1$ outside $(y_0-1, 2)$, and are $\cF$-measurable; the black curve is $\cA_1$ in $(y_0-1, 2)$, which (give $\cF$) is a Brownian bridge conditioned on staying above $\cA_2$. The red points are upper bounds and the blue curves are lower bounds.}
\label{fig:rwb-esti}
\end{figure}
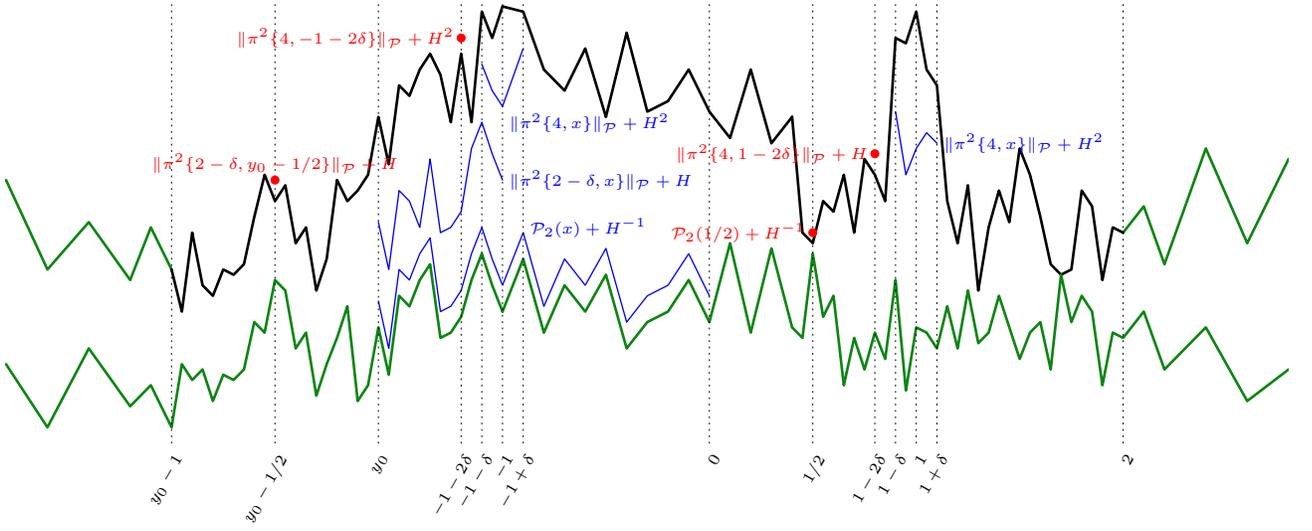

We will make use of a random variable $H>0$, which is $\cF$ measurable, and its exact definition will come up only at the end.
Consider the event $\cE_{ALE,1}'$:
\[
\|\pi^2\{4,-1-2\delta\}\|_\cA - \cA_1(-1-2\delta) > -H^2 > \sup_{-1-\delta \le x \le -1+\delta} \|\pi^2\{4,x\}\|_\cA - \cA_1(x)
\]
as illustrated in Figure \ref{fig:rwb-esti} (at $-1-2\delta$ and in the interval $[-1-\delta, -1+\delta]$).
We note that under $\cE_{ALE,1}'$, for the function $\|\pi^2\{4,\cdot\}\|_\cA - \cA_1$, its supremum in $[-1-\delta, -1+\delta]$ is smaller than its value at $-1-2\delta$.
By the above relation \eqref{eq:rel-geo-2geo} between geodesics and second geodesics, $\cE_{ALE,1}'$ implies that the $\pi\{4,-1+\delta\}$ does not jump to the second line in $[-1-\delta, -1+\delta]$, thus $\pi\{4,-1+\delta\}(-1-\delta)=1$ and $\cE_{ALE,1}$ holds.

Similarly, let $\cE_{ALE,2}'$ be the event where
\[
\|\pi^2\{4,1-2\delta\}\|_\cA - \cA_1(1-2\delta) > -H > \sup_{1-\delta \le x \le 1+\delta} \|\pi^2\{4,x\}\|_\cA - \cA_1(x),
\]
and $\cE_{ALE,3}'$ be the event where
\[
\|\pi^2\{2-\delta, y_0-1/2\}\|_\cA - \cA_1(y_0-1/2) > -H > \sup_{y_0 \le x \le -1} \|\pi^2\{2-\delta, x\}\|_\cA - \cA_1(x),
\]
which are also illustrated in in Figure \ref{fig:rwb-esti}.
We can similarly conclude that $\cE_{ALE,2}'\subset \cE_{ALE,2}$ and  $\cE_{ALE,3}'\subset \cE_{ALE,3}$.

Finally we let $\cE_{ALE,4}'$:
\[
\cA_1(1/2) - \cA_2(1/2) < H^{-1} < \inf_{y_0\le x\le 0} \cA_1(x) - \cA_2(x).
\]
It would be slightly more complicated to get $\cE_{ALE,4}$ from $\cE_{ALE,4}'$.
\begin{lemma}  \label{lem:ale4}
$\cE_{ALE,4}'\subset \cE_{ALE,4}$ when $H$ is large enough (depending on $\cF$).
\end{lemma}
The idea behind is that, using \eqref{eq:ale4-1} (in choosing $y_0$), we can show that the event $\cE_{ALE,4}'$ implies $\|\pi^2\{2+\delta, 1/2\}\|_\cA - \cA_1(1/2) > \sup_{x \le 0} \|\pi^2\{2+\delta, x\}\|_\cA - \cA_1(x)$. By the relation \eqref{eq:rel-geo-2geo}, this further implies that the geodesic $\pi\{2+\delta, 1\}$ jumps to the second line at a point $>0$, thus $\cE_{ALE,4}$ holds.

We leave the detailed proof of Lemma \ref{lem:ale4} to the end, and assume it for now.
Conditional on $\cF$, the above inequalities for the event $\cE_{ALE,1}'\cap \cE_{ALE,2}' \cap \cE_{ALE,3}' \cap \cE_{ALE,4}'$ can be summarized as follows (the inequalities are presented essentially in an order such that the points/intervals involved move from left to right):
\begin{equation}  \label{eq:A1bounds}
\begin{aligned}[c]
\cA_1(y_0-1/2) <& \|\pi^2\{2-\delta, y_0-1/2\}\|_\cA + H ,\\
\cA_1(x) >& \|\pi^2\{2-\delta, x\}\|_\cA +H,\quad \forall  y_0 \le x \le -1 ,\\
\cA_1(x) >& \cA_2(x)+H^{-1},\quad \forall y_0\le x\le 0, \\
{\cA_1(-1-2\delta)} <& \|\pi^2\{4,-1-2\delta\}\|_\cA +H^2 , \\
\cA_1(x) >& \|\pi^2\{4,x\}\|_\cA + H^2,\quad \forall -1-\delta \le x \le -1+\delta , \\
\cA_1(1/2) <& \cA_2(1/2) + H^{-1},\\
\cA_1(1-2\delta) <& \|\pi^2\{4,1-2\delta\}\|_\cA+H ,\\
\cA_1(x) >& \|\pi^2\{4,x\}\|_\cA +H, \quad \forall 1-\delta \le x \le 1+\delta,
\end{aligned}
\quad\quad
\begin{aligned}[c]
& \cE_{ALE,3}', \\
& \cE_{ALE,3}', \\
& \cE_{ALE,4}', \\
& \cE_{ALE,1}', \\
& \cE_{ALE,1}', \\
& \cE_{ALE,4}', \\
& \cE_{ALE,2}', \\
& \cE_{ALE,2}'. \\
\end{aligned}
\end{equation}

Note that by \eqref{eq:rel-geo-2geo}, $\pi^2\{2-\delta, x\}$ and $\pi^2\{4, x\}$ are uniformly bounded for $x$ in any compact interval.
Also note that $\cA_1$ on $[y_0-1, 2]$ is a Brownian bridge conditional on $>\cA_2$ (by the Brownian Gibbs property of $\cA$ mentioned in Section \ref{ss:ale}).
Then by looking up the inequalities in \eqref{eq:A1bounds} (see Figure \ref{fig:rwb-esti}), we can choose $H$ large enough, such that with positive probability, all of  them hold simultaneously, and Lemma \ref{lem:ale4} also holds.
Thus we have $$\PP[\cE_{ALE,1}'\cap \cE_{ALE,2}' \cap \cE_{ALE,3}' \cap \cE_{ALE,4}'\mid\cF]>0.$$
By averaging over $\cF$, we have that
$$\PP[\cE_{coal}'\cap \cE_{disj}]\ge \PP[\cE_{ALE,1}\cap \cE_{ALE,2} \cap \cE_{ALE,3} \cap \cE_{ALE,4}] \ge \PP[\cE_{ALE,1}'\cap \cE_{ALE,2}' \cap \cE_{ALE,3}' \cap \cE_{ALE,4}']>0.$$

It remains to prove Lemma \ref{lem:ale4}.
\begin{proof}[Proof of Lemma \ref{lem:ale4}]
We show $\cE_{ALE,4}'$ implies that
\begin{equation} \label{eq:ale4}
 \|\pi^2\{2+\delta, 1/2\}\|_\cA - \cA_1(1/2) > \sup_{x \le 0} \|\pi^2\{2+\delta, x\}\|_\cA - \cA_1(x),   
\end{equation}
which further implies $\cE_{ALE,4}$, by the relation \eqref{eq:rel-geo-2geo} between geodesics and second geodesics.
Under $\cE_{ALE,4}'$, we have that
\begin{equation}  \label{eq:ale4-15}
\begin{split}
& \|\pi^2\{2+\delta, 1/2\}\|_\cA - \cA_1(1/2) \\ > & \|\pi^2\{2+\delta, 1/2\}\|_\cA - \cA_2(1/2) - H^{-1}
\\ \ge & \sup_{y_0\le x \le 0} \|\pi^2\{2+\delta, x\}\|_\cA - \cA_2(x) - H^{-1}
\\ >& \sup_{y_0\le x \le 0} \|\pi^2\{2+\delta, x\}\|_\cA - \cA_1(x).
\end{split}
\end{equation}
Here the first and third inequality is by $\cE_{ALE,4}'$. The second inequality is essentially a quadrangle inequality. More precisely, it is due to the following reason: for any $y_0\le x \le 0$, any any second geodesic $\pi^2\{2+\delta, x\}$, we can construct a path from $2+\delta$ to $1/2$, which is the same as $\pi^2\{2+\delta, x\}$ on $(-\infty, x)$, and equals $2$ on $[x, 1/2]$. By Lemma \ref{lem:ale-weight-cal}, the weight of this path is $\|\pi^2\{2+\delta, x\}\|_\cA - \cA_2(x) + \cA_2(1/2) \le \|\pi^2\{2+\delta, 1/2\}\|_\cA$.

By \eqref{eq:ale4-1}, we can take $H$ large enough and $\cF$ measurable, such that
\[
\|\pi^2\{2+\delta, 0\}\|_\cA - \cA_2(0) - H^{-1} > \|\pi^2\{2+\delta, y_0\}\|_\cA - \cA_2(y_0).
\]
and hence by the second inequality above, we have 
\begin{equation}  \label{eq:ale4-16}
\|\pi^2\{2+\delta, 1/2\}\|_\cA - \cA_1(1/2) > \|\pi^2\{2+\delta, y_0\}\|_\cA - \cA_2(y_0).
\end{equation}
We also have
\begin{equation}  \label{eq:ale4-2}
\|\pi^2\{2+\delta, y_0\}\|_\cA - \cA_2(y_0)
\ge \sup_{x \le y_0} \|\pi^2\{2+\delta, x\}\|_\cA - \cA_1(x),
\end{equation}
which is again essentially a quadrangle inequality, and follows similar arguments.
Note that \eqref{eq:ale4} is implied by \eqref{eq:ale4-15}, \eqref{eq:ale4-16}, and \eqref{eq:ale4-2}, so the conclusion follows.
\end{proof}

\end{document}